\documentclass{article}
\pdfoutput=1
\usepackage{curves}
\usepackage{epic}
\usepackage{eepic}
\usepackage{amsmath}
\usepackage{amsfonts}
\usepackage{amsthm}
\usepackage{amssymb}
\usepackage{epsfig}

\theoremstyle{plain}
\newtheorem{Th}{Theorem}[section]
\newtheorem{Lem}[Th]{Lemma}
\newtheorem{Cor}[Th]{Corollary}
\newtheorem{Prop}[Th]{Proposition}

\theoremstyle{definition}
\newtheorem{Def}[Th]{Definition}
\newtheorem*{Not}{Notation}
\newtheorem{Ex}[Th]{Example}
\newtheorem{Assum}[Th]{Assumption}

\theoremstyle{remark}
\newtheorem*{Rem}{Remark}

\def\thm{\begin{Th}}
\def\endthm{\end{Th}}
\def\lemma{\begin{Lem}}
\def\endlemma{\end{Lem}}
\def\cor{\begin{Cor}}
\def\endcor{\end{Cor}}
\def\prop{\begin{Prop}}
\def\endprop{\end{Prop}}
\def\definition{\begin{Def}}
\def\enddefinition{\end{Def}}
\def\remark{\begin{Rem}}
\def\endremark{\end{Rem}}
\def\example{\begin{Ex}}
\def\endexample{\end{Ex}}
\def\demo{\begin{proof}}
\def\enddemo{\end{proof}}

\def\notation{\begin{Not}}
\def\endnotation{\end{Not}}
\def\assumption{\begin{Assum}}
\def\endassumption{\end{Assum}}

\def\A{\mathcal{A}}

\def\C{\mathcal{C}}

\def\E{\mathcal{E}}
\def\F{\mathcal{F}}
\def\G{\mathcal{G}}
\def\H{\mathcal{H}}
\def\L{\mathcal{L}}
\def\M{\mathcal{M}}
\def\N{\mathcal{N}}
\def\O{\mathcal{O}}
\def\P{\mathcal{P}}

\def\U{\mathcal{U}}

\def\BbR{\mathbb{R}}
\def\BbN{\mathbb{N}}

\def\BbQ{\mathbb{Q}}

\def\a{\alpha}
\def\b{\beta}
\def\c{\gamma}
\def\d{\delta}
\def\e{\epsilon}

\def\vp{\varphi}
\def\s{\sigma}
\def\GG{\Gamma}
\def\LL{\Lambda}
\def\SS{\Sigma}

\def\bp{\mathbf p}

\def\inte#1{{\rm int}(#1)}
\def\diam#1{{\rm diam}(#1)}
\def\word#1#2{{#1}_1\ldots{#1}_{#2}}
\def\sd#1#2{#1\backslash#2}
\def\word#1#2{{#1}_1\ldots{#1}_{#2}}

\def\wI{\widetilde{I}}

\def\norm#1{||#1||}
\def\supp#1{{\rm supp}({#1})}
\def\tla{\widetilde{\lambda}}
\def\ola{\overline{\lambda}}

\def\wE{\widetilde{\E}}
\def\W{\mathcal{W}}
\def\hE{\widehat{\E}}

\def\olc{\overline{c}}
\def\ulc{\underline{c}}

\def\BB{\mathbb{B}}

\usepackage{euscript}
\renewcommand{\theequation}{\arabic{section}.\arabic{equation}}

\begin{document}
\begin{center}
{\bf \Large Conductive homogeneity of compact metric spaces and construction of $p$-energy}
\par\vspace{3pt}
by
\par\vspace{3pt}
{\large  Jun Kigami}
\par\vspace{3pt}
Graduate School of Informatics\\
Kyoto University\\
e-mail: kigami@i.kyoto-u.ac.jp
\end{center}


\begin{abstract}
In the ordinary theory of Sobolev spaces on domains of $\BbR^n$, the $p$-energy is defined as the integral of $|\nabla{f}|^p$. In this paper, we try to construct a $p$-energy on compact metric spaces as a scaling limit of discrete $p$-energies on a series of graphs approximating the original space.  In conclusion, we propose a notion called conductive homogeneity under which one can construct a reasonable $p$-energy if $p$ is greater than the Ahlfors regular conformal dimension of the space. In particular, if $p = 2$, then we construct a local regular Dirichlet form and show that the heat kernel associated with the Dirichlet form satisfies upper and lower sub-Gaussian type heat kernel estimates. As examples of conductively homogeneous spaces, we present new classes of square-based self-similar sets and rationally ramified Sierpinski crosses, where no diffusions were constructed before. 
\end{abstract}

\tableofcontents

\setcounter{equation}{0}
\section{Introduction}\label{INT}

The main objective of this paper is to generalize the following elementary fact.\par
Let $I = [0, 1]$. Define
\[
\E_p^n(f) = \sum_{i = 1}^{2^n} \Big|f\Big(\frac {i - 1}{2^n}\Big) - f\Big(\frac i{2^n}\Big)\Big|^p
\]
for $n \ge 1$ and $f : I  \to \BbR$. If $f$ is smooth or more generally $f \in W^{1, p}(I)$, which is the $(1, p)$-Sobolev space, then
\[
(2^{p - 1})^n\E_p^n(f) \to \int_0^1 |\nabla{f}|^pdx
\]
as $n \to \infty$, where $\nabla{f}$ is the derivative of $f$. \par
Our naive question is what is a counterpart of this in the case of metric spaces. More precisely, our general strategy of the study is:\\
(1)\, To fix an adequate sequence of discrete graphs $\{(T_n, E_n^*)\}_{n \ge 1}$, where $T_n$ is a discrete approximation of the original metric space $(X, d)$ and $E_n^*$ is the collection of edges, i.e. pairs of points in $T_n$. For a function $f: T_n \to \BbR$, define
\[
\E_p^n(f) = \frac 12\sum_{(x, y) \in E_n^*} |f(x) - f(y)|^p,
\]
which is called the $p$-energy of the function $f$.\\
(2)\,\, To find a proper scaling constant $\s$ such that the space of functions
\[
\{f: X \to \BbR| \text{$\s^n\E_p^n(P_nf)$ is ``convergent'' as $n \to \infty$}\},
\]
where $P_nf$ is a suitable discrete approximation of $f$, is rich enough to be a ``Sobolev'' space in some sense. From our perspective, we do not care about the existence of a derivative $\nabla{f}$ but pursue the convergence of $\s^n\E_p^n(P_nf)$.\par
\begin{figure}
\hspace*{0pt}
\begin{minipage}[b]{6cm}
\centering
\includegraphics[width = 120pt]{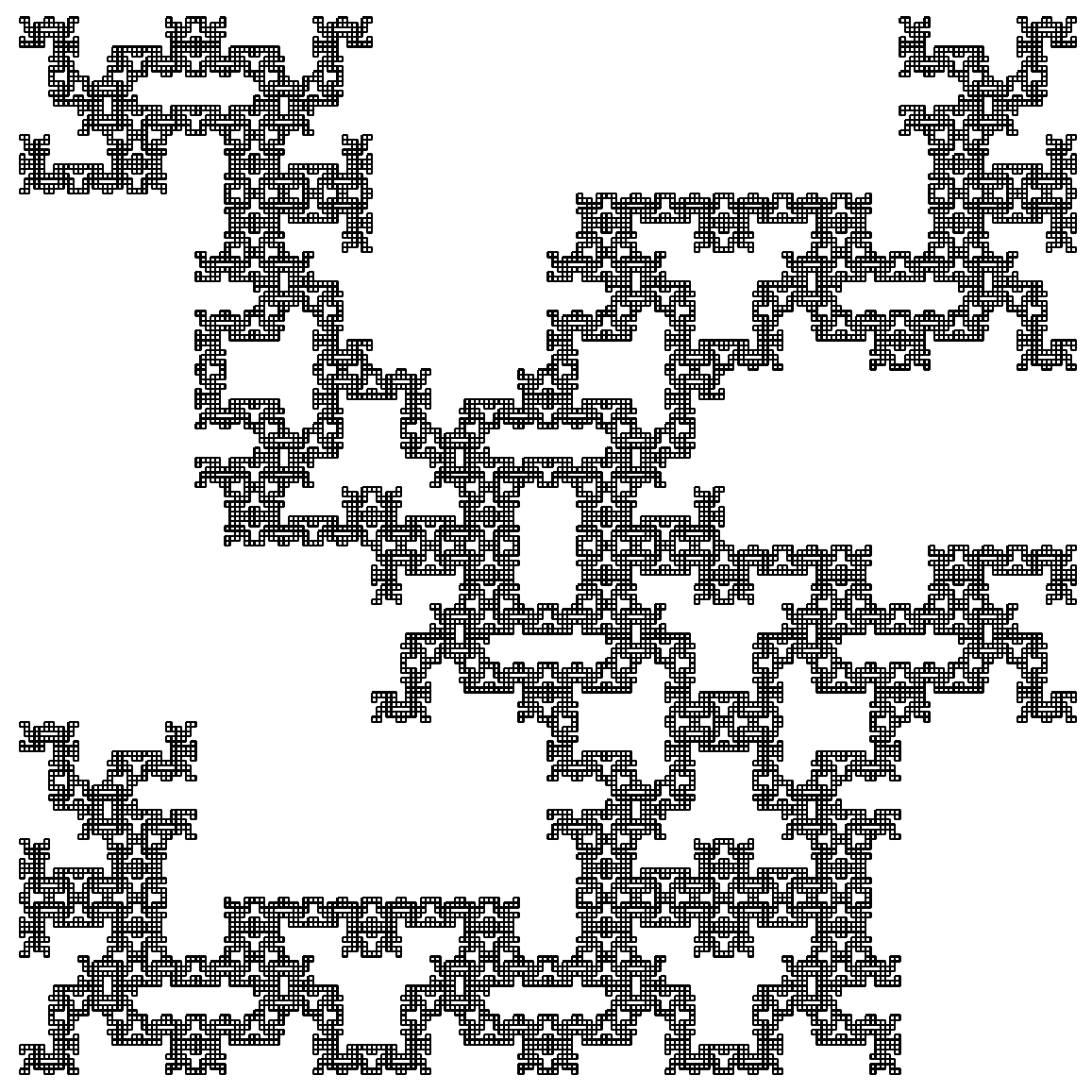}
\end{minipage}
\begin{minipage}[b]{6cm}
\hspace{-20pt}
\centering
\includegraphics[width=114pt]{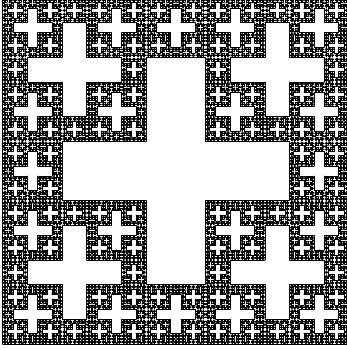}
\vspace{3pt}
\end{minipage}
\caption{Square-based self-similar sets}\label{FIGINT}
\end{figure}

Actually, in the case $p = 2$, this strategy was employed to construct Dirichlet forms inducing diffusion processes on self-similar sets like the Sierpinski gasket and the Sierpinski carpet. (See Figure~\ref{Selfsimilar}.) For the sake of simplicity, we confine ourselves to non-finitely ramified self-similar sets. (This excludes post critically finite self-similar sets represented by the Sierpinski gasket.)  Barlow and Bass constructed the Brownian motions on (generalized) Sierpinski carpets in \cite{BB1, BB2, BB3, BB4, BB5, BB6} as scaling limits of the Brownian motions on regions approximating Sierpinski carpets. Later Kusuoka and Zhou employed the above strategy for $p = 2$ and directly constructed the Dirichlet form inducing the Brownian motion on the planar Sierpinski carpet in \cite{KusZ1}. Note that all these works were done in the last century. Although more than 20 years have passed, no essential progress has been made on the construction of diffusion processes/Dirichlet forms on non-finitely ramified self-similar sets. In particular, no diffusion was constructed on square-based non-finitely ramified self-similar sets like those in Figure~\ref{FIGINT}. The right-hand one is an example of rationally ramified Sierpinski crosses treated in Section~\ref{SCR}. It has two different contraction ratios. The left-hand one is an example having no symmetry of the square.  As a by-product of our results in this paper, we will construct non-trivial self-similar local regular Dirichlet forms on classes of square-based self-similar sets including those in Figure~\ref{FIGINT}. See Sections~\ref{CYP}, \ref{ESC}, and \ref{SCR} for details.\par

 From the viewpoint of construction of Sobolev spaces on metric spaces, there have already been established theories based on upper gradients, which correspond to local Lipschitz constants of Lipschitz functions. Compared with our strategy above, this direction is to seek a counterpart of $\nabla{f}$ instead of the convergence of $\s^n\E_p^n(P_nf)$ like us. The pioneering works of this theory are  Haj{\l}asz\cite{Haj1}, Cheeger\cite{Cheeger} and Shanmugalingam\cite{Shanm1}. One can find a panoramic view of this theory in \cite{HeiKoShTy}. Recent studies by Kajino and Murugan in \cite{KajMur2} and \cite{KajMur3}, however, have suggested that they may not cover all the interesting cases. So far examples in question are higher dimensional Sierpinski gaskets, the Vicsek set, and the planar Sierpinski carpet. What they have shown in \cite{KajMur2} and \cite{KajMur3} is that the Brownian motions on those examples will not have the Gaussian heat kernel estimate under any time change by a pair $(d, \mu)$, where a metric $d$ is quasisymmetric to the Euclidean metric $d_E$ and a measure $\mu$ has the volume doubling property with respect to $d_E$. On the other hand, under the established theory, the heat kernel associated with a $(1, 2)$-Sobolev space satisfying a $(2, 2)$-Poincar{\'e} inequality should satisfy the Gaussian estimate due to the results in \cite{Gri1,Saloff1} and \cite{Sturm1}. Thus, the Dirichlet forms associated with the Brownian motions on the above-mentioned self-similar sets can hardly be one of $(1, 2)$-Sobolev spaces based on upper gradients. Note that, in these cases, there exist plenty of rectifiable curves with respect to  (the restriction of) the Euclidean metrics, which are even quasiconvex. Partly motivated by such a situation, we will try to provide an alternative theory of function spaces, which may be called Sobolev spaces or else, on metric spaces, and to construct natural diffusion processes at the same time.\par
Getting straight to the conclusion, we propose a condition called $p$-conductive homogeneity and show that under this condition, the strategy consisting of (1) and (2) succeeds for $p > \dim_{AR}(K, d)$, where $\dim_{AR}(K, d)$ is the Ahlfors regular conformal dimension of a compact metric space$(K, d)$. One can see a more precise and detailed exposition in what follows. The definition of the Ahlfors regular conformal dimension of $(K, d)$ is
\begin{multline}\label{INT.eq00}
\dim_{AR}(K, d) = \inf\{\a| \text{there exist a metric $\rho$ on $K$ which is}\\
\text{ quasisymmetric to $d$ and a Borel regular measure $\mu$}\\
\text{ which is $\a$-Ahlfors regular with respect to $\rho$}.\}
\end{multline}
\par
In the direction of our study, Shimizu has done pioneering work for the case of the planar Sierpinski carpet, PSC for short, in the very recent paper \cite{Shimizu2}. Extending Kusuoka-Zhou's method, he has constructed a $p$-energy and the corresponding $p$-energy measure for $p > \dim_{AR}({\rm PSC}, d_E)$, and done detailed analysis of those objects. In particular, he has shown that the collection of functions with finite $p$-energies is a Banach space that is reflexive and separable. His proof of reflexivity and separability can be easily extended to our general case as well. See Theorem~\ref{CPE.thm15} for details.\par

\begin{figure}
\centering
\centering
\includegraphics[width=\linewidth]{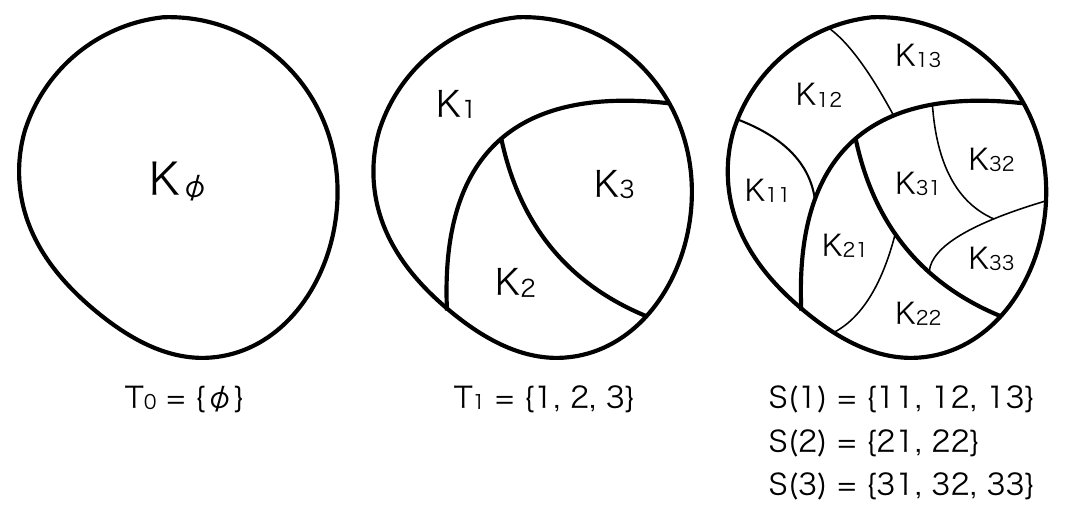}
\vspace{-30pt}
\caption{Partition}\label{FIGINT2}
\end{figure}

Our framework on metric spaces is the theory of partitions introduced in \cite{GAMS}. Let $(K, d)$ be a compact metric space. We always suppose that $(K, d)$ is connected in this paper. Roughly speaking, a partition of $K$ is a sequence of successive divisions of $K$ by some of its compact subsets. The idea is illustrated in Figure~\ref{FIGINT2}.  Let $T_0 = \{\phi\}$ and set $K_{\phi} = K$. Starting from $K$, we first divide $K$ into finite number of children $K_w$ for $w \in T_1$, i.e.
\[
K = \bigcup_{w \in T_1} K_w.
\]
$T_1$ is thought of as the collection of its children of $T_0$ and denoted by $S(\phi)$. Then we repeat this process of division, i.e. each $w \in T_1$ has a collection of its children, $S(w)$, such that
\[
K_w = \bigcup_{v  \in S(w)} K_v.
\]
Define $T_2$ as the disjoint union of the $S(w)$'s for $w \in T_1$. So repeating this process inductively, we have $\{T_n\}_{n \ge 0}$ where each $w \in T_n$ has the collection of its children $S(w) \subseteq T_{n + 1}$. Set 
\[
T = \bigcup_{n \ge 0} T_n.
\]
With several requirements described in Section~\ref{FRW}, the family $\{K_w\}_{w \in T}$ is called a partition of $K$. \par
\begin{figure}
\centering
\includegraphics[width=\linewidth]{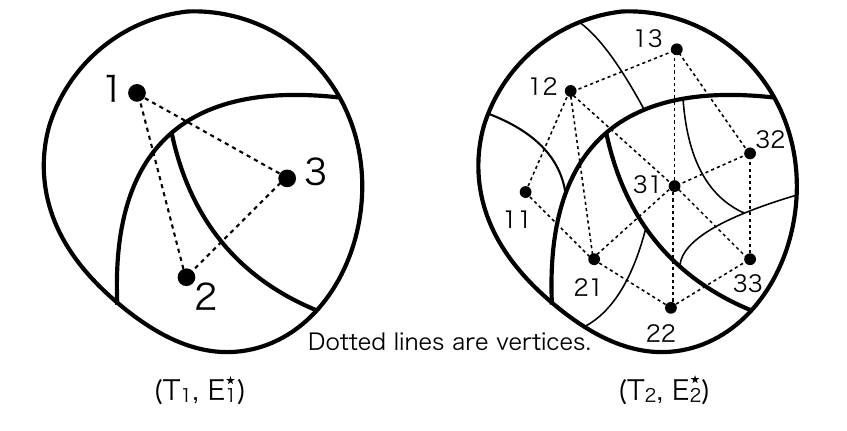}
\vspace{-30pt}
\caption{Graphs associated with a partition}\label{FIGINT3}
\end{figure}
\begin{figure}
\hspace*{0pt}
\begin{minipage}[b]{5cm}
\centering
\includegraphics[width = 120pt]{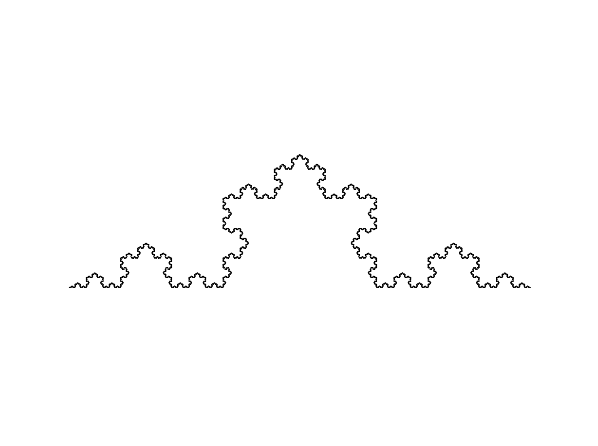}
\vspace{-25pt}
\end{minipage}
\begin{minipage}[b]{3cm}
\hspace{-30pt}
\centering
\includegraphics[width=80pt]{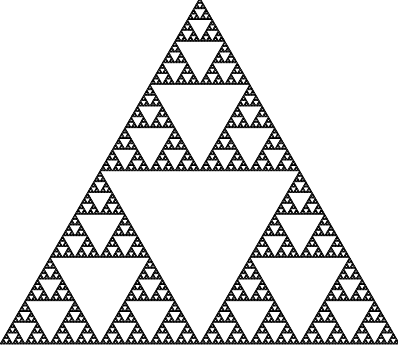}
\end{minipage}
\begin{minipage}[b]{3cm}
\centering
\includegraphics[width = 75pt]{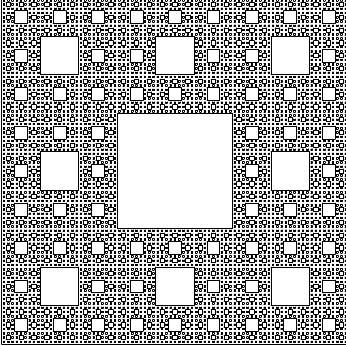}
\end{minipage}
\caption{von Koch curve, Sierpinski gasket and Sierpinski carpet}\label{Selfsimilar}
\end{figure}

For each $n \ge 1$,  $T_n$ has a natural graph structure associated with a given partition $\{K_w\}_{w \in T}$.  Namely, if
\[
E_n^* = \{(u, v) | u, v \in T_n, K_u \cap K_v \neq \emptyset\}, 
\]
then $(T_n, E_n^*)$ is a connected graph, which is illustrated in Figure~\ref{FIGINT3}. To avoid technical complexity, we are going to explain our results under  Assumption~\ref{ALFR} hereafter in the introduction. In fact, if $(K, d)$ is $\a$-Ahlfors regular for some $\a$ and the metric $d$ is $1$-adapted in the sense of \cite{GAMS} ,  then Assumption~\ref{ALFR} holds. So our setting should be broad enough.

 For $A \subseteq T_n$, we define the $p$-energy of a function on $A$ by
\[
\E_{p, A}^n(f) = \frac 12\sum_{\substack{u, v \in A\\(u, v) \in E_n^*}} |f(u) - f(v)|^p.
\] 
To carry out our strategy, we introduce two key characteristic quantities which are conductance and neighbor disparity constants. For $A_1, A_2, A \subseteq T_n$ with $A_1, A_2 \subseteq A$ and $A_1 \cap A_2 = \emptyset$, define the $p$-conductance between $A_1$ and $A_2$ in $A$ at the level $m$ by
\[
\E_{p, m}(A_1, A_2, A) = \inf\{\E_{p, A}^{n + m}(f)| f: S^m(A) \to \BbR,  f|_{S^m(A_1)} \equiv 1, f|_{S^m(A_2)} \equiv 0\},
\]
where $S^m(A) \subseteq T_{n + m}$ is the collection of the descendants in the $m$-th generation from $A$.

\remark
Attaching a resistor of resistance $1$ to each edge $(u, v) \in E^*_{n + m}$, we may consider the graph $(T_{n + m}, E_{n + m}^*)$ as an electric network. In this respect, the reciprocal of $\E_{2, m}(A_1, A_2, A)$ is the effective resistance between $A_1$ and $A_2$ within $A$ and hence $\E_{2, m}(A_1, A_2, A)$ corresponds to the effective conductance. Such an analogy has been often used in the study of random walks. See \cite{DS} for a classical reference. In potential theory, the quantity $\E_{2, m}(A_1, A_2, A)$ is called ``capacity'' as well.
\endremark

In particular, for $w \in T_n$, define 
\[
\E_{p, m}(w) = \E_{p, m}(\{w\}, \GG_1(w)^c, T_n),
\]
where $\GG_1(w)$ is the collection of neighbors of $w$ in $T_n$ given by 
\[
\GG_1(w) = \{v| v \in T_n, (w, v) \in E_n^*\}.
\]
The value $\E_{p, m}(w)$ represents the $p$-conductance between $w$ and the complement of its neighborhood $\GG_1(w)$ in the $m$-th generation from $w$.  In \cite{GAMS}, it  was shown that 
\begin{equation}\label{INT.eq100}
\limsup_{m \to \infty} \Big(\sup_{w \in T} \E_{p, m}(w)^{\frac 1m}\Big) < 1\,\,\,\,\text{if and only if}\,\,\,\, p > \dim_{AR}(K, d).
\end{equation}

The other one, the neighbor disparity constant, is defined as
\[
\s_{p, m, n} = \sup_{(w, v) \in E_n^*}\Bigg(\sup_{f: S^m(w, v) \to \BbR}\frac{ |(f)_{S^m(w)} - (f)_{S^m(v)}|^p}{\E_{p, S^m(w, v)}^{n + m}(f)}\Bigg),
\]
where $S^m(w, v) = S^m(w) \cup S^m(v)$ and $(f)_{S^m(w)}$ is the average of $f$ on $S^m(w)$ under a suitable measure $\mu$. For the case $p = 2$, this constant was introduced in \cite{KusZ1}. The neighbor disparity constant controls the difference of means of a function on neighboring cells via the $p$-energy.\par
And now, $p$-conductive homogeneity, which is the principal notion of this paper, is defined as follows.

\definition
$(K, d)$ is said to be $p$-{\it conductively homogeneous} if and only if there exists $c > 0$ such that
\[
\sup_{w \in T} \E_{p, m}(w)\sup_{n \ge 1} \s_{p, m, n} \le c
\]
for any $m \ge 1$.
\enddefinition

 At a glance, it does not quite look like ``homogeneity''. The following theorem, however, gives the legitimacy of the name.

\thm[Theorem~\ref{SUC.thm10}]\label{INT.thm10}
$(K, d)$ is $p$-conductively homogeneous if and only if there exist $\s > 0$ and $c_1, c_2 > 0$ such that
\[
c_1\s^{-m} \le \E_{p, m}(w) \le c_2\s^{-m}
\]
for any $w \in \sd{T}{\{\phi\}}$ and $m \ge 1$ and
\[
c_1\s^m \le \s_{p, m, n} \le c_2\s^m
\]
for any $m, n \ge 1$.
\endthm

The next natural question is how the conductive homogeneity is related to construction of a $p$-energy. The answer is the next theorem which follows by combining Theorems~\ref{CPE.thm01}, \ref{CPE.thm10}, \ref{CPE.thm20} and Lemma~\ref{SUC.lemma20}.

\thm\label{INT.thm20}
Suppose $p > \dim_{AR}(K, d)$ and $(K, d)$ is $p$-conductively homogeneous. Let $C(K)$ be the collection of continuous functions on $K$. Define
\[
\N_p(f) = \Big(\sup_{m \ge 0}\s^m\E_p^m(P_mf)\Big)^{\frac 1p}
\]
for $f \in L^p(K, \mu)$, where
\[
(P_mf)(w) = \frac 1{\mu(K_w)}\int_{K_w} f(x)\mu(dx),
\]
and
\[
\W^p = \{f| f \in L^p(K, \mu), \N_p(f) < \infty\}.
\]
Then $\N_p(f) = 0$ if and only if $f$ is constant on $K$, $\N_p$ is a semi-norm of $\W^p$, $(\W^p, \norm{\cdot}_{p, \mu} + \N_p(\cdot))$ is a Banach space and $\W^p$ is a dense subset of $(C(K), \norm{\cdot}_{\infty})$. Moreover, there exists $\hE_p: \W^p \to [0, \infty)$ such that $\hE_p^{\frac 1p}$ is a semi-norm of $\W^p$ which is equivalent to $\N_p(\cdot)$, $\hE_p$ satisfies the Markov property and there exist $\tau > 0$ and $c_1, c_2 > 0$ such that
\[
c_1d(x, y)^{\tau} \le \sup_{\substack {f \in \W^p\\ \hE_p(f) \neq 0}}\frac{|f(x) - f(y)|^p}{\hE_p(f)} \le c_2d(x, y)^{\tau}
\]
for any $x, y \in K$. In particular, for $p = 2$, one can choose $(\hE_2, \W^2)$ as a local regular Dirichlet form on $L^2(K, \mu)$.
\endthm

Note that by \eqref{INT.eq100}, the condition $p > \dim_{AR}(K, d)$ implies $\s > 1$. An explicit description of the constant $\tau$ is given in Lemma~\ref{SUC.lemma20}. In addition, we show a sub-Gaussian type heat kernel estimate for the diffusion process induced by the Dirichlet form $(\hE_2, \W^2)$ in Theorem~\ref{SUC.thm30}. Moreover, if $(K, d)$ is a self-similar set with rationally related contraction ratios, then a self-similar $p$-energy which is equivalent to $\N_p$ will be constructed in Section~\ref{SSF}.   \par
Another important question is how to show conductive homogeneity. The following theorem provides an equivalent and useful condition for this purpose.
\thm[Theorem~\ref{SUC.thm20}]\label{INT.thm30}
Suppose that $p > \dim_{AR}(K, d)$. $(K, d)$ is $p$-conductively homogeneous if and only if, for any $k \ge 1$,  there exists $c(k) > 0$ such that
\begin{equation}\label{INT.eq200}
\sup_{z \in T}\E_{p, m}(z) \le c(k)\E_{p, m}(u, v, S^k(w))
\end{equation}
for any $m  \ge 1$, $w \in T$ and $u, v \in S^k(w)$ with $u \neq v$.
\endthm

The condition \eqref{INT.eq200}, which is the same as \eqref{SUC.eq40} in Theorem~\ref{SUC.thm20},  is a relative of the ``Knight move'' condition in \cite{KusZ1} described in the terminology of random walks, although the word ``Knight move'' does not make sense from the appearance of \eqref{INT.eq200} any longer. The original ``Knight move'' in \cite{BB1} was the name of an argument based on the symmetry of the Sierpinski carpet to show a probabilistic counterpart of \eqref{INT.eq200}.  With certain symmetries of the space, it is possible to show the condition \eqref{INT.eq200} by the method of combinatorial modulus in \cite{BouKleiner}.  Applying Theorem~\ref{INT.thm30}, we are going to show the conductive homogeneity for examples like those in Figure~\ref{FIGINT} in Sections~\ref{ESC} and \ref{SCR}.\par
Besides applications, Theorem~\ref{INT.thm30} has a remarkable theoretical consequence; conductive homogeneity is determined only by conductance constants and independent of the neighbor disparity constants if $p > \dim_{AR}(K, d)$. This is the reason conductive homogeneity is called ``{\it conductive}'' homogeneity.
\par
 The major methodological backgrounds of this paper are Kusuoka-Zhou's arguments in \cite{KusZ1} and combinatorial moduli of path families on graphs introduced in \cite{BouKleiner}. On many occasions, we will extend Kusuoka-Zhou's results to compact metric spaces and to general values of $p$. On such occasions, we will put a reference to the original result by Kusuoka and Zhou right behind the number of the propositions and lemmas like ``Lemma 5.2\,\,(\cite[(2.12) lemma]{KusZ1})''.  Beyond Kusuoka-Zhou's arguments, the notion of combinatorial modulus will play a crucial role on several occasions. The most important one is in the proof of a sub-multiplicative inequality of conductance constants, Corollary~\ref{CMS.cor10}. Moreover, by Lemma~\ref{UFD.lemma50}, one can compare moduli of different graphs and this lemma is indispensable for showing \eqref{INT.eq200} in Sections~\ref{CYP} and \ref{SCR}.\par
Regrettably, we do not have much for the case $p \le \dim_{AR}(K, d)$. In Section~\ref{CPE2}, we will construct a function space $\W^p$ and a semi-norm $\hE_p$ on $\W^p$ under $p$-conductive homogeneity for $p \in [1, \dim_{AR}(K, d)]$. In this case, however, $\W^p$ is given as a subspace of $L^p(K, \mu)$ and we do not know whether $\W^p \cap C(K)$ is dense in $(C(K), \norm{\cdot}_{\infty})$ or not. This is due to the lack of an elliptic Harnack principle of $p$-harmonic functions on the corresponding graphs. In the case $p = 2$, using the coupling method, Barlow and Bass conquered this difficulty for higher dimensional Sierpinski carpets in \cite{BB5} and \cite{BB6}. We have little idea what is an analytic counterpart of the coupling method at this moment. It is a big open problem for future work. In particular, it is interesting to know whether the following naive conjecture is true or not.\par\vspace{5pt}\noindent
{\bf Conjecture.} \,\,$\W^p \subseteq C(K)$ if and only if $p > \dim_{AR}(K, d)$.
\par\vspace{5pt}

Now we briefly explain what happens in the cases of familiar examples.\\
{\bf 1.} Unit (hyper)cube $[-1, 1]^n$:\,\,In this case, for any $p > n$, 
\[
\W^p = W^{1, p}([-1, 1]^n)
\]
and there exists $c > 0$ such that
\[
c\hE_p(f) \le \int_{[-1, 1]^n} |\nabla{f}|^pdx \le c^{-1}\hE_p(f)
\]
for any $f \in W^{1, p}([-1, 1]^n)$. See Example~\ref{CYP.ex70} for details. Even if $p \in [1, n]$, the above results should be true but we do not have any proof for now.\\
{\bf 2.} von Koch curve (Figure~\ref{Selfsimilar}):\,\,The von Koch curve does not contain any rectifiable curve, so that the approaches using upper gradients do not work from the beginning. However, our theory does not distinguish metric spaces which are snowflake equivalent, i.e. two metric spaces $(X, d_X)$ and $(Y, d_Y)$ are snowflake equivalent if there exist  a homeomorphism $\vp: X \to Y$, $c_1, c_2 > 0$  and $\a > 0$ such that
\[
c_1d_X(x_1, x_2)^{\a} \le d_Y(\vp(x_1), \vp(x_2)) \le c_2d_X(x_1, x_2)^{\a}
\]
for any $x_1, x_2 \in X$. Since the von Koch curve is snowflake equivalent to the unit interval $[0, 1]$, we see that $\W^p$ for the von Koch curve equals $W^{1, p}([0, 1])$ for any $p > 1$.\\
{\bf 3.} Planar Sierpinski carpet (Figure~\ref{Selfsimilar}): \,\,As is mentioned above, this is one of the original motivations of this paper and it is expected that our space $\W^p$ is quite different from what one may get from the upper gradient approaches.  By Theorem~\ref{CYP.thm100}, the planar Sierpinski carpet $K$ is shown to be $p$-conductive homogeneous for any $p > \dim_{AR}(K, d_*)$, where $d_*$ is the restriction of the Euclidean metric. Moreover, let $\a_H = \frac {\log 8}{\log 3}$ and let $\beta_p = \frac{\log{8\s}}{\log3}$ , where $\s$ is the exponent appearing in Theorem~\ref{INT.thm10}.  Then by \cite[Theorem~2.19]{Shimizu2}, we have a fractional  Korevaar-Shoen type expression of $\W^p$ as follows:
\[
\W^p = \bigg\{f \bigg| f \in L^p(K, \mu), \limsup_{r \downarrow 0} \int_K\frac 1{r^{\a_H}}\int_{B_{d_*}(x, r)} \frac{|f(x) - f(y)|^p}{r^{\beta_p}}dxdy < \infty\bigg\}.
\]
Furthermore it is shown in \cite{Shimizu2} that $\b_p > p$. This fact implies that $\W^p$ should not coincide with any of the spaces obtained by approaches using upper gradients.\\
{\bf 4.} Sierpinski gasket (Figure~\ref{Selfsimilar}):\,\,Let $K$ be the standard Sierpinski gasket and let $d_*$ be the restriction of the Euclidean metric. Since $K$ is one of nested fractals and $\dim_{AR}(K, d_*) = 1$, Theorem~\ref{OEX.thm10} yields that $K$ is $p$-conductively homogeneous for any $p > 1$. Arguments analogous to those in \cite[Section~5.3]{Shimizu2} give the same fractional Korevaar-Shoen type expression of $\W^p$ as the planar Sierpinski carpet. In this case, $\a_H = \frac{\log 3}{\log 2}$ and $\b_p = \frac{\log {3\s}}{\log 2}$. We expect that $\b_p > p$ for any $p > 1$. In fact, due to \cite{BP}, we know $\b_2 = \frac{\log 5}{\log 2} > 2$. Moreover, $\b_p/p$ is monotonically decreasing by \cite[Lemma~4.7.3]{GAMS}. So at least for $p \in (1, 2]$,  $\b_p > p$ and the space $\W^p$ does not seem to be obtained by the upper gradient approaches. However in this case, if we replace the Euclidean metric with the harmonic geodesic metric and the Hausdorff measure with the Kusuoka measure, then the heat kernel associated with the new pair of the metric and the measure has the Gaussian estimate. See \cite{Ki15} for details. Consequently, the Cheeger theory \cite{Cheeger} is now in place for $\W^2$ at least. On the other hand, the replacement of the metric and the measure causes a change of the partition and, consequently, a change of the associated function space $\W^P$. So, we expect that $\W^p$ associated with the new pair may coincide with those obtained from the approaches based on upper gradients but we have no proof so far.
\par
Before the conclusion of the introduction, we mention two related works. The first one is \cite{BonkSaks}, where the authors constructed another type of  ``Sobolev spaces'' ${\dot A}_p(X)$  on a compact metric space $(Z, d)$  from its hyperbolic fillings $X$. The method is to construct a discretization  $Pf$ on $X$ of $f \in L^1(Z)$, and to consider the weak $\ell^p$-norm of the gradient of $Pf$. Their space ${\dot A}^p(Z)$ seems closely related to our space $\W^p$ but we merely know that $\W^p \subseteq {\dot A}_p(X)$ under suitable assumptions at this point. The second one is \cite{HerPeiStr}, where the authors constructed a $p$-energy on Sierpinski gasket type self-similar sets by extending the notion of harmonic structures in the case of $p = 2$ for post critically finite self-similar sets.  Their $p$-energy should be equivalent to ours, although they did not show the completeness of the domain of their $p$-energy. Despite the fact that their method can work only for finitely ramified self-similar sets even if $p =2$,  their work is the first pioneering study to construct a $p$-energy by renormalizing discrete counterparts.
\par
The organization of this paper is as follows. In Section~\ref{FRW}, we review the basics of partitions of compact metric spaces and then give a framework of this paper including standing assumptions, Assumptions~\ref{ASS.10}, \ref{ASS.15}, \ref{ASS.20} and \ref{ASS.30}. In the end, we present Assumption~\ref{ALFR}, which is stronger than the combination of all the assumptions above but more concise. In Section~\ref{CON}, we introduce the notion of conductance constant which is one of two principal quantities of this paper and we show the existence of a partition of unity associated with the conductance constant. In Section~\ref{CMS}, we introduce the notion of combinatorial moduli of path families on graphs and show a sub-multiplicative inequality for conductance constants using them. In Section~\ref{REL}, we introduce the other principal quantity, the neighbor disparity constant and show its relation with the conductance constant and a sub-multiplicative inequality of them. In Section~\ref{CPE}, we construct our function space $\W^p$ and the $p$-energy $\hE_p$ under Assumption~\ref{ID} and show Theorem~\ref{INT.thm20}. At the same time, we propose a condition called $p$-conducive homogeneity and show, in Section~\ref{SUC}, that the condition $p > \dim_{AR}(K, d)$ and $p$-conductive homogeneity imply Assumption~\ref{ID}. In Section~\ref{CPE2}, we see what we can do for $p \le \dim_{AR}(K, d)$. In Section~\ref{SUC}, we show Theorem~\ref{SUC.thm10} (= Theorem~\ref{INT.thm10}) and Theorem~\ref{SUC.thm20} (= Theorem~\ref{INT.thm30}). Moreover, in Theorem~\ref{SUC.thm30}, we give a sub-Gaussian type heat kernel estimate for the diffusion process induced by the Dirichlet form $(\E, \W^2)$ given in Section~\ref{CPE}. In Section~\ref{SSF}, we construct a self-similar $p$-energy for self-similar sets with rationally related contraction ratios. In Section~\ref{CHS}, we give a sufficient condition for the conductive homogeneity for self-similar sets. Section~\ref{CYP} is devoted to a class of self-similar sets called subsystems of cubic tiling, for which conductive homogeneity is shown through Theorem~\ref{SUC.thm20}. This class includes the Sierpinski carpets, the Menger curve, and the higher dimensional hypercubes. In Section~\ref{ESC}, we present examples of subsystems of cubic tiling having the conductive homogeneity. Also, Section~\ref{SCR} is devoted to showing conductive homogeneity of rationally ramified Sierpinski crosses. In Sections~\ref{CPO}, \ref{REX} and \ref{PTR}, we give a proof of Theorem~\ref{SUC.thm20}. In Section~\ref{FIN}, we show the conductance, Poincar{\'e} and the neighbor disparity constants are uniformly bounded from below and above. We will briefly discuss the modification of the graph structure in Section~\ref{MGS}. Finally, in Section~\ref{OPR}, we gather open problems and future directions of research. Appendices give basic facts used in this paper.

\setcounter{equation}{0}
\section{Framework}\label{FRW}

In this section, we are going to make our framework of this paper clear. It is based on the notion of partitions of compact metric spaces parametrized by rooted trees, which was introduced in \cite{GAMS}. Roughly speaking, a partition is successive divisions of a given space like the binary division of the unit interval. See \cite{GAMS} for examples. Since this notion is relatively new and unfamiliar to most of the readers, we will give a minimal but detailed account of its definition.\par
 To start with,  we present the basics of graphs and trees.

\definition\label{TWR.def10}
Let $T$ be a countable set and let $\A : T \times T \to \{0, 1\}$ which satisfies $\A(w, v) = \A(v, w)$ and $\A(w, w) = 0$ for any $w, v \in T$. We call the pair $(T, \A)$ a (non-directed) graph with the vertices $T$ and the adjacency matrix $\A$. An element $(u, v) \in T \times T$ is called an edge of $(T, \A)$ if and only if $\A(u, v) = 1$. We often identify the adjacency matrix $\A$ with the collection of edges $\{(u, v)| u, v \in T, \A(u, v) = 1\}$.\newline
(1)\,\,\,
A graph $(T, \A)$ is called locally finite if $\#(\{v | \A(w, v) = 1\}) < \infty$ for any $w \in T$, where $\#(A)$ is the number of elements of a set $A$.\\
(2)\,\,\,For $w_0, \ldots, w_n \in T$, $(w_0, w_1, \ldots, w_n)$ is called a path between $w_0$ and $w_n$ if $\A(w_i, w_{i + 1}) = 1$ for any $i = 0, 1, \ldots n - 1$. A path $(w_0, w_1, \ldots, w_n)$ is called simple if and only if $w_i \neq w_j$ for any $i, j$ with $0 \le i < j \le n$ and $|i - j| < n$. \\
(3)\,\,\,
$(T, \A)$ is called a tree if and only if there exists a unique simple path between $w$ and $v$ for any $w, v \in T$ with $w \neq v$. For a tree $(T, \A)$, the unique simple path between two vertices $w$ and $v$ is called the geodesic between $w$ and $v$ and denoted by $\overline{wv}$. We write $u \in \overline{wv}$ if $\overline{wv} = (w_0, w_1, \ldots, w_n)$ and $u = w_i$ for some $i$.
\enddefinition

Next, we define fundamental notions on trees.

\definition\label{TWR.def20}
Let $(T, \A)$ be a tree and let $\phi \in T$. The triple $(T, \A, \phi)$ is called a rooted tree with a root (or a reference point, see \cite{Woess3} e.g.) $\phi$. \\
(1)\,\,
Define $\pi:T \to T$ by
\[
\pi(w) = \begin{cases}
w_{n - 1} &\text{if $w \neq \phi$ and $\overline{\phi{w}} = (w_0, w_1, \ldots, w_{n - 1}, w_n)$,}\\
\phi &\text{if $w = \phi$}
\end{cases}
\]
 and, for $w \in T$, set 
 \[
 S(w) = \sd{\{v| \pi(v) = w\}}{\{w\}}.
 \]
An element $v \in S(w)$ is thought of as a child of $w$. Moreover,  for any $k \ge 1$, we define $S^k(w)$ inductively as
 \[
 S^{k + 1}(w) = \bigcup_{v \in S(w)} S^k(v).
 \]
which is the collection of descendants in the $k$-th generation from $w$.\\
(2)\,\,
For $w \in T$ and $m \ge 0$, we define 
\[
|w| = \min\{n| n \ge 0, \pi^n (w) = \phi\} \quad \text{and} \quad T_m = \{w| w \in T, |w| = m\}
\]
(3)\,\,
For any $w \in T$, define
\[
T(w) = \{v| \text{there exists $n \ge 0$ such that $\pi^n(v) = w$}\},
\]
which is the collection of all the descendants of $w$.\\
(4)\,\,
Define
\[
\SS = \{(w(i))_{i \ge 0}| w(i) \in T_i\,\,\text{and}\,\,w(i) = \pi(w(i + 1))\,\,\text{for any $i \ge 0$}\}.
\]
For $\omega = (\omega(i))_{i \ge 0} \in \SS$, set $[\omega]_m = \omega(m)$ for $m \ge 0$. An element $(w(i))_{i \ge 0} \in \SS$ is called a geodesic ray starting from $\phi$ in \cite{Woess3}.
\enddefinition

\remark
In \cite{GAMS}, we have used $(T)_n$ and $T_w$ in place of $T_n$ and $T(w)$ respectively.
\endremark

Throughout this paper, $T$ is a countably infinite set and $(T, \A)$ is a locally finite tree satisfying $\#(\{v| (w, v) \in \A\}) \ge 2$ for any $w \in T$.\par
Next, we define partitions.

\definition[Partition]\label{PAS.def20}
Let $(K, \O)$ be a compact metrizable topological space having no isolated point, where $\O$ is the collection of open sets.\\
(1)\,\,
A collection of non-empty compact subsets $\{K_w\}_{w \in T}$ is called a partition of $K$ parametrized by $(T, \A, \phi)$ if and only if it satisfies the following conditions (P1) and (P2). \newline
(P1)\,\,\,
$K_{\phi} = K$ and for any $w \in T$, $K_w$ has no isolated point and
\[
K_w = \bigcup_{v \in S(w)} K_v.
\]
(P2)\,\,\,
For any $\omega \in \SS$, $\cap_{m \ge 0}{K_{[\omega]_m}}$ is a single point.\newline
\enddefinition

Originally in \cite{GAMS}, we did not assume that $K$ is connected to include spaces like the Cantor set. In this paper, however, we will only deal with connected spaces. In such cases, the assumption that $K$ has no isolated point is always satisfied unless $K$ is a single point.\\
As an illustrative example of partitions, we present the case of the unit square $[-1, 1]^2$ as a self-similar set. This is an example of the general construction of partitions associated with self-similar sets discussed in Section~\ref{SSF}. 

\begin{figure}
\centering
\includegraphics[width=\linewidth]{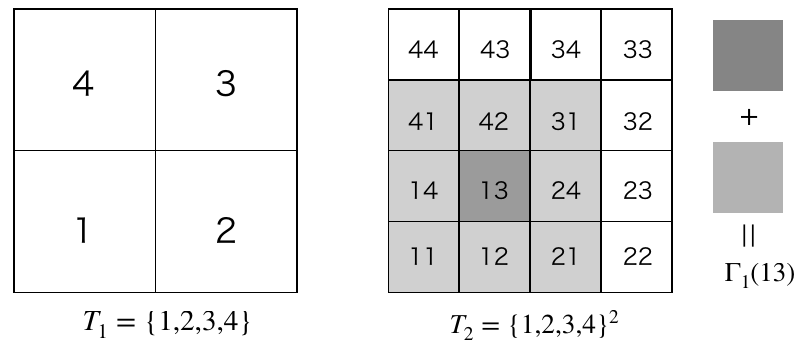}
\vspace{-25pt}
\caption{Partition of the unit square}\label{Square1}
\end{figure}

\example[the unit square]\label{FRW.ex10}
Let $K = [-1, 1]^2$ and let $S = \{1, 2, 3, 4\}$. Set $p_1 = [-1, -1], p_2 = [1, -1], p_3 = [1,1]$ and $p_4 = [-1, 1]$. For $i \in S$, define $f_i(x) = \frac 12(x - p_i) + p_i$ for any $x \in \BbR^2$. Then it is obvious that
\[
K = \bigcup_{i \in S} f_i(K).
\]
This is the expression of the unit square as the self-similar set with respect to the collection of contractions $\{f_i\}_{i \in S}$. Let 
\[
T_n = S^n = \{i_1\ldots{i_n}| i_j \in S\,\,\text{for any $j = 1, \ldots, n$}\}.
\]
 In particular let $T_0 = \{\phi\}$. Moreover define $T = \cup_{m \ge 0} T_m$ and define $\pi: T \to T$ by
\[
\pi(i_1\ldots{i_n}i_{n + 1}) = i_1\ldots{i_n}
\]
for any $i_1\ldots{i_n}i_{n + 1} \in T_{n + 1}$ for $n \ge 1$ and $\pi(\phi) = \phi$. Define $\A(w, v)$ for $w, v \in T$ as $\A(w, v) = 1$ if $\pi(w) = v$ or $\pi(v) = w$ except for $(w, v) = (\phi, \phi)$. Then $(T, \A, \phi)$ is a rooted tree. For $w = \word wn \in T_n$, define 
\[
f_w = f_{w_1}{\circ}\ldots{\circ}f_{w_n}\quad\text{and}\quad K_w = f_w(K).
\]
Then $\{K_w\}_{w \in T}$ is a partition of $K$ parametrized by $(T, \A, \phi)$. See Figure~\ref{Square1}.
\endexample

The following definition is a collection of notions concerning partitions.

\definition\label{ASS.def10}
Let $\{K_w\}_{w \in T}$ be a partition of $K$ parametrized by $(T, \A, \phi)$. 
(1)\,\,Define $O_w$ and $B_w$ for $w \in T$ by
\begin{align*}
O_w &= \sd{K_w}{\Bigg(\bigcup_{v \in \sd{T_{|w|}}{\{w\}}} K_v\Bigg)},\\
B_w &= K_w \cap \Bigg(\bigcup_{v \in \sd{T_{|w|}}{\{w\}}} K_v\Bigg).
\end{align*}
If $O_w \neq \emptyset$ for any $w \in T$, then the partition $K$ is called minimal.\\
(2)\,\,
For any $A \subseteq T_n$ and $w \in A$, define $\GG_M^A(w) \subseteq T_n$ as
\begin{multline*}
\GG_M^A(w) = \{u| u \in A, \text{there exist $u(0), \ldots, u(M) \in A$ such that}\\
\text{$u(0) = w$, $u(M) = u$ and $K_{u(i)} \cap K_{u(i + 1)} \neq \emptyset$ for any $i = 0, \ldots, M - 1$}\}.
\end{multline*}
For simplicity, for $w \in T_n$, we write $\GG_M(w) = \GG_M^{T_n}(w)$.\\
(3)\,\,$\{K_w\}_{w \in T}$ is called uniformly finite if and only if
\[
\sup_{w \in T} \#(\GG_1(w)) < +\infty.
\]
\enddefinition

If a partition is minimal, then $O_w$ is actually the interior of $K_w$ and $B_w$ is the topological boundary of $K_w$. See \cite[Proposition~2.2.3]{GAMS} for details.\par
In the case of the unit square in Example~\ref{FRW.ex10}, $K_w$ is a square and $O_w$ (resp. $B_w$) is the interior (resp. the boundary) of $K_w$. Therefore, it is minimal. Moreover, 
\[
\sup_{w \in T} \#(\GG_1(w)) \le 8,
\]
so that it is uniformly finite.\par
Now we give the first part of our framework in this paper.

 \assumption\label{ASS.10}
 $T$ is a countably infinite set, $\phi \in T$, and $(T, \A)$ is a locally finite tree satisfying $\#(\{w | (w, v) \in \A\}) \ge 2$ for any $w \in T$. $(K, \O)$ is a compact connected metrizable space. $\{K_w\}_{w \in T}$ is a partition of $K$ parametrized by $(T, \A, \phi)$ that is minimal, and uniformly finite.\\
(1)\,\,For any $w \in T$, $K_w$ is connected. \\
(2) There exist $M_*$ and $k_* \in \BbN$ such that 
\begin{equation}\label{FRW.eq100}
\pi^{k_*}(\GG_{M_* + 1}(w)) \subseteq \GG_{M_*}(\pi^{k_*}(w))
\end{equation}
for any $w \in T$.\\
(3)\,\,
There exists $M_0 \ge M_*$ such that
\begin{equation}\label{FRW.eq110}
\GG_{M_*}(u)  \cap S^k(w) \subseteq \GG_{M_0}^{S^k(w)}(u)
\end{equation}
for any $w \in T$, $k \ge 1$ and $u \in S^k(w)$.
 \endassumption
 
\begin{figure}
\centering
\includegraphics[width=\linewidth]{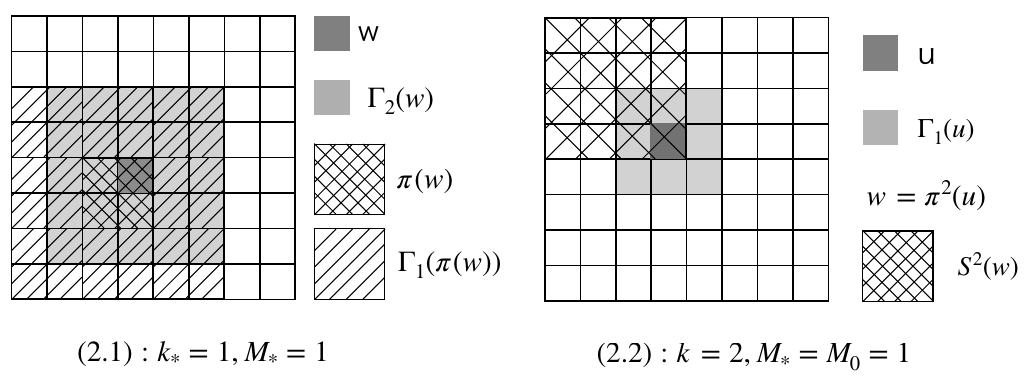}
\vspace{-25pt}
\caption{Assumption~\ref{ASS.10}: the unit square}\label{Square2}
\end{figure}

See Figure~\ref{Square2} for an illustrative exposition of Assumption~\ref{ASS.10} in the case of the unit square.
 
 \remark
As is explicitly mentioned in Proposition~\ref{ASS.prop10}, Assumption~\ref{ASS.10}-(2) is always satisfied under mild additional assumptions.
\endremark

 \remark
If $M_* = 1$, then we have $\GG_{M_*}(w) \cap A = \GG_{M_*}^A(w)$ for any $w$ and $A$. So in this case, by choosing $M_0 = M_* = 1$, Assumption~\ref{ASS.10}-(3) is always satisfied.
\endremark

Throughout this paper, we set
\begin{equation}\label{ASS.eq00}
L_* = \sup_{w \in T} \#(\GG_1(w)).
\end{equation}
Then, for any $m \in \BbN$,
\[
\sup_{w \in T} \#(\GG_m(w)) \le (L_*)^m.
\]

Under Assumption~\ref{ASS.10}-(2), if the partition $\{K_w\}_{w \in T}$ is replaced by the partition $\{K_w\}_{w \in   T^{(k_*)}}$, where $T^{(k_*)} = \cup_{i \ge 0} T_{ik_*}$, the constant $k_*$ can be regarded as $1$. So doing such replacement, we will adopt the following assumption.

\assumption\label{ASS.15}
The constant $k_*$ appearing in \eqref{FRW.eq100} is $1$.
\endassumption
 
For a given partition $\{K_w\}_{w \in T}$, we always associate the following graph structure $E_n^*$ on $T_n$.

\prop\label{FRW.prop10}
For $n \ge 0$, define
\[
E_n^* =\{(w, v)| w, v \in T_n, w \neq v, K_w \cap K_v \neq \emptyset\}.
\]
Then $(T_n, E_n^*)$ is a non-directed graph. Under Assumption~\ref{ASS.10}, $(T_n, E_n^*)$ is connected for any $n \ge 0$, and
\[
\GG_1(w) = \{v| v \in T_n, (w, v) \in E_n^*\}
\]
for any $n \ge 0$ and $w \in T_n$.
\endprop

\remark
In \cite{GAMS}, $E_n^*$ is denoted by $J^h_{1, n}$.
\endremark

\definition\label{FRW.def10}
For $w \in T_n$, define
\begin{multline*}
\partial{S^m(w)} = \{v | v \in S^m(w), \text{there exists $v' \in T_{n + m}$}\\
\text{ such that $(v, v') \in E^*_{n + m}$ and $\pi^m(v') \neq w$.}\}
\end{multline*}
\enddefinition

The set $\partial{S^m(w)}$ is a kind of a boundary of $S^m(w)$. In fact, it is easy to see 
\[
\partial{S^m(w)} = \{v | v \in S^m(w), K_v \cap B_w \neq \emptyset\},
\]
where $B_w$ is the topological boundary of $K_w$ as is mentioned above. So the next assumption means that the boundary is not the whole space.
 
\assumption\label{ASS.20}
There exists $m_0 \ge 1$ such that $\sd{S^m(w)}{\partial{S^m(w)}} \neq \emptyset$ for any $w \in T$ and $m \ge m_0$.
\endassumption

\begin{figure}
\centering
\includegraphics[width=\linewidth]{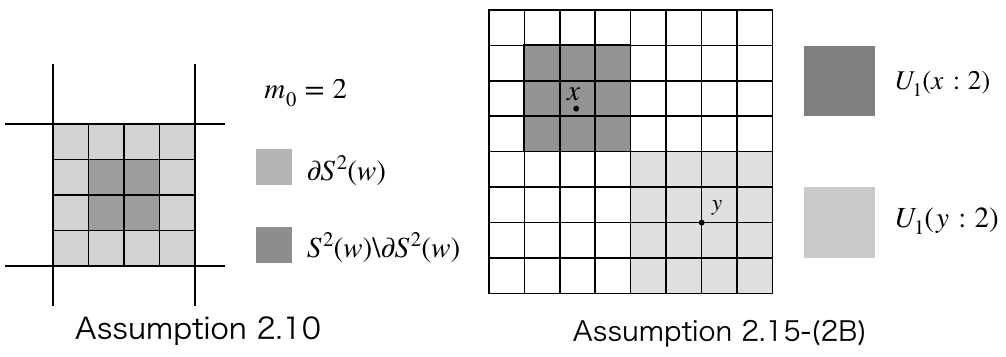}
\vspace{-25pt}
\caption{Assumptions~\ref{ASS.20} and \ref{ALFR}-(2B); the unit square}\label{Square3}
\end{figure}

In Figure~\ref{Square3}, we have an illustrative exposition of Assumption~\ref{ASS.20} in the case of the unit square.

\definition\label{MAC.def20}
For $w \in T$, $M \ge 1$ and $k \ge 1$, define
\[
B_{M, k}(w) = \{v| v \in S^k(w), \GG_{M - 1}(v) \cap \partial{S^k(w)} \neq \emptyset\}.
\]
\enddefinition

\remark
$B_{1, k}(w) = \partial{S^k(w)}$.
\endremark

The final assumption is an assumption on a measure on $K$.

\assumption\label{ASS.30}
$\mu$ is a Borel regular probability measure on $K$ satisfying
\begin{equation}\label{ASS.eq10}
\mu(K_w) = \sum_{v \in S(w)} \mu(K_v)
\end{equation}
for any $w \in T$. There exists $\c \in (0, 1)$ such that
\begin{equation}\label{ASS.eq20}
\mu(K_w) \ge \c\mu(K_{\pi(w)})
\end{equation}
for any $w \in T$. This property is called ``super-exponential'' in \cite{GAMS}. Moreover, there exists $\kappa > 0$ such that if $w, v \in T$, $|w| = |v|$ and $(w, v) \in E_{|w|}^*$, then
\begin{equation}\label{ASS.eq30}
\mu(K_w) \le \kappa\mu(K_v)
\end{equation}
\endassumption

The condition \eqref{ASS.eq30} corresponds to the gentleness of the measure $\mu$ introduced in \cite{GAMS}. Indeed, if $\mu$ has the volume doubling property, then this condition is satisfied. See Proposition~\ref{ASS.prop10} and its proof below for an exact statement.

\lemma\label{ASS.lemma20}
Under Assumptions~\ref{ASS.10}, \ref{ASS.20} and \ref{ASS.30},\\
{\rm (1)}\,\, $\mu$ is exponential, i.e. $\mu$ satisfies \eqref{ASS.eq20} and there exist $m_1 \ge 1$ and $\c_1 \in (0, 1)$ such that $\mu(K_v) \le \c_1\mu(K_w)$ for any $w \in T$ and $v \in S^{m_1}(w)$.\\
{\rm (2)}\,\,
\[
\sup_{w \in T} \#(S(w)) < \infty.
\]
\endlemma

Throughout this paper, we set
\begin{equation}\label{ASS.eq05}
N_* = \sup_{w \in T} \#(S(w)).
\end{equation}

\demo
(1)\,\,In fact, we set $m_1 = m_0$. For any $w$ with $|w| \ge 1$ and $m \ge 0$, we see that $\partial{S^m(w)} \neq \emptyset$ because $K$ is connected. Hence by Assumption~\ref{ASS.20}, $\#(S^{m_0}(w)) \ge 2$ for any $w \in T$. Let $v \in S^{m_1}(w)$. Then there exists $u \in S^{m_1}(w)$ with $v \neq u$.  By \eqref{ASS.eq20},
\[
\mu(K_w) \ge \mu(K_v) + \mu(K_u) \ge \mu(K_v) + \c^{m_1}\mu(K_w),
\]
so that  $\mu(K_v) \le (1 - \c^{m_1})\mu(K_w)$.\\
(2)\,\,
\[
\mu(K_w) = \sum_{v \in S(w)} \mu(K_v) \ge \c\sum_{v \in S(w)} \mu(K_w) = \c\#(S(w))\mu(K_w).
\]
Hence $\#(S(w)) \le 1/\c$.
\enddemo

\lemma\label{ASS.lemma30}
Under Assumptions~\ref{ASS.10}, \ref{ASS.20} and \ref{ASS.30}, 
\[
\sd{S^m(w)}{B_{M, m}(w)} \neq \emptyset
\]
for any $w \in T$, $M \ge 1$ and $m \ge Mm_0$. Moreover,
\begin{equation}\label{ASS.eq40}
\mu\Bigg(\bigcup_{v \in S^n(\sd{S^m(w)}{B_{M, m}(w)})} K_v\Bigg) \ge \c^{m_0M}\mu(K_w)
\end{equation}
for any $w \in T$, $n \ge 0$ and $m \ge Mm_0$.
\endlemma

\demo
By Assumption~\ref{ASS.20}, we can inductively choose $v_i \in S^{im_0}(w)$ for $i \ge 1$ such that $v_{i + 1} \in \sd{S^{m_0}(v_i)}{\partial{S^{m_0}(v_i)}}$ for any $i \ge 1$. At the same time, we see  $v_i  \notin B_{i, im_0}(w)$. If $m_0i < k \le m_0(i + 1)$, then $v  \notin B_{i, k}(w)$ for  $v = \pi^{m_0(i + 1) - k}(v_{i + 1})$.  So the first part of the claim has been verified. Now if $v \in \sd{S^m(w)}{B_{M, m}(w)}$, then 
\[
\mu\Bigg(\bigcup_{v \in S^n(\sd{S^m(w)}{B_{M, m}(w)})} K_v\Bigg) \ge \mu(K_v) \ge \c^{m_0M}\mu(K_w)
\]
by  Assumption~\ref{ASS.30}.
\enddemo

Until now, we have not considered any metric of $(K, \O)$, which was merely assumed to be compact and metrizable. The introduction of a metric having suitable properties enables us to integrate the above assumptions into the following one.

\assumption\label{ALFR}
 $T$ is a countably infinite set, $\phi \in T$,  and $(T, \A)$ is a locally finite tree satisfying $\#(\{w | (w, v) \in \A\}) \ge 2$ for any $w \in T$. $(K, d)$ is a compact connected metric space and $\diam{K, d} = 1$, where $\diam{A, d} = \sup_{x, y \in A} d(x, y)$ for a subset $A \subseteq B$.  $\{K_w\}_{w \in T}$ is a partition of $K$ parametrized by $(T, \A, \phi)$ that is minimal and uniformly finite.\\
(1)\,\,For any $w \in T$, $K_w$ is connected.\\
(2)\,\,There exist $M_* \ge 1$ and $r \in (0, 1)$ such that the following properties (2A), (2B), and (2C) hold;\\
(2A)\,\,Define $h_r: T \to (0, 1]$ as $h_r(w) = r^{|w|}$. Then there exist $c_1, c_2 > 0$ such that 
\[
c_1h_r(w) \le \diam{K_w, d} \le c_2h_r(w)
\]
 for any $w \in T$. \\
(2B)\,\,
For $x \in K$ and $n \ge 1$, define
\[
U_M(x: n) = \bigcup_{\substack{w \in T_n\\x \in K_w}} \bigcup_{v \in \GG_M(w)} K_v
\]
(See Figure~\ref{Square3} for examples of $U_1(\cdot: 2)$ in the case of the unit square.)
Then there exist $c_1, c_2 > 0$ such that
\[
B_d(x, c_1r^n) \subseteq U_{M_*}(x: n) \subseteq B_d(x, c_2r^n)
\]
for any $n \ge 1$ and $x \in K$, where $B_d(x, r) = \{y| d(x, y) < r\}$. \\
(2C)\,\, There exist $c > 0$ such that, for any $n \ge 1$ and $w \in T_n$, there exists $x \in K_w$ such that
\[
K_w \supseteq B_d(x, cr^n).
\]
\noindent(3)\,\,$\mu$ is a Borel regular probability measure on $K$ that is exponential and has the volume doubling property with respect to the metric $d$. Furthermore, $\mu$ satisfies \eqref{ASS.eq10} for any $w \in T$.\\
(4)\,\,There exists $M_0$ such that \eqref{FRW.eq110} holds for any $w \in T$, $k \ge 1$ and $u \in S^k(w)$.\\
(5)\,\,For any $w \in T$, $\pi(\GG_{M_* + 1}(w)) \subseteq \GG_{M_*}(\pi(w))$.
\endassumption

\remark
In the terminology of \cite{GAMS}, the condition (2A) corresponds to the bi-Lipschitz equivalence of $d$ and $h_r$, the condition (2B) says that the metric $d$ is $M_*$-adapted to $h_r$ and (2C) together with (2B) says the metric $d$ is thick. The combination of (2A), (2B) and (2C) is equivalent to that of (BF1) and (BF2) in \cite[Section~4.3]{GAMS}.
\endremark

\remark 
Modifying the original partition $\{K_w\}_{w \in T}$, we always obtain Assumption~\ref{ALFR}-(5)  from Assumption~\ref{ALFR}-(1), (2), (3) and (4).  Namely, by Proposition~\ref{ASS.prop10}, we have $k_*$ satisfying \eqref{FRW.eq100} under Assumption~\ref{ALFR}-(1), (2), (3) and (4). So, replacing the original partition $\{K_w\}_{w \in T}$ with $\{K_w\}_{w \in T^{(k_*)}}$, we may suppose $k_* = 1$.
\endremark

\prop\label{ASS.prop10}
Assumption~\ref{ALFR}-(1), (2), (3) and (4) suffice Assumptions~\ref{ASS.10}, \ref{ASS.20} and \ref{ASS.30}.
\endprop

\demo
About Assumption~\ref{ASS.10}, (1) and (3) are  included in Assumption~\ref{ALFR}. Since $d$ is $M_*$-adapted, \cite[Proposition~4.4.4]{GAMS} shows the existence of $k_*$ required in Assumption~\ref{ASS.10}-(2). By (2C) and (2B), there exists $m_0 \ge 1$ such that
\[
K_w \supseteq B_d(x, cr^n) \supseteq U_{M_*}(x: n + m_0)
\]
for any $n \ge 1$ and $w \in T_n$, where the point $x \in K_w$ is chosen as in (2C). So if $v \in T_{n + m_0}$ and $x \in K_v$, then $K_v \subseteq B_d(x, cr^n)$ and hence $K_v \cap B_w = \emptyset$. Therefore Assumption~\ref{ASS.20} is satisfied. \eqref{ASS.eq10} is included in Assumption~\ref{ALFR} and \eqref{ASS.eq20} follows from the fact that $\mu$ is exponential. Finally, \eqref{ASS.eq30} is a consequence of the volume doubling property by \cite[Theorem~3.3.4]{GAMS}.
\enddemo

Under Assumption~\ref{ALFR}, we may suppose further properties of the metric $d$ and the measure $\mu$. Namely, if $\a > \dim_{AR}(K, d)$, then by \eqref{INT.eq00}, there exist an $\a$-Ahlfors regular metric $d_*$ which is quasisymmetric to $d$ and a Borel regular measure $\nu$ which is $\a$-Ahlfors regular with respect to $d_*$, i.e. there exist $c_1, c_2 > 0$ such that
 \begin{equation}\label{ASS.eq50}
 c_1r^{\a} \le \nu(B_{d_*}(x, r)) \le c_2r^{\a}
 \end{equation}
for any $x \in K$ and $r \in (0, 2\diam{K, d}]$. Replacing $d$ and $\mu$ by $d_*$ and $\nu$ respectively, we may assume that $d$ is $\a$-Ahlfors regular.  Note that if $\mu$ is $\a$-Ahlfors regular with respect to $d$, then $\a$ is the Hausdorff dimension of $(K, d)$.

\setcounter{equation}{0}
\section{Conductance constant}\label{CON}
Hereafter in this paper, we always presume Assumptions~\ref{ASS.10}, \ref{ASS.15}, \ref{ASS.20} and \ref{ASS.30}.\par
In this section, we introduce the conductance constant $\E_{M, p, m}(w, A)$ and show the existence of a partition of unity whose $p$-energies are estimated by conductance constants from above. In the next section, using the method of combinatorial modulus, we will establish a sub-multiplicative inequality of conductance constants. \par
To begin with, we define $p$-energies of functions on graphs $(T_n, E_n^*)$ and the associated $p$-conductances between subsets. 

\notation
Let $A$ be a set. Set
\begin{equation}\label{CON.eq10}
\ell(A) = \{f| f: A \to \BbR\}.
\end{equation}
\endnotation

\definition\label{MAC.def10}
(1)\,\,Let $A \subseteq T_n$. For $f \in \ell(A)$, define $\E_{p, A}^n(u)$ by
\[
\E_{p, A}^n(f) = \frac 12\sum_{u, v \in A, (u, v) \in E_n^*} |f(u) - f(v)|^p.
\]
(2)\,\,Let $A \subseteq T_n$ and let $A_1, A_2 \subseteq A$. 
Define
\[
\E_{p, m}(A_1, A_2, A) = \inf\{\E_{p, S^m(A)}^{n + m}(f)| f \in \ell(S^m(A)), f|_{S^m(A_1)} \equiv 1, f|_{S^m(A_2)} \equiv 0\}.
\]
(3)\,\,Let $A \subseteq T_n$. For $w \in A$, define
\[
\E_{M, p, m}(w, A) \\= \E_{p, m}(\{w\}, \sd{A}{\GG_M^A(w)}, A),
\]
which is called the $p$-conductance constant of $w$ in $A$ at level $m$.
\enddefinition

For simplicity, we often denote a set consisting of a single point, $\{w\}$, by $w$. For example, if $A_1$ and $A_2$ are single points $u$ and $v$ respectively, we sometimes write $\E_{p, m}(u, v, A)$ instead of $\E_{p, m}(\{u\}, \{v\}, A)$. 

\lemma\label{MAC.lemma10}
\[
\E_{M_0, p, m}(u, S^k(w)) \le \E_{M_*, p, m}(u, T_{|w| + k}).
\]
for any $w \in T$, $k \ge 0$ and $u \in S^k(w)$.
\endlemma

\demo
This follows from the assumption \eqref{FRW.eq110}.
\enddemo

\remark
In case $M_* = 1$, we always have $\GG_1^A(w) = \GG_1(w) \cap A$. Hence even without \eqref{FRW.eq110},
\[
\E_{1, p, m}(w, S^k(w)) \le \E_{1, p, m}(w, T_{|w| + k})
\]
for any $w \in T$, $k \ge 0$ and $u \in S^k(w)$.
\endremark

The following lemma shows the existence of a partition of unity.

\lemma\label{MAC.lemma15}
Let $p \ge 1$ and let $A \subseteq T_n$. For any $w \in A$, there exists $\vp_w:S^m(A) \to [0, 1]$ such that
\[
\sum_{w \in A} \vp_w \equiv 1,\,\,\,\vp_w|_{S^m(w)} \ge (L_*)^{-M},\,\,\,\,
\vp_w|_{\sd{S^m(A)}{S^m(\GG_M^A(w))}} \equiv 0
\]
and
\[
\E_{p, S^m(A)}^{n + m}(\vp_w) \le ((L_*)^{2M + 1} + 1)^p\max_{w' \in \GG_{2M + 1}^A(w)}\E_{M, p, m}(w', A).
\]
\endlemma

\demo
Let $h_w \in \ell(S^m(A))$ satisfy $h_w|_{S^m(w)}  \equiv 1$, $h_w|_{\sd{S^m(A)}{S^m(\GG_M^A(w))}} \equiv 0$ and $\E_{M, p, m}(w, A) = \E_{p, S^m(A)}^{n + m}(h_w)$.
Define $h \in \ell(S^m(A))$ as
\[
h(v) = \sum_{w \in A} h_w(v)
\]
for any $v \in S^m(A)$. Note that $1 \le h(v) \le (L_*)^M$. Set $\vp_w = h_w/h$ and $E_{n + m}(w) = E_{n + m}^* \cap S^m(\GG_{M + 1}^A(w))^2$. It follows that $\vp_w(u) = \vp_w(v) = 0$ for any $(u, v) \notin E_{n + m}(w)$. Let $(u, v) \in E_{n + m}(w)$. Then $h_w(v)(h_{w'}(v) - h_{w'}(u)) = 0$ for any $w' \notin \GG_{2M + 1}^A(w)$. Hence
\begin{multline*}
|\vp_w(u) - \vp_w(v)| = \Big|\frac 1{h(u)h(v)}(h(v)(h_w(u) - h_w(v)) + h_w(v)(h(v) - h(u)))\Big| \\
\le |h_w(u) - h_w(v)| + \sum_{w' \in \GG_{2M + 1}^A(w)} |h_{w'}(u) - h_{w'}(v)|.
\end{multline*}
Set $C = (L_*)^{2M + 1} + 1$. Then the last inequality yields
\begin{multline*}
\E_p^{n + m}(\vp_w) = \frac 12\sum_{(u, v) \in E_{n + m}(w)} |\vp_w(u) - \vp_w(v)|^p \\
\le \frac {C^{p - 1}}2\sum_{(u, v) \in E_{n + m}(w)}\Big(|h_w(u) - h_w(v)|^p + \sum_{w' \in \GG_{2M + 1}^A(w)} |h_{w'}(u) - h_{w'}(v)|^p\Big) \\
 \le C^{p - 1}\Big(\E_{p, S^m(A)}^{n + m}(h_w) + \sum_{w' \in \GG_{2M + 1}^A(w)} \E_{p, S^m(A)}^{n + m}(h_{w'})\Big) \\
 \le C^p\max_{w' \in \GG_{2M + 1}^A(w)}\E_{M, p, m}(w', A).
 \end{multline*}
\enddemo

In particular, in the case $A = T_n$, the associated partition of unity defined below will be used to show the regularity of the $p$-energy constructed in Section~\ref{CPE}.

\definition\label{MAC.def30}
For $w \in T$, define $h_{M, w, m}^* \in \ell(T_{|w| + m})$ as the unique function $h$ satisfying $h|_{S^m(w)} = 1$, $h|_{\sd{T_{|w| + m}}{S^m(\GG_M(w))}} = 0$ and 
\[
\E_p^{|w| + m}(h) = \E_{M, p, m}(w, T_{|w|}).
\]
Moreover, define $\vp_{M, w, m}^* \in \ell(T_{|w| + m})$ by
\[
\vp_{M, w, m}^* = \frac{h_{M, w, m}^*}{\sum_{v \in T_{|w|}} h_{M, v, m}^*}.
\]
\enddefinition

By the proof of Lemma~\ref{MAC.lemma15}, 
\[
\E_p^{n + m}(\vp_{M, w, m}^*) \le ((L_*)^{2M + 1} + 1)^p\max_{v \in T_n} \E_{M, p, m}(v, T_n)
\]
for any $w \in T_n$.\par

\setcounter{equation}{0}
\section{Combinatorial modulus}\label{CMS}

Another principal tool of this paper is the notion of combinatorial modulus of a path family of a graph introduced in \cite{BouKleiner}. The general theory is briefly reviewed in Appendix~\ref{UFD}. In this section, we introduce the notion of the $p$-modulus of paths between two sets and show a sub-multiplicative inequality for them.
\definition\label{CMS.def10}
(1)\,\,Define
\[
E_{M, m}^* = \{(w, v)| w, v \in T_m, v \in \GG_M(w)\}.
\]
Note that $E^*_m = E^*_{1, m}$. Moreover, define 
\[
\theta_m(w, v) = \min\{M| v \in \GG_M(w)\}
\] 
for $w, v \in T_m$. $\theta_m(w, v)$ is the graph distance of the graph $(T_m, E^*_m)$.\\
(2)\,\,
Let $A \subseteq T_n$ and let $A_1, A_2 \subseteq A$. For $k \ge 0$, define
\begin{multline*}
\C_m^{(M)}(A_1, A_2, A) = \{(v(1), \ldots, v(l))| v(i) \in S^m(A)\,\,\text{for any $i = 1, \ldots, l$},\\
\text{there exist $v(0) \in S^m(A_1)$ and $v(l + 1) \in S^m(A_2)$}\\
\text{ such that $(v(i), v(i + 1)) \in E^*_{M, n + m}$ for any $i = 0, \ldots, l$.}\},
\end{multline*}
\begin{multline}\label{CMS.eq10}
\A_m^{(M)}(A_1, A_2, A) = \{f| f: T_{n + m} \to [0, \infty), \\
\sum_{i = 1}^l f(w(i)) \ge 1\,\,\text{for any $(w(1), \ldots, w(l)) \in \C_m^{(M)}(A_1, A_2, A)$}\}
\end{multline}
and
\begin{equation}\label{CMS.eq20}
\M_{p, m}^{(M)}(A_1, A_2, A) = \inf_{f \in \A_m^{(M)}(A_1, A_2, A)} \sum_{u \in T_{n + m}} f(u)^p.
\end{equation}

\noindent(3)\,\,
For $w \in T_n$, define
\[
\C_{N, m}^{(M)}(w) = \C_m^{(M)}(\{w\}, \GG_N(w)^c, T_{n}),\,\,\A_{N, m}^{(M)}(w) = \A_m^{(M)}(\{w\}, \GG_N(w)^c, T_{n})
\]
and
\[
\M_{N, p, m}^{(M)}(w) = \M_{p, m}^{(M)}(\{w\}, \GG_N(w)^c, T_{n}).
\]
\enddefinition 
The quantity $\M_{p, m}^{(M)}(A_1, A_2, A)$ is called the $p$-modulus of the family of paths between $A_1$ and $A_2$ inside $A$.

\remark
In \eqref{CMS.eq10} and \eqref{CMS.eq20}, the domain of $f$ is $T_{n + m}$. However, since we only use $f(u)$ for $u \in S^m(A)$ in \eqref{CMS.eq10} and the sum in \eqref{CMS.eq20} becomes smaller by setting $f(u) = 0$ for $u \in \sd{T_{n + m}}{S^m(A)}$, we may think of the domain of $f$ as $S^m(A)$.
\endremark

As is the case of conductances, if $A_1$ and $A_2$ consist of single points $u$ and $v$ respectively, then we  write $\C^{(M)}_m(u, v, A)$, $\A^{(M)}_m(u, v, A)$ and $\M^{(M)}_{p, m}(u, v, A)$ instead of $\C^{(M)}_m(\{u\}, \{v\}, A)$, $\A^{(M)}_m(\{u\}, \{v\}, A)$ and $\M^{(M)}_{p, m}(\{u\}, \{v\}, A)$ respectively.

By \cite[Proposition~4.8.4]{GAMS}, we have the following simple relation between $\E_{p, m}(A_1, A_2, A)$ and $\M^{(1)}_{p, m}(A_1, A_2, A)$. Consequently, to know $\M^{(1)}_{p, m}(A_1, A_2, A)$ is essential to know $\E_{p, m}(A_1, A_2, A)$. 

\lemma\label{CMS.lemma10}
Let $A \subseteq T_n$ and let $A_1, A_2 \subseteq A$ with $A_1 \cap A_2 = \emptyset$. Then for any $m \ge 1$ and $p > 0$,
\begin{equation}\label{CMP.eq30}
\frac 1{L_*}\E_{p, m}(A_1, A_2, A) \le \M^{(1)}_{p, m}(A_1, A_2, A) \le 2\max\{1, (L_*)^{p - 1}\}\E_{p, m}(A_1, A_2, A)
\end{equation}
\endlemma

The following theorem is the main result of this section. 

\thm[Sub-multiplicative inequality]\label{CMS.thm10}
Let $k_0, L, M \in \BbN$.  Suppose that $\pi^{k_0}(\GG_{L + 1}(u)) \subseteq \GG_{M}(\pi^{k_0}(u))$ for any $u \in T$. Then
\[
\M_{M, p, k + l}^{(1)}(w) \le c_{\ref{CMS.thm10}}\M_{M, p, k}^{(1)}(w)\max_{v \in S^k(\GG_{M}(w))}\M_{L, p, l}^{(1)}(v),
\]
for any $l \in \BbN$, $k \ge k_0$, $w \in T$ and $p > 0$, where $c_{\ref{CMS.thm10}}$ depends only $p, L_*$ and $L$.
\endthm

\remark
If $\pi^{k_0}(\GG_{L + 1}(u)) \subseteq \GG_{M}(\pi^{k_0}(u))$, then $\pi^{k}(\GG_{L + 1}(u)) \subseteq \GG_{M}(\pi^{k}(u))$ for any $k \ge k_0$.
\endremark

Similar sub-multiplicative inequalities for moduli of curve families have been shown in \cite[Proposition~3.6]{BouKleiner},
\cite[Lemma~3.8]{CarPiag} and \cite[Lemma~4.9.3]{GAMS}.\par
By Assumption~\ref{ASS.15}, the assumption $\pi^{k_0}(\GG_{L + 1}(u)) \subseteq \GG_{M}(\pi^{k_0}(u))$ is satisfied with $M = L = M_*$ and $k_0 = 1$. This fact along with Lemma~\ref{CMS.lemma10} shows the following sub-multiplicative inequality of conductance constants.

\cor\label{CMS.cor10}
For any $n, k, l  \ge 1$, $w \in T_n$ and $p \ge 1$.
\begin{equation}\label{CMP.eq40}
\E_{M_*, p, k + l}(w, T_n) \le c_{\ref{CMS.cor10}}\E_{M_*, p, k}(w, T_n)\max_{v \in S^k(\GG_{M}(w))}\E_{M_*, p, l}(v, T_{n + k}),
\end{equation}
where the constant $c_{\ref{CMS.cor10}} = c_{\ref{CMS.cor10}}(p, L_*, M_*)$ depends only on $p$, $L_*$ and $M_*$.
\endcor

The rest of this section is devoted to a proof of Theorem~\ref{CMS.thm10}.

\lemma\label{AA.lemma10}
Let $A \subseteq T_n$ and let $A_1, A_2 \subseteq A$ with $A_1 \cap A_2 = \emptyset$. Assume that $\GG_M(u) \cap S^m(A)$ is connected for any $u \in S^m(A)$. Then
\[
\M_{p, m}^{(1)}(A_1, A_2, A) \le M_{p, m}^{(M)}(A_1, A_2, A) \le   (L_*)^{(p + 1)M} \M_{p, m}^{(1)}(A_1, A_2, A).
\]
\endlemma

\demo
By definition, 
\[
\C_m^{(M)}(A_1, A_2, A) \supseteq \C_m^{(1)}(A_1, A_2, A)\,\,\text{and}\,\, \A_m^{(M)}(A_1, A_2, A) \subseteq \A_m^{(1)}(A_1, A_2, A).
\]
This shows
\[
\M_{p, m}^{(1)}(A_1, A_2, A) \le M_{p, m}^{(M)}(A_1, A_2, A).
\]
Define
\[
H_u = \GG_M(u)
\]
for any $u \in (T)_{n + m}$. Then
\[
\#(H_u) \le (L_*)^M \quad\text{and}\quad \#(\{v | u \in H_v\}) \le (L_*)^M.
\]
Let $(u(1), \ldots, u(l)) \in \C_m^{(M)}(A_1, A_2, A)$. Then there exist $u(0) \in S^m(A_1) \cap \GG_M(u(1))$ and $u(l + 1) \in S^m(A_2) \cap \GG_M(u(l))$. Since $u(0)$ and $u(1)$ is connected by a chain in $\GG_M(u(1))$ and $u(i)$ and $u(i + 1)$ is connected by a chain in $\GG_M(u(i))$ for $i = 1, \ldots, l$, we have a chain belonging to $\C^{(1)}_m(A_1, A_2, A)$  and contained in $\cup_{i = 1, \ldots, n} H_{u(i)}$. Thus Lemma~\ref{UFD.lemma50} shows
\[
 \M_{p, m}^{(M)}(A_1, A_2, A) \le  (L_*)^{(p + 1)M}\M_{p, m}^{(1)}(A_1, A_2, A).
 \]
\enddemo

\demo[Proof of Theorem~\ref{CMS.thm10}]
Let $f \in \A_{M, k}^{(L + 1)}(w)$ and let $g_v \in \A_{L, l}^{(1)}(v)$ for any $v \in (T)_{|w| + k}$. Define $h: (T)_{|w| + k + l} \to [0, \infty)$ by
\[
h(u) = \max\{f(v)g_v(u) | v \in \GG_L(\pi^l(u)) \cap S^k(\GG_{M}(w))\}\chi_{S^{k + l}(\GG_{M}(w))}(u).
\]
{\bf Claim 1.}\,\,$h \in \A_{M, k + l}^{(1)}(w)$.\\
Proof of Claim 1:\,\,
Let $(u(1), \ldots, u(m)) \in \C_{M, k + l}^{(1)}(w)$. There exist $u(0) \in S^{k + l}(w)$ and $u(m + 1) \in \sd{(T)_{|w| + k + l}}{S^{k + l}(\GG_{M}(w))}$ such that $u(0) \in \GG_{1}(u(1))$ and $u(m + 1) \in \GG_{1}(u(m))$. Set $v(i) = \pi^l(u(i))$ for $i = 0, \ldots, m + 1$. Let $v_*(0) = v(0)$ and let $i_0 = 0$. Define $n_*$, $v_*(n)$ and $i_n$ for $i = 1, \ldots, n_*$ inductively as follows:\\
If 
\[
\max\{j|i_n \le j \le m, v(j) \in \GG_{L}(v_*(n))\} = m,
\]
then $n = n_*$. If
\[
\max\{j|i_n \le j \le m, v(j) \in \GG_{L}(v_*(n))\} < m,
\]
then define
\[
i_{n + 1} = \max\{j|i_n \le j \le m, v(j) \in \GG_{L}(v_*(n))\} + 1\quad\text{and}\quad v_*(n + 1) = v(i_{n + 1}).
\]
The fact that $\pi^k(\GG_{L + 1}(v_*(0))) \subseteq \GG_{M}(\pi^k(v(0)))$ implies $n_* \ge 1$.
Since $v(i_{n + 1} - 1) \in \GG_{L}(v_*(n))$, we have $v_*(n + 1) \in \GG_{L + 1}(v_*(n))$. Hence $(v_*(1), \ldots, v_*(n_*)) \in \C_{M, k}^{(L + 1)}(w)$. Moreover, since $v_*(n - 1) \notin \GG_L(v_*(n))$ for $n = 1, \ldots, n_*$, there exist $j_n$ and $m_n$ such that  $i_{n - 1} < j_n \le m_n < i_n$ and $(u(j_n), \ldots, u(m_n)) \in \C_{L, l}^{(1)}(v_*(n))$. Since $g_{v_*(n)} \in \A_{L, l}^{(1)}(v_*(n))$, we have
\[
\sum_{i = j_n}^{m_n} h(u(i)) \ge \sum_{i = j_n}^{m_n} f(v_*(n))g_{v_*(n)}(u(i)) \ge f(v_*(n)).
\]
This and the fact that $(v_*(1), \ldots, v_*(n_*)) \in \C_{M, k}^{(L + 1)}(w)$ yield 
\[
\sum_{i = 1}^m h(u(i)) \ge \sum_{j = 1}^{n_*} f(v_*(j)) \ge 1.
\]
Thus Claim 1 has been verified. \qed\\
Set $C_0 = \max\{(L_*)^{L(p - 1)}, 1\}$. Then by Lemma~\ref{BAS.lemma10}, for $u \in S^{k + l}(\GG_{M}(w))$,
\begin{multline*}
h(u)^p \le \bigg(\sum_{v \in \GG_L(\pi^l(u)) \cap S^k(\GG_{M}(w))}f(v)g_v(u)\bigg)^p\\ \le C_0\sum_{v \in \GG_L(\pi^l(u)) \cap S^k(\GG_{M}(w))}f(v)^pg_v(u)^p.
\end{multline*}
The above inequality and Claim 1 yield
\[
\M_{M, p, k + l}^{(1)}(w) \le \sum_{u \in S^{k + l}(\GG_{M}(w))}h(u)^p
\le C_0\sum_{v \in S^k(\GG_{M}(w))}\sum_{u \in (T)_{|w| + k + l}}f(v)^pg_v(u)^p.
\]
Taking infimum over $g_v \in \A_{L, l}^{(1)}(v)$ and $f \in \A_{M, k}^{(L + 1)}(w)$, we have
\begin{align*}
\M_{M, p, k + l}^{(1)}(w) &\le C\sum_{v \in S^k(\GG_{M}(w))}f(v)^p\M_{L, p, l}^{(1)}(v)\\
 &\le C_0\sum_{v \in (T)_{|w| + k}}f(v)^p\max_{v \in S^k(\GG_{M}(w))}\M_{L, p, l}^{(1)}(v)\\
 &\le C_0\M_{M, p, k}^{(L + 1)}(w)\max_{v \in S^k(\GG_{M}(w))}\M_{L, p, l}^{(1)}(v).
 \end{align*}
 Finally, applying Lemma~\ref{AA.lemma10}, we have the desired inequality.
\enddemo

\setcounter{equation}{0}
\section{Neighbor disparity constant $\s_{p, m}(w, v)$}\label{REL}

Another important constant in this paper is $\s_{p, m}(u, v)$, which is called the neighbor disparity constant. Neighbor disparity constant controls the difference between means of a function on two neighboring cells via the $p$-energy of the function. For $p = 2$, $\s_{2, m}$ was introduced in \cite{KusZ1} for the case of self-similar sets. 

\notation
For $A \subseteq T_n$ and $f \in \ell(A)$, define
\[
(f)_A = \frac 1{\sum_{v \in A} \mu(K_w)}\sum_{v \in A} f(w)\mu(K_w).
\]
\endnotation

\definition\label{REL.def20}
For $p \ge 1$, $n \ge 1$, $m \ge 0$ and $(w, v) \in E_n^*$, define
\[
\s_{p, m}(w, v) = \sup_{f \in \ell(S^m(w) \cup S^m(v))} \frac{|(f)_{S^m(w)} - (f)_{S^m(v)}|^p}{\E_{p, S^m(w) \cup S^m(v)}^{n + m}(f)},
\]
which is called the $p$-neighbor disparity constant of $(w, v)$ at level $m$.
Moreover, define
\[
\s_{p, m, n} = \sup_{(w, v) \in E_n^*} \s_{p, m}(w, v) \quad\text{and}\quad \s_{p, m} = \sup_{n \ge 1} \s_{p, m, n}.
\]
\enddefinition

\remark
By Theorem~\ref{FIN.thm30} and Assumption~\ref{ASS.30},  $\s_{p, m, n}$ and $\s_{p, m}$ are finite. 
\endremark

One of the advantages of neighbor disparity constants is their compatibility with the integral projection $P_{n, m}$ from $\ell(T_{n + m})$ to $\ell(T_n)$ as follows.

\lemma[{\cite[(2.12) Lemma]{KusZ1}}]\label{REL.lemma20}
Let $p \ge 1$ and let $A \subseteq T_n$. Define $P_{n, m}: \ell(S^m(A)) \to \ell(A)$ by
\[
(P_{n, m}f)(w) = (f)_{S^m(w)}
\]
for any $f \in \ell(S^m(A))$ and $w \in A$. Then 
\[
\E_{p, A}^n(P_{n, m}f) \le L_*\max_{w, v \in A, (w, v) \in E_n^*}\s_{p, m}(w, v)\E_{p, S^m(A)}^{n + m}(f).
\]
In particular, if $A_1, A_2 \subseteq A$, then
\begin{equation}\label{REL.eq10}
\E_{p, 0}(A_1, A_2, A) \le L_*\max_{w, v \in A, (w, v) \in E_n^*}\s_{p, m}(w, v)\E_{p, m}(A_1, A_2, A)
\end{equation}
for any $m \ge 0$.
\endlemma

\demo
\begin{multline*}
\E_{p, A}^n(P_{n, m}(f)) = \frac 12\sum_{(w, v) \in E_n^*, w, v \in A} |(f)_{S^m(w)} - (f)_{S^m(v)}|^p\\
 \le \frac 12\sum_{(w, v) \in E_n^*, w, v \in A} \s_{p, m}(w, v)\E_{p, S^m(w) \cup S^m(v)}^{n + m}(f)\\
  \le L_*\max_{w, v \in A, (w, v) \in E_n^*}\s_{p, m}(w, v)\E_{p, S^m(A)}^{n + m}(f).
\end{multline*}
Choose $f$ such that $f|_{A_1} \equiv 1$, $f|_{A_0} \equiv 0$ and $\E_{p, m}(A_1, A_2, A) = \E_{p, S^m(A)}^{n + m}(f)$. Then
\[
\E_{p, 0}(A_1, A_2, A) \le \E_{p, A}^n(P_{n, m}f).
\]
So we have \eqref{REL.eq10}.
\enddemo

The first application of the above lemma is the following relation between conductance and neighbor disparity constants.
\lemma\label{REL.lemma25}
Let $p \ge 1$ and let $A \subseteq T_n$. For any $m, l \ge 0$ and $v \in A$,
\begin{equation}\label{REL.eq20}
\E_{M, p, m}(v, A) \le L_*\s_{p, l, n + m}\E_{M, p, m + l}(v, A).
\end{equation}
In particular, there exists $c_{\ref{REL.lemma25}}$, depending only on $M, p$ and  $L_*$, such that if $A \neq \GG_M^A(v)$, then

\begin{equation}\label{REL.eq30}
c_{\ref{REL.lemma25}}(\E_{M, p, l}(v, A))^{-1} \le \s_{p, l, n}
\end{equation}
for any $n \ge 1$ and $l \ge 0$.
\endlemma

\demo
\eqref{REL.eq20} is a special case of \eqref{REL.eq10}. To obtain \eqref{REL.eq30}, letting $m = 0$ in \eqref{REL.eq20}, we have
\[
\E_{M, p, 0}(v, A) \le L_*\s_{p, l, n}\E_{M, p, l}(v, A).
\]
By Theorem~\ref{FIN.thm10}, 
\begin{equation}
\ulc_{\E}(L_*, (L_*)^{M - 1}, p) \le \E_{M, p, 0}(v, A).
\end{equation}
This immediately implies \eqref{REL.eq30}. 
\enddemo
 
 Another important consequence of Lemma~\ref{REL.lemma20} is a sub-multiplicative inequality of neighbor disparity constants.
 
\lemma[{\cite[(2.13) Prop.-(3)]{KusZ1}}]\label{REL.lemma40}
Let $p \ge 1$.
\[
\s_{p, n + m, k} \le (L_*)^2\s_{p, n, k}\s_{p, m, k + n}
\]
for any $n, m, k \in \BbN$. 
\endlemma

\demo
Let $(w, v) \in E_k^*$. By Lemma~\ref{REL.lemma20}, for any $f \in \ell(T_{k + n + m})$,
\begin{multline*}
|(f)_{S^{n + m}(w)} - (f)_{S^{n + m}(v)}|^p  = \E_{p, \{w, v\}}^k(P_{k, n}(P_{k + n, m}f))\\
 \le L_*\s_{p, n, k}\E_{p, S^n(w) \cup S^n(v)}^{k + n}(P_{k + n, m}f)\\
  \le (L_*)^2\s_{p, n, k}\s_{p, m, k + n}\E_{p, S^{m + n}(w) \cup S^{n + m}(v)}^{n + m + k}(f).
\end{multline*}
This implies
\[
\s_{p, n + m}(w, v) \le (L_*)^2\s_{p, n, k}\s_{p, m, k + n}.
\]
\enddemo

In the rest of this section, we study an estimate of the difference $f(u) - f(v)$ for $f: T_n \to \BbR$ and $u, v \in T$ by means of the $p$-energy $\E_p^n(f)$ and neighbor disparity constants.

\lemma\label{REL.lemma80}
Let $w \in T$ and let $m \ge 1$. For any $f \in \ell(S^m(w))$ and $u \in S(w)$,
\[
|(f)_{S^m(w)} - (f)_{S^{m - 1}(u)}| \le N_*(\s_{p, m - 1, |w| + 1})^{\frac 1p}\E^{|w| + m}_{p,, S^m(w)}(f)^{\frac 1p}.
\]
\endlemma

\demo
For any $v \in S(w)$, there exist $v_0, v_1, \ldots, v_k \in S(w)$ such that $k \le N_*$, $v_0 = v$, $v_k = u$ and $(v_i, v_{i + 1}) \in E_{|w| + 1}^*$ for any $i = 0, \ldots, k - 1$. Hence
\begin{multline*}
|(f)_{S^{m - 1}(v)} - (f)_{S^{m - 1}(u)}| \le \sum_{i = 1}^{k - 1} |(f)_{S^{m - 1}(v_i)} - (f)_{S^{m - 1}(v_{i + 1})}| \\
\le N_*(\s_{p, m - 1, |w| + 1})^{\frac 1p}\E^{|w| + m}_{p, S^m(w)}(f)^{\frac 1p}.
\end{multline*}
Combining this with
\[
(f)_{S^m(w)} - (f)_{S^{m - 1}(u)} = \frac 1{\mu(w)}\sum_{v \in S(w)} \big((f)_{S^{m - 1}(v)} - (f)_{S^{m - 1}(u)}\big)\mu(v),
\]
we obtain the desired inequality.
\enddemo

\lemma\label{REL.lemma90}
Let $w \in T$ and let $n \ge m$. If $u, v \in S^n(w)$ and $\pi^{n - m}(u) \in \GG_k(\pi^{n - m}(v))$ for $k \ge 0$, then
\begin{multline}\label{REL.eq40}
|f(u) - f(v)| \le \\
\Big(2(N_*)^2\sum_{i = 1}^{n - m}(\s_{p, n - m - i, |w| + m + i})^{\frac 1p} + k(\s_{p, n - m, |w| + m})^{\frac 1p}\Big)\E_{p, S^n(w)}^{|w| + n}(f)^{\frac 1p}
\end{multline}
for any $f \in \ell(S^n(w))$.
\endlemma

\demo
Set $v(i) = \pi^{n - m - i}(v)$ for $i = 0, \ldots, n - m$. Then by Lemma~\ref{REL.lemma80},
\begin{multline}\label{REL.eq50}
|f(v) - (f)_{S^{n - m}(v(0))}| \le \sum_{i = 1}^{n - m} |(f)_{S^{n - m - i}(v(i))} - (f)_{S^{n - m - i + 1}(v(i - 1))}|\\
\le (N_*)^2\sum_{i = 1}^{n - m}(\s_{p, n - m - i, |w| + m + i})^{\frac 1p}\E_{p, S^n(w)}^{|w| + n}(f)^{\frac 1p}.
\end{multline}
The same inequality holds if we replace $v$ by $u$. Since $v(0) \in \GG_k(u(0))$, there exist $w(0), \ldots, w(l) \in \GG_k(u(0))$ such that $l \le k$, $w(0) = u(0), w(l) = v(0)$ and $(w(i), w(i - 1)) \in E_{|w| + m}^*$ for $i = 1, \ldots, l$. Then
\begin{multline}\label{REL.eq60}
|(f)_{S^{n - m}(u(0))} - (f)_{S^{n - m}(v(0))}| \le \sum_{i = 1}^l |(f)_{S^{n - m}(w(i))} - (f)_{S^{n - m}(w(i - 1))}|\\
\le k(\s_{p, n - m, |w| + m})^{\frac 1p}\E_{p, S^n(w)}^{|w| + n}(f)^{\frac 1p}.
\end{multline}
By \eqref{REL.eq50} and \eqref{REL.eq60}, we have \eqref{REL.eq40}.
\enddemo

\setcounter{equation}{0}
\section{Construction of $p$-energy: $p > \dim_{AR}(K, d)$}\label{CPE}

In this section, we are going to construct a $p$-energy on $K$ as a scaling limit of the discrete counterparts $\E_p^n$'s step by step under Assumption~\ref{ID}, which consists of the following two requirements \eqref{ID1} and \eqref{ID2}:\\
\eqref{ID1}\,Neighbor disparity constants and $(\text{conductance constants})^{-1}$ have the same asymptotic behavior,\\
\eqref{ID2}\,Conductance constants have exponential decay.\par
Under these assumptions, the $p$-energy $\hE_p$ is constructed in Theorem~\ref{CPE.thm10}. Furthermore, in the case $p = 2$, we construct a local regular Dirichlet form in Theorem~\ref{CPE.thm20}.\par
The  question when Assumption~\ref{ID} is fulfilled will be addressed in Section~\ref{SUC}.\par
As in the previous sections, we continue to suppose that Assumptions~\ref{ASS.10}, \ref{ASS.15}, \ref{ASS.20} and \ref{ASS.30} hold.  Moreover, throughout this section, we fix $p \ge 1$.

\definition\label{CPE.def00}
For $M \ge 1, m \ge 0$ and $n \ge 1$, define
\[
\E_{M, p, m, n} = \max_{v \in T_n}\E_{M, p, m}(v, T_n).
\]
\enddefinition

\remark
Theorem~\ref{FIN.thm10} shows that $\E_{M, p, m, n}$ is finite.
\endremark

\assumption\label{ID}
There exist $c_1, c_2 > 0$ and $\a \in (0, 1)$ such that
\begin{equation}\label{ID1}
c_1 \le \E_{M_*, p, m, n} \s_{p, m, n} \le c_2
\end{equation}
and
\begin{equation}\label{ID2}
\E_{M_*, p, m, n} \le c_2\a^m
\end{equation}
for any $m \ge 0, n \ge 1$.
\endassumption

By \cite[Theorems~4.7.6 and 4.9.1]{GAMS},  we have the following characterization of the condition \eqref{ID2} under Assumption~\ref{ALFR}.

\prop\label{CPE.prop10}
Under Assumption~\ref{ALFR}, \eqref{ID2} holds if and only if
\[
p > \dim_{AR}(K, d).
\]
\endprop

Note that since $K$ is assumed to be connected, we have $\dim_{AR}(K, d) \ge 1$,  so that $p > 1$.

In the following definition, we introduce the principal notion of this paper called conductive homogeneity. Due to Theorem~\ref{CPE.thm01}, conductive homogeneity yields \eqref{ID1}.

\definition[{\bf Conductive Homogeneity}]\label{CPE.def05}
Define
\[
\E_{M, p, m} = \sup_{w \in T, |w|  \ge 1}\E_{M, p, m}(w, T_{|w|}).
\]
A compact metric space $K$ (with a partition $\{K_w\}_{w \in T}$ and a measure $\mu$) is said to be $p$-conductively homogeneous if and only if  
\begin{equation}\label{ID3}
\sup_{m \ge 0} \s_{p, m}\E_{M_*, p, m} < \infty.
\end{equation}
\enddefinition

\remark
As in the case of $\E_{M, p, m, n}$, $\E_{M, p, m}$ is always finite due to Theorem~\ref{FIN.thm10}.
\endremark

\thm\label{CPE.thm01}
If $K$ is $p$-conductively homogeneous, then \eqref{ID1} holds.
\endthm
A proof of Theorem~\ref{CPE.thm01} will be provided in Section~\ref{SUC}. \par
Under conductive homogeneity, it will be shown in Theorem~\ref{SUC.thm10} that there exist $c_1, c_2 > 0$ and $\s > 0$ such that
\[
c_1\s^m \le \s_{p, m, n} \le c_2\s^m
\]
and 
\[
c_1\s^{-m} \le \E_{M_*, p, m}(v, T_n) \le c_2\s^{-m}
\] 
for any $m \ge 1, n \ge 0$ and $v \in T_n$. This is why we have given the name ``homogeneity'' to this notion.\par
Now we start to construct a $p$-energy under Assumption~\ref{ID}. An immediate consequence of Assumption~\ref{ID} is the following multiplicative property of $\s_{p, m, n}$.

\lemma\label{CPE.lemma00}
There exist $c_1, c_2 > 0$ such that
\[
c_1\s_{p, m, n + k}\s_{p, n, k} \le \s_{p, n + m, k} \le c_2\s_{p, m, n + k}\s_{p, n, k}
\]
for any $k \ge 1$, and $m, n \ge 0$.
\endlemma

\demo
By \eqref{CMP.eq40}, we have
\[
\E_{M_*, p, n + m, k} \le c\E_{M_*, p, m, n + k}\E_{M_*, p, n, k}.
\]
This along with \eqref{ID1} show
\[
c_1\s_{p, m, n + k}\s_{p, n, k} \le \s_{p, n + m, k}.
\]
The other half of the desired inequality follows from Lemma~\ref{REL.lemma40}.
\enddemo

Next, we study some geometry associated with the partition $\{K_w\}_{w \in T}$.

\definition\label{CPE.def10}
Let $L \ge 1$. Define
\begin{multline*}
n_L(x, y) = \max\{n| \text{there exist $w, v \in T_n$ such that}\\
\text{$x \in K_w, y \in K_v$ and $v \in \GG_L(w)$}\}.
\end{multline*}
 Furthermore, fix $r \in (0, 1)$ and define
 \begin{equation}\label{CPE.eq00}
\d_L(x, y) = r^{n_L(x, y)}.
\end{equation}
\enddefinition

Recall that $h_r: T \to (0, 1]$ is given as $h_r(w) = r^{|w|}$. Since $\LL_s^{h_r} = T_n$ if $r^{n - 1} > s \ge r^n$, where
\[
\LL_s^{h_r} = \{w| w \in T, h_r(\pi(w)) > s \ge  h_r(w)\},
\]
 $\d_L$ is nothing but $\d_L^{h_r}$ defined in \cite[Definition~2.3.8]{GAMS}.\par
By \cite[Proposition~2.3.7]{GAMS} and the discussions in its proof, we have the following fact.

\prop\label{CPE.prop20}
Suppose that $d$ is a metric on $K$ giving the original topology $\O$ of $K$. Let $L \ge 1$. There exists a monotonically non-decreasing function $\eta_L: [0, 1] \to [0, 1]$ satisfying $\lim_{t \downarrow} \eta_L(t) = 0$ and 
\[
\d_L(x, y) \le \eta_L(d(x, y))
\]
for any $x, y \in K$.
\endprop

\demo
Define
\[
\LL_{s, 0}^{h_r}(x) = \{v| v \in \LL^{h_r}_s, x \in K_v\},\,\, U^{h_r}_0(x, s) = \bigcup_{v \in \LL^{h_r}_{s, 0}(x)} K_v
\]
and 
\[
U^{h_r}_1(x, s) = \bigcup_{y \in U^{h_r}_0(x, s)} U^{h_r}_0(y, s).
\]
for $s \in (0, 1]$ and $x \in K$.
First we show that for any $\e > 0$, there exists $\c_{\e} > 0$ such that $\d_L(x, y) \le \e$ whenever $d(x, y) \le \c_{\e}$. If this is not the case, then there exist $\e_0 > 0$, $\{x_n\}_{n \ge 1}$ and $\{y_n\}_{n \ge 1}$ such that $d(x_n, y_n) \le \frac 1n$ and $\d_L(x_n, y_n) > \e_0$. Since $K$ is compact, choosing an adequate subsequence $\{n_k\}_{k \to \infty}$, we see that there exists $x \in K$ such that $x_{n_k} \to x$ and $y_{n_k} \to x$ for $k \to \infty$. By \cite[Proposition~2.3.7]{GAMS}, $U_0^{h_r}(x, \e_0/2)$ is a neighborhood of $x$. Hence both $x_{n_k}$ and $y_{n_k}$ belong to $U_0^{h_r}(x, \e_0/2)$ for sufficiently large $k$. So, there exist $w, v \in \LL_{\e_0/2, 0}^{h_r}(x)$ such that $x_{n_k} \in K_w$ and $y_{n_k} \in K_v$. Since $x \in K_w \cap K_v$, we see that $y \in U_1^{h_r}(x, \e_0/2)$, so that $\d_L(x_{n_k}, y_{n_k})  \le \e_0/2$. This contradicts the assumption that $\d_L(x_n, y_n) \ge \e_0$.  Thus our claim at the beginning of this proof is verified. Note that with a modification if necessary, we may assume that $\c_{\e}$ is monotonically non-decreasing as a function of $\e$ and $\lim_{\e \downarrow 0} \c_{\e} = 0$. Define
\[
\eta_L(t) = \inf\{\e| \e > 0, t \le \c_{\e}\}.
\]
Now it is routine to see that $\eta$ is the desired function.
\enddemo

Let $T_n = \{w(1), \ldots, w(l)\}$, where $l = \#(T_n)$. Inductively we define $\widetilde{K}_w$ by
\[
\widetilde{K}_{w(1)} = K_{w(1)}
\]
and
\[
\widetilde{K}_{w(k + 1)} = \sd{K_{w(k + 1)}}{\Big(\bigcup_{i = 1, \ldots, k} \widetilde{K}_{w(i)}\Big)}.
\]
Note that \eqref{ASS.eq10} implies that $\mu(B_w) = 0$ for any $w \in T_n$ and hence we have
\[
\widetilde{K}_w \supseteq O_w\quad\text{and}\quad\mu(\sd{K_w}{\widetilde{K}_w}) = 0
\]
for any $w \in T_n$. The latter equality is due to \eqref{ASS.eq10}. Now define $J_n: \ell(T_n) \to \BbR^K$ by
\begin{equation}\label{CPE.eq500}
J_nf = \sum_{w \in T_n} f(w)\chi_{\widetilde{K}_w}.
\end{equation}
Since $\widetilde{K}_w$ is a Borel set, $J_nf$ is $\mu$-measurable for any $f \in \ell(T_n)$. The definitions of $\widetilde{K}_w$ and $J_n$ depend on an enumeration of $T_n$ but $J_nf$ stays the same in the $\mu$-a.e. sense regardless of an enumeration.\par
Define
\begin{equation}\label{CPE.eq200}
\wE_p^m(f) = \s_{p, m - 1, 1}\E_p^m(f).
\end{equation}

The next lemma yields the control of the difference of values of $J_nf$ through $\wE_p^n(f)$.

\lemma\label{CPE.lemma10}
Suppose that Assumption~\ref{ID} holds. There exists $C > 0$ such that for any $n \ge 1$,$f \in \ell(T_n)$ and $x, y \in K$,
\begin{equation}\label{CPE.eq10}
|(J_nf)(x) - (J_nf)(y)| \le C\a^{\frac mp}\wE^n_p(f)^{\frac 1p},
\end{equation}
where $m = \min\{n_L(x, y), n\}$.
\endlemma

\demo
Let $m = \min\{n_L(x, y), n\}$. Then there exist $w, w' \in T_m$, $v \in S^{n - m}(w)$ and $u \in S^{n - m}(w')$ such that $x \in K_v$, $y \in K_u$, $(J_nf)(x) = f(v)$, $(J_nf)(y) = f(u)$ and $w' \in \GG_{L + 2}(w)$. By \eqref{REL.eq40},
\begin{equation}\label{CPE.eq15}
|f(u) - f(v)| \le  c\sum_{i = 0}^{n - m}(\s_{p, n - m - i, m + i})^{\frac 1p} \E_{p}^n(f)^{\frac 1p},
\end{equation}
where $c = \max\{2(N_*)^2, L_*(L + 2)\}$. Lemma~\ref{CPE.lemma00} shows that
\[
c_1\s_{p, m + i - 1, 1}\s_{p, n - m - i, m + i} \le \s_{p, n - 1, 1}.
\]
Combining this with Assumption~\ref{ID}, we obtain
\[
\s_{p, n - m - i, m + i} \le c_3\a^{m + i}\s_{p, n - 1, 1}.
\]
Using \eqref{CPE.eq15}, we see
\[
|f(u) - f(v)| \le c_4\a^{\frac mp}\wE_p^n(f)^{\frac 1p}.
\]
\enddemo

By this lemma, the boundedness of $\wE_p^n(f_n)$ gives a kind of equi-continuity to the family $\{f_n\}_{n \ge 1}$ and hence an analogue of Arzel{\`a}-Ascoli theorem, which we present in Appendix~\ref{APP1},  shows the existence of a uniform limit as follows.

\lemma\label{CPE.lemma20}
Suppose that Assumption~\ref{ID} holds.  Define $\tau = \frac{\log{\a}}{\log r}$. Let $f_n \in \ell(T_n)$ for any $n \ge 1$. If 
\[
\sup_{n \ge 1} \wE_p^n(f_n) < \infty\quad\text{and}\quad \sup_{n \ge 1}|(f_n)_{T_n}| < \infty
\]
then there exist a subsequence $\{n_k\}_{k \ge 1}$ and $f \in C(K)$ such that $\{J_{n_k}f_{n_k}\}$ converges uniformly to $f$ as $k \to \infty$, $\wE_p^{n_k}(f_{n_k})$ is convergent as $k \to \infty$ and
\begin{equation}\label{CPE.eq20}
|f(x) - f(y)|^p \le C\eta_L(d(x, y))^{\tau}\lim_{k \to \infty}\wE_p^{n_k}(f_{n_k}),
\end{equation}
where $\eta_L$ was introduced in Proposition~\ref{CPE.prop20}.
\endlemma

\demo
Set $C_* = \sup_{n \ge 1} \wE_p^n(f_n)$. By Lemma~\ref{CPE.lemma10}, if $n \ge n_L(x, y)$, then
\begin{equation}\label{CPE.eq100}
|J_nf_n(x) - J_nf_n(y)| \le C\a^{\frac{n_L(x, y)}p}(C_*)^{\frac 1p} \le C\eta_L(d(x, y))^{\frac{\tau}p}(C_*)^{\frac 1p}.
\end{equation}
In case $n < n_L(x, y)$, then there exist $w, w' \in T_n$ such that $x \in K_w$, $J_nf_n(x) = f(w)$, $y \in K_{w'}$, $J_nf_n(w') = f(w')$ and $w' \in \GG_{L + 2}(w)$. So there exists an $E_n^*$-path $(w(0), \ldots, w(L + 2))$ satisfying $w(0) = w$ and $w' = w(L + 2)$. By Lemma~\ref{BAS.lemma10}, 
\[
|f(w) - f(w')|^p \le (L + 2)^{p - 1}\sum_{i = 0}^{L + 1}|f(w(i)) - f(w(i + 1))|^p \le (L + 2)^{p - 1}\E_p^n(f_n).
\]
On the other hand, since $\wE_p^n(f_n) \le C_*$, Assumption~\ref{ID} implies
\[
\E_p^n(f_n) \le (\s_{p, n - 1, 1})^{-1}C_* \le c_2\E_{M_*, p, n - 1, 1}C_* \le (c_2)^2\a^{n - 1}C_*.
\]
Thus we have
\begin{equation}\label{CPE.eq110}
|J_nf_n(x) - J_nf_n(y)| \le c\a^{\frac np}(C_*)^{\frac 1p}.
\end{equation}
Making use of \eqref{CPE.eq100} and \eqref{CPE.eq110}, we see that
\[
|J_nf_n(x) - J_nf_n(y)| \le C\eta_L(d(x, y))^{\frac{\tau}p}(C_*)^{\frac 1p} + c\a^{\frac np}(C_*)^{\frac 1p}
\]
for any $x, y \in K$. Applying Lemma~\ref{GHT.lemma10} with $X = K, Y = \BbR, u_i = J_if_i$, we obtain the desired result.
\enddemo

\definition\label{CPE.def20}
Define $P_n: L^1(K, \mu) \to \ell(T_n)$ by
\[
(P_nf)(w) = \frac 1{\mu(w)}\int_{K_w} fd\mu
\] 
for any $n, m \ge 1$.  For $f \in \ell(T_k),$ we define
\[
P_nf = P_n{J_kf}.
\]
\enddefinition

The next lemma is one of the keys to the construction of a $p$-energy. A counterpart of this fact has already been used in Kusuoka-Zhou's construction of Dirichlet forms on self-similar sets in \cite{KusZ1}.
 
\lemma\label{CPE.lemma30}
Under Assumption~\ref{ID}, there exists $C > 0$ such that for any $n, m \ge 1$ and $f \in L^1(K, \mu) \cup (\cup_{k \ge 1}\ell(T_k))$,
\begin{equation}\label{CPE.eq30}
C\wE^n_p(P_nf) \le \wE^{n + m}_p(P_{n + m}f).
\end{equation}
In particular,
\begin{equation}\label{CPE.eq40}
C\sup_{n \ge 0} \wE^n_p(P_nf) \le \liminf_{n \to \infty} \wE^n_p(P_nf) \le \limsup_{n \to \infty} \wE^n_p(P_nf) \le \sup_{n \ge 0} \wE^n_p(P_nf)
\end{equation}
for any $f \in L^1(K, \mu)$.
\endlemma
\remark
This lemma holds without \eqref{ID2}.
\endremark

\demo
By Lemma~\ref{REL.lemma20}, 
\[
\E_p^n(P_nf) \le L_*\s_{p, m, n}\E_p^{n + m}(P_{n + m}f).
\]
Hence
\[
\frac 1{\s_{p, n - 1, 1}}\wE_p^n(P_nf) \le L_*\frac{\s_{p, m, n}}{\s_{p, n + m - 1, 1}}\wE_p^{n + m}(P_{n + m}f).
\]
By Lemma~\ref{CPE.lemma00}, we have \eqref{CPE.eq30}.
\enddemo

By virtue of the last lemma, we have a proper definition of the domain $\W^p$ of a $p$-energy given in Theorem~\ref{CPE.thm10}  and its semi-norm $\N_p$.

\lemma\label{CPE.lemma200}
Define
\[
\W^p = \{f| f \in L^p(K, \mu), \sup_{n \ge 1} \wE_p^n(P_nf) < +\infty\},
\]
and
\[
\N_p(f) = \sup_{n \ge 1} \wE_p^n(P_nf)^{\frac 1p}
\]
for $f \in \W^p$. Then $\W^p$ is a normed linear space with norm $\norm{\cdot}_{p, \mu} + \N_p(\cdot)$, where $\norm{\cdot}_{p, \mu}$ is the $L^p$-norm. Moreover, for any $f \in \W^p$, there exists $f_* \in C(K)$ such that $f(x) = f_*(x)$ for $\mu$-a.e. $x \in K$. In this way, $\W^p$ is regarded as a subset of $C(K)$ and 
\begin{equation}\label{CPE.eq60}
|f(x) - f(y)|^p \le C\eta_L(d(x, y))^{ \tau}\N_p(f)^p
\end{equation}
for any $f \in \W^p$ and $x, y \in K$, where $\eta_L$ was introduced in Proposition~\ref{CPE.prop20}.  In particular, $\N_p(f) = 0$ if and only if $f$ is constant on $K$.
\endlemma

If no confusion may occur, we write $\norm{\cdot}_p$ in place of $\norm{\cdot}_{p, \mu}$ hereafter.\par
In fact,  $(\W^p, \norm{\cdot}_{p} + \N_p(\cdot))$ turns out to be a Banach space in Lemma~\ref{CPE.lemma60}.

\demo
Note that 
\begin{equation}\label{CPE.eq50}
\wE_p^n(f + g)^{\frac 1p} \le \wE_p^n(f)^{\frac 1p} + \wE_p^n(g)^{\frac 1p}
\end{equation}
and so $\wE_p^n(\cdot)^{\frac 1p}$ is a semi-norm. This implies that $\N_p(\cdot)$ is a semi-norm of $\W^p$.\par
For $f \in \W^p$, by Lemma~\ref{CPE.lemma20}, there exist $\{n_k\}_{k \ge 1}$ and $f_* \in C(K)$ such that
\[
\norm{J_{n_k}P_{n_k}f - f_*}_{\infty} \to 0
\]
as $k \to \infty$ and
\[
|f_*(x) - f_*(y)|^p \le C\eta_L(d(x, y))^{\tau}\limsup_{n \to \infty} \E_p^n(P_nf).
\]
Since $\int_{K_w} P_{n_k}fd\mu \to \int_{K_w} f_*d\mu$ as $k \to \infty$,  it follows that $\int_{K_w}fd\mu = \int_{K_w} f_*d\mu$ for any $w \in T$. Hence $f = f_*$ for $\mu$-a.e. $x \in K$. Thus we identify $f_*$ with $f$ and so $f \in C(K)$.  Moreover, \eqref{CPE.eq60} holds for any $x, y \in K$. By \eqref{CPE.eq60}, $\N_p(f) = 0$ if and only if $f$ is constant on $K$.
\enddemo

We now examine the properties of the normed space $(\W^p, \norm{\cdot}_{p} + \N_p(\cdot))$.  The intermediate goals are to show its completeness (Lemma~\ref{CPE.lemma60}) and that it is dense in $C(K)$ with respect to the supremum norm (Lemma~\ref{CPE.lemma80}).

\lemma\label{CPE.lemma40}
Suppose that Assumption~\ref{ID} holds. The identity map $I: (\W^p, \norm{\cdot}_{p} + \N_p(\cdot)) \to (C(K), \norm{\cdot}_{\infty})$ is continuous.
\endlemma

\demo
Let $\{f_n\}_{n \ge 1}$ be a Cauchy sequence in $(\W^p, \norm{\cdot}_{p} + \N_p(\cdot))$. Fix $x_0 \in K$ and set $g_n(x) = f_n(x) - f_n(x_0)$. Then
\begin{multline*}
|g_n(x) - g_m(x)| = |(f_n(x) - f_m(x)) - (f_n(x_0) - f_m(x_0)| \\
\le C\eta_L(d(x, x_0))^{\frac {\tau}p}\N_p(f_n - f_m)
\end{multline*}
for any $x \in K$ and $n, m \ge 1$. Thus $\{g_n\}_{n \ge 1}$ is a Cauchy sequence in $C(K)$ with the norm $\norm{\cdot}_{\infty}$, so that  there exists $g \in C(K)$ such that $\norm{g - g_n}_{\infty} \to 0$ as $n \to \infty$.  On the other hand, since $\{f_n\}_{n \ge 1}$ is a Cauchy sequence of $L^p(X, \mu)$, there exists $f \in L^p(X, \mu)$ such that $\norm{f_n - f}_p \to 0$ as $n \to \infty$. Thus $f_n(x_0) = f_n - g_n$ converges as $n \to \infty$ in $L^p(K, \mu)$. Let $c$ be its limit. Then $f = g + c$ in $L^p(K, \mu)$. Therefore, $f \in C(K)$ and $\norm{f_n - f}_{\infty} \to 0$ as $n \to \infty$. 
\enddemo

Define $\overline{\W}^p$ as the completion of $(\W^p, \norm{\cdot}_{p} + \N_p(\cdot))$. Then the map $I$ is extended to a continuous map from $\overline{\W}_p \to C(K)$, which is denoted by $I$ as well for simplicity. 

\lemma[Closability]\label{CPE.lemma50}
Suppose that Assumption~\ref{ID} holds. The extended map $I: \overline{\W}^p \to C(K)$ is injective. In particular, $\overline{\W}^p$ is identified with a subspace of $C(K)$.
\endlemma

\demo
Let $\{f_n\}_{n \ge 1}$ be a Cauchy sequence in $(\W^p, \norm{\cdot}_{p} + \N_p(\cdot))$. Suppose $\lim_{n \to \infty}\norm{f_n}_{\infty} = 0$. Note that
\[
\wE_p^k(P_kf_n - P_kf_m) \le \sup_{l \ge 1}\wE_p^l(P_lf_n - P_lf_m) = \N_p(f_n - f_m)^p
\]
for any $k, n, m \ge 1$. Hence, for any $\e > 0$, there exists $N \in \BbN$ such that 
\[
\wE_p^k(P_kf_n - P_kf_m) \le \e
\]
for any $n, m \ge N$ and $k \ge 1$. As $\norm{f_m}_{\infty} \to 0$ as $m \to \infty$, we see that
\[
\wE_p^k(P_kf_n) \le \e
\]
for any $n \ge N$ and $k \ge 1$ and hence $\N_p(f_n)^p \le \e$ for any $n \ge N$. Therefore, $\N_p(f_n) \to 0$ as $n \to \infty$, so that $f_n \to 0$ in $\W^p$ as $n \to \infty$.
\enddemo

\lemma\label{CPE.lemma60}
Suppose that Assumption~\ref{ID} holds.
\[
\overline{\W}^p = \W^p.
\]
\endlemma

\demo
Let $\{f_n\}_{n \ge 1}$ be a Cauchy sequence of $\W^p$ and let $f$ be its limit in $\overline{\W^p}$. It follows that $\norm{f - f_n}_{\infty} \to 0$ as $n \to \infty$. Using the same argument as in the proof of Lemma~\ref{CPE.lemma50}, we see that for sufficiently large $n$,
\[
C\wE_p^k(P_kf_n - P_kf) \le \e
\]
for any $k \ge 1$. Since
\[
\wE_p^k(P_kf)^{\frac 1p} \le \wE_p^k(P_kf - P_kf_n)^{\frac 1p} + \wE_p^k(P_kf_n)^{\frac 1p},
\]
it follows that $\sup_{k \ge 1}\wE_p^k(P_kf) < \infty$ and hence $f \in \W^p$.
\enddemo

\lemma\label{CPE.lemma70}
Suppose that Assumption~\ref{ID} holds.\\
{\rm (1)}\,\,Let $\{n_k\}_{k \ge 1}$ be a monotonically increasing sequence of $\BbN$. Suppose that $f_{n_k} \in \ell(T_{n_k})$ for any $k \ge 1$, that $\sup_{k \ge 1}\wE_p^{n_k}(f_{n_k}) < \infty$ and that there exists $f \in C(K)$ such that $\norm{J_{n_k}f_{n_k} - f}_{\infty} \to 0$ as $n \to \infty$. Then $f \in \W^p$.\\
{\rm (2)}\,\,Let $f, g \in \W^p$. Then $f{\cdot}g \in \W^p$.
\endlemma

\demo
(1)\,\,Set $C_1 = \sup_{k \ge 1}\wE_p^{n_k}(f_{n_k})$.
By \eqref{CPE.eq30}, if $n \le n_l$, then
\[
C\wE_p^{n}(P_{n}f_{n_l}) \le \wE_p^{n_l}(f_{n_l}) \le C_1.
\]
Letting $l \to \infty$, we obtain
\[
C\wE_p^{n}(P_{n}f) \le C_1
\]
for any $k \ge 1$.  This implies $f \in \W^p$.\\
(2)\,\,For any $\vp, \psi \in \ell(T_n)$, 
\begin{multline*}
\E_p^n(\vp{\cdot}\psi) = \frac 12\sum_{(w, v) \in E_n^*} |\vp(w)\psi(w) - \vp(v)\psi(v)|^p\\
\le 2^{p - 1}\frac 12\sum_{(w, v) \in E_n^*} \big(|\vp(w)|^p|\psi(w) - \psi(v)|^p + |\vp(w) - \vp(v)|^p|\psi(v)|^p\big)\\
\le 2^{p - 1}\big(\norm{\vp}_{\infty}\E_p^n(\vp) + \norm{\psi}_{\infty}\E_p^n(\psi)\big).
\end{multline*}
Hence if $h_n = P_nf{\cdot}P_ng$, then
\[
\wE_p^n(h_n) \le 2^{p - 1}\big(\norm{f}_{\infty}\wE_p^n(P_nf) + \norm{g}_{\infty}\wE_p^n(P_ng)\big).
\]
Since $f, g \in \W_p$, we see that $\sup_{n \ge 1} \wE_p^n(h_n) < \infty$. Moreover, $\norm{J_nh_n - fg}_{\infty} \to 0$ as $n \to \infty$. Using (1), we conclude that $fg \in \W^p$.
\enddemo

\lemma\label{CPE.lemma75}
Suppose that Assumption~\ref{ID} holds. There exist a monotonically increasing sequence $\{m_j\}_{j \in \BbN}$ and $h_{M_*, w}^*, \vp_{M_*, w}^* \in \W^p$ for $w \in T$ such that\\
{\rm (a)}\,\,For any $w \in T$,
\[
\lim_{j \to \infty}\norm{J_{m_j}h^*_{M_*, w, m_j - |w|} - h_{M_*, w}^*}_{\infty}  = \lim_{j \to \infty} \norm{J_{m_j}\vp^*_{M_*, w, m_j - |w|} - \vp_{M_*, w}^*}_{\infty} = 0,
\]
where $h^*_{M_*, w, m}$ and $\vp^*_{M_*, w, m}$ are defined in Definition~\ref{MAC.def30}. For negative vales of $m$,  we formally define $h^*_{M_*. w, k - |w|} = P_{k}h^*_{M_*, w, 0}$  and $\vp^*_{M_*, w, k - |w|} = P_k\vp^*_{M_*, w, 0}$ for $k = 0, 1, \ldots, |w|$.\\
{\rm (b)}\,\,$\{\wE_p^{m_j}(h_{M_*, w, m_j - |w|}^*)\}_{j \ge 1}$ and $\{\wE_p^{m_j}(\vp_{M_*, w, m_j - |w|}^*)\}_{j \ge 1}$ converge as $j \to \infty$,\\
{\rm (c)}
Set $U_M(w) = \cup_{v \in \GG_M(w)} K_w$. For any $w \in T$, $h_{M_*, w}^*: K \to [0, 1]$ and
\[
h_{M_*, w}(x) = \begin{cases}
1 \quad&\text{if $x \in K_w$,}\\
0 \quad&\text{if $x \notin U_{M_*}(w)$.}
\end{cases}
\]
{\rm (d)}\,\,
For any $w \in T$, $\vp_{M_*, w}^*: K \to [0, 1]$, $\supp{\vp_{M_*, w}^*} \subseteq U_{M_*}(w)$, and
\[
\vp_{M_*, w}^*(x) \ge (L_*)^{-M_*}
\]
for any $x \in K_w$.  Moreover, for any $n \ge 1$, 
\[
\sum_{w \in T_n} \vp_{M_*, w}^* \equiv 1.
\]
{\rm (e)}\,\,For any $w \in T$ and $x \in K$,
\[
\vp_{M_*, w}^*(x) =\frac{ h_{M_*, w}^*(x)}{\sum_{v \in T_{|w|}} h_{M_*, v}^*(x)}.
\]

\endlemma

Note that $\{\vp_{M_*, w}^*\}_{w \in T_n}$ is a partition of unity subordinate to the covering $\{U_{M_*}(w)\}_{w \in T_n}$.

\demo
For ease of notation, write  $\vp_{w, m}^* = \vp_{M_*, w, m}^*$ and $h_{w, m}^* = h_{M_*, w, m}^*$. By Lemma~\ref{MAC.lemma15}, \eqref{ID1} and Lemma~\ref{CPE.lemma00}, we see that
\begin{multline*}
\wE_p^{|w| + m}(\vp_{w, m}^*) \le ((L_*)^{2M + 1} + 1)^p\s_{p, |w| + m - 1, 1}\E_{M, p, m}(w, T_{|w|})\\
\le C\s_{p, |w| + m - 1, 1}\s_{p, m, |w|}^{-1} \le C'\s_{p, |w| - 1, 1}
\end{multline*}
for any $w \in T$ and $m \ge 0$. Similarly, 
\[
\wE_p^{|w| + m}(h_{w, m}^*) \le C'\s_{p, |w| - 1, 1}.
\]
Hence Lemma~\ref{CPE.lemma20} shows that, for each $w$, there exists $\{n_k\}_{k \to \infty}$ such that $\{J_{|w| + n_k}h^*_{w, n_k}\}_{k \ge 1}$ (resp. $\{J_{||w| + m_j}\vp^*_{w, n_k}\}_{k \ge 1}$) converges uniformly as $k \to \infty$. Let
$h_w^*$ (resp. $\vp_w^*$) be its limit. Lemma~\ref{CPE.lemma70}-(1) implies that $h_w^* \in \W^p$ and $\vp_w^* \in \W^p$. By the diagonal argument, we choose $\{m_j\}_{j \ge 1}$ such that (a) and (b) hold. The statements (c), (d) and (e) are straightforward from the properties of $h_{w, m}^*$ and $\vp_{w, m}^*$.
\enddemo

\lemma\label{CPE.lemma80}
Suppose that Assumption~\ref{ID} holds. Then $\W^p$ is dense in $(C(K), \norm{\cdot}_{\infty})$.
\endlemma

\demo
Choose $x_w \in K_w$ for each $w \in T$. For $f \in C(K)$, define
\[
f_n = \sum_{w \in T_n} f(x_w)\vp_{M_*, w}^*.
\]
Then by Lemma~\ref{CPE.lemma75}, it follows that $\norm{f_n - f}_{\infty} \to 0$ as $n \to \infty$. Hence $\W^p$ is dense in $C(K)$.
\enddemo

\definition\label{CPE.def100}
For $f \in L^P(K, \mu)$, define $\overline{f}$ by
\[
\overline{f}(x) = \begin{cases}
1\quad&\text{if $f(x) \ge 1$},\\
f(x) \quad&\text{if $0 < f(x) < 1$,}\\
0\quad&\text{if $f(x) \le 0$}
\end{cases}
\]
for $x \in K$.
\enddefinition

Now we construct the $p$-energy $\hE_p$ as a $\Gamma$-cluster point of $\wE_p^n(P_n\,\cdot\,)$. The use of $\Gamma$-convergence in construction of Dirichlet forms on self-similar sets has been around for some time. See \cite{Griyang} and \cite{CaoQui} for example.

\thm\label{CPE.thm10}       
Suppose that Assumption~\ref{ID} holds. Then there exist $\hE_p:\W^p \to [0, \infty)$,  and $c  > 0$ such that \\
{\rm (a)}\,\,$(\hE_p)^{\frac 1p}$ is a semi-norm on $\W^p$ and 
\begin{equation}\label{CPE.eq70}
c\N_p(f) \le \hE_p(f)^{\frac 1p} \le \N_p(f)
\end{equation}
for any $f \in \W^p$.\\
{\rm(b)}\,\,For any $f \in \W^p$, $\overline{f} \in \W^p$ and
\[
\hE_p(\overline{f}) \le \hE_p(f).
\]
{\rm (c)}\,\,For any $f \in \W^p$,
\[
|f(x) - f(y)|^p \le c\eta_L(d(x, y))^{\tau}\hE_p(f).
\]
 In particular, for $p = 2$, $(\hE_2, \W^2)$ is a regular Dirichlet form on $L^2(K, \mu)$ and the associated non-negative self-adjoint operator has compact resolvent.
\endthm

The property (c) in the above theorem is called the Markov property.\par

\thm[Shimizu \cite{Shimizu2}]\label{CPE.thm15}
Suppose that Assumption~\ref{ID} holds.  Then the Banach space $(\W^p, \norm{\cdot}_p + \hE_p(\cdot))$ is reflexive and separable.
\endthm

\remark
In \cite{Shimizu2}, the reflexivity and separability are shown in the case of the planar Sierpinski carpet. His method, however, can easily be extended to our general case and one has the above theorem.
\endremark

\demo[Proof of Theorem~\ref{CPE.thm10}]
Define $\hE_p^n: L^p(K, \mu) \to [0, \infty)$ by $\hE_p^n(f) = \wE_p^n(P_nf)$ for $f \in L^p(K, \mu)$. Then by \cite[Proposition~2.14]{Braides1}, there exists a $\Gamma$-convergent subsequence $\{\hE_p^{n_k}\}_{k \ge 1}$. Define $\hE_p$ as its limit. Let $f \in \W^p$. Then 
\[
\hE_p(f) \le \liminf_{k \to \infty} \hE_p^{n_k}(f) \le \sup_{n \ge 1} \wE_p^n(P_nf) = \N_p(f)^p.
\]
Let $\{f_{n_k}\}_{k \ge 1}$ be  a recovering sequence for $f$, i.e. $\norm{f - f_{n_k}}_p \to 0$ as $k \to \infty$ and $\lim_{k \to \infty} \hE_p^{n_k}(f_{n_k}) = \E_p(f)$. By \eqref{CPE.lemma30}, if $n_k \ge n$, then
\[
C\wE_p^n(P_nf_{n_k}) \le \wE_p^{n_k}(P_{n_k}f_{n_k}) = \hE_p^{n_k}(f_{n_k}).
\]
Letting $k \to \infty$, we obtain
\[
C\wE_p^n(P_nf) \le \hE_p(f),
\]
so that
\[
C\N_p(f)^p \le \hE_p(f).
\]
The semi-norm property of $\hE_p(\cdot)^{\frac 1p}$ is straightforward from basic properties of $\Gamma$-convergence.\par
Next we show that $\hE_p(\overline{f}) \le \hE_p(f)$ for any $f \in \W^p$. Define 
\begin{equation}\label{CPE.eq600}
Q_nf = \sum_{w \in T_n} (P_nf)(w)\chi_{K_w}.
\end{equation}
Then
\begin{multline*}
\int_K|f(y) - Q_nf(y)|^p\mu(dy) \le \sum_{w \in T_n} \int_{K_w}\Big(\frac1{\mu(w)}\int_{K_w} |f(y) - f(x)|\mu(dx)\Big)^p\mu(dy)\\
\le \sum_{w \in T_n} \frac 1{\mu(w)}\int_{K_w \times K_w} |f(y) - f(x)|^p\mu(dx)\mu(dy).
\end{multline*}
This shows that if $f \in C(K)$, then $\norm{f - Q_nf}_p \to 0$ as $n \to \infty$. Let $\{f_{n_k}\}_{k \ge 1}$ be a recovering sequence for $f$. Since
\begin{multline*}
\norm{\overline{f} - \overline{Q_ng}}_p \le \norm{\overline{f} - \overline{Q_nf}}_p + \norm{\overline{Q_nf} - \overline{Q_ng}}_p \\
\le \norm{f - Q_nf}_p + \norm{Q_nf - Q_ng}_P \le \norm{f - Q_nf}_p + \norm{f - g}_p,
\end{multline*}
it follows that $\norm{\overline{f} - \overline{Q_{n_k}f_{n_k}}}_p \to 0$ as $n \to \infty$. Then
\begin{multline*}
\E_p(\overline{f}) \le \liminf_{k \to \infty} \hE_p^{n_k}(\overline{Q_{n_k}f_{n_k}}) = \liminf_{k \to \infty} \wE_p^{n_k}(\overline{P_{n_k}f_{n_k}})\\
\le \liminf_{k \to \infty} \wE_p^{n_k}(P_{n_k}f_{n_k}) = \lim_{k \to \infty} \hE_p^{n_k}(f_{n_k}) = \E_p(f).
\end{multline*}
Finally for $p = 2$, since a $\Gamma$-limit of quadratic forms is a quadratic form, we see that $(\hE_2, \W^2)$ is a regular Dirichlet form on $L^2(K, \mu)$. Since the identity map from $(\W^2, \norm{\cdot}_2 + \N_p(\cdot))$ to $(C(K), \norm{\cdot}_{\infty})$ is a compact operator, by \cite[Exercise 4.2]{Dav2}, the non-negative self-adjoint operator associated with $(\E_2, \W^p)$ has compact resolvent.
\enddemo

For the case $p = 2$, due to the above theorem,  $\W^2$ is separable.  Hence, we may replace $\Gamma$-convergence with point-wise convergence as seen in the following theorem. This enables us to obtain the local property of our Dirichlet form, which turns out to be a resistance form as well.

\thm\label{CPE.thm20}
Suppose that Assumption~\ref{ID} holds for $p = 2$. Then there exists a sub-sequence $\{m_k\}_{k \ge 1}$ such that $\{\E_2^{m_k}(P_{m_k}f, P_{m_k}g)\}_{k \ge 1}$ converges as $k \to \infty$ for any $f, g \in \W^2$. Furthermore, define $\E(f, g)$ as its limit. Then $(\E, \W^2)$ is a local regular Dirichlet form on $L^2(K, \mu)$, and there exist $c_1, c_2, c_3 > 0$ such that
\begin{equation}\label{CPE.eq80}
c_1\N_2(f) \le \E(f, f)^{\frac 12} \le c_2\N_2(f)
\end{equation}
and
\begin{equation}\label{CPE.eq90}
|f(x) - f(y)|^2 \le c_3\eta_L(d(x, y))^{\tau}\E(f, f)
\end{equation}
for any $f \in \W^2$ and $x, y \in K$. In particular, $(\E, \W^2)$ is a resistance form on $K$ and the associated resistance metric $R$ gives the original topology $\O$ of $K$.
\endthm

\demo
{\bf Existence of $\{m_k\}_{k \ge 1}$}:\,\,
By Lemma~\ref{CPE.thm10}, the non-negative self-adjoint operator $H$ associated with the regular Dirichlet form $(\hE_2, \W^2)$ has compact resolvent. Hence there exist a complete orthonormal basis $\{\vp_i\}_{i \ge 1}$ of $L^2(K, \mu)$ and $\{\lambda_i\}_{i \ge 1} \subseteq [0, \infty)$ such that
\[
H\vp_i = \lambda_i\vp_i\quad\text{and}\quad \lambda_i \le \lambda_{i + 1}
\]
for any $i \ge 1$ and $\lim_{i \to \infty} \lambda_i = \infty$. Note that $\{\frac{\vp_i}{\sqrt{1 + \lambda_i}}\}_{i \ge 1}$ is a complete orthonormal system of $ (\W^2, (\cdot, \cdot)_{2, \mu} + \hE_p(\cdot, \cdot))$. Hence setting
\[
\F = \{a_{i_1}\psi_{i_1} + \cdots + a_{i_m}\psi_{i_m}| m \ge 1, i_1, \ldots, i_m \ge 1, a_{i_1}, \ldots, a_{i_m} \in \BbQ\},
\]
we see that $\F$ is a dense subset of $\W^p$. For any $f, g \in \F$, since
\[
|\wE_2^n(P_nf, P_ng)| \le \wE_2^n(P_nf)^{\frac 12}\wE_2^n(P_ng)^{\frac 12} \le \N_2(f)\N_2(g),
\]
some sub-sequence of $\{\wE_2^n(P_nf, P_ng)\}_{n \ge 1}$ is convergent.  Since $\F \times \F$ is countable, the standard diagonal argument shows the existence of a sub-sequence $\{m_k\}_{k \ge 1}$ such that $\wE_2^{m_k}(P_{m_k}f, P_{m_k}g)$ converges as $k \to \infty$ for any $f, g \in \F$. Define $\E_2(f, g)$ as its limit. For $f, g \in \W^2$, choose $\{f_i\}_{i \ge 1} \subseteq \F$ and $\{g_i\}_{i \ge 1} \in \F$ such that $f_i \to f$ and $g_i \to g$ as $i  \to \infty$ in $\W^2$. Write $\wE_k(u, v) = \wE_2^{m_k}(P_{m_k}u, P_{m_k}v)$ for ease of notation. Then
\begin{multline*}
|\wE_k(f, g) - \wE_l(f, g)|  \le |\wE_k(f, g) - \wE_k(f_i, g)| + |\wE_k(f_i, g) - \wE_k(f_i, g_i)|\\
 + |\wE_k(f_i. g_i) - \wE_l(f_i, g_i)| + |\wE_l(f_i, g_i) - \wE_l(f_i, g)|  + |\wE_l(f_i, g) - \wE_l(f, g)| \\
 \le |\wE_k(f_i, g_i) - \wE_l(f_i, g_i)| + 2\N_2(f_i)\N_2(g - g_i)
 + 2\N_2(f - f_i)\N_2(g).
\end{multline*}
This shows that $\{\wE_k(f, g)\}_{k \ge 1}$ is convergent as $k \to \infty$. \eqref{CPE.eq80} is straightforward.\\
{\bf Strongly local property}:\,\,Let $f, g \in \W^p$. Assume that there exists an open set $U \subseteq K$ such that $\supp{f} \subseteq U$ and $g|_U$ is a constant. Then for sufficiently large $k$, $\wE_k(f, g) = 0$, so that $\E(f, g) = 0$.\\
{\bf Markov property}:\,\,By \eqref{CPE.eq70} and \eqref{CPE.eq80}, 
\[
0 \le \E(f, f) \le \hE_2(f, f)
\]
for any $f \in \W^2$. Since $(\hE_2, \W^2)$ is a regular Dirichlet form, by \cite[Theorem~2.4.2]{ChenFuku}, we see that $\E(f, g) = 0$ whenever $f, g \in \W^2$ and $f(x)g(x) = 0$ for $\mu$-a.e. $x \in K$. Now by the same argument as in the proof of \cite[Theorem~2.1]{BBKT}, we have the Markov property.\\
{\bf Resistance form}:\,\,Among the conditions for a resistance form in \cite[Definition~3.1]{Ki16}, (RF1), (RF2), (RF3) and (RF5) are immediate from what we have already shown. (RF4) is deduced from \eqref{CPE.eq90}. In fact, \eqref{CPE.eq90} yields that
\[
R(x, y) \le c\eta_L(d(x, y))^{\tau}
\]
for any $x, y \in K$. Assume that $R(x_n, x) \to 0$ as $n \to \infty$ and $\varlimsup_{n \to \infty}d(x, x_n) > 0$. Note that the collection of
\[
U_L^{h_r}(x, r^n) = \bigcup_{w \in T_n: x \in K_w}\Bigg(\bigcup_{v \in \GG_L(w)} K_v\Bigg)
\]
for $n \ge 1$ is a fundamental system of neighborhoods of $x$ by \cite[Proposition~2.3.9]{GAMS}. Therefore there exist $n \ge 1$ and $\{x_{m_k}\}_{k \ge 1}$ such that $x_{m_k} \notin U_L^{h_r}(x, r^n)$ for any $k \ge 1$. Choose $w \in T_n$ such that $x \in K_w$. Then $x_{m_k}$ belongs to $K_v$ for some $v \in \GG_L(w)^c$. So, $h^*_{L, w}(x) = 1$ and $h^*_{L, w}(x_{m_k}) = 0$. Hence
\[
R(x_{m_k}, x) \ge \frac 1{\E(h^*_{L, w})}
\]
for any $k \ge 1$. This contradicts the fact that $R(x, x_{m_k}) \to 0$ as $k \to \infty$. Thus we have shown $d(x_n, x) \to 0$ as $n \to \infty$. Hence the topology induced by the resistance metric $R$ is the same as the original topology $\O$.
\enddemo

\setcounter{equation}{0}
\section{Construction of $p$-energy: $p \le \dim_{AR}(K, d)$}\label{CPE2}

In this section, we will consider how much we can salvage the results in the previous section if $p \le \dim_{AR}(K, d)$. Honestly, what we will have in this section is far from satisfactory mainly because we have no proof of the conjecture saying that  $\W^P \cap C(K)$ is dense in $C(K)$ with respect to the supremum norm.  In spite of this, we present what we have now for future study. \par
Throughout this section, we assume \eqref{ID1}. Then, Lemma~\ref{CPE.lemma30} still holds. Replacing $(C(K), \norm{\cdot}_{\infty})$ by $(L^p(K, \mu), \norm{\cdot}_{p})$ in the statements and proofs of  Lemma~\ref{CPE.lemma50} and \ref{CPE.lemma60}, we have the following statement.\par

\lemma\label{CPE2.lemma10}
$\W^p$ is a Banach space with the norm $\norm{\cdot}_p + \N_p(\cdot)$.
\endlemma

\lemma\label{CPE2.lemma20}
Let $p > 1$. If $\{f_n\}_{n \ge 1}$ is a bounded sequence in the Banach space $\W^p$, then there exist $\{n_k\}_{k \ge 1}$ and $f \in \W^p$ such that $f$ is the weak limit of $\{f_{n_k}\}_{k \ge 1}$ in $L^p(K, \mu)$,
\[
\norm{f}_p \le \sup_{n \ge 1} \norm{f_n}_p\quad\text{and}\quad \N_p(f) \le \sup_{n \ge 1} \N_p(f_n).
\]
\endlemma

\demo
Since $L^p(K, \mu)$ is reflexive, $\{f_n\}$ has a weakly convergent sub-sequence $\{f_{n_k}\}_{k \ge 1}$. (See \cite[Section V.2]{Yosida}.) Let $f \in L^p(K, \mu)$ be its weak limit.   Since the map $f \to (P_mf)(w)$ is continuous, we see that $P_mf_{n_k} \to P_mf$ as $k \to \infty$ and hence
\[
\wE_p^{m}(P_mf) = \lim_{k \to \infty} \wE_p^m(P_mf_{n_k}) \le \sup_{k \ge 1} \N_p(f_{n_k})^{\frac 1p}.
\]

\enddemo

\lemma\label{CPE2.lemma30}
Let $p > 1$. Suppose that $f_n \in \ell(T_n)$ for any $n \ge 1$ and that
\[
\sup_{n \ge 1} \norm{J_nf_n}_p < \infty\quad \text{and}\quad \sup_{n \ge1} \wE_p^n(f_n) < \infty.
\] 
Then there exist a subsequence $\{n_k\}_{k \ge 1}$ and $f \in \W^p$ such that $f$ is the weak limit of $\{J_{n _k}f_{n_k}\}_{k \ge 1}$ in $L^p(K, \mu)$ and
\[
\norm{f}_p \le \sup_{n \ge 1} \norm{J_nf_n}_p\quad\text{and}\quad  C\N_p(f)^p \le \sup_{n \ge 1} \wE_p^n(f_n).
\]
\endlemma

\demo
Since $L^p(K, \mu)$ is reflexive, $\{J_nf_n\}$ has a weak convergent sub-sequence $\{J_{n_k}f_{n_k}\}_{k \ge 1}$. (See \cite[Section V.2]{Yosida}.)  Let  $f \in L^p(K, \mu)$ be its weak limit.  By Lemma~\ref{CPE.lemma30}, if $n_k \ge m$, then
\[
C\wE_p^m(P_mJ_{n_k}f_{n_k}) \le \wE_p^{n_k}(P_{n_k}J_{n_k}f_{n_k}) = \wE_p^{n_k}(f_{n_k}) \le \sup_{n \ge 1} \wE_p^n(f_n).
\]
Letting $k \to \infty$, we see
\[
C\wE_p^m(P_mf) \le \sup_{n \ge 1} \wE_p^n(f_n)
\]
for any $m \ge 1$. Thus $f \in \W^p$ and $C\N_p(f)^p \le \sup_{n \ge 1} \wE_p^n(f_n)$.
\enddemo

Using this lemma, we have a counterpart of Lemma~\ref{CPE.lemma75} as follows.

\lemma\label{CPE2.lemma40}
There exist $\{h_w^*\}_{w \in T}$ and $\{\vp_w^*\}_{w \in T} \subseteq \W^p$ such that\\
{\rm (a)}
Set $U_{M_*}(w) = \cup_{v \in \GG_{M_*}(w)} K_v$. For any $w \in T$, $h_w^*: K \to [0, 1]$ and
\[
h_{w}^*(x) = \begin{cases}
1\quad&\text{if $x \in K_w$,}\\
0 \quad&\text{if $x \notin U_{M_*}(w)$}.
\end{cases}
\]
{\rm (b)}\,\,
For any $w \in T$, $\vp_w^*: K \to [0, 1]$, $\supp{\vp_w^*} \subseteq U(w)$,  and
\[
\vp_w^*(x) \ge (L_*)^{-M_*}
\]
for any $x \in K_w$. Moreover, for any $n \ge 1$,
\[
\sum_{w \in T_n} \vp_w^* \equiv 1.
\]
{\rm (c)}\,\,For any $w \in T$ and $x \in K$,
\[
\vp_w^*(x) =\frac{ h_w^*(x)}{\sum_{v \in T_{|w|}} h_v^*(x)}.
\]
\endlemma

By the above lemma, 

\lemma\label{CPE2.lemma50}
$\W^p$ is dense in $L^p(K, \mu)$.
\endlemma

Finally, we have the following result on the construction of a $p$-energy.

\lemma\label{CPE2.lemma60}
There exist $\hE_p: \W^p \to [0, \infty)$ and $c_1, c_2 > 0$ such that $\hE_p^{\frac 1p}$ is a semi-norm, 
\[
c_1\N_p(f)^p \le \hE_p(f) \le c_2\N_p(f)^p\quad\text{and}\quad \hE_p(\overline{f}) \le \hE_p(f)
\]
for any $f \in \W^p$. In particular, for $p = 2$, $(\hE_2, \W^2)$ is a Dirichlet form on $L^2(K, \mu)$.

\endlemma

\setcounter{equation}{0}
\section{Conductive homogeneity}\label{SUC}
In this section, we study the notion of conductive homogeneity, namely, its consequence and how one can show it. \par
Throughout this section, we suppose that Assumptions~\ref{ASS.10}, \ref{ASS.15}, \ref{ASS.20} and \ref{ASS.30} hold.\par
The first theorem explains the reason why it is called ``homogeneity''.

\thm\label{SUC.thm10}
$K$ is $p$-conductively homogeneous if and only if there exist $c_1, c_2 > 0$ and $\s > 0$ such that
\begin{equation}\label{SUC.eq10}
c_1\s^{-m} \le \E_{M_*, p, m}(v, T_n) \le c_2\s^{-m},
\end{equation}
and
\[
c_1\s^m \le \s_{p, m, n} \le c_2\s^m
\]
for any $m \ge 0$, $n \ge 1$ and $v \in T_n$. 
\endthm

An immediate corollary of this theorem is Theorem~\ref{CPE.thm01}.

\cor[Theorem~\ref{CPE.thm01}]\label{SUC.cor10}
If $K$ is $p$-conductively homogeneous, then \eqref{ID1} holds.
\endcor

\demo[Proof of Theorem~\ref{SUC.thm10}]
Assume that $K$ is $p$-conductively homogeneous. Then by \eqref{REL.eq30}, there exists $c_1 > 0$ such that 
\[
c_1 \le \s_{p, m}\E_{M_*, p, m}.
\]
Also by Lemma~\ref{REL.lemma40}, there exists $c_2 > 0$ such that
\begin{equation}\label{SUC.eq20}
\s_{p, m + n} \le c_2\s_{p, m}\s_{p, n}
\end{equation}
for any $n, m \ge 0$. Moreover by \eqref{CMP.eq40}, there exists $c_3 > 0$ such that
\[
\E_{M_*, p, m + n} \le c_3\E_{M_*, p, m}\E_{M_*, p, n}
\]
for any $n, m \ge 0$. These inequalities along with \eqref{ID3} shows that there exist $c_4, c_5 > 0$ such that
\[
c_4\s_{p, m}\s_{p, n} \le \s_{p, m + n} \le c_5\s_{p, m}\s_{p, n}
\]
and
\[
c_4 \le \s_{p, m}\E_{M_*, p, m} \le c_5
\]
for any $m, n \ge 0$. From these, there exist $c_6,  c_7 > 0$ and $\s > 0$ such that
\[
c_6\s^m \le \s_{p, m} \le c_7\s^m\quad\text{and}\quad c_6\s^{m} \le (\E_{M_*, p, m})^{-1} \le c_7\s^{m}
\]
for any $m \ge 0$. Hence for any $w \in T$ and $n \ge 1$,
\[
c_6\s^m \le (\E_{p, m})^{-1} \le (\E_{M_*, p, m}(w, T_n))^{-1} \quad\text{and}\quad \s_{p, m, n} \le c_7\s^m.
\]
Making use of \eqref{REL.eq30}, we see that there exists $c_8 > 0$ such that
\[
c_6\s^m \le (\E_{M_*, p, m}(w, T_n))^{-1} \le c_8\s_{p, m, n} \le c_8c_7\s^m
\]
for any $m \ge 0$, $n \ge 1$ and $w \in T_n$.\par
The converse direction is straightforward.
\enddemo

Next we show another consequence of conductive homogeneity. For simplicity, we set $\E_{p, m}(u, v, S^k(w)) = \E_{p, m}(\{u\}, \{v\}, S^k(w))$. (In other words, we deliberately confuse $u$ with $\{u\}$.)

\lemma\label{SUC.lemma10}
 If $K$ is $p$-conductively homogeneous, then there exists $c_{\ref{SUC.lemma10}} > 0$, depending only on $p, L_*, N_*, M_*, k$,  such that
\[
\E_{M_*, p, m} \le c_{\ref{SUC.lemma10}}\E_{p, m}(u, v, S^k(w)) 
\]
for any $m \ge 0$, $w \in T$ and $u, v \in S^k(w)$ with $u \neq v$.
\endlemma

\demo
By \eqref{REL.eq10}, we see that
\[
\E_{p, 0}(u, v, S^k(w)) \le L_*\s_{p, m}\E_{p, m}(u, v, S^k(w)).
\]
Using Theorem~\ref{FIN.thm10}, it follows that
\[
\ulc_{\E}(L_*,  (N_*)^k , p) \le \E_{p, 0}(u, v, S^k(w)) \le L_*\s_{p, m}\E_{p, m}(u, v, S^k(w)).
\]
Now Theorem~\ref{SUC.thm10} suffices.
\enddemo

When $p > \dim_{AR}(K, d)$, the converse direction of the above lemma is actually true.

\thm\label{SUC.thm20}
Assume that there exist $c > 0$ and $\a \in (0, 1)$ such that
\begin{equation}\label{SUC.eq30}
\E_{M_*, p, m} \le c\a^m
\end{equation}
for any $m \ge 0$. Then $K$ is $p$-conductively homogeneous if and only if for any $k \ge 1$, there exists $c(k) > 0$ such that
\begin{equation}\label{SUC.eq40}
\E_{M_*, p, m} \le c(k)\E_{p, m}(u, v, S^k(w)) 
\end{equation}
for any $m \ge 0$, $w \in T$ and $u, v \in S^k(w)$ with $u \neq v$. In particular, under Assumption~\ref{ALFR}, if $p > \dim_{AR}(K, d)$, then whether $K$ is $p$-conductively homogeneous or not is independent of neighbor disparity constants.
\endthm

The last part of the theorem justifies the name ``{\it conductive}'' homogeneity.\par
The condition \eqref{SUC.eq30} is the same as \eqref{ID2}. Recall that, by Proposition~\ref{CPE.prop10}, \eqref{SUC.eq30} holds if and only if $p > \dim_{AR}(K, d)$ under Assumption~\ref{ALFR}.\par
The condition \eqref{SUC.eq40} is an analytic relative of the  ``Knight move'' condition described in probabilistic terminologies in \cite{KusZ1}. The name ``Knight move''  originated from the epoch-making paper \cite{BB1} where Barlow and Bass constructed the Brownian motion on the Sierpinski carpet. \par
The proof of the ``only if'' part of the above theorem is Lemma~\ref{SUC.lemma10}. A proof of the ``if'' part will be given in Sections~\ref{CPO}, \ref{REX} and \ref{PTR}.  \par
In Sections~\ref{CYP}, \ref{ESC} and \ref{SCR}, we are going to give examples for which one can show $p$-conductive homogeneity by Theorem~\ref{SUC.thm20}.\par
In the rest of this section, we study asymptotic behaviors of the heat kernel associated with the diffusion process induced by the Dirichlet form $(\E, \W^2)$ under Assumption~\ref{ALFR}. The next lemma shows that the associated resistance metric is bi-Lipschitz equivalent to a power of the original metric.

\lemma\label{SUC.lemma20}
Suppose that Assumption~\ref{ALFR} holds, $p > \dim_{AR}(K, d)$ and $K$ is $p$-conductively homogeneous. Let $\s$ be the same as in Theorem~\ref{SUC.thm10} and set $\tau_p = -\frac{\log{\s}}{\log r}$. Then there exist $c_1, c_2 > 0$ such that
\begin{equation}\label{SUC.eq45}
c_1d(x, y)^{\tau_p} \le \sup_{f \in \W^p, \hE_p(f) \neq 0} \frac{|f(x) - f(y)|^p}{\hE_p(f)} \le c_2d(x, y)^{\tau_p}
\end{equation}
for any $x, y \in K$. In particular, if $2 > \dim_{AR}(K, d)$, then
\begin{equation}\label{SUC.eq50}
c_1d(x, y)^{\tau_2} \le R(x, y) \le c_2d(x, y)^{\tau_2}
\end{equation}
for any $x, y \in K$, where $R(x, y)$ is the resistance metric associated with the resistance form $(\E, \W^2)$.
\endlemma

\demo
Since $\E_p^m(h^*_{M_*, w, m - |w|}) = \E_{M_*, p, m - |w|}(w, T_{|w|})$, we have
\[
c_1\s^{-m + |w|} \le \E_p^m(h^*_{M_*, w, m - |w|}) \le c_2\s^{-m + |w|}.
\]
by \eqref{SUC.eq10}. This shows 
\[
c_1\s^{|w|} \le \hE_p(h^*_{M_*, w}) \le c_2\s^{|w|}.
\]
Note that $d$ is $M_*$- adapted to $h_r$ by Assumption~\ref{ALFR}.  Hence by \cite[(2.4.1)]{GAMS}, 
\begin{equation}\label{SUC.eq60}
c_1d(x, y)  \le \d_{M_*}(x, y) \le c_2d(x, y)
\end{equation}
for any $x, y \in K$. Choose $n = n_{M_*}(x, y) + 1$. Let $w \in T_n$ satisfying $x \in K_w$. Since $n > n_{M_*}(x, y)$, it follows that if $v \in T_n$ and $y \in K_v$, then $v \notin \GG_{M_*}(w)$. Hence $h^*_{M_*, w}(x) = 1$ and $h^*_{M_*, w}(y) = 0$. Therefore \eqref{CPE.eq00} and \eqref{SUC.eq60} yield
\begin{multline*}
\sup_{f \in \W^p, \hE_p(f) \neq 0} \frac{|f(x) - f(y)|^p}{\hE_p(f)} \ge  \frac 1{\hE_p(h^*_{M_*, w})} \\
\ge  c(\s_p)^{-n} \ge c'r^{n_{M_*}(x, y)\tau_p} \ge c''d(x, y)^{\tau_p}.
\end{multline*}
On the other hand in this case, $\eta_{M_*}(t) = t$ by \eqref{SUC.eq60}. Hence Theorem~\ref{CPE.thm10}-(c) implies the other side of the desired inequality.
\enddemo

Due to the general theory of resistance forms in \cite{Ki16}, once we have \eqref{SUC.eq50}, it is straightforward to obtain asymptotic estimates of the heat kernel.

\thm\label{SUC.thm30}
Suppose that Assumption~\ref{ALFR} holds, $2 > \dim_{AR}(K, d)$ and $K$ is $2$-conductively homogeneous. Set $\tau_* = \tau_2$. Then there exists a jointly continuous hear kernel $p_{\mu}(t, x, y)$ on $(0, \infty) \times K \times K$ associated with the diffusion process induced by the local regular Dirichlet form $(\E, \W^2)$ on $L^2(K, \mu)$. Moreover\\
{\rm (1)}\,\,There exist $\b \ge 2$, a metric $\rho$, which is quasisymmetric to $d$,  and constants $c_1, c_2, c_3, c_4 > 0$ such that
\begin{equation}\label{SUC.eq70}
p_{\mu}(t, x, y) \le \frac{c_1}{\mu(B_{\rho}(x, t^{\frac 1{\b}}))}\exp\bigg(-c_2\Big(\frac{\rho(x, y)^{\b}}t\Big)^{\frac 1{\b - 1}}\bigg)
\end{equation}
for any $(t, x, y) \in (0, \infty) \times K \times K$ and
\begin{equation}\label{SUC.eq80}
\frac{c_3}{\mu(B_{\rho}(x, t^{\frac 1{\b}}))} \le p_{\mu}(t, x, y)
\end{equation}
for any $y \in B_{\rho}(x, c_4t^{\frac 1{\b}})$.\\
{\rm (2)}\,\,Suppose that $\mu$ is $\a_H$-Ahlfors regular with respect to the metric $d$. Set $\b_* = \tau_* + \a_H$. Then $\b_* \ge 2$ and there exist $c_7, c_8, c_9, c_{10} > 0$ such that
\begin{equation}\label{SUC.eq100}
p_{\mu}(t, x, y) \le c_6t^{-\frac {\a_H}{\b_*}}\exp\bigg(-c_7\Big(\frac{d(x, y)^{\b_*}}t\Big)^{\frac 1{\b_* - 1}}\bigg)
\end{equation}
for any $(t, x, y) \in (0, \infty) \times K \times K$ and
\begin{equation}\label{SUC.eq110}
c_9t^{-\frac {\a_H}{\b_*}} \le p_{\mu}(t, x, y)
\end{equation}
for any $y \in B_d(x, c_{10}t^{\frac {\a_H}{\b_*}})$. In addition, suppose that $d$ has the chain condition, i.e. for any $x, y \in K$ and $n \in \BbN$, there exist $x_0, \ldots, x_n \in K$ such that $x_0 = x, x_n = y$ and $d(x_i, x_{i + 1}) \le Cd(x, y)/n$, where the constant $C > 0$ is independent of $x, y$ and $n$. Then there exist $c_{11}$, $c_{12} > 0$ such that
\begin{equation}\label{SUC.eq120}
 c_{11}t^{-\frac {\a_H}{\b_*}}\exp\bigg(-c_{12}\Big(\frac{d(x, y)^{\b_*}}t\Big)^{\frac 1{\b_* - 1}}\bigg) \le p_{\mu}(t, x, y).
\end{equation}
\endthm

The exponent $\a_H$ above is in fact the Hausdorff dimension of $(K, d)$. The exponents $\b$ and $\b_*$ are called the walk dimensions.

\demo
We make use of \cite[Theorems~15.10 and 15.11]{Ki16}. Since $\mu$ has the volume doubling property with respect to $d$, \eqref{SUC.eq50} shows that $\mu$ has the volume doubling property with respect to $R$ as well. Since $K$ is connected, $(K, R)$ is uniformly perfect. Moreover, since $(\E, \W^2)$ has the local property, the annulus comparable condition (ACC) holds by \cite[Proposition~7.6]{Ki16}. Thus, the condition (C1) of \cite[Theorem~15.11]{Ki16} is verified and so is the condition (C3) of \cite[Theorem~15.11]{Ki16}. Using \cite[Theorem~15.11]{Ki16},  we have \eqref{SUC.eq70}. Consequently, by \cite[Theorem~15.10]{Ki16}, we see \eqref{SUC.eq80}.  Thus we have shown the first part of the statement. The fact that $\b \ge 2$, which is beyond the reach of  \cite[Theorem~15.10]{Ki16}, is due to \cite{Hino0}. See also \cite[Theorem~22.2]{Ki21}.\par
About the second part, assuming $\a_H$-Ahlfors regularity, i.e.\eqref{ASS.eq50}, we see that
\[
h_d(x, s) = s^{\tau_* + \a_H} = s^{\b_*},
\]
where $h_d(x, s)$ is defined as
\[
h_d(x, s) = \sup_{y \in B_d(x, s)}R(x, y)\cdot\mu(B_d(x, s)).
\]
Hence following the flow of exposition of \cite[Theorem~15.10]{Ki16}, we have
\[
g(s) = s^{\b_*}\quad \text{and}\quad \Phi(s) = s^{\b_* -1},
\]
where $g$ and $\Phi$ appear in the statement of \cite[Theorem~15.10]{Ki16}. Consequently, by \cite[Theorem~15.10]{Ki16}, we obtain \eqref{SUC.eq100}, \eqref{SUC.eq110} and \eqref{SUC.eq120}. The fact that $\b_* \ge 2$ can be shown in the same way as we did for $\b$ above.
\enddemo

\setcounter{equation}{0}
\section{Self-similar sets and self-similarity of energy}\label{SSF}

In this section, we consider the case where $K$ is a self-similar set with rationally related contraction ratios and construct self-similar energies under conductive homogeneity. Throughout this section, we fix a self-similar structure $\L = (K, S, \{f_s\}_{s \in S})$. The notion of the self-similar structure was introduced to give a purely topological description of self-similar sets. See \cite[Section~1.3]{AOF} for details.

\definition\label{SSF.def00}
Let $K$ be a compact metrizable space, let $S$ be a finite set, and let $\{f_s\}_{s \in S}$ be a family of continuous injective maps from $K$ to itself. \\
(1)\,\,The triple $(K, S, \{f_s\}_{s \in S})$ is called a self-similar structure if there exists a continuous surjective map $\chi: S^{\BbN} \to K$ such that 
\begin{equation}\label{SSF.eq200}
\chi(s_1s_2\ldots) = f_{s_1}(\chi(s_2s_3\ldots))
\end{equation}
for any $s_1s_2\ldots \in S^{\BbN}$, where $S^{\BbN}$ is equipped with the product topology.\\
(2)\,\,Define $W_* = \cup_{n \ge 0} S^n$, where $S^0 = \{\phi\}$. We use $\word wn$ to denote $(w_1, w_2, \ldots, w_n) \in S^n$. For $\word wn \in S^n$, set
\[
f_w = f_{w_1}{\circ}\ldots{\circ}f_{w_n}\quad\text{and}\quad K_w = f_w(K).
\]
In particular, $f_{\phi}$ is an identity map and $K_{\phi} = K$.
\enddefinition

By \cite[Proposition~3.3]{AOF}, if $(K, S, \{f_s\}_{s \in S})$ is a self-similar structure, $\chi: S^{\BbN} \to K$ is uniquely given by
\[
\{\chi(s_1s_2\ldots)\} = \bigcap_{m \ge 0} K_{\word sm}
\]
for any $s_1s_2\ldots \in S^{\BbN}$.\par
Typically, an example of self-similar structures is given by a self-similar set with respect to a family of  contractions. Let $(X, d)$ be a complete metric spaces and let $\{f_i\}_{i = 1, \ldots, N}$ be a family of contractions of $(X, d)$, i.e. $f_i: X \to X$ and 
\[
\sup_{x, y \in X, x \neq y} \frac{d(f_i(x), f_i(y))}{d(x, y)} < 1
\]
for any $i \in \{1, \ldots, N\}$. Then there exists a unique non-empty compact subset $K$ of $X$ satisfying
\begin{equation}\label{SSF.eq210}
K = \bigcup_{i = 1}^N f_i(K).
\end{equation}
See \cite[Theorem~1.1.4]{AOF} for example. The set $K$ is called a self-similar set with respect to $\{f_i\}_{i = 1, \ldots, N}$. By \cite[Theorem~1.2.3]{AOF}, if $S = \{1, \ldots, N\}$, then $(K, S, \{f_i\}_{i \in S})$ is a self-similar structure.\par
Let $r \in (0, 1)$ and let $j_s \in \BbN$ for $s \in S$. Define 
\begin{equation}\label{SSF.eq100}
j(w) = \sum_{i = 1}^m j_{w_i}\quad \text{and} \quad g(w) = r^{j(w)}
\end{equation}
for $w = \word wm \in S^m$. (In particular, $j(\phi) = 0$ and $g(\phi) = 1$.)  Define $\tilde{\pi}(\word wm) = \word w{m - 1}$ for $w = \word wm \in S^m$ and 
\begin{equation}\label{SSF.eq110}
\LL^g_{r^n} = \{w| w = \word wm \in W_*, g(\tilde{\pi}(w)) > r^n \ge g(w)\}.
\end{equation}
Note that $\LL^g_{r^n} \cap \LL^g_{r^{n + 1}}$ can be non-empty. (See Section~\ref{SCR} for example.)  So to distinguish $w \in \LL^g_{r^n}$ and $w \in \LL^g_{r^{n + 1}}$, we set
\[
T_n = \{(n, w)| w \in \LL_{r^n}\}
\]
and define $T = \cup_{n \ge 0} T_n$. There is a natural map $\iota: T \to W_*$ given by $\iota(n, w) = w$. Define
\[
\A = \{((n, v), (n + 1, w))| n \ge 0, \text{$v = w$ or  $v = \tilde{\pi}(w)$}\}
\]
Then $(T, \A, \phi)$ is a rooted tree and $\{K_w\}_{w \in T}$ is a partition of $K$ parametrized by $(T, \A, \phi)$. \par
Next, define $\a_H$ to be the unique number satisfying
\[
\sum_{s \in S} r^{j_s\a_H} = 1
\]
and let $\mu$ be the self-similar measure on $K$ with weight $\{r^{j_s\a_H}\}_{s \in S}$.\par
In the rest of this section, we presume the following assumption.
\assumption\label{SSF.ass10}
There exists a metric $d$ on $K$ giving the original topology of $K$ and Assumption~\ref{ALFR} holds with the metric $d$.  
\endassumption
Under this assumption, in particular, due to Assumption~\ref{ALFR}-(3), there exist $c_1, c_2 > 0$ such that
\[
c_1r^{j(w)} \le \diam{K_w, d} \le c_2r^{j(w)}
\]
for any $w \in T$. This enable us to regard the contraction ratio of $f_s$ as $r^{j_s}$. From this fact, we say that the contraction ratios of $\{f_s\}_{s \in S}$ are rationally related.\par

Under our assumptions, let  $\s$  be the same constant as in Theorem~\ref{SUC.thm10}. In this case, $\mu$ is $\a_H$-Ahlfors regular with respect to the metric $d$ and $\a_H$ coincides with the Hausdorff dimension of $(K, d)$. Note that even if we replace the definition \eqref{CPE.eq200} of $\wE_p^m(u)$ by
\begin{equation}\label{SSF.eq00}
\wE^m_p(u) = \s^m\E^m_p(u),
\end{equation}
all the arguments in Section~\ref{CPE} work and the results are unchanged. Our goal of this section is the next theorem.

\thm\label{SSF.thm10}
Suppose that $p > \dim_{AR}(K, d)$ and that $K$ is $p$-conductively homogeneous.\\
{\rm (1)}\,\,For any $w \in W_*$ and $f \in \W^p$, 
\[
f{\circ}f_w \in \W^p.
\]
{\rm (2)}\,\,There exists $\E_p: \W^p \to [0, \infty)$ satisfying\\
{\rm (a)}\,\,$(\E_p)^{\frac 1p}$ is a semi-norm on $\W^p$ and there exist $c_1, c_2 > 0$ such that
\begin{equation}\label{SSF.eq10}
c_1\N_p(f) \le \E_p(f)^{\frac 1p} \le c_2\N_p(f)
\end{equation}
and
\[
c_1d(x, y)^{\tau_p} \le \sup_{f \in \W^2, \E_p(f) \neq 0}\frac{|f(x) - f(y)|^p}{\E_p(f)} \le c_2d(x, y)^{\tau_p}
\]
for any $f \in \W^p$ and $x, y \in K$.\\
{\rm(b)}\,\,For any $f \in \W^p$, $\overline{f} \in \W^p$ and
\[
\E_p(\overline{f}) \le \E_p(f).
\]
{\rm (c)}
For any $f \in \W^p$,
\[
\E_p(f) = \sum_{s \in S} \s^{j_s}\E_p(f{\circ}f_s).
\]
In particular, for $p = 2$, $(\E_2, \W^2)$ is a local regular Dirichlet form on $L^2(K, \mu)$.
\endthm

\demo

Define
\begin{multline*}
\U = \{A(\cdot)| \text{$A(\cdot)$ is a semi-norm on $\W^p$, there exist $c_1, c_2 > 0$ such that}\\
\text{ $c_1\N_p(f) \le A(f) \le c_2\N_p(f)$ for any $f \in \W^p$}\}.
\end{multline*}
For $A_1, A_2 \in \U$, we write $A_1 \le A_2$ if and only if $A_1(f) \le A_2(f)$ for any $f \in \W^p$. We give $\U$ the point-wise convergence topology, i.e. $\{A_n\}_{n \ge 1} \subseteq \U$ is convergent to $A \in \U$ as $n \to \infty$ if and only if $A_n(f) \to A(f)$ as $n \to \infty$ for any $f \in \W^p$. Then due to the separability of $\W^p$ described in Theorem~\ref{CPE.thm15}, $\U$ is an ordered topological cone in the sense of \cite{Ki10}. \par
Let $w \in W_*$. For any $v = \word vk \in \LL_{r^{n - j(w)}}$, since
\[
g(w\word v{k - 1}) = g(w)g(\word v{k - 1}) > g(w)r^{n - j(w)} = r^n \ge g(wv),
\]
it follows that $wv \in \LL_{r^n}$. This shows that $\{(n, wv)| v \in  \LL_{r^{n - j(w)}}\} \subseteq T_n$.  In fact, $T_n = \cup_{w \in S^m} \{(n, wv)| v \in  \LL_{r^{n - j(w)}}\}$, which is a disjoint union. This yields
\[
\sum_{w \in S^m} \E^{n - j(w)}_p(P_{n - j(w)}(f{\circ}f_w)) \le \E^n_p(P_nf)
\]
for any $f \in L^p(K, \mu)$. Therefore,
\[
\sum_{w \in S^m} \s^{j(w)}\hE^{n - j(w)}(f{\circ}f_w) \le \hE^n_p(f).
\]
This inequality implies that $\s^{j(w)}\sup_{n \ge j(w)}\hE^{n - j(w)}(f{\circ}f_w) \le \N_p(f)^{p} < \infty$ for any $f \in \W^p$, so that $f{\circ}f_w \in \W^p$. Thus we have verified the statement (1). Again by the above inequality,
\begin{multline}\label{SSF.eq20}
c\sum_{w \in S^m} \s^{j(w)}\N_p(f{\circ}f_w)^{p} \le \sum_{w \in S^m}\s^{j(w)}\varliminf_{n \to \infty}\hE^{n - j(w)}(f{\circ}f_w) \\\le \sup_{n \ge 0}\hE_p^n(f) = \N_p(f)^{p}.
\end{multline}
Note that
\[
\sum_{(n, v) \in T_n} \s^{j(v)}\hE_p^{k - j(v)}(f{\circ}f_v) \le \sum_{w \in S^m} \s^{j(w)}\hE^{n + k  - j(w)}(f{\circ}f_w).
\]
By \eqref{CPE.eq40}, taking $\varliminf$ in the left-hand side and $\sup$ in the right-hand side, we see that
\begin{equation}\label{SSF.eq30}
c\sum_{(n, v) \in T_n} \s^{j(v)}\N_p(f{\circ}f_v)^{p} \le \sum_{w \in S^m} \s^{j(w)}\N_p(f{\circ}f_w)^{p}.
\end{equation}
On the other hand, for any $(n, v) \in T_n$ and $x \in K_v$, the self-similarity of $\mu$ and \eqref{SUC.eq45} show
\begin{multline*}
|(P_nf)(v) - f(x)| \le \int_K |f{\circ}f_v(y) - f{\circ}f_v(x_0)|\mu(dy) \\
\le c\int_K d(x_0, y)^{\frac {\tau_*}p}\mu(dy)\N_p(f{\circ}f_v) \le c'\N_p(f{\circ}f_v),
\end{multline*}
where $x_0 = (f_v)^{-1}(x)$. Hence if $((n, v), (n, u)) \in E_n^*$, then 
\[
|(P_nf)(v) - (P_nf)(u)| \le c'(\N_p(f{\circ}f_v) + \N_p(f{\circ}f_w)).
\]
This along with \eqref{SSF.eq30} yields
\begin{multline*}
\hE^n_p(f) = \frac {\s^n}2\sum_{((n, v), (n, u)) \in E_n^*} |(P_nf)(v) - (P_nf)(u)|^p \\
\le C\sum_{(n, v) \in T_n}\s^{j(w)}\N_p(f{\circ}f_v)^p \le C'\sum_{w \in S^m} \s^{j(w)}\N_p(f{\circ}f_w)^{p}.
\end{multline*}
Taking $\sup$ in the right-hand side, we have
\begin{equation}\label{SSF.eq40}
\N_p(f)^p \le C'\sum_{w \in S^m} \s^{j(w)}\N_p(f{\circ}f_w)^{p}.
\end{equation}
Now for $A \in \U$, define $\F(A)$ by
\[
\F(A)(f) = \Big(\sum_{s \in S} \s^{j_s}A(f{\circ}f_s)^p\Big)^{\frac 1p}.
\]
For any $A \in \U$, since $A \le c_2\N_p$, \eqref{SSF.eq20} implies
\[
\F(A) \le c_2\F(\N_p) \le c'\N_p.
\]
On the other hand, the fact $c_1\N_p \le A$ and \eqref{SSF.eq40} yield
\[
\F(A) \ge c_1\F(\N_p) \ge c''\N_p.
\]
Thus $\F(A) \in \U$ and $\F: \U \to \U$. It is easy to see that $\U$ is continuous and $\F(A + B) \le \F(A) + \F(B)$. Combining \eqref{SSF.eq20} and \eqref{SSF.eq40}, we see that there exist $C_1, C_2 > 0$ such that
\[
c_1\N_p  \le \F^j(\N_p) \le c_2\N_p
\]
for any $j \ge 1$. So, by \cite[Theorem~1.5]{Ki10}, there exists $\E_* \in \U$ such that $\F(\E_*) = \E_*$. Define
\[
\U_M = \{A | A \in \U, A(\overline{f}) \le A(f)\,\,\text{for any $f \in \W^p$}\}.
\]
Then $\hE_p \in \U_M$ and $\U_M$ is a closed subset of $\U$.  Hence by \cite[Corollary~1.6]{Ki10}, we see there exists $\E' \in \U_M$ such that $\F(\E') = \E'$.  Letting $\E = (\E')^p$, we have the desired $\E$. In the case $p = 2$, define
\begin{multline*}
\U_{DF} = \{A |  A \in \U, \text{$A$ satisfies the parallelogram law,}\\
\text{ the resulting quadratic form has both Markov and local property}\}.
\end{multline*}
Then $\U_{DF}$ is a closed subspace of $\U$ and Theorem~\ref{CPE.thm20} ensures that $\U_{DF} \neq \emptyset$. So again by \cite[Corollary~1.6]{Ki10}, we have the desired local regular Dirichlet form.
\enddemo

\section{Conductive homogeneity of self-similar sets}\label{CHS}
\def\IT{\mathcal{I}\mathcal{T}}
\def\I{\mathcal{I}}

In this section, we present a sufficient condition for conductive homogeneity of self-similar sets. The idea is originated from \cite{BouKleiner}, where the authors used symmetries of the spaces to show the combinatorial Loewner property of the Sierpinski carpet and the Menger curve, also known as the Menger sponge. Our sufficient condition, Theorem~\ref{CHS.thm10}, will be used in Sections~\ref{CYP} and \ref{OEX}.\par
In this section,  we assume that  $(K, S, \{f_s\}_{s \in S})$ is a self-similar structure  and adopt the setting in Section~\ref{SSF}. For simplicity, we also assume that $j_s = 1$ for any $s \in S$, so that $g(w) = r^{|w|}$ and $T_m = S^m$.
\definition\label{CHS.def10}
(1)\,\,
For any $e = (w, v) \in E_m^*$, define
\[
X(e) = (f_w)^{-1}(f_w(K) \cap f_v(K))
\]
and $\vp_e: X(e) \to X(e^r)$ by $(f_v)^{-1}{\circ}f_w|_{X(e)}$, where $e^r = (v, w)$ for $e = (w, v)$.
Furthermore, define
\[
\IT(K, T) = \{(X(e), X(e^r), \vp_e) | m \ge 1, e \in E_m^*\}.
\]
An element of $\IT(K, T)$ is called an intersection type of $(K, T)$. \\
(2)\,\,
A homeomorphism $g: K \to K$ is said to be a symmetry of $(K, T)$ if there exists $g^*: T \to T$ such that $|g^*(w)| = |w|$ and $g(K_w) = K_{g^*(w)}$ for any $w \in T$. Define $\G_{(K, T)}$ as the collection of symmetries of $(K, T)$. \\
(3)\,\,For any $n \ge 0$, define $\psi_n: \cup_{m \ge 0} T_{n + m} \to T$ by $\psi_n(v) = u$ if $v \in T_{n + m}$ and $v = \pi^m(v)u$. 
\enddefinition

\remark
The notion of intersection types and the set $\IT(K, T)$ were introduced in \cite{Ki13}.
\endremark

Note that $\psi_n(T_{n + m}) = T_m$ and $(f_{\pi^m(v)})^{-1}(K_v) = K_{\psi_n(v)}$ for any $v \in T_{n + m}$.

\notation
For $A \subseteq T$, set
 \begin{equation}\label{SCR.eq10}
 K(A) = \bigcup_{v \in A} K_v. 
 \end{equation}
\endnotation

\thm\label{CHS.thm10}
Suppose that there exist a finite subset $\I \subseteq \IT(K, T)$ and finite subgroups  $\G_0$ and $\G_1$ of $\G_{(K, T)}$satisfying the following properties {\rm (a)}, {\rm (b)} and {\rm (c)}{\rm:}\\
{\rm (a)}\,\,$(T_m, E_m^{\I})$ is connected for any $m \ge 1$, where 
\[
E_m^{\I} = \{e| e \in E_m^*, (X(e), X(e^r), \vp_e)  \in \I\}.
\]
{\rm (b)}\,\,For any $(X, Y, \vp) \in \I$ and $x \in X$, there exists $g \in \G_0$ such that $g(x) = \vp(x)$.\\
{\rm (c)}\,\,For any $n \ge 1$, $w \in T_n$ and $\bp \in \C^{(1)}_{M, m}(w)$, there exists $\U_{\bp} \subseteq \cup_{g \in \G_1}g^*(\psi_n(\bp))$ such that $K(\U_p)$ is connected and $g(K(\U_p)) \cap X \neq \emptyset$ for any $(X, Y, \vp) \in \I$ and $g \in \G_0$.\\
Then for any $p \ge 1$, $n, k \ge 1$, $m \ge 1$, $u_*, v_* \in T_k$, and $w \in T_n$,
\begin{equation}\label{CHS.eq100}
\M^{(1)}_{M, p, m}(w) \le (L_*)^M\#(\G_1)^{p + 1}\#(T_k)^p\M_{p, m}^{(1)}(u_*, v_*, T_k).
\end{equation}
Furthermore, if Assumption~\ref{SSF.ass10} holds with $M_*  = M$, then $K$ is $p$-conductively homogeneous for any $p > \dim_{AR}(K, d)$.
\endthm
\remark
Strictly,  a path $\bp = (w(1), \ldots, w(k))$ of a graph is not a subset of vertices but a sequence of them. However,  we use $\bp$ to denote a subset $\{w(1), \ldots, w(k)\}$ if no confusion may occur. For example, in the expression $\psi_n(\bp)$ above, we regard $\bp$ as a subset of $T_{n + m}$. 
\endremark
\demo
For $u \in S^m(\GG_1(w))$, define $H_u \subseteq T_{k + m}$ by
\[
H_u = \{vg^*(\psi_n(u))| g \in \G_1, v \in T_k\}.
\]
Then $\#(H_u) \le \#(T_k)\#(\G_1)$ for any $u \in S^m(\GG_1(w))$ and $\#(\{u| v \in H_u\}) \le \#(\GG_M(w))\#(\G_1)$ for any $v \in T_{k + m}$.\par
Now, since $(T_k, E_k^{\I})$ is connected, there exists $(w(0), w(1), \ldots, w(l), w(l + 1)) \in (T_k)^{l + 2}$ such that $w(0) = u_*, w(l + 1) = v_*$, $(w(i), w(i + 1)) \in E_k^{\I}$ for any $i = 0, 1, \ldots, l$. Set $e_i = (w(i), w(i + 1))$. Then $(X(e_i), X(e_i)^r, \vp_{e_i}) \in \I$. \\
{\bf Claim}: There exist $\A_i \subseteq T_m, x_i \in K$  and $g_i, h_i \in \G_0$ for $i = 1, 2, \ldots, l $ such that\\
(i)\,\,$\A_i = (h_i)^*(\U_{\bp})$ and $K(\A_i) \cap X(e_i) \neq \emptyset$,\\
(ii)\,\,$x_i \in K(\A_i) \cap X(e_i)$ and $g_i(x_i) = \vp_{e_i}(x_i)$,\\
and\\
(iii)\,\,$\A_{i + 1} = (g_i)^*(\A_i)$.\\
Proof of Claim:\,\,For $i = 1$, let $h_1$ be the identity map. Then $\A_1 = \U_{\bp}$. Since $K(\A_1) \cap X(e_1) \neq \emptyset$ by (c), we may choose $x_1 \in K(\A_1) \cap X(e_1)$. By (b), there exists $g_1 \in \G_0$ such that $g_1(x_1) = \vp_{e_1}(x_1)$.\par
Assume that we have the desired objects for $i  \in \{1, \ldots, l - 1\}$. Letting $h_{i + 1} = g_i{\circ}h_i \in \G_0$ and $\A_{i + 1} = (g_i)^*(\A_i)$, we obtain \[
\A_{i + 1} = (g_i)^*(h_i)^*(\U_{\bp}) = (h_{i + 1})^*(\U_{\bp}).
\]
 Using (c), we see that $K(\A_{i + 1}) \cap X(e_{i + 1}) \neq \emptyset$. Choose $x_{i + 1}\in K(\A_{i + 1}) \cap X(e_{i + 1})$. By (b), there exists $g_{i + 1} \in \G_0$ such that $g_{i + 1}(x_{i + 1}) = \vp_{e_{i + 1}}(x_{i + 1})$. \\
Thus by induction, the claim has been proven.\qed\par
Now, by (c), $X(e_0) \cap K(\A_1) \neq \emptyset$. This implies 
\begin{equation}\label{CHS.eq10}
f_{w(1)}(K(\A_1)) \cap K_{w(0)} \neq \emptyset. 
\end{equation}
Next, Claim-(ii) yields $f_{w(i + 1)}(g_i(x_i)) = f_{w(i)}(x_i)$. Moreover, since $g_i(x_i) \in K((g_i)^*(\A_i)) = K(\A_{i + 1})$, we have
\begin{equation}\label{CHS.eq20}
f_{w(i)}(K(\A_i)) \cap f_{w(i + 1)}(K(\A_{i + 1})) \neq \emptyset
\end{equation}
for $i = 1, \ldots, l$. Since $\A_i = (h_i)^*(\U_{\bp}) \subseteq \cup_{g \in \G_1} g^*(\psi_n(\bp))$,  it follows that $\cup_{i = 1}^l  w(i)\A_i  \subseteq \cup_{u \in \bp}H_u$. Note that $K(\cup_{i = 1}^l w(i)\A_i) = \cup_{i = 1}^l f_{w(i)}(\A_i)$. By \eqref{CHS.eq20} and \eqref{CHS.eq100}, we see that $K(\cup_{i = 1}^l w(i)\A_i)$ is connected and intersects with $K_{w(0)}$.  Thus there exists $\bp_0 \in \C_{m}^{(1)}(u_*, v_*, T_k)$ included in $\cup_{i = 1}^l w(i)\A_i \subseteq \cup_{u \in \bp} H_u$. Consequently, Lemma~\ref{UFD.lemma50} shows \eqref{CHS.eq100}. The conductive homogeneity follows from Lemma~\ref{CMS.lemma10} and Theorem~\ref{SUC.thm20}.
\enddemo

\setcounter{equation}{0}
\section{Subsystems of (hyper)cubic tiling}\label{CYP}

In this section, we present three classes of hypercube-based self-similar sets as examples of conductively homogeneous spaces. The first one given in Theorem~\ref{CYP.thm100} includes generalized  Sierpinski carpets studied in the series of papers \cite{BB1, BB2, BB3, BB4, BB5, BB6} by Barlow and Bass, the Menger curves (also known as the Menger sponge), and the hypercubes $[-1, 1]^L$ for $L \ge 1$. Unlike those examples, however, our examples also contain self-similar sets with fewer, or even no, symmetries of a hypercube. See Section~\ref{ESC}, where we present examples of self-similar sets having conductive homogeneity.\par
We start with basic notations on the hypercube $[-1, 1]^L$ and its symmetry group.

\definition\label{SHD.def00}
Let $L \in \BbN$ and let $C_*^L = [-1, 1]^L$. Moreover let $\BB_L$ be the $L$-dimensional hyperoctahedral group, that is,
\[
\BB_L = \{g| g \in O(L), g(C_*^L) = C_*^L\},
\]
where $O(L)$ is the collection of orthogonal transformations of $\BbR^L$.
\noindent Define
\[
B_{j, i} = \{(x_1, \ldots, x_L)| (x_1, \ldots, x_I) \in [-1, 1]^L,  x_j = i\}
\]
for $j = \{1, \ldots, L\}$ and $i \in \{-1. 0 , 1\}$. Then the boundary of $[-1, 1]^L$ consists of $\{B_{j, i}\}_{j \in \{1, \ldots, L\}, i \in \{1, -1\}}$. For $s = (s_1, \ldots, s_L) \in \{1, \ldots, N\}^L$, define
\[
C_{s}^{L, N} = \prod_{i = 1}^L \Big[\frac{2s_i - 2 - N}N, \frac{2s_i- N}N\Big],
\]
and
\[
c_s^{L, N} = \Big(\frac{2s_1 - 1 - N}N, \ldots, \frac{2s_L - 1 - N}N\Big).
\]
\enddefinition

If no confusion may occur, we use $C_*$, $C_s$ and $c_s$ instead of $C_*^L$, $C_s^{L, N}$ and $c_s^{L, N}$ respectively hereafter. \par
In the course of this section, we are going to deal with particular elements of $\BB_L$. 

\definition\label{SHD.def10}
Define $R_j \in \BB_L$ as the reflection in the hyperplane $B_{j, 0}$ for $j = \{1, \ldots, L\}$.  Furthermore, define $R_{j_1, j_2}^{i}$ as the reflection in the hyperplane
\[
\H_{j_1, j_2}^i = \{(x_1, \ldots, x_L)| x_{j_1} = ix_{j_2}\}
\]
for $j_1, j_2 \in \{1, \ldots, L\}$ with $j_1 \neq j_2$ and $i \in \{1, -1\}$.
\enddefinition

In the next definition, we introduce key notions of this section. \par
Throughout this section, we fix $L \ge 1$ and $N \ge 2$.

\definition\label{CYP.def10}
(1)\,\, A self-similar structure $(K, S, \{f_s\}_{s \in S})$ is called a subsystem of $L$-dimensional hypercubic tiling, or a subsystem of cubic tiling for short,  if $K \subseteq C_*$, $S \subseteq \{1, \ldots, N\}^L$ and, for any $s \in S$, $f_s$ is a restriction of a similitude from $\BbR^L$ to itself satisfying $f_s(C_*) = C_s$, i.e. there exists  $\Phi_s \in \BB_L$ such that
\begin{equation}\label{CYP.eq10}
f_s(x) = \frac 1N\Phi_sx + c_s
\end{equation}
for any $x \in \BbR^L$. A subsystem of cubic tiling $(K, S, \{f_s\}_{s \in S})$ is called non-degenerate if $K \cap B_{j, i} \neq \emptyset$ for any $j \in \{1, \ldots, L\}$ and $i \in \{1, -1\}$.\\
(2)\,\, A continuous map $\vp: C_* \to C_*$ is called an $N$-folding map if and only if, for any $s \in \{1, \ldots, N\}^L$,  there exists $A_s \in \BB_L$ such that
\begin{equation}\label{CYP.eq20}
\vp(x) = NA_s(x - c_s)
\end{equation}
for any $x \in C_s$. If no confusion may occur, we omit $N$ in the expression of an ``$N$-folding'' map and say a ``folding map'' for simplicity.\\
(3)\,\,Let $\L = (K, S, \{f_s\}_{s \in S})$ be a subsystem of cubic tiling. We use the framework of Section~\ref{SSF} to define $(T, \A, \phi)$ with $r = \frac 1N$ and $j_s = 1$ for any $s \in S$. In this case, $T_n = S^n$ for any $n \ge 1$. Define a graph $(T_n, E_n^{\ell})$ by 
\begin{multline*}
E_n^{\ell} = \{(w, v)| w, v \in T_n,  w \neq v, f_w(C_*) \cap f_v(C_*) = f_w(B_{j, i})\,\\
\text{for some $j \in \{1, \ldots, L\}$ and $i \in \{1, -1\}$}\}.
\end{multline*}
$\L$ is said to be strongly connected if and only if $(T_n, E_n^{\ell})$ is connected for any $n \ge 1$.\\
(4)\,\,Let $\L = (K, S, \{f_s\}_{s \in S})$ be a subsystem of cubic tiling. $\L$ is called locally symmetric if and only if $K_w \cup K_v$ is invariant under the reflection in the hyperplane including $f_w(C_*) \cap f_v(C_*)$ for any $n \ge 1$ and $(w, v) \in E_n^{\ell}$.
\enddefinition

\remark
Let $\L$ be a subsystem of cubic tiling which is non-degenerate and locally symmetric.  Then $E_n^{\ell} \subseteq E_n^*$ by the following arguments. Assume that $(w, v) \in E_n^{\ell}$. Set 
\begin{equation}\label{CYP.eq25}
\ell_{w, v} = f_w(C_*) \cap f_v(C_*).
\end{equation}
 By non-degeneracy, $K_w \cap \ell_{w, v} \neq \emptyset$ and by local symmetry, $K_w \cap \ell_{w, v} = K_v \cap \ell_{w, v} \neq \emptyset$. Hence $(w, v) \in E_n^*$. Note that even if $(w, v) \in T_n$ and $f_w(C_*) \cap f_v(C_*) \neq \emptyset$, it may happen that $K_w \cap K_v = \emptyset$.
\endremark

By properties of cubic tiling, it is easy to see that Assumption~\ref{ALFR} holds. In summary, we have the next proposition. Recall that the edges of $T_n$ is given not by $E_n^{\ell}$ but by $E_n^*$ as it has always been in the previous sections.

\prop\label{CYP.prop10}
Let $\L = (K, S, \{f_s\}_{s \in S})$ be a subsystem of cubic tiling. Then $\{K_w\}_{w \in T}$ is a partition of $K$ parametrized by the tree $(T, \A, \phi)$. Let $d_*$ be the restriction of the Euclidean metric on $K$ and let $\mu$ be the self-similar measure satisfying $\mu(K_w) = (\#(S))^{-|w|}$ for any $w \in T$.  Then Assumption~\ref{ALFR} is satisfied with $d = d_*$, $r = \frac 1N$, $M_* = 1$, $M_0 = 1$, $N_* = \#(S)$ and $L_*  \le 3^L - 1$. In this case, $\mu$ is $\a_H$-Ahlfors regular with respect to $d_*$, where $\a_H = \frac{\log {\#(S)}}{\log N}$.
\endprop

The exponent $\a_H$ coincides with the Hausdorff dimension of $(K, d_*)$. Note that $\#(S) \le N^L$. Since $\#(S) = N^L$ implies $K = C_*$, we see that $\a_H < L$ unless $K = C_*$. \par
The following theorems are the main results of this section.

\thm\label{CYP.thm100}
Let  $\L = (K, S, \{f_s\}_{s \in S})$ be a subsystem of cubic tiling. Assume that $\L$ is non-degenerate, locally symmetric, and strongly connected. Moreover, suppose that the following condition {\rm (SDR)} is satisfied{\rm:}\\
{\rm (SDR)}\,\,For any $j_1, j_2 \in \{1, \ldots, L\}$ with $j_1 \neq j_2$, there exists $i \in \{1, -1\}$ such that $R_{j_1, j_2}^i \in \G_{(K, T)}$.\\
Then $K$ is $p$-conductively homogeneous for any $p > \dim_{AR}(K, d_*)$.
\endthm

The name (SDR) represents ``symmetric with respect to diagonal reflections'' as $R_{j_1, j_2}^i$ is the reflection in the diagonal hyperplane $\H_{j_1, j_2}^i$. For generalized Sierpinski carpets, the Menger curve and the hypercube, it follows that $\G_{(K, T)} = \BB_L$ and the condition (SDR) is satisfied. However, $\G_{(K, T)}$ does not necessarily coincide with $\BB_L$ to satisfy (SDR). For example, the group generated by $\{R_{j_1, j_2}^1| j_1, j_2 \in \{1, \ldots, L\}, j_1 \neq j_2\}$ is (isomorphic to) the symmetric group of order $L$, $\mathcal{S}_L$, that is a proper subgroup of $\BB_L$, and if $\mathcal{S}_L \subseteq \G_{(K, T)}$, then the condition (SDR) is satisfied. See Example~\ref{CYP.ex60}.\par
In the case $L = 2$, the advantage of being planar gives another two classes having conductive homogeneity.

\thm\label{CYP.thm20}
Let $L = 2$ and let  $\L = (K, S, \{f_s\}_{s \in S})$ be a subsystem of $2$-dimensional cubic tiling. Assume that $\L$ is non-degenerate, locally symmetric, and strongly connected. Moreover, assume one of the following two conditions {\rm (RS)} or {\rm (NS)}.\\
{\rm (RS)}\,\,  $\Theta_{\pi/2} \in \G_{(K, T)}$, where $\Theta_{\pi/2}$ is the rotation by $\pi/2$ around $(0, 0)$.\\
{\rm (NS)}\,\, For each $i, j \in \{1, \ldots, N\}$, there exist $i_1, j_1 \in \{1, \ldots, N - 1\}$ such that 
\[
\{(i_1, j), (i_1 + 1, j), (i, j_1), (i, j_1 + 1)\} \cap S = \emptyset.
\]
Then $K$ is $p$-conductively homogeneous for any $p > \dim_{AR}(K, d_*)$.
\endthm

The expressions (RS) and (NS) represent ``rotational symmetry'' and ``no symmetry'' respectively.\par
At a glance at definitions, it may look difficult to verify the conditions like ``non-degenerate'', ``strongly continuous'', and  ``locally symmetric''. In the course of discussion, however, we will show useful criteria concerning only the first iteration $\{f_s(C_*)\}_{s \in S}$ to check those conditions. \par
Proofs of the above theorems will be given later in this section after necessary preparations. The main idea of the proof is to construct a family of paths required in the condition (c) of Theorem~\ref{CHS.thm10}  by using local symmetry and an additional geometric condition (SDR), (RS) or (NS). Such an idea was used in \cite{BouKleiner} and can be traced back to the ``Knight move'' argument by Barlow-Bass \cite{BB1}. In those previous works, however, the full $\BB_L$-symmetry of the space was required but we find that weaker (or even no) symmetry is good enough under the presence of local symmetry.\par
Now we start to study the conditions ``non-degenerate'', ``strong continuous'', and ``locally symmetric''.  First, we study the nature of folding maps, which turns out to be closely related to the local symmetry.

\lemma\label{CYP.lemma10}
Let $\vp: C_* \to C_*$ be a folding map characterized as \eqref{CYP.eq20}. Then for any $s, t \in \{1, \ldots, N\}^L$, 
\[
A_s = A_tR_j \quad\text{if $C_s \cap C_t = \frac 1NB_{j, i} + c_s$ for some $i \in \{1, -1\}$}.
\]
\endlemma

\demo
 Assume that $C_s \cap C_t = \frac 1NB_{j, i} + c_s$. Then $C_s \cap C_t = \frac 1NB_{j, -i} + c_t$ as well and $x - c_t = R_j(x - c_s)$ for any $x \in C_s \cap C_t$. On the other hand, as $\vp$ is a folding map, we see that 
\[
NA_s(x - c_s) = NA_t(x - c_t)
\]
for any $x \in C_s \cap C_t$. Hence $A_s(x - c_s) = A_tR_j(x - c_s)$ for any $x \in C_s \cap C_t$. This immediately implies $A_s = A_tR_j$.
\enddemo

Note that $R_{j_1}R_{j_2} = R_{j_2}R_{j_1}$ for any $j_1, j_2 \in \{1, \ldots, L\}$. So, by the above lemma, we can determine all the folding maps as follows.
\lemma\label{CYP.lemma20}
Fix $s^* = (s^*_1, \ldots, s^*_L) \in \{1, \ldots, N\}^L$. For $A \in \BB_L$, define $\vp_{s^*, A}: C_* \to C_*$ by 
\[
\vp_{s^*, A}(x) = NA\prod_{j = 1}^L(R_j)^{|s^*_j - s_j|}(x - c_{(i, j)}^N)
\]
for any $x \in C_{(s_1, \ldots, s_L)}$.  Then $\vp_{s_0, A}$ is a folding map. Moreover, $\{\vp_{s^*, A}| A \in \BB_L\}$ is the totality of folding maps for any $s^* \in \{1, \ldots, N\}^L$.
\endlemma

\begin{figure}
\centering
\includegraphics[width=\linewidth]{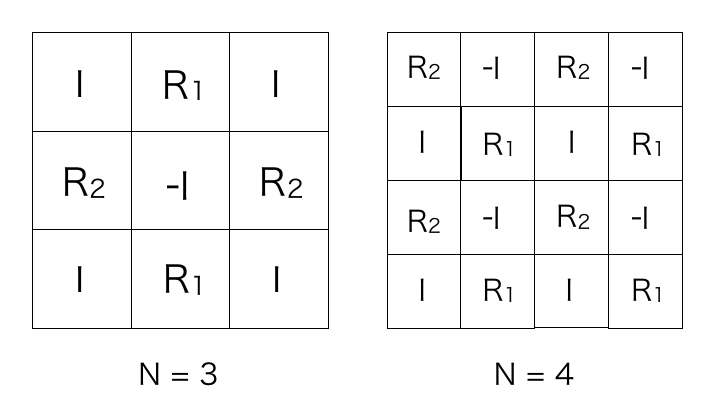}
\caption{Folding maps}\label{FigFold}
\end{figure}

Examples of folding maps in the case of $L = 2$ are given in Figure~\ref{FigFold}.  In each example, $s^* = (1, 1)$ and $A = I$. The element of $\BB_2$ in each square  indicates the corresponding $A(R_1)^{|s_1 - s^*_1|}(R_2)^{|s_2 - s^*_2|}$.

\notation
Let $\L = (K, S, \{f_s\}_{s \in S})$ be a subsystem of cubic tiling. Set
\[
K^{(m)} = \bigcup_{w \in T_m} f_w(C_*).
\]
\endnotation

Due to the next lemma, one can easily determine non-degeneracy of $K$ by examining $K^{(1)}$.

\lemma\label{CYP.lemma30}
Let $\L = (K, S, \{f_s\}_{s \in S})$ be a  subsystem of cubic tiling. Then $\L$ is non-degenerate if and only if $K^{(1)} \cap B_{j, i} \neq \emptyset$  for any $j \in \{1, \ldots, L\}$ and $i \in \{1, -1\}$.
\endlemma

\demo
Since $K \subseteq K^{(1)}$, the ``only if'' part is obvious. Assume that $K^{(1)} \cap B_{j, i} \neq \emptyset$  for any $j \in \{1, \ldots, L\}$ and $i \in \{1, -1\}$. We are going to show that $K^{(k)} \cap B_{j, i} \neq \emptyset$  for any $j \in \{1, \ldots, L\}$, $i \in \{1, -1\}$, and $k \in \{1, \ldots, n\}$ by induction on $n$. Assume that the claim holds for $n$. Let $w \in T_n$ satisfying $f_w(C_*) \cap B_{j, i} \neq \emptyset$. Since $(f_w)^{-1}(f_w(C_*) \cap B_{j, i}) = B_{j_1, i_1}$ for some $j_1 \in \{1, \ldots, L\}$ and $i_1 \in \{1, -1\}$, there exists $s \in T_1$ such that $f_s(C_*) \cap (f_w)^{-1}(f_w(C_*) \cap B_{j, i}) \neq \emptyset$.  This implies that $f_{ws}(C_*) \cap B_{j, i} \neq \emptyset$. Thus we have shown the desired statement for $n + 1$. Now by induction, $K^{(k)} \cap B_{j, i} \neq \emptyset$  for any $j \in \{1, \ldots, L\}$, $i \in \{1, -1\}$. Since $K^{(n)}$ is monotonically decreasing and $K = \cap_{n \ge 1} K^{(n)}$, it follows that $K \cap B_{j, i} \neq \emptyset$ for any $j \in \{1, \ldots, L\}$ and $i \in \{1, -1\}$.
\enddemo

The locally symmetric property can also be determined by the first step of the iteration as follows.

\lemma\label{CYP.lemma40}
Let $\L = (K, S, \{f_s\}_{s \in S})$ be a  subsystem of cubic tiling. Then $\L$ is locally symmetric if and only if $K_s \cup K_t$ is invariant under the reflection in $\ell_{s, t}$ for any $(s, t) \in E_1^{\ell}$.  
\endlemma
\demo
The ``only if'' part is obvious. We show the following statement by induction on $n \ge 1$. \\
For any $k \in \{1, \ldots, n\}$,  and $(w, v) \in E_k^{\ell}$, $K_w \cup K_v$ is invariant under the reflection in $\ell_{w, v}$.\\
The case $n = 1$ is exactly the assumption of the lemma. Suppose that the statement holds for $n$. Let $(w, v) \in E_{n + 1}^{\ell}$.  In case $\pi^n(w) = \pi^n(v)$, let $s = \pi^n(w)$. Then $w = sw'$ and $v = sv'$ for some $w', v' \in T_n$. Since $f_w(C_*) = f_s(f_{w'}(C_*))$ and $f_v(C_*) = f_s(f_{v'}(C_*))$, we see $\ell_{w', v'} \in E_n^{\ell}$. By induction hypothesis,  $K_{w'} \cap K_{v'}$ is invariant under the reflection in $\ell_{w', v'}$. Applying $f_s$, we see that $K_w \cup K_v$ is invariant under the reflection in $\ell_{w,v}$. In case $\pi^n(w) \neq \pi^n(v)$, let $s = \pi^n(w)$ and let $t = \pi^n(v)$. Since $\ell_{w, v} \subseteq \ell_{s, t} = f_s(B_{j, i})$ for some $j \in \{1, \ldots, L\}$ and $i \in \{1, -1\}$, we obtain $(s, t) \in E_1^{\ell}$.  So, $K_s \cup K_t$ is invariant under the reflection in $\ell_{s, t}$. Denoting this reflection by $R$, we see that $R$ coincides with the reflection in $\ell_{w, v}$. Since $R(f_w(C_*)) = f_v(C_*)$, it follows that $R(K_w) = R(K_s \cap f_w(C_*)) =K_t \cap f_v(C_*) = K_v$.  So we have verified the statement for $n + 1$. Thus by induction, we have the desired result.
\enddemo

Next, we consider the strongly connectedness.

\lemma\label{CYP.lemma50}
Let $\L = (K, S, \{f_s\}_{s \in S})$ be a locally symmetric subsystem of cubic tiling. If $\L$ is non-degenerate and $(T_1, E_1^{\ell})$ is connected, then $\L$ is strongly connected.
\endlemma

\demo
By the non-degeneracy, we see that $K^{(n)} \cap B_{j, i} \neq \emptyset$ for any $j \in \{1, \ldots, L\}$ and $i \in \{1, -1\}$.\par
We are going to show that $(T_k, E_k^{\ell})$ is connected for any $ k \in \{1, \ldots, n\}$ by induction on $n \ge 1$. Assume that $w, v \in T_{n + 1}$. If $\pi^n(w) = \pi^n(v)$, then there exist $w', v' \in T_n$ such that $w = sw'$ and $v = sv'$, where $s = \pi^n(w)$. Since $w'$ and $v'$ are connected by an $E_n^{\ell}$-path, $w$ and $v$ are connected by an $E_{n + 1}^{\ell}$-path. In case $\pi^n(w) \neq \pi^n(v)$, let $s = \pi^n(w)$ and let $t = \pi^n(v)$. Then $w = sw'$ and $v = tv'$ for some $w', v' \in T_n$. Since $(T_1, E_1^{\ell})$ is connected, there exists an $E_1^{\ell}$-path $(s(0), \ldots, s(m))$ such that $s(0) = s$, $s(m) = t$ and $(s(i), s(i + 1)) \in E_1^{\ell}$ for any $i = 0, \ldots, m - 1$.  For each $i = 0, \ldots, m - 1$, since $\cup_{w' \in T_n} f_{w'}(C_*) \cap B_{j, i} \neq \emptyset$ for any $j = \{1, \ldots, L\}$ and $i \in \{1, -1\}$, there exists $u(i) \in T_n$ such that $f_{s(i)u(i)}(C_*) \cap \ell_{s(i), s(i + 1)} \neq \emptyset$. Since $\L$ is locally symmetric, there exists $v(i) \in T_n$ such that $f_{s(i + 1)v(i)}(C_*)$ is the image of $f_{s(i)u(i)}(C_*)$ by the reflection in $\ell_{s(i), s(i + 1)}$. Define $v(-1) = w'$ and $u(m) = v'$. Then $w = s(0)v(-1)$ and $v = s(m)u(m)$. Since $(T_n, E_n^{\ell})$ is connected, $v(i - 1)$ and $u(i)$ are connected by an $E_n^{\ell}$-path for any $i = 0, \ldots, m - 1$. Adding $s(i)$ at the top, we obtain an $E_{n + 1}^{\ell}$-path between $s(i)v(i - 1)$ and $s(i)u(i)$. Combining all these $E_{n + 1}^{\ell}$-paths, we obtain an $E_{n + 1}^{\ell}$-path between $w$ and $v$. Thus $(T_{n + 1}, E_{n + 1}^{\ell})$ is connected. By induction, we see that $\L$ is strongly connected.
\enddemo

\lemma\label{CYP.lemma60}
Let  $\L = (K, S, \{f_s\}_{s \in S})$ be a subsystem of cubic tiling. Assume that $K \cap \inte{C_*} \neq \emptyset$. For any $s \in \{1, \ldots, N^m\}^L$, if $K \cap \inte{C_{s}^{L, N^m}} \neq \emptyset$, then there exists $w \in T_m$ such that $f_w(C_*) = C_{s}^{L, N^m}$.
\endlemma

\demo
Suppose that $f_w(C_*) \neq C_{s}^{L, N^m}$ for all $w \in T_m$. Then $f_w(C_*) \cap C_{s}^{L, N^m}$ is included in the boundary of $C_s^{L, N^m}$ and hence $f_w(C_*) \cap \inte{C_{s}^{L, N^m}} = \emptyset$. So, 
\[
K^{(m)} \cap \inte{C_{s}^{L, N^m}} = \bigcup_{w \in T_m} \big(f_w(C_*) \cap \inte{C_{s}^{L, N^m}}\big) = \emptyset.
\]
Since $K \subseteq K^{(m)}$, it follows that $K \cap \inte{C_{s}^{L, N^m}} = \emptyset$.
\enddemo

The following relation between a folding map and a subsystem of cubic tiling will be used to characterize local symmetry.

\lemma\label{CYP.lemma70}
Let  $\L = (K, S, \{f_s\}_{s \in S})$ be a subsystem of cubic tiling. Assume that $K \cap \inte{C_*} \neq \emptyset$. Let $\vp$ be a folding map. Then the following four statements are equivalent:\\
{\rm (a)}\,\,$\vp(K) = K$.\\
{\rm (b)}\,\,$\vp{\circ}f_s(K^{(m)}) = K^{(m)}$ for any $s \in S$ and $m \ge 0$.\\
{\rm (c)}\,\,$\vp{\circ}f_s(K) = K$ for any $s \in S$.\\
{\rm (d)}\,\,$\vp(K^{(m + 1)}) = K^{(m)}$ for any $m \ge 0$.
\endlemma

\demo
(a) $\Rightarrow$ (b):\,\,Let $s \in S$. Then $\vp{\circ}f_s(K) \subseteq K$. For any $w \in T_m$, there exists $\tau = (\tau_1, \ldots, \tau_L) \in \{1, \ldots, N^m\}^L$ such that $\vp{\circ}f_s(f_w(C_*)) = C_{\tau}^{L, N^m}$. Now
\[
K \supseteq \vp{\circ}f_s(f_w(K \cap \inte{C_*})) = \vp{\circ}f_s{\circ}f_w(K) \cap \inte{C_{\tau}^{L, N^m}}.
\]
Since $K \cap \inte{C_*} \neq \emptyset$, this implies $K \cap \inte{C_{\tau}^{L, N^m}} \neq \emptyset$. Lemma~\ref{CYP.lemma60} shows that $\vp{\circ}f_s(f_w(C_*)) = C_{\tau}^{L, N^m} \subseteq K^{(m)}$, so that 
\[
\vp{\circ}f_s(K^{(m)}) = \bigcup_{w \in T_m} \vp{\circ}f_s(f_w(C_*)) \subseteq K^{(m)}.
\]
Note that $\vp{\circ}f_s \in \BB_L$ preserves the Lesbegue measure of a set. Hence we see $\vp{\circ}f_s(K^{(m)}) = K^{(m)}$. \\
(b) $\Rightarrow$ (c):\,\, Since $\cap_{m \ge 0} K^{(m)} = K$, 
\[
\vp{\circ}f_s(K) = \vp{\circ}f_s\Big(\bigcap_{m \ge 0} K^{(m)}\Big) = \bigcap_{m \ge 0} K^{(m)} = K.
\]
(c) $\Rightarrow$ (a):\,\,Since $K = \cup_{s \in S} f_s(K)$, 
\[
\vp(K) = \vp\Big(\bigcup_{s \in S} f_s(K)\Big) = K.
\]
(b) $\Rightarrow$ (d):\,\, Since $\cup_{s \in S} f_s(K^{(m)}) = K^{(m + 1)}$,
\[
\vp(K^{(m + 1)}) = \vp\Big(\bigcup_{s \in S} f_s(K^{(m)})\Big) = K^{(m)}.
\]
(d) $\Rightarrow$ (a):\,\, Since $\cap_{m \ge 0} K^{(m)} = K$, 
\[
\vp(K) = \vp\Big(\bigcap_{m \ge 0} K^{(m + 1)}\Big) = \bigcap_{m \ge 0} K^{(m)} = K.
\]
\enddemo

The next theorem tells that a locally symmetric subsystem of cubic tiling is almost an inverse of a folding map.
\thm\label{CYP.thm10}
Let  $\L = (K, S, \{f_s\}_{s \in S})$ be a subsystem of cubic tiling. \\
{\rm (1)}\,\,
If $\L$ is strongly connected and locally symmetric, then there exists a folding map satisfying
\begin{equation}\label{LS}
\vp^n\circ{f_w}(K^{(m)}) = K^{(m)}
\end{equation}
for any $n \ge 1, m \ge 0$ and $w \in T_n$. In particular, 
\[
\vp^n(K^{(n + m)}) = K^{(m)}
\]
for any $n \ge 1$, $m \ge 0$ and
\[
\vp^n(K) = K
\]
for any $n \ge 1$. Furthermore, define $F_s: C_* \to C_s$ by $F_s  = (\vp|_{C_s})^{-1}$ for each $s \in S$. Then
\[
K = \bigcup_{s \in S} F_s(K)
\]
and $(K, S, \{F_s\}_{s \in S})$ is a self-similar structure.\\
{\rm (2)}\,\,
Suppose that $K \cap \inte{C_*} \neq \emptyset$. If there exists a folding map $\vp$ such that $\vp(K) = K$, then $\L$ is locally symmetric.
\endthm

\demo
(1)\,\,Fix $s \in S$. Recall that there exists $\Phi_s \in \BB_L$ such that
\[
f_{s}(x) = \frac 1N\Phi_{s}x + c_{s}
\]
for any $x \in C_*$. Set $A_{s} = (\Phi_{s})^{-1}$ and define $\vp = \vp_{s_0, A_{s}}$. Since $\vp{\circ}f_{s} = I$, it follows that $\vp^n{\circ}(f_{s})^n = I$ for any $n \ge 1$. Thus letting $s_n = \underset{\text{$n$-times}}{ss\cdots{s}}$, we see that $\vp^n{\circ}f_{s_n}(K) = K$. Choose $\tau = (\tau_1, \ldots, \tau_L) \in \{1, \ldots, N^n\}^L$ such that $C_{\tau}^{L, N^n} = f_{s_n}(C_*)$. Let $w \in T_n$. Choose $\xi = (\xi_1, \ldots, \xi_L) \in \{1, \ldots, N^n\}^L$ such that $C_{\xi}^{L, N^n} = f_w(C_*)$.  Since $\L$ is strongly connected, there exists  an $E_n^{\ell}$-path $(w(0), \ldots, w(m))$ between $s_n$ and $w$. Following this path and applying the reflections in $\ell_{w(i), w(i + 1)}$, we see that
\[
K_w - c_{\xi}^{L, N^n} = R(K_{s_n} - c_{\tau}^{L, N^n}),
\]
where $R = \prod_{j = 1}^L (R_j)^{|\tau_j - \xi_j|}$. Note that $\vp^n$ is an $N^n$-folding map. For any $\c \in \{1, \ldots, N^n\}^L$, there exists $A_{\c} \in \BB_L$ such that
\[
\vp^n(x) = N^nA_{\c}(x - c_{\c}^{L, N^n})
\]
for any $x \in C_{\c}^{L, N^n}$.  Applying Lemma~\ref{CYP.lemma20} to $\vp^n$, we see that
\begin{multline*}
\vp^n{\circ}f_w(K) = \vp^n(K_w) = N^nA_{\xi}(K_w - c_{\xi}^{L, N^n}) \\
=N^nA_{\tau}RR(K_{s_n} - c_{\tau}^{L, N^n}) = \vp^n(K_{s_n}) = K.
\end{multline*}
Hence 
\[
\vp^n{\circ}f_w(K) = K
\]
for any $n \ge 1$ and $w \in T_n$. Since $K \subseteq K^{(m)}$, it follows that $\vp^n{\circ}f_w(K^{(m)}) \supseteq K$. Note that $\vp^n{\circ}f_w(K^{(m)}) = \cup_{\c\in B} C_{\c}^{L, N^n}$  for some subset $B \subseteq \{1, \ldots, N^n\}^L$ and $K^{(m)}$ is the minimal of such unions containing $K$. This shows $\vp^n{\circ}f_w(K^{(m)}) \supseteq K^{(m)}$. Since $\vp^n{\circ}f_w$ preserves the Lebesgue measure of a set, we conclude that $\vp^n{\circ}f_w(K^{(m)}) = K^{(m)}$. Since $K^{(m + n)} = \cup_{w \in T_n} f_w(K^{(m)})$, we obtain $\vp^n(K^{(n + m)}) = K^{(m)}$. Note that $K = \cup_{w \in T_n} f_w(K)$. Hence $\vp^n(K) = K$. Moreover, if $\vp(x) = NA_s(x - c_s)$ for $x \in C_s$, then by Lemma~\ref{CYP.lemma70}-(c), we have $K = NA_s(K_s - c_s)$. This implies $K_s = \frac 1N(A_s)^{-1}K + c_s$. Hence letting $F_s(x) = \frac 1N(A_s)^{-1}x + c_s$, we see $K = \cup_{s \in S} F_s(K)$.
\\
(2)\,\,
  Suppose that $(s, t)  \in E_1^{\ell}$. Then by Lemma~\ref{CYP.lemma20},  there exist $A_s \in \BB_L$ and $j \in \{1, \ldots, L\}$ such that 
  \[
\vp(x) = NA_s(x - c_s)
\]
for any $x \in C_s$ and
\[
\vp(x) = NA_sR_j(x - c_t)
\]
for any $x \in C_t$. Since $\vp{\circ}f_s(K) = K$ and $\vp{\circ}f_t(K) = K$ by Lemma~\ref{CYP.lemma70}, it follows that
\[
K_s - c_s = \frac 1N(A_s)^{-1}K\quad \text{and}\quad K_t - c_t = \frac 1NR_j(A_s)^{-1}K.
\]
Therefore,
\[
R(K_s - c_s) = R\frac 1N(A_s)^{-1}K = K_t - c_t,
\]
so that $K_t \cup K_s$ is invariant under the reflection in $\ell_{s, t}$. Thus Lemma~\ref{CYP.lemma40} shows that $\L$ is locally symmetric.
\enddemo

By (2) of the above theorem, we immediately have the following sufficient condition for the local symmetry.

\cor\label{CYP.cor20}
Let  $S \subseteq \{1, \ldots, N\}^L$. Assume that $B_{j, i} \cap (\cup_{s \in S} C_s) \neq \emptyset$ for any $j \in \{1, \ldots, L\}$ and $i \in \{1, -1\}$. Let $\vp$ be an $N$-folding map. Define
\[
f_s = (\vp|_{C_s})^{-1}
\]
for any $s \in S$. Let $K$ be the unique non-empty compact set satisfying
\[
K = \bigcup_{s \in S} f_s(K).
\]
Then, $\L = (K, S, \{f_s\}_{s \in S})$ is non-degenerate and locally symmetric.
\endcor
\demo
Since $B_{i, i} \cap (\cup_{s \in S} C_s) \neq \emptyset$ for any $j \in \{1, \ldots, L\}$ and $i \in \{1, -1\}$, Lemma~\ref{CYP.lemma30} shows that $\L$ is non-degenerate and hence $K \cap \inte{C_*} \neq \emptyset$. Moreover it is immediate to see that $\vp(K) = K$. Now Theorem~\ref{CYP.thm10}-(2) suffices.
\enddemo

Note that by Theorem~\ref{CYP.thm10}-(1), any subsystem of cubic tiling that is locally symmetric and strongly continuous is given by an inverse of a folding map described in Corollary~\ref{CYP.cor20}.\par
Now we are ready to give a proof of Theorem~\ref{CYP.thm100}.

\demo[Proof of Theorem~\ref{CYP.thm100}]
By Theorem~\ref{CYP.thm10}, we may assume that $\L$ is given by an inverse of a folding map described in Corollary~\ref{CYP.cor20} without loss of generality.  Note that 
\begin{equation}\label{CYP.eq30}
(\vp^m|_{f_w(C_*)})^{-1} = f_w
\end{equation}
for any $m \ge 1$ and $w \in T_m$. For any $m \ge 1$ and $e = (w, v) \in E_m^{\ell}$, by \eqref{CYP.eq30}, 
\[
\vp^m|_{f_w(C_*) \cap f_v(C_*)} = (f_w)^{-1}|_{f_w(C_*) \cap f_v(C_*)} = (f_v)^{-1}|_{f_w(C_*) \cap f_v(C_*)}.
\]
Hence $X(e) = X(e^r)$ and $\vp_e = I$, where $I$ is the identity map. Now let
\[
\I = \{(X(e), X(e^r), \vp_e)| e \in \cup_{m \ge 1} E_m^{\ell}\},
\]
and set $\G_0 = \{I\}$ and $\G_1 = \G_{(K, T)} \cap \BB_L$. We are going to make use of Theorem~\ref{CHS.thm10}. The condition (a) of Theorem~\ref{CHS.thm10} follows because $\L$ is strongly connected. Since $\vp_e = I$ for any $e \in \cup_{m \ge 1} E_m^{\ell}$, the condition (b) of Theorem~\ref{CHS.thm10} is obvious.\par
Now it only remains to show the condition (c) of Theorem~\ref{CHS.thm10}. Let $w \in T_n$. Suppose that $f_w(C_*) = \prod_{i = 1}^L [\a_i, \a_i + 2/N^n]$. Then every path $\bp \in \C^{(1)}_{1, m}(w)$ contains a path between hyperplanes
\[ 
\{(x_1, \ldots, x_L)| x_{j} = \a_{j}\}\quad\text{and}\quad \{(x_1, \ldots, x_L)| x_{j} = \a_{j} - 2/N^n\}
\]
or
\[ 
\{(x_1, \ldots, x_L)| x_{j} = \a_{j} + 2/N^n\}\quad\text{and}\quad \{(x_1, \ldots, x_L)| x_{j} = \a_{j} + 4/N^n\}
\]
for some $j \in \{1, \ldots, L\}$. This implies that there exists $j_* \in \{1, \ldots, L\}$ such that $\vp^n(K(\bp)) \cap B_{j_*, i} \neq \emptyset$ for any $i \in \{1, -1\}$. Note that $\vp^m(K(\bp)) = K(\psi_n(\bp))$. Hence there exists a path $\bp_{j_*}  \subseteq \psi_n(\bp)$ between $B_{j_*, -1}$ and $B_{j_*, 1}$. By the condition (SDR), for any $j_1 \neq j_*$, there exists $i_* \in \{1, -1\}$ such that $R_{j_*, j_1}^{i_*} \in \G_{(K, T)}$. Set $\bp_{j_1} = (R_{j_*, j_1}^{i_*})^*(\bp_{j_*})$. Then $K(\bp_{j_1}) \cap B_{j_1, i} \neq \emptyset$ for any $i \in \{1, -1\}$. Moreover $K(\bp_{j_*})$ and $K(\bp_{j_1})$ intersects at $H_{j_*, j_1}^{i_*}$. Thus set $\bp_* = \cup_{k = 1}^L\bp_k$. Then $\bp_*$ is connected and $K(\bp_*) \cap B_{k, i} \cap K \neq \emptyset$ for any $k \in \{1, \ldots, L\}$ and $i \in \{1, -1\}$. Moreover, $\bp_* \subseteq \cup_{g \in \G_{(K, T)} \cap \BB_L} g^*(\psi_n(\bp))$. Thus we have verified the condition (c) of Theorem~\ref{CHS.thm10}.
\enddemo

\demo[Proof of Theorem~\ref{CYP.thm20}]
The arguments are the same as in the proof of Theorem~\ref{CYP.thm100} except the deduction of the condition (c) of Theorem~\ref{CHS.thm10}. \par
In the case of (RS), to construct $\bp_{j_1}$ from $\bp_{j_*}$, we use $\Theta_{\pi/2}$ in place of $R_{j_*, j_1}^{i_*}$. Then the advantage of being planar yields $K(\bp_{j_*}) \cap K(\bp_j) \neq \emptyset$. The rest is the same as in the proof of Theorem~\ref{CHS.thm10}.\par
Next assume (NS). Let $w \in T_n$ and let $\bp = (w(1), \ldots, w(k)) \in \C^{(1)}_{M, m}(w)$ with $M = 4N - 3$.  Note that
\[
\#(\{\pi^m(w(1)), \ldots, \pi^m(w(k))\}) \ge M.
\]
We are going to show that 
\begin{equation}\label{CYP.eq40}
K(\psi_n(\bp)) \cap B_{j, i} \neq \emptyset
\end{equation}
for any $j \in \{1, 2\}$ and $i  \in \{1, -1\}$.  Suppose $K(\psi_n(\bp)) \cap B_{1, 1} = \emptyset$. Since $\vp^{-n}(B_{1, 1})$ forms vertical lines at intervals of $\frac{2}{N^{n}}$,  we see that $K(\bp)$ is contained in the interior of a vertical strip $\cup_{j = 1, \ldots, N^{n}} C_{(i_*, j)}^{2, N^{n}} \cup C_{(i_* + 1, j)}^{2, N^{n}}$, which is denoted by $Z_{i_*}$, for some $i_*$. Let $C_1, \ldots, C_l$ be the collection of connected components of
\[
\Big(\bigcup_{w \in T_{n}} f_w(Q) \Big) \cap  Z_{i_*}
\]
and set
\[
D_i = \{v | v \in T_{n}, f_v(C_*) \subseteq C_i\}
\]
for $ i = 1, \ldots, l$. Then by (NS), we see that
\[
\#(D_i) \le 2(2N - 2).
\]
Note that $\cup_{i = 1}^k  f_{\pi^m(w(i))}(C_*) \subset C_{i_*}$ for some $i_*$. Hence
\[
4N - 4 \ge \#(D_{i_*}) \ge \#(\{\pi^m(w(i))| i = 1, \ldots, k\}) \ge M = 4N - 3.
\]
This contradiction shows \eqref{CYP.eq40}. Thus setting $\U_p = \psi_n(\bp)$, we have the condition (c) of Theorem~\ref{CHS.thm10}.
\enddemo

To conclude this section, we present a useful criterion to determine whether $g \in \BB_L$ is a symmetry of $(K, T)$ or not.

\lemma\label{CYP.lemma80}
Let $\L = (K, S, \{f_s\}_{s \in S})$ be a subsystem of cubic tiling. Assume that $\L$ is non-degenerate, locally symmetric and strongly connected. Let $\vp$ be the folding map satisfying the condition of Theorem~\ref{CYP.thm10}-(1). Then for $g \in \BB_L$, if there exists a map $g_*: S \to S$ such that, for any $s \in S$,  $g(C_s) = C_{g_*(s)}$ and $A_{g_*(s)}g(A_s)^{-1} = g^k$ for some $k \ge 0$, then $g \in \G_{(K, T)}$.
\endlemma

Recall that $A_s \in \BB_2$ is given in Definition~\ref{CYP.def10}-(2).

\demo
We are going to show that $g(K^{(n)}) = K^{(n)}$ for any $n \ge 1$ by induction. For $n = 1$, Since $g(C_s) = C_{g_*(s)}$, it follows $g(K^{(1)}) = K^{(1)}$. Next assume that $g(K^{(n)}) = K^{(n)}$. Then by Theorem~\ref{CYP.thm10}, $\vp{\circ}f_s(K^{(n)}) = K^{(n)}$,  so that $A_s\Phi_s(K^{(n)}) = K^{(n)}$. Hence
\[
f_s(K^{(n)}) = \frac 1N(A_s)^{-1}(K^{(n)}) + c_s.
\]
Set $t = g_*(s)$. Then
\begin{multline*}
g(f_s(K^{(n)})) = \frac 1Ng(A_s)^{-1}(K^{(n)}) + c_{t} = \frac 1N(A_t)^{-1}A_tg(A_s)^{-1}(K^{(n)}) + c_t\\
= \frac 1N(A_t)^{-1}g^k(K^{(n)}) + c_t = f_t(K^{(n)}).
\end{multline*}
Since $K^{(n + 1)} = \cup_{s \in S} f_s(K^{(n)})$, this yields $g(K^{(n + 1)}) = K^{(n + 1)}$. Thus using induction, we see that $g(K^{(n)}) = K^{(n)}$ for any $n \ge 1$. Since $\cap_{n \ge 1} K^{(n)} = K$, we obtain $g(K) = K$. Now, since $g(K^{(n)}) = K^{(n)}$, it follows that, for any $w \in T_n$, there exists $v \in T_n$ such that $g(f_w(C_*)) = f_v(C_*)$. Set $v = g_*(w)$. Then $g_*: T_n \to T_n$. Since $g(f_w(C_*)) = f_{g_*(v)}(C_*)$ and $g(K_w) \subseteq K$, we see that
\[
g(K_w) \subseteq g(f_w(C_*)) \cap K = f_{g_*(w)}(C_*) \cap K = K_{g_*(w)}.
\]
Using $g^{-1}$ in place of $g$ in the arguments above, we obtain $g^{-1}(K_{g_*(w)}) \subseteq K_w$ as well. Thus we have shown $g(K_w) = K_{g_*(w)}$, so that $g \in \G_{(K, T)}$.
\enddemo

\setcounter{equation}{0}
\section{Examples: subsystems of (hyper)cubic tiling}\label{ESC}

In this section, we present examples of subsystems of cubic tiling having conductive homogeneity.\par
 We begin with planar examples where $\dim_{AR}(K, d_*) \le \dim_H(K, d_*) < 2$, so that they are $2$-conductively homogeneous and have self-similar local regular Dirichlet forms constructed in Theorem~\ref{SSF.thm10}.

\example[Chipped Sierpinski carpet]\label{CYP.ex10}

\begin{figure}
\hspace*{0pt}
\begin{minipage}[b]{6cm}
\centering
\includegraphics[width = 115pt]{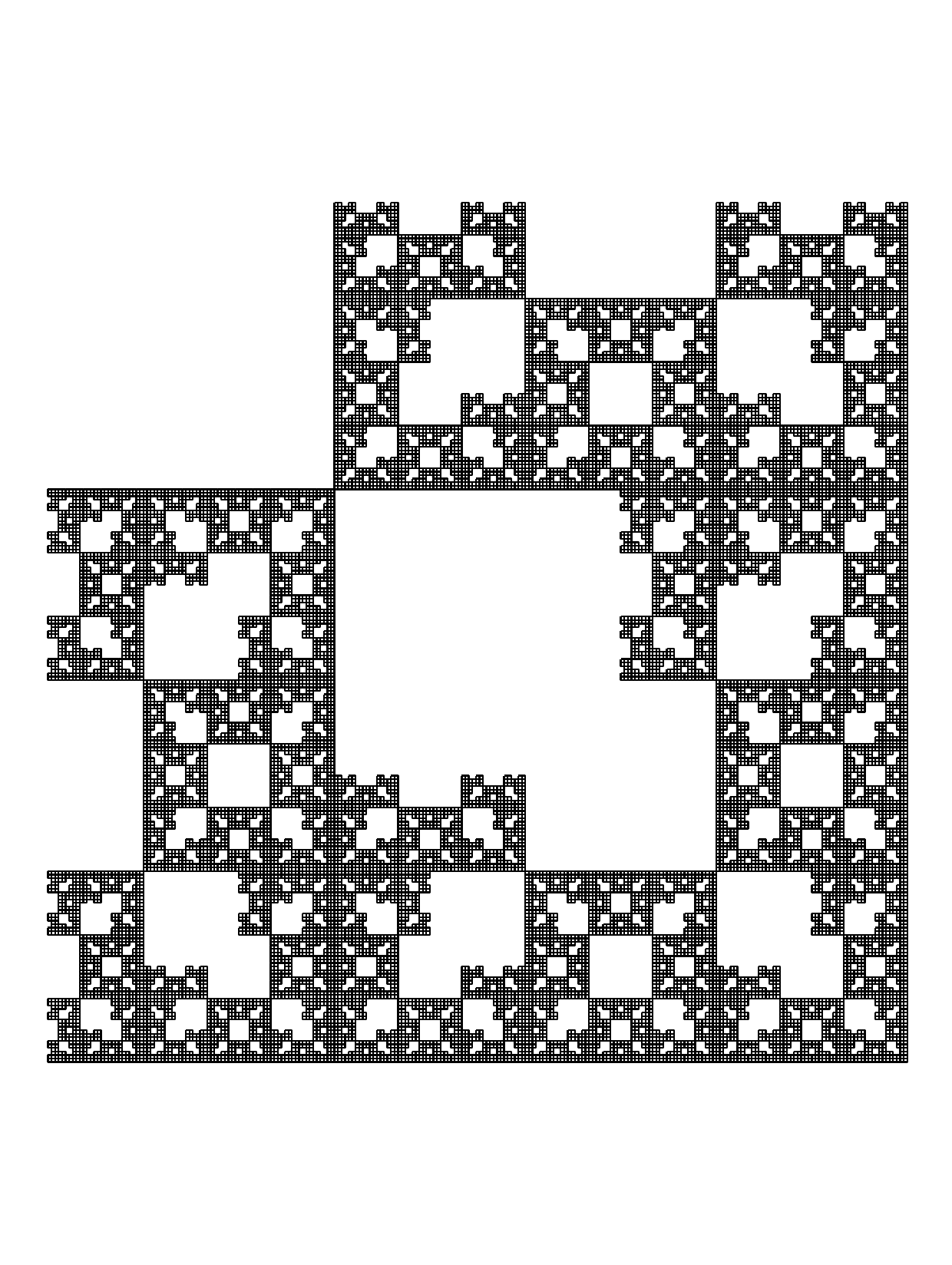}
\vspace{-17pt}
\end{minipage}
\begin{minipage}[b]{6cm}
\hspace{-27pt}
\centering
\includegraphics[width=125pt]{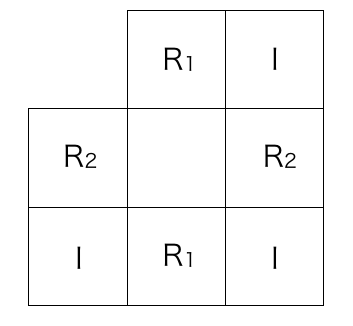}
\end{minipage}
\caption{Chipped Sierpinski carpet}\label{FIG1}
\vspace{15pt}
\end{figure}

Let $L = 2$ and let $N = 3$. Let $S$ be the set of squares in the right figure of Figure~\ref{FIG1} where one of $R_1, R_2$ or $I$ is written. The corresponding $f_s$ is given by
\[
f_s(x) = \frac 1N{\Phi_s}x + c_s^3,
\]
where $\Phi_s  \in \BB_2$ is indicated in Figure~\ref{FIG1}. Note that if the upper-left square belonged to $S$ as well, then $K$ would be the Sierpinski carpet.  Lemma~\ref{CYP.lemma30} and Corollary~\ref{CYP.cor20} show that $\L$ is non-degenerate and locally symmetric respectively. Then using Lemma~\ref{CYP.lemma50}, we see that $\L$ is strongly connected. Finally Lemma~\ref{CYP.lemma80} shows that $R_{1, 2}^{-1} \in \G_{(K, T)}$, so that (SDR) is satisfied. Thus we have confirmed all the assumptions in Theorem~\ref{CYP.thm100}.  Note that $K \cap \partial{C_*}$ has two different ingredients, the line segment, and the Cantor set. The lack of a rotational symmetry enables such a phenomenon. Another unique feature is the ``countably ramified'' property, that is,  after removing a certain countable set, every remaining point becomes a connected component. In this example, since there are enough number of straight lines inside $K$, $(K, d_*)$ has the chain condition and hence the heat kernel associated with $(\E, \W^2)$ satisfies \eqref{SUC.eq100} and \eqref{SUC.eq120}.
\endexample

\example\label{CYP.ex20}
\begin{figure}
\hspace*{0pt}
\begin{minipage}[b]{6cm}
\centering
\includegraphics[width = 114pt]{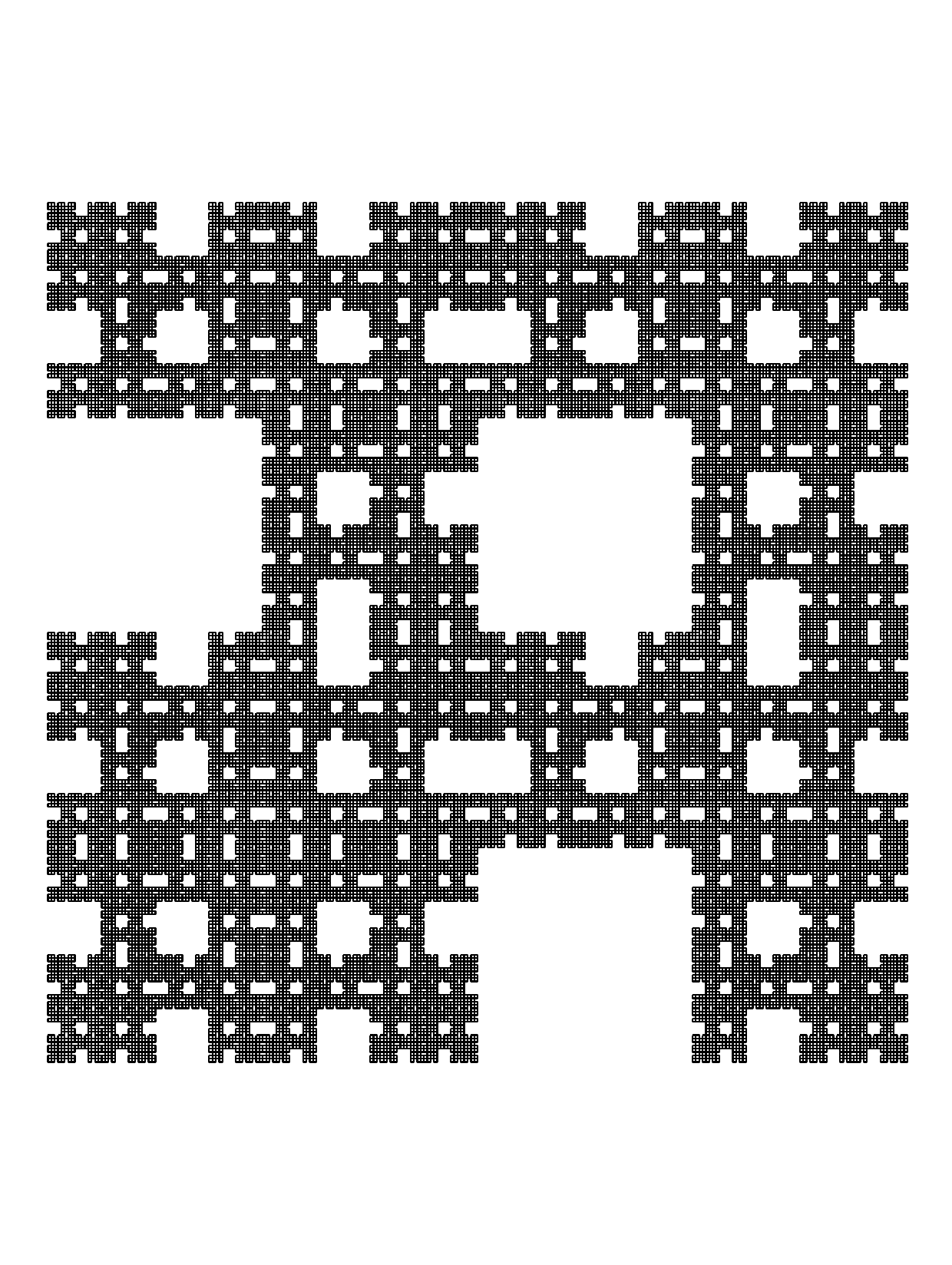}
\vspace{-17pt}
\end{minipage}
\begin{minipage}[b]{6cm}
\hspace{-20pt}
\centering
\includegraphics[width=122pt]{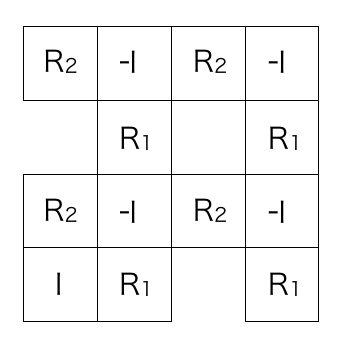}
\vspace{-3pt}
\end{minipage}
\caption{Non-countably ramified example}\label{FIG2}
\end{figure}
Let $L = 2$ and let $N = 4$. As in Example~\ref{CYP.ex10},  $S$ and $\{\Phi_s\}_{s \in S}$ are indicated in the right figure of Figure~\ref{FIG2}.  It is easy to see that the corresponding self-similar structure is non-degenerate, locally symmetric, and strongly connected in the same way as Example~\ref{CYP.ex10}. Moreover, Lemma~\ref{CYP.lemma80} shows that $R_{1, 2}^{1} \in \G_{(K, T)}$, so that (SDR) is satisfied. Thus we have confirmed all the assumptions of Theorem~\ref{CYP.thm100}. Unlike the chipped Sierpinski carpet, this example is not ``countably ramified''. In this example, like the chipped Sierpinski carpet, $K$ contains enough straight lines. This implies that $(K, d_*)$ has the chain condition, so that the heat kernel associated with $(\E, \W^2)$ satisfies \eqref{SUC.eq100} and \eqref{SUC.eq120}
\endexample

\example[Moulin/Pinwheel]\label{CYP.ex30}
\begin{figure}[ht]
\hspace*{0pt}
\begin{minipage}[b]{6cm}
\centering
\includegraphics[width = 114pt]{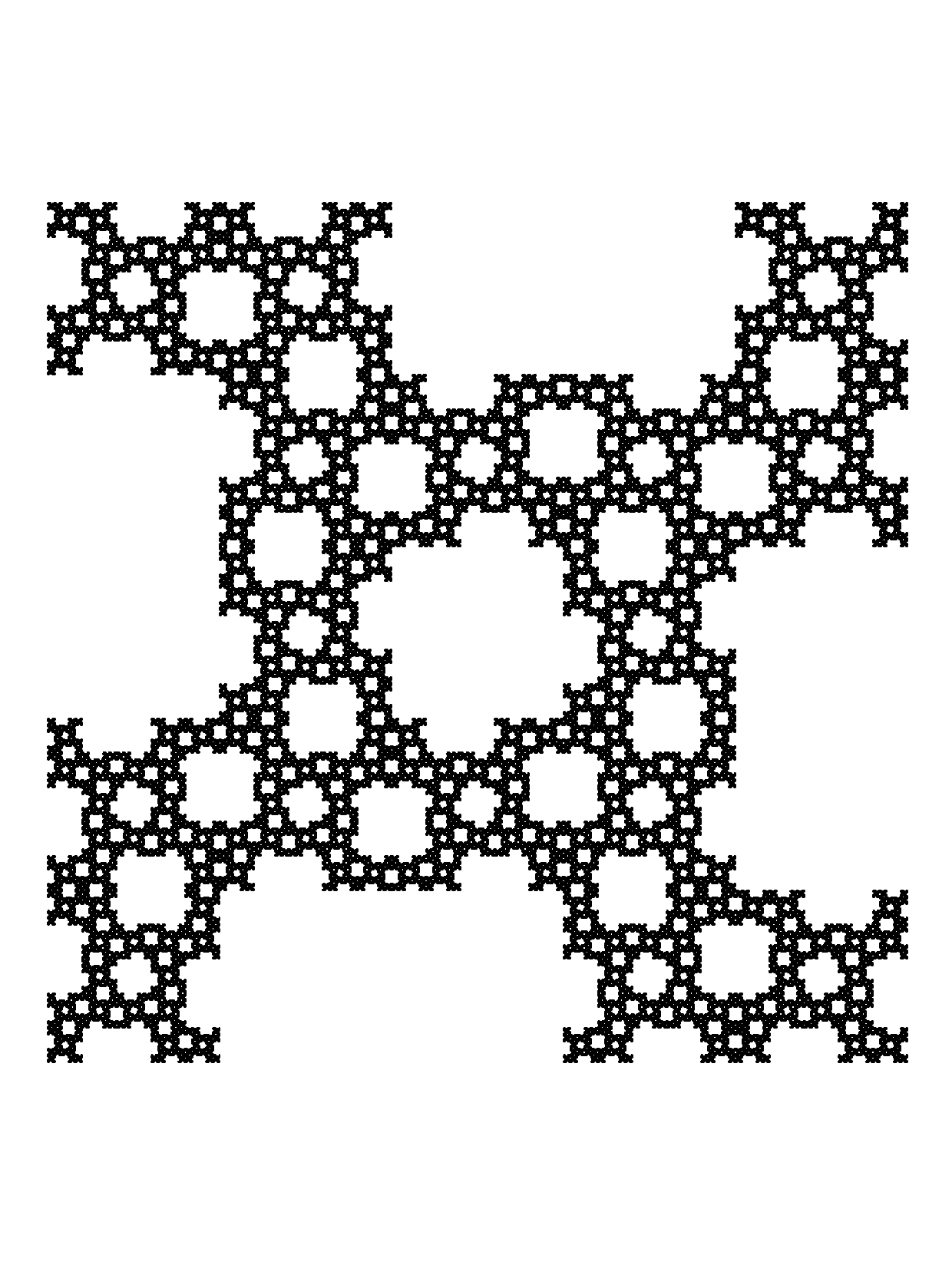}
\vspace{-20pt}
\end{minipage}
\begin{minipage}[b]{6cm}
\hspace{-20pt}
\centering
\includegraphics[width=123pt]{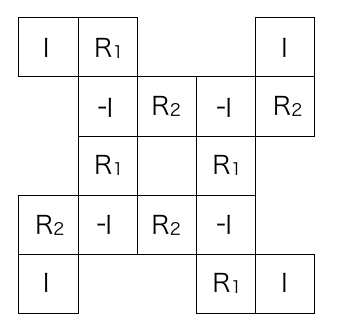}
\vspace{-5pt}
\end{minipage}\caption{Moulin/Pinwheel}\label{FIG3}
\end{figure}

Let $L = 2$ and let $N = 5$. As in the above examples,  $S$ and $\{\Phi_s\}_{s \in S}$ are indicated in the right figure of Figure~\ref{FIG3}. The assumptions of Theorem~\ref{CYP.thm20} are verified in exactly the same way as before including (RS), i.e. $\Theta_{\pi/2} \in \G_{(K, T)}$. In this example, unlike previous ones, $(K, d_*)$ does not have the chain condition and hence we have \eqref{SUC.eq100} and \eqref{SUC.eq110}.
\endexample

The next two examples satisfy (NS).

\begin{figure}
\hspace*{0pt}
\begin{minipage}[b]{6cm}
\centering
\includegraphics[width = 113pt]{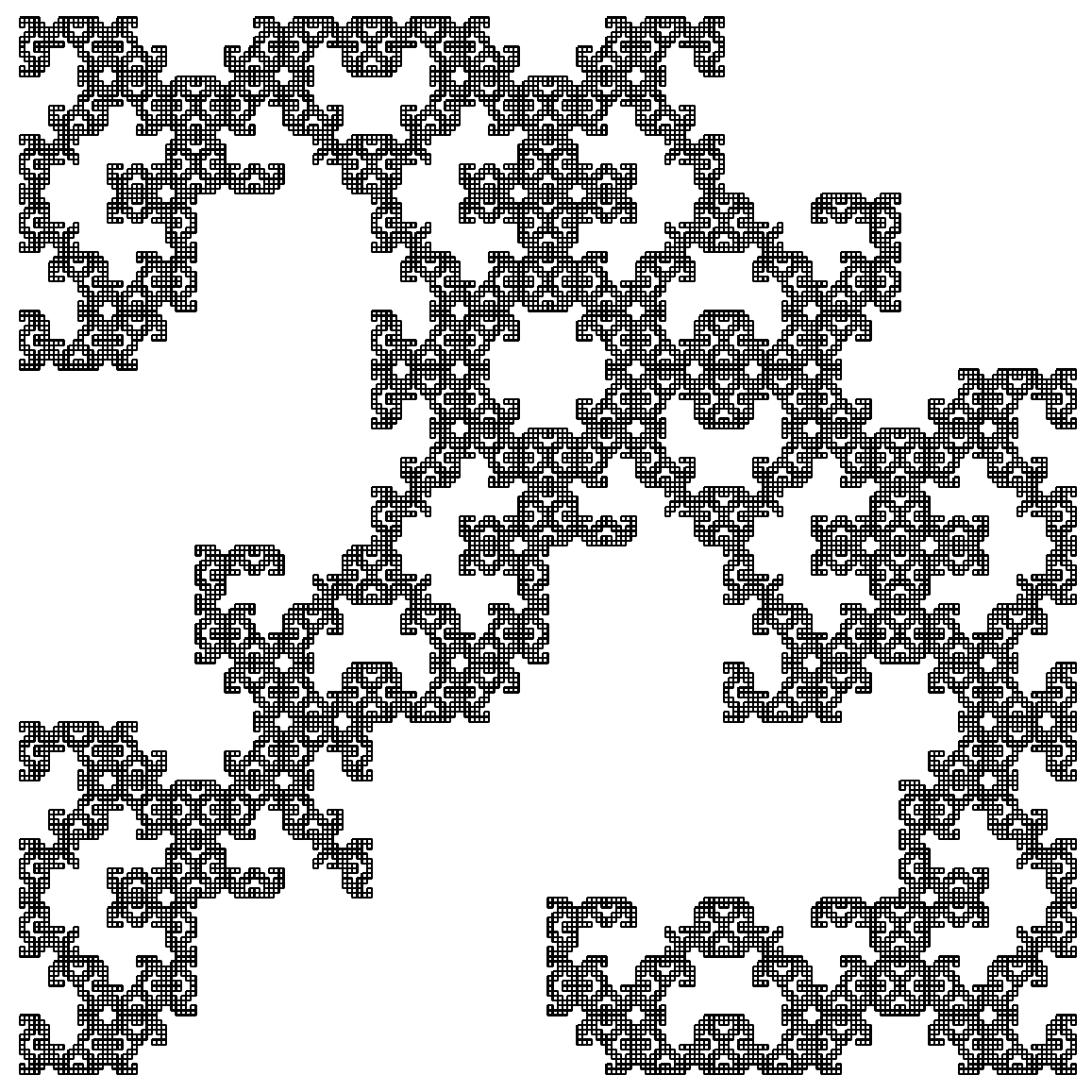}
\vspace{-5pt}
\end{minipage}
\begin{minipage}[b]{6cm}
\hspace{-20pt}
\centering
\includegraphics[width=116pt]{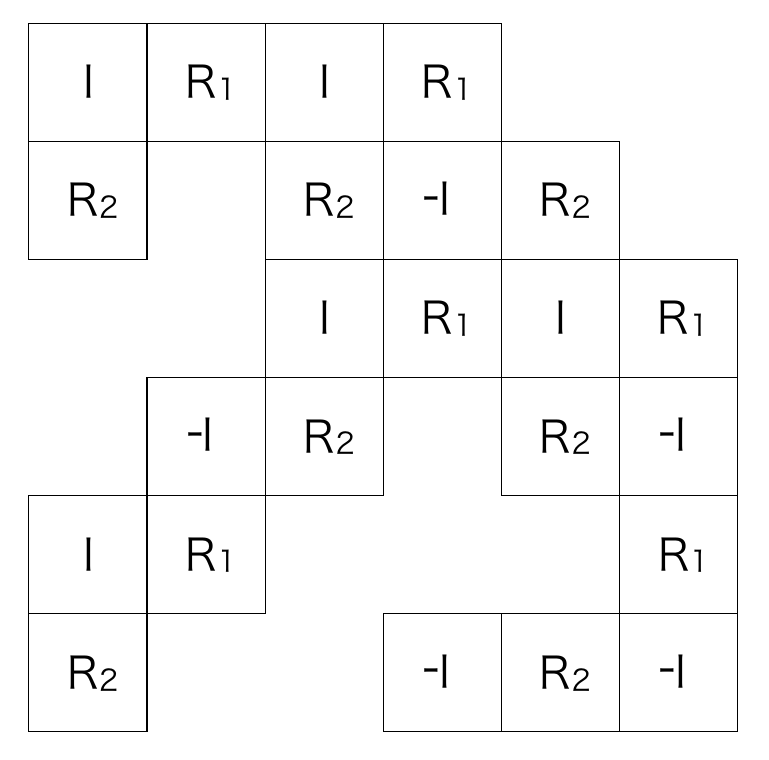}
\vspace{-6pt}
\end{minipage}
\caption{Non-symmetric example 1}\label{FIG4}
\end{figure}

\example\label{CYP.ex40}
Let $L = 2$ and let $N = 6$. As in the previous examples, $S$ and $\{\Phi_s\}_{s \in S}$ are indicated in the right figure of Figure~\ref{FIG4}. In the same manner as before, we verify local symmetry, non-degeneracy and strongly connectedness. The condition (NS) is immediate from the right figure of Figure~\ref{FIG4}. We have $\#(S) = 23$, so that $ \dim_H(K, d_*) = \log{23}/\log 6$. Obviously there is no $\BB_2$-symmetry.

\endexample

\begin{figure}
\hspace*{0pt}
\begin{minipage}[b]{6cm}
\centering
\includegraphics[width = 113pt]{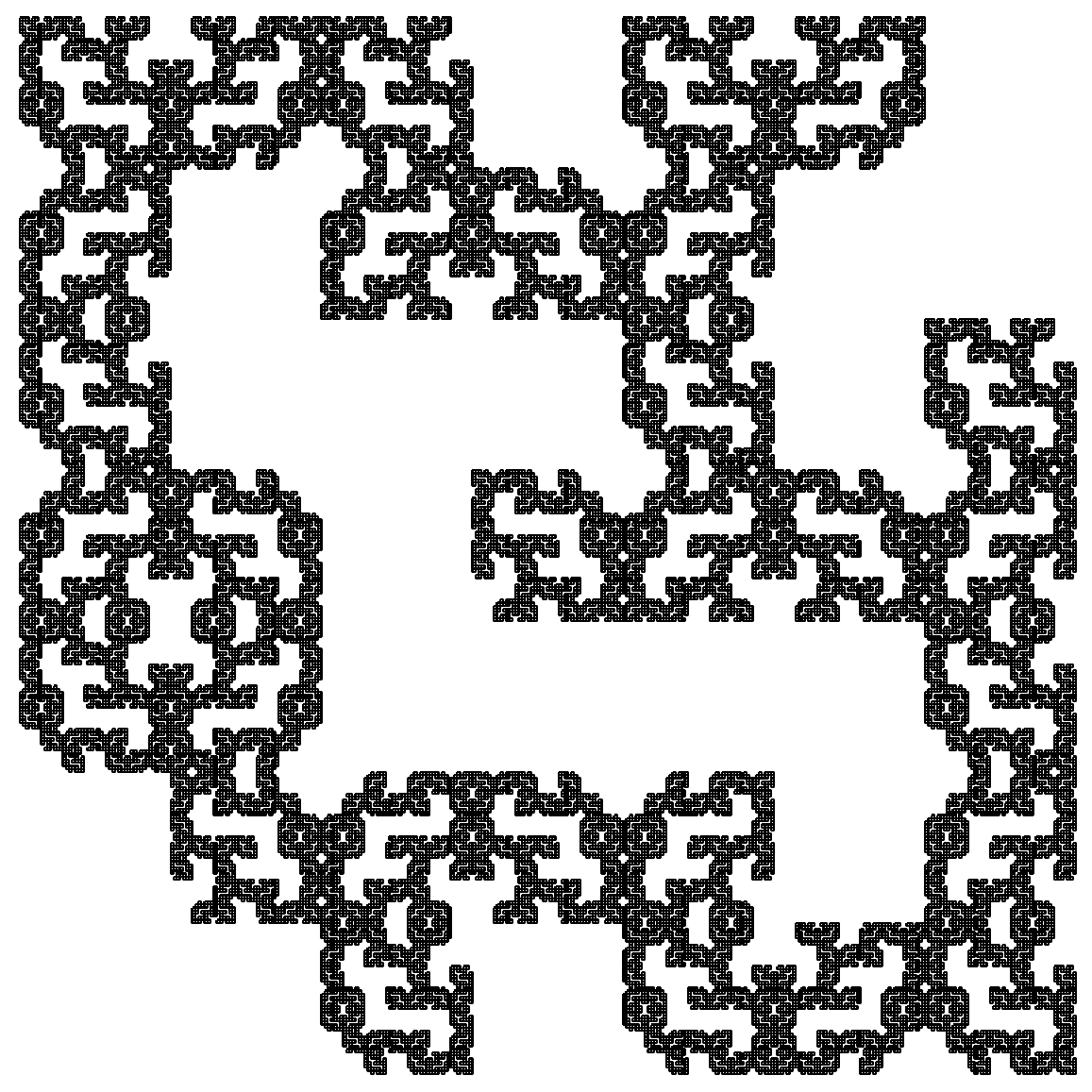}
\vspace{-5pt}
\end{minipage}
\begin{minipage}[b]{6cm}
\hspace{-20pt}
\centering
\includegraphics[width=116pt]{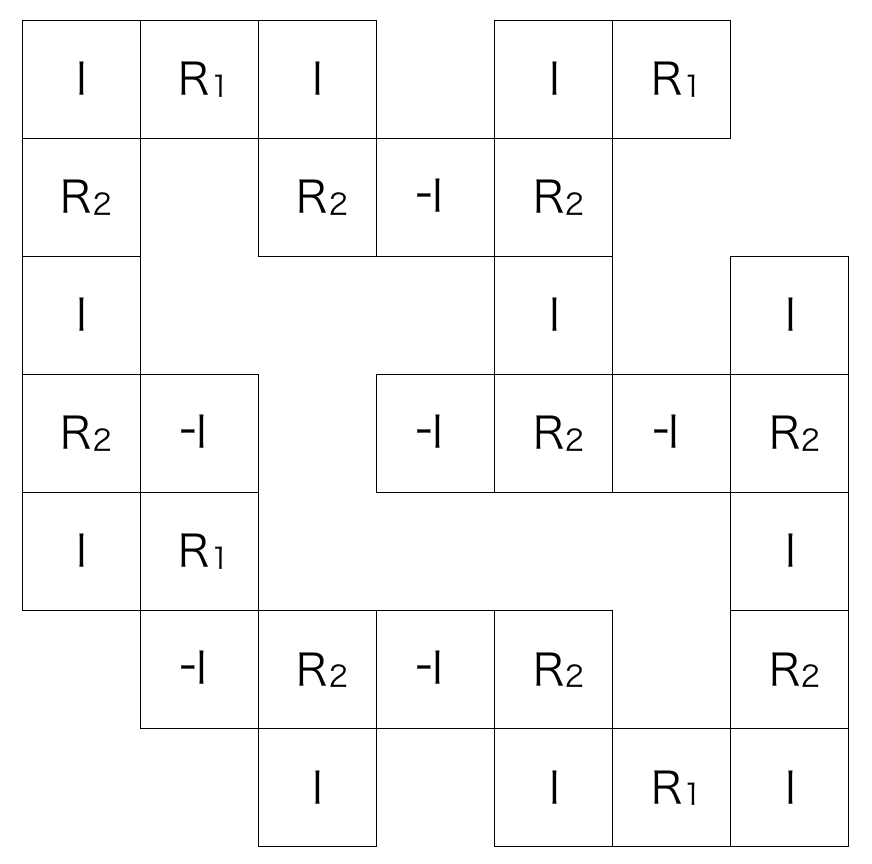}
\vspace{-6pt}
\end{minipage}
\caption{Non-symmetric example 2}\label{FIG5}
\end{figure}

\example\label{CYP.ex50}
Let $L = 2$ and let $N = 7$. As in the previous examples, $S$ and $\{\Phi_s\}_{s \in S}$ are indicated in the right figure of Figure~\ref{FIG5}. In the same manner as before, we verify local symmetry, non-degeneracy and strongly connectedness. The condition (NS) is immediate from the right figure of Figure~\ref{FIG5}. In this example $\#(S) = 30$, so that $\dim_H(K, d_*) = \log{30}/\log 7$. Note that 
\[
\dim_H (K \cap R_{2, 1}) = \frac {\log 5}{\log 7}\,\,\text{while}\,\,\dim_H(K \cap R_{2, -1}) = \frac{\log 4}{\log 7}.
\]

\endexample

In the following examples, we may choose an arbitrary $L \ge 2$.

\example\label{CYP.ex60}
Let $S = \sd{\{1, \ldots, N\}^L}{\{s_*\}}$, where $s_* = (1, \ldots, 1)$.  Also let $\vp = \vp_{s_*, I}$, i.e. $\vp$ is a folding map given by
\[
\vp(x) = NA_s(x - c_s)
\]
for any $s = (s_1, \ldots, s_L) \in \{1, \ldots, N\}^L$ and $x \in Q_s$, where $A_s = \prod_{j = 1}^L (R_j)^{|s_i - 1|}$. Note that $(A_s)^{-1} = A_s$. Define
\[
f_s(x) = \frac 1NA_sx + c_s
\]
and let $K$ be the unique non-empty compact set satisfying
\[
K = \bigcup_{s \in S} f_s(K).
\]
Then $\L = (K, S, \{f_s\}_{s \in S})$ is a self-similar structure. By Corollary~\ref{CYP.cor20}, $\L$ is non-degenerate and locally symmetric. Moreover, Lemma~\ref{CYP.lemma50} shows that $\L$ is strongly connected. Additionally, using Lemma~\ref{CYP.lemma80}, we see that $R_{j_1, j_2}^1 \in \G_{(K, T)}$ for any $j_1, j_2 \in \{1, \ldots, L\}$ with $j_1 \neq j_2$. In fact, $\G_{(K, T)}$ is generated by $\{R_{j_1, j_2}^1| j_1, j_2 \in \{1, \ldots, L\}, j_1 \neq j_2\}$ and it is isomorphic to the symmetric group of order $L$. Hence by Theorem~\ref{CYP.thm100}, $K$ is $p$-conductively homogeneous for any $p > \dim_{AR}(K, d_*)$. Note that $\G_{(K, T)}$ is a proper subgroup of $\BB_L$ in this case.
\endexample

\example[Hypercube]\label{CYP.ex70}
Let $S = \{1, \ldots, N\}^L$ and let $f_s(x) = \frac 1Nx + c_s$ for any $s \in S$ and $x \in [-1, 1]^L$. Set $K = [-1, 1]^L$. Then $(K, S, \{f_s\}_{s \in S})$ is a self-similar structure. Obviously, $\L$ is non-degenerate, strongly connected and locally symmetric. Moreover, $\G_{(K, T)} = \BB_L$. By Theorem~\ref{CYP.thm100} $K$ is $p$-conductively homogeneous for any $p > L$. In fact,  for any $p > L$, we see that $W^{1, p}(K) = \W^p$ and there exist $c  > 0$ such that
\begin{equation}\label{CYP.eq50}
c\E_p(f) \le \int_{K} |\nabla{f}|^pdx \le c^{-1}\E_p(f)
\end{equation}
 for any $f \in W^{1, p}(K)$, where $\E_p$ is the self-similar $p$-energy constructed in Section~\ref{SSF}. The rest of this example is devoted to showing these facts. Choose $A = \{w(1), w(2), w(3)\} \subseteq T_n$ such that $K_{w(1)}, K_{w(2)}$ and $K_{w(3)}$ are three consecutive cubes in $x_1$-direction, i.e. $K_{w(1)} \cap K_{w(2)} = f_{w(1)}(B_{1, 1}) = f_{w(2)}(B_{1, -1})$ and $K_{w(2)} \cap K_{w(3)} = f_{w(2)}(B_{1, 1}) = f_{w(3)}(B_{1, -1})$. Let $A_1 = \{w(1)\}$ and let $A_2 = \{w(3)\}$. Then, the function attaining the infimum in the definition of $\E_{p, m}(A_1, A_2, A)$ depends only on the first variable $x_1$ and is a piecewise linear function in the direction of $x_1$. Consequently, we see that
\[
\E_{p, m}^{\ell}(A_1, A_2, A) \ge 2^{m(L - p) -1}.
\]
On the other hand, the comparison of moduli shows
\[
\M_{p, m}^{(1)}(A_1, A_2, A) \le \M_{1, p, m}^{(1)}(w)
\]
for any $w \in T$. Therefore, there exists $c_2 > 0$ such that
\[
c_22^{m(L - p)} \le \E_{1, p, m}(w, T_{|w|})
\]
for any $m \ge 1$ and $w \in T$. \par
Now, for $f: K \to \BbR$, we define $\tilde{f}_m: T_m \to T$ by $\tilde{f}_m(w) = f(f_w(0))$.  Then there exists $c >0$ such that
\begin{equation}\label{CYP.eq60}
2^{m(p - L)}\E_{p, T_m}(\tilde{f}_m) \to c\int_{K} |\nabla{f}|^pdx.
\end{equation}
as $m \to \infty$ for any $f \in C^{\infty}(K)$. This shows that there exists $c_3 > 0$ such that $\E_{1, p, m}(w, T_{|w|}) \le c_32^{m(L - p)}$ for any $w \in T$. Thus the scaling exponent of $\s$ appearing in \eqref{SUC.eq10} is $2^{L - p}$. Combining this fact and arguments analogous to those in \cite[Section~5.3]{Shimizu2}, we have the following Korevaar-Shoen type expression of $\W^p$:
\[
\W^p = \bigg\{f\bigg| f \in L^p(K, dx), \limsup_{r \downarrow 0}\int_K\frac 1{r^L}\int_{B_{d_*}(x, r)}\frac{|f(x) - f(y)|^p}{r^p}dydx < \infty\bigg\}.
\]
This expressing enable us to identify $\W^p$ with $W^{1, p}(K)$. By \eqref{CYP.eq60}, we see that \eqref{CYP.eq50} holds for any $f \in C^{\infty}(K)$. Since $C^{\infty}(K)$ is dense in $W^{1, p}(K)$, \eqref{CYP.eq50} holds for any $f \in \W^p$.
\endexample

\setcounter{equation}{0}
\section{Rationally ramified Sierpinski crosses}\label{SCR}
\begin{figure}[h]
\hspace*{0pt}
\begin{minipage}[b]{6cm}
\centering
\includegraphics[width = 140pt]{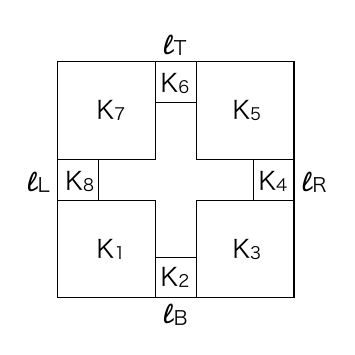}

\end{minipage}
\begin{minipage}[b]{6cm}
\hspace{-40pt}
\centering
\includegraphics[width=100pt]{cr7.pdf}
\vspace{20pt}
\end{minipage}
\caption{the $\rho_*$-Sierpinski cross: $\rho_* = \sqrt 2 - 1$}\label{FIGSCR}
\end{figure}

In this section, we present another class of conductively homogeneous spaces called rationally ramified Sierpinski crosses. This example is a planar square-based self-similar set as those in the last section but the sizes of the squares constituting it are not one but two. See Figure~\ref{FIGSCR}. Consequently, although it has full $\BB_2$-symmetry, we should make a little more complicated discussion than that of the previous section to show the conductive homogeneity.\par
The family of Sierpinski crosses was introduced in \cite[Example~1.7.5]{Ki13}.
\definition\label{SCR.def00}
Let $r_1, r_2 \in (0, 1)$ satisfying $2r_1 + r_2 = 1$ and $r_1 \ge r_2$. Let $p_1 = (-1, -1)$, $p_2 = (0, -1)$, $p_3 = (1, -1)$, $p_4 = (1, 0)$, $p_5 = (1, 1)$, $p_6 = (0, 1)$, $p_7 = (-1, 1)$ and $p_8 = (-1, 0)$. Set $S = \{1, \ldots, 8\}$. For $s \in S$, define $F_s: C_* \to C_*$ as
\[
F_s(x) = \begin{cases}
 r_1(x - p_s) + p_s \quad&\text{if $s$ is odd,}\\
 r_2(x - p_s) + p_s \quad&\text{if $s$ is even.}
\end{cases}
 \]
 The self-similar set $K$ with respect to the family of contractions $\{F_s\}_{s \in S}$ is called the $(r_1)$-Sierpinski cross. Define $\ell_L = \{-1\} \times [-1, 1]$, $\ell_R = \{1\} \times [-1, 1]$, $\ell_B = [-1, 1] \times \{-1\}$, and $\ell_T = [-1, 1] \times \{1\}$, where the symbols, L, R, B, and T correspond to left, right, bottom, and top respectively.
 \enddefinition

In this section, we will show that if an $(r_1)$-Sierpinski cross $K$ is rationally ramified, then it is $p$-conductively homogeneous for any $p > \dim_{AR}(K, d_*)$. Roughly speaking an $(r_1)$-Sierpinski cross is rationally ramified if  $\cup_{v \in \GG_1( w)} K_v$, which represents the local geometry around $w \in T$, has finite types of variety up to the isometries when $w \in T$ varies.  See \cite{Ki13} for the exact definition. In fact, in \cite[Proposition~1.7.6]{Ki13}, it is shown that an $(r_1)$-Sierpinski cross is rationally ramified if and only if $1 - r_1 = (r_1)^m$ for some $m \ge 2$. For simplicity of arguments, we confine ourselves to the case $m = 2$ hereafter in this section. The generalization to other values of $m$ is a little complicated but the essential idea is the same.\par
In the case $m = 2$, the value of $r_1$ equals $\sqrt{2} -1$. Set $\rho_* = \sqrt{2} -1 $. Our main object of study is now the $\rho_*$-Sierpinski cross. We take advantage of the framework of Section~\ref{SSF} with  $r = \rho_*$ and
 \[
 j_s = \begin{cases}
 1 \quad&\text{if $s$ is odd,}\\
 2 \quad&\text{if $s$ is even.}
 \end{cases}
 \]
to define $(T, \A, \phi)$ and the associated  partition of $K$. In this case, $g(w)$ is the contraction ratio of the map $F_w = F_{w_1}\circ\ldots\circ{F_{w_m}}$ for $w = \word wm \in S^m$. Note that $g(w) = (\rho_*)^n$ or $(\rho_*)^{n + 1}$ for any $(n, w) \in T_n$. For example $\LL^g_{\rho_*} = S$ and
\begin{multline*}
\LL^g_{(\rho_*)^2} = \{1s, 3s, 5s, 7s| s \in S, s: \text{even}\} \cup  \{1s, 3s, 5s, 7s| s \in S, s: \text{odd}\}\\ \cup \{2, 4, 6, 8\}.
\end{multline*}
Note that $g(1s) = (\rho_*)^3$ if $s$ is even and $g(1s) = (\rho_*)^2$ if $s$ is odd. Moreover, $\LL^g_{\rho_*} \cap \LL^g_{(\rho_*)^2} \neq \emptyset$ in this case. Let $d_*$ be the restriction of the Euclidean metric to $K$. Let $h_{\rho_*}(n, w) = (\rho_*)^n$ for $(n, w) \in T_w$.  It is straightforward to see that $d_*$ is $1$-adapted to the weight function $h_{\rho_*}$, i.e. Assumption~\ref{ALFR}-(2B) holds with $M_* = 1$.\par
For simplicity, to denote an element in $T_n$, we use $w$ in place of $(n, w)$ hereafter as long as no confusion may occur.\par
The Hausdorff dimension of $(K, d_*)$ is given by the unique number $\a_H$ satisfying
\[
4(\rho_*)^{2\a_H} + 4(\rho_*)^{\a_H} = 1. 
\] 
Consequently, we see that 
\[
\a_H = 1 + \frac{\log 2}{\log{(1 + \sqrt{2})}}.
\]
Let $\mu$ be the self-similar measure with weight $(\mu_i)_{i \in S}$, where
\[
\mu_i = \begin{cases}
(\rho_*)^{\a_H}\quad&\text{if $i$ is odd,}\\
(\rho_*)^{2\a_H}\quad&\text{if $i$ is even.}
\end{cases}
\]
Then $\mu$ is the normalized $\a_H$-dimensional Hausdorff measure and is $\a_H$-Ahlfors regular with respect to $d_*$. After those observations, it is easy to see that Assumption~\ref{ALFR} is satisfied with $M_* = M_0 = 1$, $N_* = 8$. Moreover we see that $L_* \le 8$.  \par
The main result of this section is as follows.

 \thm\label{SCR.thm10}
 For any $p > 0$, $n, m, k \ge 1$, $w \in T_n$ and $u, v \in T_k$,
 \[
 \M_{1, p, m}^{(1)}(w) \le 8(24)^{p + 1}\#(T_{k + 1})^p\M_{p, m}^{(1)}(u, v, T_k).
 \]
  \endthm
  
  An immediate consequence of the above theorem is the conductive homogeneity of the Sierpinski cross.
  
\cor\label{SCR.cor10}
The $\rho_*$-Sierpinski cross $K$ is $p$-conductively homogeneous for any $p > \dim_{AR}(K, d_*)$. Moreover, there exists a self-similar $p$-energy $\E_p$ on $\W^p$. In particular, there exists a local regular Dirichlet form $(\E, \W^2)$ on $L^2(K, \mu)$ whose associated heat kernel satisfies \eqref{SUC.eq100} and \eqref{SUC.eq120}.
\endcor
 Note that due to the two different values of $j_s$, the self-similarity of the $p$-energy $\E_p$ is given as
\[
\E_p(f) = \s\sum_{s:\text{odd}}\E_p(f{\circ}F_s) + \s^2\sum_{s:\text{even}} \E_p(f{\circ}F_s)
\]
for any $f \in \W^p$.
\demo[Proof of Corollary~\ref{SCR.cor10}]
By \eqref{CMP.eq30}, it follows that
\[
\E_{1, p, m}(w, T_n) \le c_p\#(T_{k + 1})^p\E_{p, m}(u, v, T_k)
\]
for any $n, m, k \ge 1$, $w \in T_n$ and $u, v \in T_k$.  Moreover since $p > \dim_{AR}(K, d_*)$, there exist $c > 0$ and $\a \in (0, 1)$ such that
\[
\E_{1, p, m} \le c\a^m
\]
for any $m \ge 1$. Thus we have obtained \eqref{SUC.eq30} and \eqref{SUC.eq40}, so that $K$ is $p$-conductively homogeneous by Theorem~\ref{SUC.thm20}. In particular, since $\a_H < 2$, $K$ is $2$-conductively homogeneous and we have $(\E, \W^2)$. Since $(K, d_*)$ has the chain condition, we have \eqref{SUC.eq100} and \eqref{SUC.eq120} by Theorem~\ref{SUC.thm30}.
\enddemo
  
  To show Theorem~\ref{SCR.thm10}, we need to prepare several notions.
  
 \definition\label{SCR.def10}
(1)\,\, Set
\[
U = \{(2, 13), (2, 31), (4, 35), (4, 53), (6, 57), (6, 75), (8, 17), (8, 71)\}.
\]
For $(i, jk) \in U$, define $R_{i, jk}: K_i \to K_{jk}$ as the reflection in the line segment $K_i \cap K_{jk}$. Moreover, define $R^*_{i, jk}(w)$ for $w \in T(i) \cup T(jk)$ as  the unique $v \in T(i) \cup T(jk)$ satisfying $R_{i, jk}(K_w) = K_v$.  $R_{i, jk}^*$ is a map from $T(i) \cup T(jk)$ to itself.\\
(2)\,\,For $g \in \BB_2$, define $g^*:T \to T$ by
\[
g^*(w) = v,
\]
where $v$ is the  the unique $v \in T$ satisfying $g(K_w) = K_v$. Note that $g^*|_{T_n}: T_n \to T_n$. \\
(3)\,\,
For $w \in T$, if $w \notin T(2) \cup T(4) \cup T(6) \cup T(8)$, then define
\[
\H_w = \{g^*(v)| g \in \BB_2\}.
\]
Otherwise if $w \in T(i)$ for $i = 2, 4, 6, 8$, then define
\[
\H_w = \{g^*(v)| g \in \BB_2\} \cup \{g_*(R^*_{i, jk}(v))| g \in \BB_2, (i, jk) \in U\}.
\]
 \enddefinition
 
 Note that $\#(\H_w) \le 24$ for any $w \in T_n$.\par
 By the construction of $T_n$, we see that $g(w) = (\rho_*)^n$ or $g(w) = (\rho_*)^{n + 1}$ for any $w \in T_n^n$. In fact, we immediately obtain the following lemma.
 
 \lemma\label{SCR.lemma05}
 Set
 \[
 T_n^n = \{w | w \in T_n, g(w) =(\rho_*)^n\}\quad\text{and}\quad T_n^{n + 1} = \{w | w \in T_n, g(w) = (\rho_*)^{n + 1}\}.
 \]
 {\rm (1)}\,\,For any $w \in T^n_n$, $wv \in T_{n + m}$ if and only if $v \in T_m$.\\
 {\rm (2)}\,\,For any $w \in T_n^{n + 1}$, $wv \in T_{n + m}$ if and only if $v \in T_{m - 1}$.\\
 {\rm (3)}\,\,$w \in T_{n + 1}^{n + 1}$ if and only if $w \in T_n^{n + 1}$ or $w = \tau{j}$ for some $\tau \in T_n^n$ and $j \in \{1, 3, 5, 7\}$. \\
 {\rm (4)}\,\,$w \in T_{n + 1}^{n + 2}$ if and only if $w = \tau{j}$ for some $\tau \in T_n^n$ and $j \in \{2, 4, 6, 8\}$.
 \endlemma
 
 \definition\label{SCR.def20}
(1)\,\, 
Define $\psi^*_{n,m}: S^m(T_n^n) \to T_m$ by 
\[
\psi^*_{n, m}(wv) = v
\]
for $w \in T_n^n$ and $v \in T_m$.\\
(2)\,\,For $w \in T$, define $\H^0_w \subseteq T$ by
\[
\H^0_w = \begin{cases}
\{w, R^*_{i, jk}(w)\} \quad&\text{if $w \in T(jk)$ for some $(i, jk) \in U$,}\\
\{w\}\quad&\text{otherwise}.
\end{cases}
\]
For $w \in T_{n + 1}^{n + 1}$ and $u \in T$, define
\[
\H_{wu}^n = \begin{cases}
\{\tau{v}| v \in \H_{ju}^0\}&\text{if $w = \tau{j}$ for some $\tau \in T_n^n$ and $j \in \{1, 3, 5, 7\}$},\\
\{wu\} &\text{if $w \in T_n^{n + 1}$.}
\end{cases}
\]
(3)\,\,Define 
 \begin{equation}\label{SCR.eq20}
 K_{\%} = \bigcup_{s \in S, K_s \cap \ell_{\%} \neq \emptyset} K_s
 \end{equation}
 for $\% \in \{T, B, R, L\}$. For example, $K_B = K_1 \cup K_2 \cup K_3$.
\enddefinition

Note that if $w \in T_n$, then $\H^0_w \in T_n$ and that if $w \in T_{n + 1}^{n + 1}$ and $u \in T_{m - 1}$, then $\H_{wu}^n \subseteq T_{n + m}$.

\lemma\label{SCR.lemma10}
Assume that there exists a path $\bp = (w(1), \ldots, w(l))$ of $T_{m - 1}$ contained in $K_L$ such that $K_{w(1)} \cap \ell_B \neq \emptyset$, $K_{w(l)} \cap \ell_T \neq \emptyset$, and $\bp$ is $R_2^*$-invariant. Set 
\[
\H^*_u = \bigcup_{w \in T_{k + 1}^{k + 1}}\bigcup_{v \in \H_u} \H_{wv}^{k + 1}
\]
for $u \in T_{m - 1}$. Then for any $u_1, u_2 \in T_k$, there exists $\bp_0 \in \C_m^{(1)}(\{u_1\}, \{u_2\}, T_k)$ such that
\begin{equation}\label{SCR.eq25}
\bp_0 \subseteq \bigcup_{i = 1}^l \H^*_{w(i)}.
\end{equation}
\endlemma

\remark
Strictly, $\bp_0$ is not a subset but a sequence of points. However, in \eqref{SCR.eq25}, we use $\bp_0$ to denote a subset consisting of the points in the sequence. We use such abuse of notations if no confusion may occur.
\endremark

\begin{figure}
\vspace{-10pt}
\centering
\includegraphics[width= \linewidth]{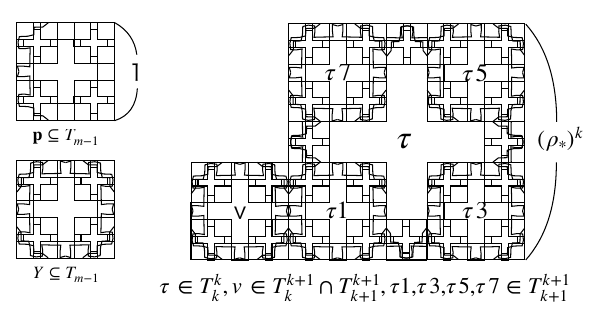}
\vspace{-25pt}
\caption{$\bp$, $Y$ and a part of $\H^*(Y)$}\label{SCRc}
\end{figure}

\demo
Set
\[
Y = \bp \cup \Theta_{\pi/2}^*(\bp) \cup \Theta_{\pi}^*(\bp) \cup \Theta_{3\pi/2}^*(\bp).
\]
Then $Y = g^*(Y)$ for any $g \in \BB_2$. Let 
\[
\H^*(Y) = \bigcup_{w \in T^{k + 1}_{k + 1}}\bigcup_{v \in Y} \H^{k + 1}_{wv}.
\]
See Figure~\ref{SCRc} for an illustration of $\bp$, $Y$ and a part of $\H^*(Y)$. It follows that $K(\H^*(Y))$ is a connected set intersecting $K_u$ for any $u \in T_k$. Therefore, we can choose a path $\bp_0$ connecting $K_{u_1}$ and $K_{u_2}$ from $\H^*(Y)$. Since $\H^(Y) \subseteq \cup_{i = 1}^l \H^*_{w(i)}$, we have the desired statement.
\enddemo

\begin{figure}
\centering
\includegraphics[width= \linewidth]{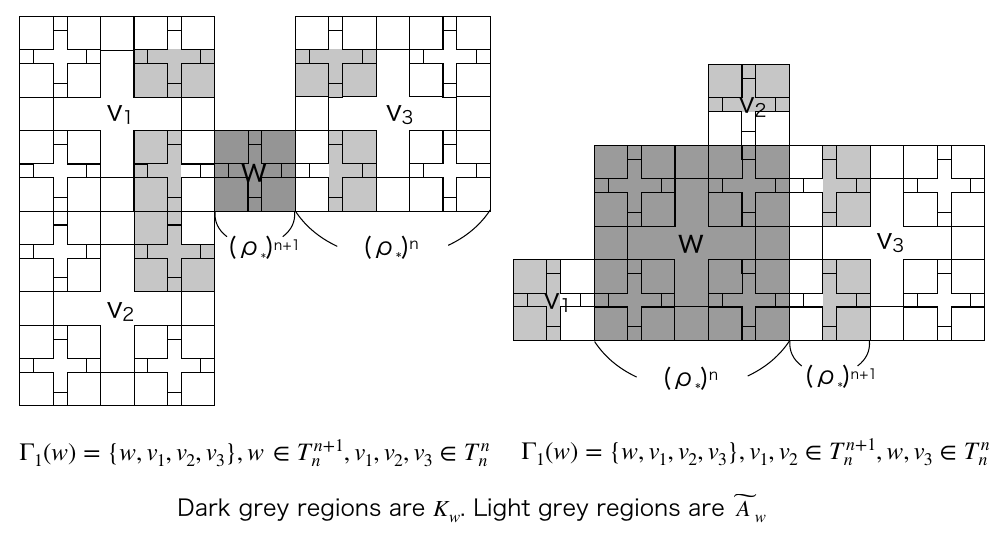}
\vspace{-25pt}
\caption{Two examples of $\widetilde{A}_w$}\label{SCrnbd}
\end{figure}

 \demo[Proof of Theorem~\ref{SCR.thm10}]
 Let $w \in T_n$ and let $u_1, u_2 \in T_k$.  For any $\bp \in \C^{(1)}_{1, m}(w)$, set
 \[
\H_{m - 1}(\bp) = \bigcup_{u \in \psi_{n + 1, m - 1}^*(\bp \cap S^{m - 1}(T_{n + 1}^{n + 1}))} \H_u.
 \]
Then $\H_{m - 1}(\bp) \subseteq T_{m - 1}$ and $g^*(\H_{m - 1}(\bp)) = \H_{m - 1}(\bp)$ for any $g \in \BB_2$. \\
 {\bf Claim 1}\,\,There exists a path $\bp^*$ contained in $\H_{m - 1}(\bp)$ such that one of the following four statements is true:\\
(a)\,\,$K(\bp^*) \cap \ell_B \neq \emptyset$ and $K(\bp^*) \cap K_T \neq \emptyset$,\\
(b)\,\,$K(\bp^*) \cap \ell_T \neq \emptyset$ and $K(\bp^*) \cap K_B \neq \emptyset$,\\
(c)\,\,$K(\bp^*) \cap \ell_L \neq \emptyset$ and $K(\bp^*) \cap K_R \neq \emptyset$,\\
(d)\,\,$K(\bp^*) \cap \ell_R \neq \emptyset$ and $K(\bp^*) \cap K_L \neq \emptyset$.\\
{\bf  Proof of Claim 1}:
Let $F_w(C_*) = [a, a + h] \times [b, b + h]$, where $h = (\rho_*)^n$ if $w \in T_n^n$ and $h = (\rho_*)^{n + 1}$ if $w \in T_n^{n + 1}$. Define 
\[
A_{w, \c} =   [a - \c, a + h + \c] \times [b - \c, b + h +\c]
\]
and $\widetilde{A}_w = K \cap (A_{w, (\rho_*)^{n + 1}}\backslash{A_{w, (\rho_*)^{n + 2}}})$. Two typical examples of $\widetilde{A}_w$ is illustrated in Figure~\ref{SCrnbd}. Since $K_{w(1)} \cap  K_w \neq \emptyset$ and $K_{w(l)} \cap A_{w, (\rho_*)^{n + 1}} = \emptyset$, a part of $\bp$ contained in $\widetilde{A}_w$ connects\\
$\{(a - (\rho_*)^{n + 1}, y)| y \in [-1, 1]\}$ and $\{(a -  (\rho_*)^{n + 2}, y)| y \in [-1, 1]\}$,\\
$\{(a + h, y +  (\rho_*)^{n + 2})| y \in [-1, 1]\}$ and $\{(a + h + (\rho_*)^{n + 1}, y)| y \in [-1, 1]\}$,\\
$\{(x, b - (\rho_*)^{n + 1})| x \in [-1, 1]\}$ and $\{(x, b -  (\rho_*)^{n + 2})| x \in [-1, 1]\}$,\\ or\\
$\{(x, b + h +  (\rho_*)^{n + 2})| x \in [-1, 1]\}$ and $\{x, b + h + (\rho_*)^{n + 1})| x \in [-1, 1]\}$.\\
According to the four possibilities above, we have (a), (b), (c) or (d), where the exact correspondence depends on $w$.\qed\par\vspace{20pt}
\begin{figure}
\vspace{-30pt}
\centering
\includegraphics[width= \linewidth]{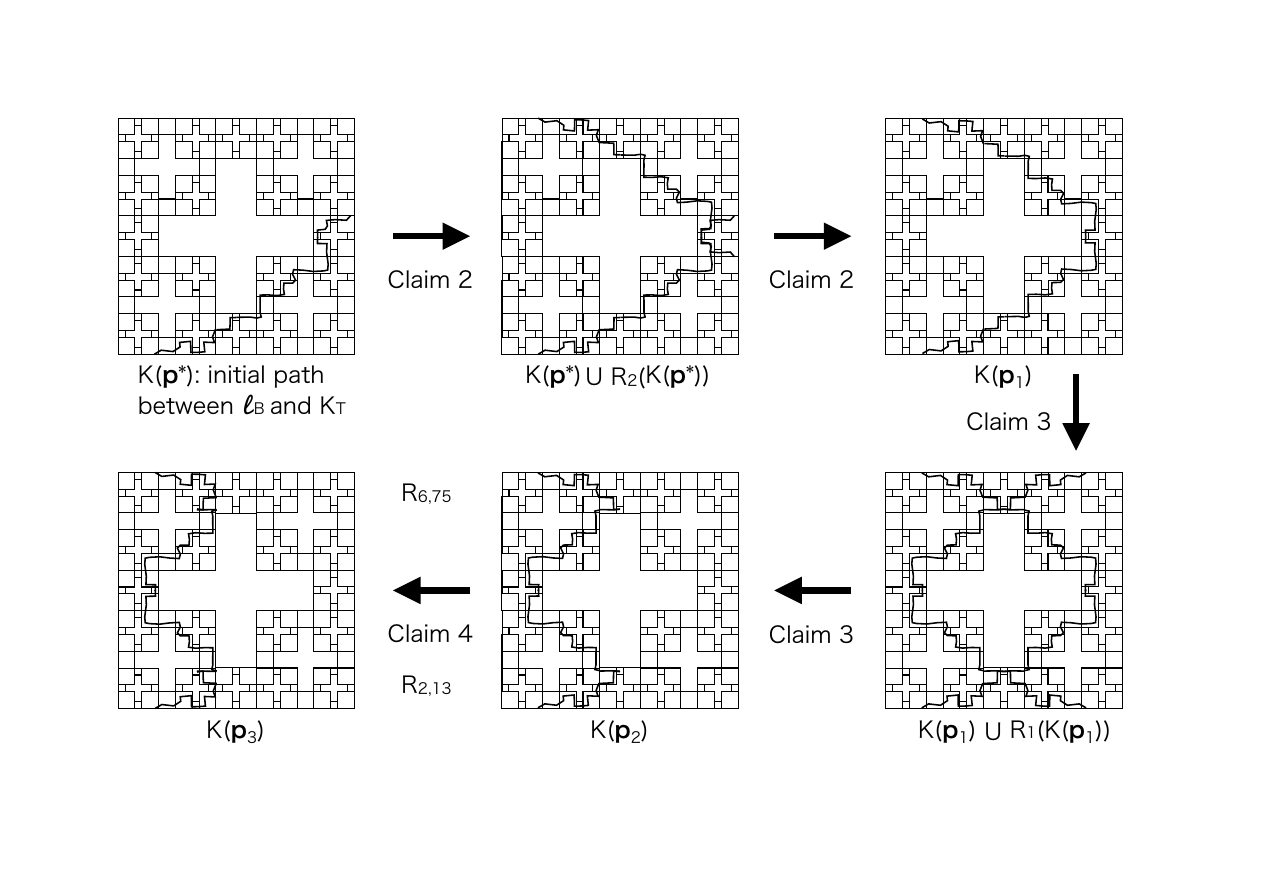}
\vspace{-50pt}
\caption{Modifications of a path}\label{SCrpathm}
\end{figure}

\vspace{-20pt}
 Hereafter we assume the case (a) in Claim 1 in the course of discussion. Other cases may be treated exactly in the same manner. In the following claims, we are going to modify the initial path $\bp^*$ step by step. This process of modification is illustrated in Figure~\ref{SCrpathm}.\\
 
\vspace{-10pt}
\noindent {\bf Claim 2}\,\,
 $\bp^* \cup R_2^*(\bp^*)$ contains an $R_2$-symmetric path $\bp_1 = (v(0), \ldots, v(l_1))$ between $\ell_B$ and $\ell_T$, i.e. $K_{v(0)} \cap \ell_B \neq \emptyset$, $R_2^*(v(i)) = v(l_1 - i)$ for $i = 1, \ldots, l_1$.\\
 {\bf Proof of Claim 2}:\,\,Let $p_* = (w(1), \ldots, w(l))$. By (a), $K(\bp^*)$ intersects with the line segment $[-1, 1] \times \{0\}$. Set $i_* = \min\{i | w(i) \cap [-1, 1] \times \{0\} \neq \emptyset\}$. Then connecting $(w(1), \ldots, w(i_*))$ and its image by $R_2^*$, we obtain a desired path.\qed\\
{\bf Claim 3}\,\,$R_1^*(\bp_1) \cup \bp_1$ contains an $R_2$-symmetric path $\bp_2$ such that $K(\bp_2) \subseteq [-1, 0] \times [-1, 1]$.\\
{\bf Proof of Claim 3}:\,\,If $\bp_1$ or $R_1^*(\bp_1)$ is contained in the left half of $C_*$, then choose $\bp_1$ or $R_1^*(\bp_1)$ accordingly as our path. Otherwise, applying $R_1$ to $K(\bp_1) \cap [0, 1] \times [-1, 1]$, we obtain a desired path. \qed\\
{\bf Claim 4}\,\, Set $\H_{\bp^*} = \cup_{u \in \bp^*} \H_u$. Then there exists an $R_2^*$-symmetric path $\bp_3 \subseteq \H_{\bp^*}$ contained in $K_L$ such that $K(\bp_3) \cap \ell_T \neq \emptyset$ and $K(\bp_3) \cap \ell_B \neq \emptyset$.\\
{\bf Proof of Claim 4}:\,\,If $K(\bp_2) \subseteq K_L$, then we set $\bp_2 = \bp_3$. Otherwise, use $R_{2, 13}^*$ (resp. $R_{6, 75}^*$) to reflect the part $K(\bp_2) \cap K_2$ (resp. $K(\bp_2) \cap K_6$) into $K_{13}$ (resp. $K_{75}$). Then we obtain a desired path. \qed \par
\noindent Now we have a path $p_3$ satisfying all the assumptions of Lemma~\ref{SCR.lemma10}. Applying Lemma~\ref{SCR.lemma10} with $\bp = \bp_3$, we obtain a path $\bp_0 \in \C_m^{(1)}(\{u_1\}, \{u_2\}, T_k)$. For $u \in S^m(\GG_1(w))$, define
\[
H_u = \begin{cases}
\displaystyle\bigcup_{v \in \H_{\psi^*_{n + 1, m - 1}(u)}} \H_v^*\quad&\text{if $u \in S^{m - 1}(T_{n + 1}^{n + 1})$,}\\
\displaystyle\quad\quad\quad\quad\emptyset\quad&\text{otherwise.}
\end{cases}
\]
Then it follows that
\[
\bp_0 \subseteq \bigcup_{v \in \bp} H_{v}.
\]
Since $\#(\H_u) \le 24$ and $\#(\GG_1(w)) \le 8$, 
\[
\#(H_w) \le 48\#(T_{k + 1})\quad \text{and} \quad \#(\{v| u \in H_v\}) \le 24\cdot8.
\]
So, Lemma~\ref{UFD.lemma50} suffices.
\enddemo

\setcounter{equation}{0}
\section{Nested fractals}\label{OEX}
In this section, we show conductive homogeneity of a class of self-similar sets, called strongly symmetric self-similar sets, that are highly symmetric and finitely ramified. This class is a natural extension of nested fractals introduced by Lindstr{\o}m\cite{Li1}, where Brownian motions were constructed on them. In \cite[Section~3.8]{AOF}, Lindsr{\o}m's results were extended to strongly symmetric self-similar sets. Typical examples of strongly symmetric self-similar sets are the Sierpinski gasket, the pentakun (``Kun'' means ``Mr.'' in Japanese), and the snowflake, whose definitions are given below.\par

Let $\rho \in (0, 1)$ and let $S$ be a finite subset of $\BbR^L$ for some $L  \in \BbN$. For each $q \in S$, let $f_q: \BbR^L \to \BbR^L$ be a $\rho$-similitude whose fixed point is $q$, i.e. there exists $U_q \in O(L)$ such that 
\[
f_q(x) = \rho{U_q}(x - q) + q
\]
for any $x \in \BbR^L$. Let $K$ be the self-similar set with respect to the family of contractions $\{f_q\}_{q \in S}$.  Then the triple $(K, S, \{f_q\}_{q \in S})$ is a self-similar structure as is explained in Section~\ref{SSF}.

\assumption\label{OEX.ass10}
(1)\,\,If $p, q \in S$ and $p \neq q$, then $p \notin f_q(K)$. \\
(2)\,\,There exists $U \subseteq S$ such that
\[
\bigcup_{\substack{ q_1, q_2 \in S\\ q_1 \neq q_2}} f_{q_1}^{-1}(f_{q_1}(K) \cap f_{q_2}(K)) = U.
\]
(3)\,\,
$K$ is connected.
\endassumption
For purpose of normalization, we assume
\[
\sum_{q \in U} q = 0
\]
hereafter. 
\prop\label{OEX.prop10}
Under Assumption~\ref{OEX.ass10}, $(K, S, \{f_q\}_{q \in S})$ is a p.c.f. self-similar structure with 
\begin{equation}\label{OEX.eq40}
V_0 = U.
\end{equation}
Moreover, define $\{V_m\}_{m \ge 1}$ inductively by
\[
V_{m + 1} = \bigcup_{i \in S} f_i(V_m).
\]
Then
\begin{equation}\label{OEX.eq50}
V_m \subseteq V_{m + 1}
\end{equation}
for any $m \ge 0$.
\endprop

The definitions of p.c.f. self-similar structures and $V_0$ along with the proof of \eqref{OEX.eq40} is given in Appendix~\ref{SNS}. \eqref{OEX.eq50} is due to \cite[Lemma~1.3.11]{AOF}.\\
For the self-similar structure $(K, S, \{f_q\}_{q \in S})$, we adopt the framework in Section~\ref{SSF} with $r = \rho$ and $j_q = 1$ for any $q \in S$. In this case,
\[
T_m = S^m = \{\word wm| w_i \in S\,\,\text{for any $i = 1, \ldots, m$}\}.
\]
Then we see that
\[
V_0 = \bigcup_{e \in E_1^*} X(e).
\]
Moreover, by \cite[Proposition~1,3,5-(2)]{AOF}, it follows that
\begin{equation}\label{OEX.eq100}
K_w \cap K_v = f_w(V_0) \cap f_v(V_0) \subseteq V_m
\end{equation}
for any $w, v \in T_m$ with $w \neq v$. This implies that
\begin{equation}\label{OEX.eq110}
V_0 = \bigcup_{(X, Y, \vp) \in \IT(K, T)} X.
\end{equation}
Let $\a_H = -\log N/\log{\rho}$. Note that $N\rho^{\a_H} = 1$. Let $\mu$ be the self-similar measure with weight $(\rho^{\a_H}, \ldots, \rho^{\a_H})$. Basic properties of $\mu$ is given in Appendix~\ref{SNS}. Also, let $d_*$ be the restriction of the Euclidean metric to $K$.  \par
The following assumption is an equivalent condition of Assumption~\ref{ALFR}-(2B) when $d$ is the (restriction of) Euclidean metric. Essentially the same assumptions have been around from time to time for almost 30 years.  See \cite[Assumption~2.2]{Kum3} and \cite[Assumption (P)]{PPStos}. The assumption is believed to be true for nested fractals but we have no proof so far. In \cite{PPStos}, it was shown that this assumption is true if $U_q$ is the same for any $q \in S$. In Appendix~\ref{SNS}, we show this assumption is true if $U_q$ is the identity map for any $q \in V_0$.

\assumption\label{OEX.ass30}
There exists $c > 0$ such that $d(K_w, K_v) \ge c\rho^{|w|}$ for any $n \ge 1$,  and $(w, v) \in E_n^*$, where $d(A, B) = \inf_{x \in A, y \in B} |x - y|$ for subsets $A, B \subseteq \BbR^L$.
\endassumption

\prop\label{OEX.prop20}
Under Assumptions~\ref{OEX.ass10} and \ref{OEX.ass30}, Assumption~\ref{ALFR} is satisfied with $d = d_*$, $r = \rho$, and  $M_* = M_0 = 1$.
\endprop

The above proposition is proven in Appendix~\ref{SNS}.

\definition\label{OEX.def10}
(1)\,\,Let $m_* = \#\{|x - y|\, |\,x, y \in V_0, x \neq y\}$, where $|x|$ is the Euclidean length of $x \in \BbR^L$. Define 
\[
l_0 = \min\{|x - y|\, |\,x, y \in V_0, x \neq y\}.
\]
 Moreover, define $l_i$ for $i = 0, 1, \ldots m_* - 1$ inductively by 
 \[
 l_{i + 1} = \min\{|x - y|\, |\,x, y \in V_0, x \neq y, |x - y| > l_i\}.
 \]
(2)\,\,
A sequence $(x_i)_{i = 1, \ldots k} \subseteq V_m$ is called an $m$-walk if there exists $w(i) \in T_m$ such that $x_i, x_{i + 1} \in f_{w(i)}(V_0)$ for any $i = 1, \ldots, k - 1$.\\
(3)\,\,
A $0$-walk $(x_i)_{i = 1, \ldots, k}$ is called a strict $0$-walk (between $x_1$ and $x_k$) if $|x_i - x_{i + 1}| = l_0$ for any $i = 1, \ldots, k - 1$.\\
(4)\,\,Define
\begin{multline*}
\G = \{g | g \in O(L), g(V_0) = V_0\\
\text{there exists $g*: T \to T$ such that $g(f_w(V_0)) = f_{g^*(w)}(V_0)$ for any $w \in T$.}\}
\end{multline*}
(5)\,\,
For any $x, y \in \BbR^L$ with $x \neq y$, define 
\[
H_{xy} = \{z| z \in \BbR^L, |x - z| = |y - z|.\}
\]
($H_{xy}$ is the hyperplane bisecting the line segment $xy$.) Also let $g_{xy}: \BbR^L \to \BbR^L$ be reflection in $H_{xy}$.
\enddefinition

\definition\label{OEX.ass20}
$(K, S, \{f_q\}_{q \in S})$ is said to be strongly symmetric if Assumption~\ref{OEX.ass10} is satisfied and there exists a finite subgroup $\G_*$ of $\G$ such that the following properties hold:\\
(1)\,\,
For any $x, y \in V_0$ with $x \neq y$, there exists a strict $0$-walk between $x$ and $y$.\\
(2)\,\,If $x, y, z \in V_0$ and $|x - y| = |x - z|$, then there exists $g \in \G_*$ such that $g(x) = x$ and $g(y) = z$.\\
(3)\,\,For any $i = 1, \ldots, m_* - 2$, there exist $x, y$ and $z \in V_0$ such that $|x - y| = l_i$, $|x - z| = l_{i + 1}$ and $g_{yz} \in \G_*$.\\
(4)\,\,$V_0$ is $\G_*$-transitive, i.e. for any $x, y \in V_0$, there exists $g \in \G_*$ such that $g(x) = y$.
\enddefinition

\remark
By Definition~\ref{OEX.ass20}-(4), $|q_1| = |q_2|$ for any $q_1, q_2 \in V_0$.
\endremark

\definition\label{OEX.def20}
A self-similar structure $(K, S, \{f_q\}_{q \in S})$ is called a nested fractal if Assumption~\ref{OEX.ass10} holds and $g_{xy} \in \G$ for any $x, y \in V_0$ with $x \neq y$.
\enddefinition

By \cite[Proposition~3.8.7]{AOF}, we have the following proposition.

\prop\label{OEX.prop30}
A nested fractal is strongly symmetric.
\endprop

We give three examples of strongly symmetric self-similar sets. Note that Assumption~\ref{OEX.ass30} is satisfied for all the three examples because of Lemma~\ref{SNS.lemma30}. The first two are nested fractals.

\begin{figure}[h]
\hspace*{0pt}

\begin{minipage}[b]{6cm}

\centering
\includegraphics[width = 100pt]{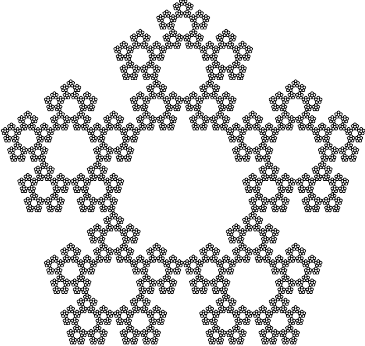}
\caption{Pentakun}\label{Penta}
\end{minipage}
\begin{minipage}[b]{6cm}
\hspace{-10pt}
\centering
\includegraphics[width=100pt]{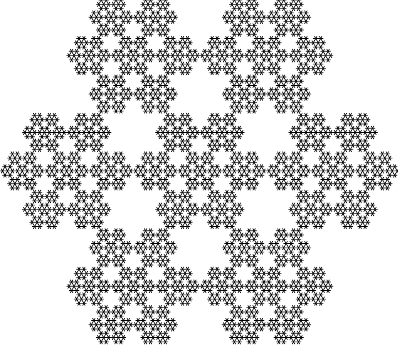}
\caption{Snowflake}\label{Snow}
\end{minipage}
\end{figure}

\example\label{OEX.ex10}[Pentakun: Figure~\ref{Penta}]
Let $S = \{p_1, \ldots, p_5\}$ be a collection of vertices of a regular pentagon satisfying $\sum_{i = 1}^5 p_i = 0$ and let $\rho = \frac{3 - \sqrt 5}2$. Then the associated self-similar set $K$, called pentakun, is strongly symmetric. (See \cite[Example~3.8.11]{AOF}.) In this case $\G = \G_* = D_5$, which is the group of symmetries of a regular pentagon, and $V_0 = \{p_1, \ldots, p_5\}$. 

\endexample

\example\label{OEX.ex20}[Snowflake: Figure~\ref{Snow}]
Let $\{p_1, \ldots, p_6\}$ be a collection of vertices of a regular hexagon satisfying $\sum_{i = 1}^6 p_i = 0$ and let $S = \{p_1, \ldots, p_7, 0\}$. Furthermore let $\rho = \frac 13$. Then the associated self-similar set, called snowflake, is strongly symmetric. (See \cite[Example~3.8.12]{AOF}.) In this case $\G = \G_* = D_6$, which is the group of symmetries of a regular hexagon and $V_0 = \{p_1, \ldots, p_6\}$. 
\endexample

The last example is not a nested fractal.

\example\label{OEX.ex30}
Let
\[
S = \{-1, 0, 1\}^3 \cup \Big\{-\frac 12, \frac 12\Big\}^3,\,\,U = \{1, -1\}^3,
\]
and $\rho = 1/5$.  Note that $U$ is the collection of vertices of the cube $[-1, 1]^3$ and 
\[
f_{q}([-1, 1]^3) = \Big[\frac{4q_1 - 1}{5}, \frac{4q_1 +1}{5}\Big] \times \Big[\frac{4q_2 -1}{5}, \frac{4q_2 +1}{5}\Big] \times \Big[\frac{4q_3 - 1}{5}, \frac{4q_3 +1}{5}\Big].
\]
for any $q = (q_1, q_2, q_3) \in S$.
It is straight forward to see that the associated self-similar set is strongly symmetric with $V_0 = U$ and $\G = \G_* = \BB_3$. This self-similar set is not a nested fractal because $g_{xy} \notin \G$ if $x = (-1, -1, -1)$ and $y = (1, 1, 1)$.
\endexample

Using Theorem~\ref{CHS.thm10}, we have the following theorem.

\thm\label{OEX.thm10}
 Suppose that $(K, S, \{f_i\}_{i \in S})$ is strongly symmetric and that Assumption~\ref{OEX.ass30} holds. Then $(K, d_*)$ is $p$-conductively homogeneous for any $p > \dim_{AR}(K, d_*)$.
\endthm

As for $\dim_{AR}(K, d_*)$, it was shown in \cite{TysonWu1} that $\dim_{AR}(K, d_*) = 1$ if $(K, d_*)$ is the Sierpinski gasket. In general, we have the following fact.

\prop\label{OEX.prop40}
Suppose that $(K, S, \{f_i\}_{i \in S})$ is strongly symmetric and that Assumption~\ref{OEX.ass30} holds. Then $\dim_{AR}(K, d_*) < 2$.
\endprop

\demo
For $m \ge 0$, define $\widetilde{E}_m = \{(f_w(x), f_w(y))| w \in T_m, x, y \in V_0,x \neq y\}$. Then $\{(V_m, \widetilde{E}_m)\}_{m \ge 0}$ is a proper system of horizontal networks in the sense of \cite[Definition~4.6.5]{GAMS}.  Define
\begin{multline*}
\L_s(V_0) = \{(D_{xy})_{x, y \in V_0}| \text{there exists $(D_0, \ldots, D_{m_* - 1}) \in [0, \infty)^{m_*}$ such that}\\
\text{$D_0 = 1$, $D_{xy} = D_i$ if $|x - y| = l_i$, and $\sum_{y \in V_0} D_{xy} = 0$ for any $x \in V_0$.}\}
\end{multline*}
In particular, let $D^1 \in \L_s(V_0)$ satisfy $(D^1)_{xy} = 1$ for any $x, y \in V_0$ with $x \neq y$.
For $D = (D_{xy})_{x, y \in V_0} \in \L_s(V_0)$, define
\[
\E_{2, m}^{D}(f) = \frac 12\sum_{w \in T_m, x, y \in V_0} D_{xy}(f(f_w(x)) - f(f_w(y)))^2
\]
for $f \in \ell(V_m)$ and
\[
\E_{2, m, w}^D = \inf\{\E_{2, n + m}^D(f)| f \in \ell(V_{n + m}), f|_{V_{n + m} \cap K_w} = 1, f|_{V_{n + m} \cap (\cup_{v \notin \GG_1(w)} K_v)} = 0\}
\]
for any $w \in T_n$. Then by \cite[Theorem~3.8.10 and Corollary~3.1.9]{AOF}, there exist $D_* \in \L_s(V_0)$ and $\s > 1$ such that $(D_*, (\s^{-1}, \ldots, \s^{-1}))$ is a harmonic structure, that is, for any $f \in \ell(V_m)$,
\[
\s^m\E_{2, m}^{D_*}(f) = \min\{\s^{m + 1}\E_{2, m + 1}^{D_*}(g)| g \in \ell(V_{m + 1}), g|_{V_m} = f\}.
\]
This implies that there exist $c_1, c_2 > 0$ and $k \ge 1$ such that
\[
c_1\s^{-m} \le \sup_{w \in \sd{T}{T_k}} \E_{2, m, w}^{D_*} \le c_2\s^{-m}.
\]
On the other hand, there exist $c_3, c_4 > 0$ such that
\[
c_3\E_{2, m}^{D_*}(f) \le \E_{2, m}^{D^1}(f)  \le c_4E_{2, m}^{D_*}(f)
\]
for any $m \ge 0$ and $f \in \ell(V_m)$. Thus we see that $\sup_{w \in T}\E_{2, m, w}^{D_1} \le C\s^{-m}$ for any $m \ge 0$. Therefore, by \cite[Theorems~4.6.9 and 4.9.1]{GAMS}, it follows that $\dim_{AR}(K, d_*) < 2$.
\enddemo

The rest of this section is devoted to proving Theorem~\ref{OEX.thm10}. We suppose that $(K, S, \{f_i\}_{i \in S})$ is strongly symmetric hereafter in this section. By \cite[Proposition~3.8.19]{AOF}, we have the following theorem.

\lemma\label{OEX.lemma05}
If $(K, S, \{f_i\}_{i \in S})$ is strongly symmetric, then $g(K_w) = K_{g^*(w)}$ for any $g \in \G$ and $w \in T$. In particular, $\G \subseteq \G_{(K, T)}$.
\endlemma

\lemma\label{OEX.lemma10}
If $(K, S, \{f_i\}_{i \in S})$ is strongly symmetric,  $x_1, x_2, y_1, y_2 \in V_0$ and $|x_1 - x_2| = |y_1 - y_2|$, then there exists $g \in \G_*$ such that $g(x_1) = y_1$ and $g(x_2) = y_2$.
\endlemma

\demo
By Definition~\ref{OEX.ass20}-(4), there exists $g_1 \in \G$ such that $g_1(x_1) = y_1$. Let $g_1(x_2) = z$. Then $|y_1 - y_2| = |y_1 - z|$.  Hence by Definition~\ref{OEX.ass20}-(2), there exists $g_2 \in \G$ such that $g_2(y_1) = y_1$ and $g_2(z) = y_2$. Thus letting $g = g_2{\circ}g_1$, we see that $g(x_1) = g_2(y_1) = y_1$ and $g(x_2) = g_2(z) = y_2$.
\enddemo

\definition\label{OEX.def30}
A path $(w(1), \ldots, w(k))$ of $(T_m, E_m^*)$ is said to connect $x \in K$ and $y \in K$ if $x \in K_{w(1)}$ and $y \in K_{w(k)}$.
\enddefinition

\lemma\label{OEX.lemma20}
Let $\bp$ be a path of $(T_m, E_m^*)$ connecting $x_1 \in V_0$ and $x_2 \in V_0$. Suppose $|x_1 - x_2| = l_i$ for some $i = 1, \ldots, m_* -1 $. Then there exist a path $\bp_1$ of  $(T_m, E_m^*)$, $x \in V_0$ and $y \in V_0$ such that $\bp_1$ connects $x$ and $y$, $\bp_1 \subseteq \cup_{g \in \G_*} g^*(\bp)$ and $|x - y| = l_{i - 1}$.
\endlemma


 \notation
 For a path $\bp = (w(1), \ldots, w(k))$ and $g \in \G$, set
 \[
 g^*(\bp) = (g^*(w(1)), \ldots, g^*(w(k))).
 \]
 \endnotation
\demo
By Definition~\ref{OEX.ass20}-(2), there exist $x, y, z \in V_0$ such that $|x - y| = l_{i - 1}$, $|x - z| = l_i$ and $g_{yz} \in \G_*$. Also, Lemma~\ref{OEX.lemma10} shows that there exists $h \in \G_*$ such that $h(x_1) = x$ and $h(x_2) = z$. Since $|x - y| < |x - z|$, $x$ and $z$ belong to different sides of $H_{yz}$. Hence the path $h^*(\bp)$ intersects with $H_{yz}$. Therefore, $h^*(\bp)$ and $(g_{yz})^*\circ{h^*}(\bp)$ has an intersection at $H_{yz}$. Since $(g_{yz})^*\circ{h^*}(\bp)$ connects $g_{yz}(x)$ and $y = g_{yz}(z)$, we can extract a path $\bp_1$ from $h^*(\bp) \cup (g_{yz})^*\circ{h^*}(\bp)$ connecting $x$ and $y$, and included in $\cup_{g \in \G_*} g^*(\bp)$. Since $|x - y| = l_{i - 1}$, $\bp_1$ is a desired path.
\enddemo

\lemma\label{OEX.lemma30}
Let $\bp$ be a path of $(T_m, E_m^*)$ connecting two distinct points in $V_0$. Then for any $x, y \in V_0$, there exists a path $\bp'$ of  $(T_m, E_m^*)$ connecting $x$ and $y$ such that $\bp' \subseteq \cup_{g \in \G_*} g^*(\bp)$.
\endlemma

\demo
Using Lemma~\ref{OEX.lemma20} inductively, we see that there exists a path $\bp_0$ of $(T_m, E_m^*)$ connecting two distinct points $z_1$ and $z_2$ in $V_0$ such that $|z_1 - z_2| = l_0$ and $\bp_0 \subseteq \cup_{g \in \G_*} g^*(\bp)$. By Definition~\ref{OEX.ass20}-(1), there exists a strict $0$-walk $(x_1, \ldots, x_{j_0})$ satisfying $x_1 = x$ and $x_{j_0} = y$. By Lemma~\ref{OEX.lemma10}, for any $j = 1, \ldots, j_0 - 1$, there exists $g_j \in \G_*$ such that $g_j(z_1) = x_j$ and $g_j(z_2) = x_{j + 1}$. Concatenating $(g_1)^*(\bp_0), \ldots, (g_{j_0 - 2})^*(\bp_0)$ and $(g_{j_0 - 1})^*(\bp_0)$, we obtain a desired path connecting $x$ and $y$.
\enddemo

\demo[Proof of Theorem~\ref{OEX.thm10}]
We are going to use Theorem~\ref{CHS.thm10}. Let $\I = \IT(K, T)$ and let $\G_0 = \G_1 = \G_*$. By \eqref{OEX.eq110} and the fact that $\I = \IT(K, T)$, we see that $E_m^{\I} = E_m^*$. Hence the condition (a) of Theorem~\ref{CHS.thm10} is satisfied. The condition (b) is also satisfied due to the fact that $\G_*$ is transition on $V_0$.\par
Let $w \in T_n$, let $u, v \in T_k$ and let $\bp \in \C^{(1)}_{1, m}(w)$.  Then $\bp$ contains a path connecting two distinct points in $\cup_{w' \in T_n} f_{w'}(V_0)$. Thus $\psi_n(\bp)$ contains a path between two distinct points in $V_0$. By Lemma~\ref{OEX.lemma30}, for any $x, y \in V_0$, there exists a path $\bp_{xy} \subseteq \cup_{g \in \G_*} g^*(\psi_n(\bp))$ connecting $x$ and $y$. Set $\U_{\bp} = \cup_{x, y \in V_0} \bp_{xy}$. Then since $K(\U_{\bp}) \supseteq V_0$, it follows that $g(K(\U_{\vp})) \supseteq V_0$ for any $g \in \G_*$. Moreover, $K(\U_{\bp})$ is connected and $\U_{\bp} \subseteq \cup_{g \in \G_*} g^*(\psi_n(\bp))$. Thus we have verified the condition (c) of Theorem~\ref{CHS.thm10}. Now, Theorem~\ref{CHS.thm10} suffices.
\enddemo

\setcounter{equation}{0}
\section{Conductance and Poincar{\'e} constants}\label{CPO}

From this section, we start preparations for a proof of Theorem~\ref{SUC.thm20}. To begin with, we will introduce Poincar{\'e} constants and study a relationship between Poincar{\'e} and conductance constants in this section.\par
The next lemma concerns an extension of functions on $T_n$ to those on $T_{n + m}$ by means of the partition of unity $\{\vp_w\}_{w \in T_n}$ given in Lemma~\ref{MAC.lemma15}.

\lemma[{\cite[(2.8) Lemma]{KusZ1}}]\label{MAC.lemma20}
Let $p \ge 1$, let $A \subseteq T_n$ and let $\{\vp_w\}_{w \in A}$ be the partition of unity given in Lemma~\ref{MAC.lemma15}.  Define $\hat{I}_{A, m}: \ell(A) \to \ell(S^m(A))$ by 
\[
(\hat{I}_{A, m}f)(u) = \sum_{w \in A} f(w)\vp_w(u).
\]
Then
\[
\E_{p, A}^{n + m}(\hat{I}_{A, m}f) \le c_{\ref{MAC.lemma20}}\Big(\max_{w \in A} \E_{M, p, m}(w, A)\Big)\E_{p, A}^n(f),
\]
where the constant $c_{\ref{MAC.lemma20}} = c_{\ref{MAC.lemma20}}(p, L_*, M)$ depends only on $p$, $L_*$ and $M$.
\endlemma

\def\tf{\tilde{f}} 
\demo
Let $(a_k(u, v))_{u, v \in T_k}$ be the adjacency matrix of $(T_k, E^*_k)$. Set $\tilde{f} = \hat{I}_{A, m}f$. Then
\begin{equation}\label{FRW.eq10}
\E_{p}^{n + m}(\tilde{f}) = \frac 12\sum_{w \in A}\sum_{v \in S^m(w)}\sum_{u \in S^m(\GG_1^A(w))} a_{n + m}(u, v)|\tilde{f}(u) - \tilde{f}(v)|^p.
\end{equation}
Suppose $v \in S^m(w), u \in S^m(\GG_1^A(w))$ and $(u, v) \in E_{n + m}^*$. Then $\vp_{w'}(u) = \vp_{w'}(v) = 0$ for any $w' \notin \GG_{M + 1}^A(w)$. Hence
\[
\sum_{w' \in \GG_{M + 1}^A(w)} \vp_{w'}(u) = \sum_{w' \in \GG_{M + 1}^A(w)} \vp_{w'}(v) = 1.
\]
Using this, we see
\begin{multline*}
\tf(u) - \tf(v) = \sum_{w' \in \GG_{M + 1}^A(w)} f(w')(\vp_{w'}(u) - \vp_{w'}(v)) \\
= \sum_{w' \in \GG_{M + 1}^A(w)} (f(w') - f(w))(\vp_{w'}(u) - \vp_{w'}(v)).
\end{multline*}
Let $q \ge 1$ be the conjugate of $p$, i.e. $\frac 1p + \frac 1q = 1$. Then by Lemma~\ref{BAS.lemma20}
\begin{multline*}
|\tf(u) - \tf(v)|^p \le \sum_{w' \in \GG_{M + 1}^A(w)}|f(w') - f(w)|^p\Big(\sum_{w' \in \GG_{M + 1}^A(w)} |\vp_{w'}(u) - \vp_{w'}(v)|^{q}\Big)^{p/q}\\
\le C_1\sum_{w' \in \GG_{M + 1}^A(w)}|f(w') - f(w)|^p\sum_{w' \in \GG_{M + 1}^A(w)} |\vp_{w'}(u) - \vp_{w'}(v)|^{p},
\end{multline*}
where $C_1 = \max\{1, (L_*)^{(M + 1)(p - 2)}\}$. If $w \in A$ and $w' \in \GG_{M + 1}^A(w)$, then there exist $w(0), \ldots, w(M + 1) \in A$ such that $w(0) = w$, $w(M + 1) = w'$, $(w(j), w(j + 1)) \in E_n^*$ for  any $j = 0, \ldots, M$. Then
\[
|f(w') - f(w)|^p \le (M + 1)^{p - 1}\sum_{j = 0}^{M} |f(w(j)) - f(w(j + 1))|^p.
\]
Since $\#(\GG_{M + 1}^A(w)) \le (L_*)^{M + 1}$, it follows that
\[
\sum_{w' \in \GG_M^A(w)} |f(w') - f(w)|^p \le C_2\sum_{w', w'' \in \GG_M^A(w), (w', w'') \in E_n^*} |f(w') - f(w'')|^p,
\]
where $C_2 = (M + 1)^{p - 1}(L_*)^{M}$. On the other hand,
\begin{multline*}
\sum_{v \in S^m(w)}\sum_{u \in S^m(\GG_1^A(w))} a_{n + m}(u, v)\sum_{w' \in \GG_{M + 1}^A(w)} |\vp_{w'}(u) - \vp_{w'}(v)|^{p} \\
\le 2\sum_{w' \in \GG_{M + 1}^A(w)} \E_{p, S^m(A)}^{ n + m}(\vp_{w'}, \vp_{w'}) \le 2(L_*)^{M + 1}\max_{w' \in A} \E_{p, S^m(A)}^{n + m}(\vp_{w'}).
\end{multline*}
Hence, by \eqref{FRW.eq10}, 
\begin{multline*}
\E_{p, S^m(A)}^{m + n}(\tf) \le C_1C_2(L_*)^{M + 1}\max_{w \in A} \E_{p, S^m(A)}^{n + m}(\vp_{w}) \times \\\sum_{w \in A}\Big(\sum_{w', w'' \in \GG_{M + 1}^A(w), (w', w'') \in E_n^*} |f(w') - f(w'')|^p\Big) \\
\le C_1C_2(L_*)^{2(M + 1)}\max_{w \in A} \E_{p, S^m(A)}^{n + m}(\vp_{w})\E_{p, A}^{ n}(f).
\end{multline*}
So, Lemma~\ref{MAC.lemma15} suffices.
\enddemo

There is another simple way of extension of functions on $T_n$ to those on $T_{n + k}$.

\lemma\label{MAC.lemma40}
Let $p \ge 1$ and let $A \subseteq T_n$. Define $\wI_{A, k}: \ell(A) \to \ell(S^k(A))$ by
\[
\wI_{A, k}f = \sum_{w \in A} f(w)\chi_{S^k(w)}.
\]
Then
\[
\E_{p, S^k(A)}^{n + k}(\wI_{A, k}f) \le \max_{w \in A} \#(\partial{S^k(w)})\E_{p, A}^n(f).
\]
\endlemma

\demo
Let $\hat{f} = \wI_{A, k}f$. Then $\hat{f}(u) = \hat{f}(v)$ if $\pi^k(u) = \pi^k(v)$. So if $(u, v) \in E_{n + k}^*$ and $\hat{f}(u) \neq \hat{f}(v)$, then $(\pi^k(u), \pi^k(v)) \in E_n^*$. Fix $(w, w') \in E_n^*$. Then
\[
\#\{(u, v)| (u, v) \in E_{n + k}, \pi^k(u) = w, \pi^k(v) = w'\} \le \#(\partial{S^k(w)}).
\]
This immediately implies the desired statement.
\enddemo

Combining two previous extensions, we have the following estimate.

 \lemma[{\cite[(2.9) Lemma]{KusZ1}}]\label{MAC.lemma50}
Let $p \ge 1$ and let $A \subseteq T_n$. Then, there exists $I_{A, k, m}: \ell(A) \to \ell(S^{ k + m}(A))$ such that for any $f \in \ell(A)$,
\begin{multline}\label{MAC.eq10}
\E_{p, S^{k + m}(A)}^{n + k + m}(I_{A, k, m}f)\\
 \le c_{\ref{MAC.lemma50}}\max_{w \in A} \#(\partial{S^k(w)})\max_{v \in S^k(A)}\E_{M, p, m}(v, S^k(A))\E_{p, A}^n(f),
\end{multline}
where the constant $c_{\ref{MAC.lemma50}} = c_{\ref{MAC.lemma50}}(p, L_*, M)$ depends only on $p$, $L_*$ and $M$, and
\begin{equation}\label{MAC.eq20}
(I_{A, k, m}f)(u) = f(w)
\end{equation}
for any $w \in A$ and $u \in S^m(\sd{S^k(w)}{B_{M, k}(w)})$.
\endlemma

\demo
Define $I = \hat{I}_{S^k(A), m}\circ\wI_{A, k}$. Combining Lemmas~\ref{MAC.lemma20} and \ref{MAC.lemma40}, we immediately obtain \eqref{MAC.eq10}.  Let $u \in S^{m + k}(A)$. Set $v = \pi^m(u)$ and $w = \pi^k(w)$.  If $\GG_M^{S^k(A)}(v) \subseteq S^k(w)$, then
\begin{multline*}
(If)(u) = \sum_{v' \in S^k(A)} f(\pi^k(v'))\vp_{v'}(u) = \sum_{v' \in \GG_M^{S^k(A)}(v)} f(\pi^k(v'))\vp_{v'}(u)\\ = \sum_{v' \in \GG_M^{S^k(A)}(v)} f(w)\vp_{v'}(u) = f(w).
\end{multline*}
If $v \in \sd{S^k(w)}{B_{M, k}(w)}$, then $\GG_M^{S^k(A)}(v) \subseteq \GG_M(v) \subseteq S^k(w)$. So the above equality suffices for \eqref{MAC.eq20}.
\enddemo

Next we introduce $p$-Poincar{\'e} constants. In fact, there are two kinds of Poincar{\'e} constants $\lambda_{p, m}(A)$ and $\tla_{p, m}(A)$ but they are almost the same in view of \eqref{MAC.eq100}.

\definition\label{MAC.def120}
Define $\mu(w) = \mu(K_w)$ for $w \in T$.  For $A \subseteq T_n$, define $\mu(A) = \sum_{w \in A} \mu(w)$ and $\mu_A: A \to [0, \infty)$ by
\[
\mu_A(w) = \frac{\mu(w)}{\mu(A)}
\]
for $w \in A$. For $f \in \ell(A)$, define
\[
(f)_{A} = \sum_{u \in A} f(u)\mu_A(u)
\]
and
\[
\norm{f}_{p, \mu_A} = \Big(\sum_{u \in A} |f(u)|^p\mu_A(u)\Big)^{\frac 1p}.
\]
Moreover, define
\[
\lambda_{p, m}(A) = \sup_{f \in \ell(S^m(A))}\frac{\inf_{c \in \BbR}(\norm{f - c\chi_{S^m(A)}}_{p, \mu_{S^m(A)}})^p}{\E_{p, S^m(A)}^{n + m}(f)}
\]
and
\[
\tla_{p, m}(A) = \sup_{f \in \ell(S^m(A))}\frac{(\norm{f - (f)_{S^m(A)}}_{p, \mu_{S^m(A)}})^p}{\E_{p, S^m(A)}^{n + m}(f)}.
\]
\enddefinition

\remark
By Lemma~\ref{FPE.lemma05}, it follows that
\begin{equation}\label{MAC.eq100}
\Big(\frac12\Big)^p\,\tla_{p, m}(A) \le \lambda_{p, m}(A) \le \tla_{p, m}(A).
\end{equation}
\endremark

Using the previous lemmas, we have a relation between Poincar{\'e} and conductance constants as follows.

\lemma[{\cite[(2.10) Prop.]{KusZ1}}]\label{MAC.lemma60}
Let $p \ge 1$ and let $A \subseteq T_n$. For any $m \ge 1$ and $k \ge Mm_0$, 
\[
\max_{w \in A} \#(\partial{S^k(w)})\max_{v \in S^k(A)}\E_{M, p, m}(v, S^k(A))\lambda_{p, k + m}(A) \ge c_{\ref{MAC.lemma60}}\lambda_{p, 0}(A),
\]
where the constant $c_{\ref{MAC.lemma60}} = c_{\ref{MAC.lemma60}}(\c, m_0, p, L_*, M)$ depends only on $\c, m_0, p, L_*$ and $M$.
\endlemma

\demo
Choose $f_0 \in \ell(A)$ such that $\E_{p, A}^n(f_0) = 1$ and
\[
\big(\min_{c \in \BbR}\norm{f_0 - c\chi_A}_{p, \mu_A}\big)^p = \lambda_{p, 0}(A).
\]
 Letting $f = I_{A, k, m}f_0$, by Lemma~\ref{MAC.lemma50}, we see that
\begin{equation}\label{MAC.eq50}
\E_{p, S^{m + k}(A)}^{n + k + m}(f) \le c_{\ref{MAC.lemma50}}\max_{w \in A} \#(\partial{S^k(w)})\max_{v \in S^k(A)}\E_{M, p, m}(v, S^k(A)).
\end{equation}
On the other hand, by \eqref{MAC.eq20} and \eqref{ASS.eq40},
\begin{multline*}
\frac 1{\mu(A)}\sum_{v \in S^{k + m}(A)} |f(v) -c|^p\mu(v)  = \frac 1{\mu(A)}\sum_{w \in A}\sum_{v \in S^m(S^k(w))}  |f(v) - c|^p\mu(v) \\
\ge \frac 1{\mu(A)}\sum_{w \in A}\sum_{v \in S^m(\sd{S^k(w)}{B_{M, k}(w)})}  |f_0(w) - c|^p\mu(v) \\
\ge \c^{m_0M}\frac 1{\mu(A)}\sum_{w \in A}|f_0(w) - c|^p\mu(w) \ge \c^{m_0M}\lambda_{p, 0}(A).
\end{multline*}
This and \eqref{MAC.eq50} yield the desired inequality.
\enddemo

 \setcounter{equation}{0}
\section{Relations of constants}\label{REX}
In this section, we will establish relations between conductance, neighbor disparity, and Poincar{\'e} constants towards a proof of Theorem~\ref{SUC.thm20}.

\definition\label{POC.def10}
For $w \in T$ and $n \ge 0$, define 
\[
\xi_n(w) = \max_{v \in S^n(w)} \frac{\mu(v)}{\mu(w)}
\]
\enddefinition

First, we consider a relation between Poincar{\'e} and neighbor disparity constants.

\lemma[{\cite[(2.13) Prop.-(1)]{KusZ1}}]\label{REL.lemma30}
Let $p \ge 1$. For any $w \in T$ and $n, m \ge 1$, 
\[
\tla_{p, n + m}(w) \le 
2^{p -1}\Big(\xi_n(w)\max_{v \in S^n(w)}\tla_{p, m}(v) + L_*\tla_{p, n}(w)\s_{p, m, n + |w|}\Big).
\]
\endlemma

\demo
By Theorem~\ref{BAS.thm05}, for any $f \in \ell(S^{n + m}(w))$,
\begin{multline*}
\frac 1{\mu(w)}\sum_{u \in S^{n + m}(w)} |f(u) - (f)_{S^{n + m}(w)}|^p\mu(v) \\
\le \frac{C_p}{\mu(w)}\sum_{v \in S^n(w)}\sum_{u \in S^m(v)}\Big( |f(u) - (f)_{S^m(v)}|^p + |(f)_{S^m(v)} - (f)_{S^{n + m}(w)}|^p\Big)\mu(u),
\end{multline*}
where $C_p = 2^{p - 1}$ for $p \neq 2$ and $C_2 = 1$.
Examining the first half of the above inequality, we obtain
\begin{multline*}
\frac 1{\mu(w)}\sum_{v \in S^n(w)}\sum_{u \in S^m(v)} |f(u) - (f)_{S^m(v)}|^p\mu(u) \\
\le \sum_{v \in S^n(w)} \frac{\mu(v)}{\mu(w)}\tla_{p, m}(v)\E^{|w| + n + m}_{p, S^m(v)}(f)
\le \xi_n(w)\max_{v \in S^n(w)}\tla_{p, m}(v)\E^{|w| + n + m}_{p, S^{n + m}(w)}(f).
\end{multline*}
For the other half, by Lemma~\ref{REL.lemma20},
\begin{multline*}
\frac 1{\mu(w)}\sum_{v \in S^n(w)}\sum_{u \in S^m(v)}|(f)_{S^m(v)} - (f)_{S^{n + m}(w)}|^p\mu(u)\\
= \sum_{v \in S^n(w)} \frac{\mu(v)}{\mu(w)}|(P_{n + |w|, m}f)(v) - (P_{n + |w|, m}f)_{S^n(w)}|^p\\
\le \tla_{p, n}(w)\E^{|w| + n}_{p, S^n(w)}(P_{n + |w|, m}f)\\
\le L_*\tla_{p, n}(w)\max_{v, v' \in S^n(w), (v, v') \in E^*_{n + |w|}} \s_{p, m}(v, v')\E^{n + m  + |w|}_{p, S^{n + m}(w)}(f).
\end{multline*}
Combining all, we see
\begin{multline*}
\tla_{p, n + m}(w) \\
\le C_p\Big(\xi_n(w)\max_{v \in S^n(w)}\tla_{p, m}(v) + L_*\tla_{p, n}(w)\max_{v, v' \in S^n(w), (v, v') \in E^*_{n + |w|}} \s_{p, m}(v, v')\Big).
\end{multline*}
\enddemo

\definition\label{REX.def10}
Define
\[
\ola_{p, m} = \sup_{w \in T} \tla_{p, m}(w).
\]
\enddefinition

By Theorem~\ref{FIN.thm20}, $\ola_{p, m}$ is finite for any $m \ge 1$.\par
Making use of Lemma~\ref{REL.lemma30}, we have the following inequality.

\lemma\label{REX.lemma20}
Define
\[
\xi_n = \sup_{w \in T} \xi_n(w).
\]
Then
\begin{equation}\label{REX.eq30}
\ola_{p, n + m} \le  2^{p - 1}\big(\xi_n\ola_{p, m} + L_*\ola_{p, n}\s_{p, m}\big),
\end{equation}
for any $n, m \ge 1$.
\endlemma

\remark
By Lemma~\ref{ASS.lemma20}, $\mu$ is exponential, so that there exist $\xi \in (0, 1)$  and $c > 0$ such that
\[
\xi_n \le c\xi^n
\]
for any $n \ge 1$.
\endremark

Next, we examine the relationship between the conductance and Poincar{\'e} constants.

\lemma\label{REX.lemma10}
For any $w \in T$, $l, m \ge 1$ and $k \ge m_0M_0$,
\begin{equation}\label{REX.eq10}
\overline{D}_k\E_{M_*, p, m, |w| + k + l}\tla_{p, k + m + l}(w) \ge c_{\ref{REX.lemma10}}\tla_{p, l}(w),
\end{equation}
where $\overline{D}_k = \max_{v \in \sd{T}{\{\phi\}}}\#(\partial{S^k(v)})$ and the constant $c_{\ref{REX.lemma10}} = 2^{-p}c_{\ref{MAC.lemma10}}$ depends only on $\c, m_0, p, L_*$ and $M_0$. In particular,
\begin{equation}\label{REX.eq20}
\overline{D}_k\E_{M_*, p, m}\ola_{p, k + m + l} \ge c_{\ref{REX.lemma10}}\ola_{p, l}
\end{equation}
\endlemma

\demo
Applying Lemma~\ref{MAC.lemma60} with $M = M_0$ and $A = S^l(w)$, we obtain
\[
\overline{D}_k\max_{v \in S^{k + l}(w)} \E_{M_0, p, m}(v, S^{k + l}(w))\lambda_{p, k + m}(S^l(w)) \ge c_{\ref{MAC.lemma60}}\lambda_{p, 0}(S^l(w)).
\]
Lemma~\ref{MAC.lemma10} shows
\[
\E_{M_0, p, m}(v, S^{k + l}(w)) \le \E_{M_*, p, m}(v, T_{|w| + k + l}) \le \E_{M_*, p, m, |w| + k + l}.
\]
Moreover, $\lambda_{p, k + m}(S^l(w)) = \lambda_{p, k + m + l}(w)$ and $\lambda_{p, 0}(S^l(w)) = \lambda_{p, l}(w)$ by definition. So letting $c_{\ref{REX.lemma10}} = 2^{-p}c_{\ref{MAC.lemma10}}$, we obtain \eqref{REX.eq10}.
\enddemo

The next theorem is one of the main results of this section.

\thm\label{REX.thm20}
Assume that $p  > 1$. If either
\begin{equation}\label{REX.eq100}
\lim_{n \to \infty} \xi_n\E_{p, n - m_0M_0} = 0
\end{equation}
or
\begin{equation}\label{REX.eq110}
\lim_{n \to \infty} \xi_n\overline{D}_{n - 1} = 0,
\end{equation}
then there exists $C > 0$ such that
\begin{equation}\label{REX.eq120}
\ola_{p, m} \le C\s_{p, m},
\end{equation}
\begin{equation}\label{REX.eq130}
\ola_{p, m + n} \le C\ola_{p, n}\s_{p, m}
\end{equation}
and
\begin{equation}\label{REX.eq140}
(\E_{M_*, p, n})^{-1}\ola_{p, m} \le C\ola_{p, m + n}
\end{equation}
for any $n, m \ge 1$.
\endthm

\remark
\eqref{REX.eq130} corresponds to \cite[(2.4)]{KusZ1} and \eqref{REX.eq140} corresponds to \cite[(2.3)]{KusZ1}.
\endremark

Unlike \eqref{REX.eq100}, \eqref{REX.eq110} does not depend on $p$. So, once \eqref{REX.eq110} holds, then we have \eqref{REX.eq120}, \eqref{REX.eq130} and \eqref{REX.eq140} for any $p > 1$. See Proposition~\ref{REL.prop10} after the proof for more discussion on \eqref{REX.eq110}.

\demo
For ease of notation, we write $\ola_m = \ola_{p, m}$, $\s_m = \s_{p, m}$ and $\E_{M_*, p, m} = \E_m$. By \eqref{REX.eq20}, if $n > k \ge m_0M_0$, then
\begin{equation}\label{REX.eq35}
\overline{D}_k\E_{n - k}\ola_{n + m} \ge c_{\ref{REX.lemma10}}\ola_m.
\end{equation}
This and \eqref{REX.eq30} show 
\begin{equation}\label{REX.eq40}
\ola_{n + m} \le 2^{p - 1}((c_{\ref{REX.lemma10}})^{-1}\overline{D}_k\E_{n - k}\xi_n\ola_m + L_*\ola_n\s_m).
\end{equation}
Suppose that \eqref{REX.eq100} holds. Let $k = m_0M_0$. Then there exists $n_0$ such that, for any $n \ge n_0$,
\[
2^{p - 1}(c_{\ref{REX.lemma10}})^{-1}\overline{D}_{m_0M_0}\E_{n - m_0M_0}\xi_n \le \frac 12
\]
and hence by \eqref{REX.eq40},
\begin{equation}\label{REX.eq50}
\ola_{n + m} \le 2^pL_*\ola_n\s_m.
\end{equation}
Next suppose that \eqref{REX.eq110} holds. Then there exists $n_0$ such that, for any $n \ge n_0$,
\[
2^{p - 1}(c_{\ref{REX.lemma10}})^{-1}\overline{D}_{n - 1}\E_{1}\xi_n \le \frac 12,
\]
so that we have \eqref{REX.eq50} as well. Thus we have seen that if either \eqref{REX.eq100} or \eqref{REX.eq110} holds, then there exists $n_0$ such that \eqref{REX.eq50} holds for any $n \ge n_0$.

Now, let $n_* = \max\{m_0M_0 + 1, n_0\}$. Then by \eqref{REX.eq35} and \eqref{REX.eq50}, 
\[
c_{\ref{REX.lemma10}}(\overline{D}_{m_0M_0})^{-1}(\E_{p, n_* - m_0M_0})^{-1}\ola_{m} \le \ola_{n_* + m} \le 2^pL_*\ola_{n_*}\s_{m}
\]
for any $m \ge 1$. This immediately implies \eqref{REX.eq120}. Using this and \eqref{SUC.eq20}, we have
\[
\ola_{m + n} \le \s_{m + n} \le C\s_{m}\s_{n}.
\]
Therefore, for any $m \ge 1$ and $n \in \{1, \ldots, n_0\}$,
\[
\frac{\ola_{m + n}}{\ola_{n}\s_{p, m}} \le C\frac{\s_{n}}{\ola_{n}} \le C\max_{n = 1, \ldots, n_0}\frac{\s_{n}}{\ola_{n}}.
\]
So we have verified \eqref{REX.eq130} for any $n, m \ge 1$. Letting  $k = m_0M_0$ in \eqref{REX.eq20} and using \eqref{REX.eq130}, we obtain \eqref{REX.eq140} as well.
\enddemo

The following proposition gives a geometric sufficient condition for \eqref{REX.eq110}.

\prop\label{REL.prop10}
Suppose that Assumption~\ref{ALFR} holds. Assume that $\mu$ is $\a_H$-Ahlfors regular with respect to the metric $d$.  If there exist $\tilde{\a} < \a_H$ and $c > 0$ such that
\[
\#(\partial{S^m(w)}) \le cr^{-m\tilde{\a}}
\]
for any $w \in T$ and $m \ge 0$, then \eqref{REX.eq110} holds.
\endprop

Under the assumptions of Proposition~\ref{REL.prop10}, $\a_H = \dim_H(K, d)$, which is the Hausdorff dimension of $(K, d)$, while $\dim_H(B_w, d) \le \tilde{\a}$ for any $w \in T$. So, roughly speaking,  Proposition~\ref{REL.prop10} says that if 
\[
\dim_H(K, d) > \sup_{w \in T} \dim_H(B_w, d), 
\]
then \eqref{REX.eq110} is satisfied. By this proposition, one can verify \eqref{REX.eq110} for generalized Sierpinski carpets for example.

\demo
By \cite[Theorem~3.1.21]{GAMS}, there exist $c_1, c_2 > 0$ such that
\[
c_1r^{\a_H|w|} \le \mu(K_w) \le c_2r^{\a_H|w|}
\]
for any $w \in T$. Hence $\xi_n \le cr^{\a_Hn}$, while $\overline{D}_n \le r^{-\tilde{\a}{n}}$.
\enddemo

To conclude this section, we present a lemma providing a control of the difference of a function on $T_n$ through $\E_p^n(f)$ and the Poincar{\'e} constant.

\lemma\label{REL.lemma70}
For any $w \in T$, $n \ge m \ge 1$, $f \in \ell(S^n(w))$, and $u, v \in S^n(w)$, if $\pi^{n - m}(u) = \pi^{n - m}(v)$, then
\[
|f(u) - f(v)| \le 2\c^{-\frac 1p}\E^{n + |w|}_{p, S^n(w)}(f)^{\frac 1p}\sum_{i = 1}^{n - m}\big(\ola_{p, i}\big)^{\frac 1p}.
\]
\endlemma

\demo
Let $u \in S^n(w)$. Set $S_i(u) = S^i(\pi^i(u))$ for $u \in S^n(w)$ and $i = 0, 1, \ldots, n$. By Lemma~\ref{FPE.lemma10} and \eqref{ASS.eq20},  for any $k = 1, \ldots, n$, 
\begin{multline*}
|f(u) - (f)_{S_k(u)}| \le \sum_{i = 1}^{k} |(f)_{S_{i - 1}(u)} - (f)_{S_i(u)}| \\
\le \sum_{i = 1}^{k}\bigg(\frac{\mu(\pi^i(u))}{\mu(\pi^{i - 1}(u))}\bigg)^{\frac 1p}\Big(\tla_{s, p, i}(\pi^i(u))\E_{p, S_i(u)}^{n + w}(f)\Big)^{\frac 1p}
\\\le \c^{-\frac 1p}\E_{p, S^n(w)}^{n + |w|}(f)^{\frac 1p}\sum_{i = 1}^{k}\big(\tla_{p, i}(\pi^i(u))\big)^{\frac 1p}.
\end{multline*}
Hence
\begin{multline*}
|f(u) - f(v)| 
\le |f(u) - (f)_{S_{n - m}(u)}| + |(f)_{S_{n - m}(v)} - f(v)|\\
\le \c^{-\frac 1p}\E_{p, S^n(w)}^{n + |w|}(f)^{\frac 1p}\Big(\sum_{i = 1}^{n - m}\Big(\big(\tla_{p, i}(\pi^i(v))\big)^{\frac 1p} + \big(\tla_{p, i}(\pi^i(w))\big)^{\frac 1p}\Big) \Big).
\end{multline*}
\enddemo

\setcounter{equation}{0}
\section{Proof of Theorem~\ref{SUC.thm20}}\label{PTR}
Finally, we are going to give a proof of the ``if'' part of Theorem~\ref{SUC.thm20}. Recall that by \eqref{SUC.eq30},  there exist $c > 0$ and  $\a \in (0, 1)$ such that
\[
\E_{M_*, p, m} \le c\a^{m}
\]
for any $m \ge 0$. Then since $\xi_n \le 1$, the condition \eqref{REX.eq100} is satisfied and hence \eqref{REX.eq120}, \eqref{REX.eq130} and \eqref{REX.eq140} turn out to be true.

\lemma\label{PTR.lemma10}
Set $\rho = \a^{\frac 1p}$. There exists $C > 0$ such that for any $w \in T$, $k, m \ge 1$ with $m \ge k$ and $f \in \ell(S^m(w))$,  if $u, v \in S^m(w)$ and $\pi^{m - k}(u) = \pi^{m - k}(v)$, then
\begin{equation}\label{PTR.eq20}
|f(u) - f(v)| \le C\rho^k(\ola_{p, m})^{\frac 1p}\E_{p, S^m(w)}^{|w| + m}(f)^{\frac 1p}.
\end{equation}
\endlemma

\demo
By \eqref{REX.eq140}, 
\begin{equation}\label{PTR.eq30}
\ola_{p, i} \le C\ola_{p, m}\E_{p, m - i} \le C\ola_{p, m}\rho^{p(m - i)}.
\end{equation}
Using this and applying Lemma~\ref{REL.lemma70}, we have
\begin{multline*}
|f(u) - f(v)| \le C\E^{|w| + m}_{p, S^m(w)}(f)^{\frac 1p}\sum_{i = 1}^{m - k}( \ola_{p, i})^{\frac 1p} \le C\E^{|w| + m}_{p, S^m(w)}(f)^{\frac 1p} (\ola_{p, m})^{\frac 1p}\sum_{i = k}^{m - 1} \rho^i.
\end{multline*}
\enddemo

\lemma\label{PTR.lemma20}
There exist $n_* \ge 1$ and $m_* \ge n_*$ such that if $m \ge m_*$, then there exist $w \in T$ and $f \in \ell(S^m(w))$ such that
\[
\min_{u \in S^{m - n_*}(y_1)} f(u) - \max_{u \in S^{m - n_*}(y_2)} f(u) \ge \frac 18
\]
for some $y_1, y_2 \in S^{n_*}(w)$ and
\[
\E_{p, S^m(w)}^{|w| + m}(f) \le \frac 2{\s_{p, m}}.
\]
\endlemma

\demo
Choose $(w_1, w_2) \in T \times T$ with $|w_1| = |w_2|$ such that $(w_1, w_2) \in E^*_{|w_1|}$, $\s_{p, m}(w_1, w_2) \ge \frac12\s_{p, m}$ and choose $f \in \ell(S^m(w_1) \cup S^m(w_2))$ such that $(f)_{S^m(w_1)} - (f)_{S^m(w_2)} = 1$ and 
\begin{equation}\label{PTR.eq40}
\E_{p, S^m(w_1) \cup S^m(w_2)}^{|w_1| + m}(f) = \frac 1{\s_{p, m}(w_1, w_2)}.
\end{equation}
{\bf Claim 1}:\,\,There exists $c_1 > 0$, which is independent of $m, w_1$ and $w_2$, such that if $u_1, u_2 \in S^m(w_1) \cup S^m(w_2)$ and $(u_1, u_2) \in E_{|w_1| + m}^*$, then
\begin{equation}\label{PTR.eq50}
|f(u_1) - f(u_2)| \le c_1\rho^m.
\end{equation}
Proof of Claim 1:\,\,By \eqref{REX.eq120}, \eqref{PTR.eq30} and \eqref{PTR.eq40}, we have
\begin{multline*}
|f(u_1) - f(u_2)|^p \le \E_{p, S^m(w_1) \cup S^m(w_2)}^{|w_1| + m}(f) = \frac 1{\s_{p, m}(w_1, w_2)}\\ \le \frac 2{\s_{p, m}} \le \frac C{\ola_{p, m}} \le C\rho^{pm}.
\end{multline*}
\qed\\
{\bf Claim 2}:\,\,There exists $c_2 > 0$, which is independent of $m, w_1$ and $w_2$, such that if $u_1, u_2 \in S^m(w_1) \cup S^m(w_2)$ and $\pi^{m - k}(u_1) = \pi^{m - k}(u_2)$ for $1 \le k \le m$, then
\[
|f(u_1) - f(u_2)| \le c_2\rho^k.
\]
Proof of Claim 2:\,\,It follows that $u_1, u_2 \in S^m(w_i)$ for $i = 1$ or $2$. Using Lemma~\ref{PTR.lemma10}, we obtain
\begin{multline*}
|f(u_1) - f(u_2)| \le C\rho^k(\ola_{p, m})^{\frac 1p}\E_{p, S^m(w_i)}^{|w| + m}(f)^{\frac 1p}\\ \le C\rho^k(\ola_{p, m})^{\frac 1p}\E_{p, S^m(w_1) \cup S^m(w_2))}^{|w| + m}(f)^{\frac 1p} \le C\rho^k(\ola_{p, m})^{\frac 1p}(\s_{p, m})^{-\frac 1p}.
\end{multline*}
Now \eqref{REX.eq120} immediately shows the claim.\qed\\
Define
\[
n_* = \inf\{n| n \in \BbN, 1 \ge 16c_2\rho^n\}\,\,\text{and}\,\, m_* = \max\{n_*, \inf\{m| m \in \BbN, 1 \ge 2c_1\rho^m\}\} .
\]
Hereafter we assume that $m \ge m_*$.\\
{\bf Claim 3}:\,\,For $i = 1$ or $2$, there exist $u_1, u_2 \in S^m(w_i)$ such that $u_2 \in \partial{S^m(w_i)}$ and
\[
|f(u_1) - f(u_2)| \ge \frac 14.
\]
Proof of Claim 3:\,\,Choose $v_{11}, v_{12} \in S^m(w_1)$ and $v_{21}, v_{22} \in S^m(w_2)$ such that $f(v_{11}) \ge (f)_{S^m(w_1)}, f(v_{22}) \le (f)_{S^m(w_2)}$ and $(v_{12}, v_{21}) \in E_{|w_1| + m}^*$. Since 
\[
f(v_{11}) - f(v_{12}) + f(v_{12}) - f(v_{21}) + f(v_{21}) - f(v_{22}) = f(v_{11}) - f(v_{22}) \ge 1,
\]
\eqref{PTR.eq50} shows that, for either $i = 1$ or $2$, 
\[
|f(v_{i1}) - f(v_{i2})| \ge \frac 12(1 - c_1\rho^m) \ge \frac 14.
\]
Letting $u_1 = v_{i1}$ and $u_2 = v_{i2}$, we have the claim. \qed\par
Let $w = w_i$ where $i$ is chosen in Claim 3. Exchanging $f$ by $-f$ if necessary, we see that there exists $u_1 \in S^m(w)$ and $u_2 \in \partial{S^m(w)}$ such that
\[
f(u_1) - f(u_2) \ge \frac 14.
\]
Set $y_i = \pi^{m - n_*}(u_i)$ for $i = 1, 2$. Note that $y_i \in S^{n_*}(w)$. By Claim 2, 
\begin{equation}\label{PTR.eq60}
\min_{u \in S^{m - n_*}(y_1)}f(u) - \max_{u \in S^{m - n_*}(y_2)}f(u) \ge \frac 14 - 2c_2\rho^{n_*} \ge \frac 18
\end{equation}

\enddemo

\demo[Proof of Theorem~\ref{SUC.thm20}]
Let $m \ge m_*$. Then there exist $w \in T$ and  $f \in S^m(w)$ satisfying the conclusions of Lemma~\ref{PTR.lemma20}. Set $c_0 = \max_{u \in S^{m - n_*}(y_2)} f(u)$. Define 
\[
h(v) = \begin{cases}
 1\quad&\text{if $8(f(v) -c_0) \ge 1$,}\\
 8(f(v) - c_0)\quad&\text{if $0 < 8(f(v) - c_0) < 1$,}\\
 0\quad&\text{if $8(f(v) - c_0) < 0$}
 \end{cases}
\]
for any $v \in S^m(w)$. Then $h|_{S^{n_*}(y_1)} \equiv 1$, $h|_{S^{n_*}(y_2)} \equiv 0$ and 
\[
\E_{p, m - n_*}(y_1, y_2, S^{n_*}(w)) \le \E_{p, S^m(w)}^{|w| + m}(h) \le 8^p\E_{p, S^m(w)}^{|w| + m} \le \frac{2^{3p + 1}}{\s_{p, m}}.
\]
By \eqref{SUC.eq40},
\[
\E_{M_*, p, m - n_*} \le c(n_*)\E_{p, m - n_*}(y_1, y_2, S^{n_*}(w)) \le \frac{c(n_*)2^{3p + 1}}{\s_{p, m}}.
\]
Making use of the sub-multiplicative property of $\E_{M_*, p, n}$, we have
\[
\E_{M_*, p, m} \le C\E_{M_*, p, n_*}\E_{M_*, p, m - n_*}.
\]
Finally, the last two inequalities show
\[
\E_{M_*, p, m}\s_{p, m} \le C\E_{M_*, p, n_*}c(n_*)2^{3p + 1}
\]
for any $m \ge m_*$, where the right-hand side is independent of $m$. Thus $K$ is $p$-conductively homogeneous.
\enddemo

\setcounter{equation}{0}
\section{Uniformity of constants}\label{FIN}

In this section, we study the uniformity of conductance, Poincar{\'e} and neighbor disparity constants with respect to the structure of graphs.

\definition\label{FIN.def05}
(1)\,\,A pair $(V, E)$ is called a (non-directed) graph if and only if $V$ is a countable set and $E \subseteq V \times V$ such that $(u, v) \in V$ if and only if $(v, u) \in V$. For a  graph $(V, E)$, $V$ is called the vertices and $E$ is called the edges.\\
(2)\,\,
Let $(V, E)$ and $(V', E')$ be graphs. A bijective map $\iota: V \to V'$ is called an isomorphism between $(V, E)$ and $(V', E')$ if and only if ``$(w, v) \in E$'' is equivalent to ``$(\iota(w), \iota(v)) \in E'$'' for any $u, v \in V$.\\
(3)\,\,
Let $(V, E)$ be a graph. For $p > 0$ and $f \in \ell(V)$, define $\E_P^{(V, E)}(f) \in [0, \infty]$ by
\[
\E_p^{(V, E)}(f) = \frac 12\sum_{(u, v) \in E} |f(u) - f(v)|^p.
\]
(4)\,\,
Let $(V, E)$ be a graph and let $A, B \subseteq V$ with $A \cap B = \emptyset$. Define
\[
\E^{(V, E)}_p(A, B) = \inf\{\E^{(V, E)}_p(f)| f \in \ell(V), f|_A \equiv 1, f|_B \equiv 0\}.
\]
\enddefinition

 In this section, we always identify isomorphic graphs.\par
 First, we study the uniformity of conductance constants.
 
\definition\label{FIN.def10}
For $L, N \ge 1$, define
\begin{multline*}
\G_{\E}(L, N) = \{(V, E)| \text{$(V, E)$ is a connected graph,}\\
V = \{{\bf t}, {\bf b}\} \cup V_*, \text{where the union is  a disjoint union and ${\bf t} \neq {\bf b}$},\\
1 \le \#(V_*) \le LN, \#(\{v| v \in E, (w, v) \in E\}) \le L\,\,\text{for any $w \in V_*$}\}.
\end{multline*}
\enddefinition

Since $\G_{\E}(L, N)$ is a finite set up to graph isomorphisms, we have the following theorem.

\thm\label{FIN.thm10}
For any $L, N \ge 1$ and $p > 0$,
\[
0 < \inf_{(V, E) \in \G_{\E}(L, N)} \E_p^{(V, E)}(\{{\bf t}\}, \{{\bf b}\}) \le \sup_{(V, E) \in \G_{\E}(L, N)} \E_p^{(V, E)}(\{{\bf t}\}, \{{\bf b}\}) < \infty
\]
\endthm

\definition\label{FIN.def20}
Define
\[
\ulc_{\E}(L, N, p) = \inf_{(V, E) \in \G_{\E}(L, N)} \E_p^{(V, E)}(\{{\bf t}\}, \{{\bf b}\})
\]
and
\[
 \olc_{\E}(L, N, p) = \sup_{(V, E) \in \G_{\E}(L, N)} \E_p^{(V, E)}(\{{\bf t}\}, \{{\bf b}\}).
 \]
\enddefinition

Next we consider Poincar{\'e} constants.

\definition\label{FIN.def30}
For $L \ge 1$ and $N \ge 2$, define
\begin{multline*}
\G(L, N) = \{(V, E)| \text{$(V, E)$ is a connected graph}, \\ 2 \le \#(V) \le N, \text{$\#(\{v| v \in V, (w, v) \in E, \}) \le L$ for any $w \in V$}\}.
\end{multline*}
For a connected graph $(V, E)$, define 
\[
\P(V, E) = \Big\{\mu\Big| \mu \in V \to [0, 1], \sum_{v \in V} \mu(v) = 1\Big\}.
\]
For $\mu \in \P(V, E)$, define
\[
(f)_{\mu} = \sum_{v \in V} f(v)\mu(v)
\]
for $f \in \ell(V)$ and
\[
\tla_{p, \mu}^{(V, E)} = \sup_{f \in \ell(V)} \frac{\sum_{v \in V} |f - (f)_{\mu}|^p\mu(v)}{\E_p^{(V, E)}(f)}
\]
for $p > 0$.
\enddefinition

\lemma\label{FIN.lemma10}
Let $(V, E)$ be a connected finite graph. Then for any $p \ge 1$, 
\[
0 < \inf_{\mu \in \P(V, E)} \tla_{p, \mu}^{(V, E)} \le \sup_{\mu \in \P(V, E)} \tla_{p, \mu}^{(V, E)} < \infty.
\]
\endlemma
\demo
 Write $\E_p = \E_p^{(V, E)}$. For any $p \ge 1$,
\[
 | (f)_{\mu}| + \E_p(f)^{\frac 1p}
\]
is a norm on $\ell(V)$. Therefore if
\[
\F_{\mu} = \{f| f \in \ell(V), \E_p(f) = 1, (f)_{\mu} = 0\},
\]
then $\F_{\mu}$ is a compact subset of $\ell(V)$. Fix $\mu_* \in \P(V, E)$ and set $\F = \F_{\mu_*}$. For any $f \in \ell(V)$ with $\E_p(f) \neq 0$, define $f_* = \E_p(f)^{-\frac 1p}(f - (f)_{\mu_*})$. Then $f_* \in \F$ and
\[
\frac{\sum_{v \in V} |f(v) - (f)_{\mu}|^p\mu(v)}{\E_p(f)} = \sum_{v \in V} |f_*(v) - (f_*)_{\mu}|^p\mu(v).
\]
Hence letting $F(\mu, f_*) = \sum_{v \in V} |f_*(v) - (f_*)_{\mu}|^p\mu(v)$, we see that
\[
\tla_{p, \mu}^{(V, E)} = \sup_{f_* \in \F} F(\mu, f_*).
\]
Since $\P(V, E) \times \F$ is compact and $F(\mu, f_*)$ is continuous on $\P(V, E) \times \F$, it follows that
\begin{multline*}
0 < \inf_{\mu \in \P(V, E), f_* \in \F} F(\mu, f_*) \le \inf_{\mu \in \P} \tla_{p, \mu}^{(V, E)} \\
\le \sup_{\mu \in \P(V, E)} \lambda_{p, \mu}^{(V, E)} < \sup_{\mu \in \P(V, E), f_* \in \F} F(\mu, f_*) < \infty.
\end{multline*}
\enddemo

Since $\G(L, N)$ is a finite set, the above lemma implies the following theorem.

\thm\label{FIN.thm20}
For $p \ge 1$,
\[
0 < \inf_{(V, E) \in \G(L, N), \mu \in \P(V, E)} \tla_{p, \mu}^{(V, E)} \le \sup_{(V, E) \in \G(L, N), \mu \in \P(V, E)} \tla_{p, \mu}^{(V, E)} < \infty.
\]
\endthm

\definition\label{FIN.def40}
Define
\[
\ulc_{\lambda}(p, L, N) = \inf_{(V, E) \in \G(L, N), \mu \in \P(V, E)} \tla_{p, \mu}^{(V, E)}
\]
and
\[
\olc_{\lambda}(p, L, N) = \sup_{(V, E) \in \G(L, N), \mu \in \P(V, E)} \tla_{p, \mu}^{(V, E)}.
\]
\enddefinition

Finally, we study neighbor disparity constants.

\definition\label{FIN.def50}
Define
\begin{multline*}
\G_{\s}(L, N) = \{(V, E, V_1, V_2)| (V, E) \in \G(L, N), V_1, V_2 \subseteq V, V_1 \neq \emptyset, V_2 = \emptyset, \\ 
V =  V_1 \cup V_2, V_1 \cap V_2 \neq \emptyset\}.
\end{multline*}
Let $(V, E)$ be a graph and let $\mu \in \P(V, E)$. For $U \subseteq V$ and $f \in \ell(V)$, define
\[
\mu(U) = \sum_{v \in U} \mu(v)\quad\text{and}\quad (f)_{U, \mu} = \frac 1{\mu(U)}\sum_{v \in U} f(v)\mu(v)
\]
if $\mu(U) > 0$.
For $(V, E, V_1, V_2) \in \G_{\s}(L, N)$, $\mu \in \P(V, E)$ and $p \ge 1$, define
\[
\s_{p, \mu}^{(V, E)}(V_1, V_2) = \sup_{f \in \ell(V), \E_p^{(V, E)}(f) \neq 0} \frac{|(f)_{V_1, \mu} - (f)_{V_2, \mu}|^p}{\E_p^{(V, E)}(f)},
\]
and
\[
\P(V, E, V_1, V_2, \kappa) = \{\mu | \mu \in \P(V, E), \mu(V_1) \ge \kappa\mu(V_2)\,\,\text{and}\,\,\mu(V_2) \ge \kappa\mu(V_1)\}.
\]
for $\kappa \in (0, 1]$.
\enddefinition

\thm\label{FIN.thm30}
For any $p \ge 1$, $L, N \ge 1$ and $\kappa \in (0, 1]$, 
\begin{multline*}
0 < \inf\{\s_{p, \mu}^{G}(V_1, V_2)| (G, V_1, V_2) \in \G_{\s}(L, N), \mu \in \P(G, V_1, V_2, \kappa)\}\le\\
 \sup\{\s_{p, \mu}^{G}(V_1, V_2)| (G, V_1, V_2) \in \G_{\s}(L, N), \mu \in \P(G, V_1, V_2, \kappa)\} < \infty.
\end{multline*}
\endthm

\demo
First fix $(G, V_1, V_2) \in \G_{\s}(L, N)$ and fix $\mu_* \in \P(G, V_1, V_2, \kappa)$. Define $\F$ as in the proof of Lemma~\ref{FIN.lemma10}. For any $f \in \ell(V)$, setting $f_* = \E_p(f)^{-\frac 1p}(f - (f)_{\mu_*})$, we see that $f_* \in \F$ and
\[
\frac{|(f)_{V_1, \mu} - (f)_{V_2, \mu}|^p}{\E_p(f)} = |(f_*)_{V_1, \mu} - (f_*)_{V_2, \mu}|^p.
\]
for any $\mu \in \P(G, V_1, V_2, \kappa)$. Let $F: \F \times \P(G, V_1, V_2, \kappa) \to \BbR$ by $F(f) = |(f)_{V_1, \mu} - (f)_{V_2, \mu}|$. Since $F$ is continuous and $\F \times \P(G, V_1, V_2, \kappa)$ is compact, 
\begin{multline*}
0 < \inf_{\mu \in \P(G, V_1, V_2, \kappa), f \in \F} F(f, \mu)  \le \inf_{\mu \in \P(G, V_1, V_2)} \s_{p, \mu}^{G}(V_1, V_2) \\
\le \sup_{\mu \in \P(G, V_1, V_2, \kappa), f \in \F} F(f, \mu) = \sup_{\mu \in \P(G, V_1, V_2, \kappa)} \s_{p, \mu}^{G}(V_1, V_2) < \infty.
\end{multline*}
Now the desired statement follows by the fact that $\G_{\s}(L, N)$ is a finite set up to graph isomorphisms.
\enddemo

\definition\label{FIN.def60}
Define
\[
\ulc_{\s}(L, N, \kappa) = \inf\{\s_{p, \mu}^{G}(V_1, V_2)| (G, V_1, V_2) \in \G_{\s}(L, N), \mu \in \P(G, V_1, V_2, \kappa)\}
\]
and
\[
\olc_{\s}(L, N, \kappa) =  \sup\{\s_{p, \mu}^{G}(V_1, V_2)| (G, V_1, V_2) \in \G_{\s}(L, N), \mu \in \P(G, V_1, V_2, \kappa)\}.
\]
\enddefinition

\setcounter{equation}{0}
\section{Modification of the structure of a graph}\label{MGS}

In the original work of Kusuoka-Zhou\cite{KusZ1}, they used a subgraph of $(T_n, E_n^*)$ to define their version of $\E_2^m$ in the case of the Sierpinski carpet. Namely, in our terminology, their subgraph is 
\[
E_n^1 = \{(u, v)| (u, v) \in E_1^*, \dim_H{(K_v \cap K_u)} = 1\}
\]
and their energy is
\[
\E_{p}^{1, n}(f) = \frac 12\sum_{(u, v) \in E_n^1} |f(u) - f(v)|^p
\]
for $f \in \ell(T_n)$. (They only consider the case $p = 2$.) Our theory in this paper works well if we replace our energy $\E_p^n$ with Kusuoka-Zhou's energy $\E_p^{1, n}$ because they are uniformly equivalent, i.e. there exist $c_1, c_2 > 0$ such that
\[
c_2\E_p^n(f) \le \E_p^{1, n}(f) \le c_2\E_p^n(f)
\]
for any $n \ge 1$ and $f \in \ell(T_n)$. More generally, if we replace our graph $(T_n, E_n^*)$ with a subgraph $(T_n, E_n)$ satisfying the conditions (A) and (B) below, all the results in this paper remain true except for changes in the constants.\par
The conditions (A) and (B) are;\\
\noindent (A)\,\,$G_n = (T_n, E_n)$ is a connected graph for each $n$ having the following properties:\\
(i)\,\, If $(w, v) \in E_n$, then $K_w \cap K_v \neq \emptyset$.\\
(ii)\,\, If $(w, v) \in E_n$ for $n \ge 1$, then $\pi(w) = \pi(v)$ or $(\pi(w), \pi(v)) \in E_{n - 1}$.\\
(iii)\,\,If $(w, v) \in E_n$ for $n \ge 1$, then there exist $w_1 \in S(w)$ and $v_1 \in S(v)$ such that $(w_1, w_2) \in E_{n + 1}$.\\
(iv)\,\,For any $n \ge 0$ and $w, v \in T_n$ with $K_w \cap K_v \neq \emptyset$, there exist $w(0), \ldots, w(k) \in \GG_1(w)$ satisfying $w(0) = w, w(k) = v$ and $(w(i), w(i + 1)) \in E_n$ for any $i = 0, \ldots, k - 1$. \\
(B)\,\,For any $w \in T$, the graphs $(S^n(w), E_{n + |w|}^{S^n(w)})$ associated with the partition $T(w)$ of $K_w$ satisfies the counterparts of conditions (i), (ii), (iii) and (iv) of (3).\par
Naturally, the graph $(T_n, E_n^*)$ satisfied the conditions (A) and (B).

\section{Open problems}\label{OPR}
In the final section, we gather some of open problems and future directions of research. \\
{\bf 1.}\,\,Regularity of $\W^p$ for $p \in [1, \dim_{AR}(K, d)]$;\,\,As we have already mentioned, it is not known whether or not $C(K) \cap W^p$ is dense in $\W^p$ for $p \in [1, \dim_{AR}(K, d)]$. The first step should be to establish an elliptic Harnack  principle for $p$-harmonic functions on approximating graphs and/or the limiting object $(\W^p, \hE_p(\cdot) + \norm{\cdot}_{p, \mu})$. Even in the case of $p = 2$, this problem is open except  for the case of generalized Sierpinski carpets. The conjecture presented in the introduction is closely related to this problem as well.\\
{\bf 2.}\,\,Construction of $p$-form and $p$-Laplacian:\,\, In this paper, we have constructed a $p$-energy $\hE_p(f)$ but not a $p$-form $\hE_p(f, g)$. On a graph $G = (V, E)$, if we define
\[
\E_p(f, g) = -\sum_{x \in V} (\Delta_pf)(x)g(x)
\]
for $f, g \in \ell(V)$, where $\Delta_p$ is the $p$-Laplacian defined by
\begin{align*}
\Phi_p(t) &= \begin{cases}
|t|^{p - 2}t\quad&\text{if $t \neq 0$}\\
0 \quad&\text{if $t = 0$,}
\end{cases}\\
(\Delta_pf)(x) &= \sum_{y \in V, (x, y) \in E} \Phi_p(f(y) - f(x)),
\end{align*}
then it follows that
\[
\E_p(f) = \frac 12\sum_{(x, y) \in E} |f(x) - f(y)|^p = \E_p(f, f).
\]
As a natural counterpart, we expect to have a $p$-form $\hE_p(f, g)$ which is linear in $g$, satisfies $\hE_p(f) = \hE_p(f, f)$ for any $f \in \W^p$, and has an expression such as
\[
\E_p(f, g) = -\int_K (\Delta_pf)(x)g(x)\mu(dx).
\]
{\bf 3.} $p$-energy measure:\,\,In the case $p = 2$, there is the notion of energy measures associated with a strongly local regular Dirichlet form $(\E, \F)$, where $\E$ is the form and $\F$ is the domain.  Roughly speaking, the energy measure $\mu_f$ associated with $f \in \F$ is a positive Radon measure satisfying
\[
\int_X u(x)d\mu_f(dx) = 2\E(uf, f) - \E(f^2, u)
\]
for any $u  \in \F \cap C_0(X)$.  See \cite{FOT} for details. So, what is a counterpart of this in the case of $\hE^p$? Is there any natural measure $\mu_f$ for $f \in \W^p$ such that
\[
\int_K d\mu_f(dx) = \hE_p(f)\text{?}
\]
For  $\BbR^n$, the answer is yes and $\mu_f = |\nabla{f}|^pdx$. For the planar Sierpinski carpet, this problem has already been studied in \cite{Shimizu2}. However, we know almost nothing beyond those examples.\\
{\bf 4.} Fractional Korevaar-Shoen type expression:\,\,As we have already mentioned, a fractional Korevaar-Shoen type expression of $\W^p$ has already shown in \cite{Shimizu2} in the case of the planar Sierpinski carpet. Namely, we have
\[
\W^p = \bigg\{f \bigg| f \in L^p(K, \mu), \limsup_{r \downarrow 0} \int_K\frac 1{r^{\a_H}}\int_{B_{d_*}(x, r)} \frac{|f(x) - f(y)|^p}{r^{\beta_p}}dxdy < \infty\bigg\}.
\]
and it is shown in \cite{Shimizu2} that $\b_p > p$ for any $p > 1$. How about other cases? Suppose that Assumption~\ref{ALFR} holds and $\mu$ is $\a_H$-Ahlfors regular with respect to the metric $d$. Then we expect that $\b_p = \a_H + \tau_p$ and we know $\a_H + \tau_p \ge p$ by \cite[(4.6.14)]{GAMS}. Now our questions are:\\
$\cdot$\,\, Do we have a fractional Korevaar-Shoen type expression as above?\\
$\cdot$\,\,When does $\beta_p > p$ hold? (Apparently, if $K = [-1, 1]^L$, then $\beta_p = p$.)\\
Related question is\\
 If $\b_p = p$, then does $\W^p$ coincide with any of the Sobolev type spaces given by approaches using upper gradients?\\
 {\bf 5.} Without local symmetry:\,\,In Sections~\ref{CYP}, \ref{ESC}, \ref{SCR} and \ref{OEX}, we have shown the conductive homogeneity of self-similar sets having local symmetry, which helped us to extend a path from one piece of $K_w$ to neighbors by the reflection in its boundaries.  However, the local symmetry does not seem indispensable for having conductive homogeneity. Intuitively the essence should be the balance of conductances in different directions, for example, the vertical and the horizontal directions for square-based self-similar sets. Unfortunately, we have not had any example without local symmetry yet except for finitely ramified cases.

\par\vspace{20pt}\par
\renewcommand{\theequation}{\thesection.\arabic{equation}}
\setcounter{equation}{0}
\appendix
\noindent{\Large \bf Appendices}
\par\vspace{-10pt}
\section{Basic inequalities}\label{BAS}

The next two lemmas can be deduced from the H{\"o}lder inequality.

\lemma\label{BAS.lemma10}
For $p \in (0, \infty)$, 
\[
\Big|\sum_{i = 1}^n a_i\Big|^p \le \max\{1, n^{p - 1}\}\sum_{i = 1}^n |a_i|^p
\]
for any $n \ge 1$ and $a_1, \ldots, a_n \in \BbR$.
\endlemma

\lemma\label{BAS.lemma20}
Let $p, q \in [1, \infty]$ satisfying $\frac 1p + \frac 1q = 1$. Then for any $n \in \BbN$ and $a_1, \ldots, a_n \in \BbR$,
\[
\Big(\sum_{i = 1}^n |a_i|^q\Big)^{1/q} \le \max\{1, n^{(p - 2)/p}\}\Big(\sum_{i = 1}^n |a_i|^p\Big)^{1/p}.
\]
\endlemma

The following fact implies the comparison \eqref{MAC.eq100} of two types of Poincar{\'e} constants, $\lambda_{p, m}$ and $\tla_{p, m}$.

\thm[{\cite[Lemma~4.17]{Bjorn2}}]\label{BAS.thm05}
Let $\mu$ be a finite measure on a set $X$. Then for any $f \in L^p(X, \mu)$ and $c \in \BbR$,
\[
\norm{f - c}_{p, \mu} \ge \frac 12\norm{f - (f)_{\mu}}_{p, \mu},
\]
where $\norm{\cdot}_{p, \mu}$ is the $L^p$-norm with respect to $\mu$ and $(f)_{\mu} = \mu(X)^{-1}\int_X fd\mu$.
\endthm

The following lemma is a discrete version of the above theorem.

\cor\label{BAS.cor10}
Let $(\mu_i)_{i = 1, \ldots, n} \in (0, 1)^n$ with $\sum_{i = 1}^n \mu_i = 1$. Then
\[
\sum_{i = 1}^n |x - a_i|^p\mu_i \ge \Big(\frac{1}2\Big)^p\sum_{i = 1}^n \Big|\sum_{j = 1}^n \mu_ja_j - a_i\Big|^p\mu_i
\]
for any $x, a_1, \ldots, a_n \in \BbR$.
\endcor

\setcounter{equation}{0}
\section{Basic facts on $p$-energy}\label{FPE}
Let $G = (V, E)$ be a finite graph. For $A \subseteq V$, set $E_A = \{(x, y)| x, y \in A, (w, y) \in E\}$ and $G_A = (A, E_A)$.

\definition\label{FPE.def10}
Let $\mu: V \to (0, \infty)$ and let $A \subseteq V$. Define $\supp{\mu} = \{x| x \in V, \mu(x) > 0\}$.
Let $p  > 0$. For $u \in \ell(V)$, define
\begin{align*}
\E^G_{p}(u) &= \frac 12\sum_{(x, y) \in E}|u(x) - u(y)|^p\\
\norm{u}_{p, \mu} &= \Big(\sum_{x \in V} |u(x)|^p\mu(x)\Big)^{1/p}\\
(u)_{\mu} &= \frac 1{\sum_{y \in V} \mu(y)}\sum_{x \in V} \mu(x)u(x)\\
\end{align*}
and
\[
\lambda^G_{p, \mu} = \sup_{u \in \ell(V), u \neq 0} \frac{(\min_{c \in \BbR}\norm{u - c\chi_{V}}_{p, \mu})^p}{\E^G_{p}(u)},
\]
where $\chi_V \in \ell(V)$ is the characteristic function of the set $V$.
\enddefinition

For $A \subseteq U$, set
\[
\E_p^A = \E_p^{G_A} \quad\text{and}\quad \lambda_{p, \mu}^A = \lambda_{p, \mu|_A}^{G_A}.
\]

\lemma\label{FPE.lemma05}
Define
\[
\tla_{p, \mu}^G = \sup_{u \in \ell(V), u \neq 0} \frac{(\norm{u - (u)_{\mu}\chi_{V}}_{p, \mu})^p}{\E^G_{p}(u)}.
\]
Then
\[
\Big(\frac{1}2\Big)^p\,\tla_{p, \mu}^G \le \lambda_{p, \mu}^G \le \tla_{p, \mu}^G.
\]

\endlemma

\demo
By Corollary~\ref{BAS.cor10},
\[
\sum_{x \in V} |u(x) - (u)_{\mu}|^p\mu(x) \ge \min_{c \in \BbR}\sum_{x \in V} |u(x) - c|^p\mu(x) \ge \Big(\frac{1}2\Big)^p\sum_{x \in V} |u(x) - (u)_{\mu}|^p\mu(x).
\]
\enddemo

\lemma[{\cite[(1.5) Prop.-(2)]{KusZ1}}]\label{FPE.lemma10}
Let $p \in [1, \infty)$ and let $\mu: V \to (0, \infty)$. Assume that $A \subseteq B \subseteq V$. Then for any $u \in \ell(B)$, 
\[
|(u)_A - (u)_B|  \le \frac 1{\mu(A)^{\frac 1p}}\Big(\tla_{p, \mu}^B\E_p^B(u)\Big)^{\frac 1p}.
\]

\endlemma

\demo
By the H{\"o}lder inequality, 
\[
|(u)_A - (u)_B|  \le \frac 1{\mu(A)}\int_B\chi_A|u - (u)_B|d\mu \le \frac 1{\mu(A)^{\frac 1p}}\Big(\int_B |u - (u)_B|^pd\mu\Big)^{\frac 1p}.
\]
\enddemo

\setcounter{equation}{0}
\section{Useful  facts on combinatorial modulus}\label{UFD}
In this appendix, we have useful facts on combinatorial modulus. In particular, the last lemma, Lemma~\ref{UFD.lemma50}, is a result on the comparison of moduli in two different graphs. This lemma plays a key role on several occasions in this paper.\par
Let $V$ be a countable set and let $\P(V)$ be the power set of $V$. For $\rho: V \to [0, \infty)$ and $A \subseteq V$, define
\[
L_{\rho}(A) = \sum_{x \in A} \rho(x).
\]
For $\U \subseteq \P(V)$, define
\[
\A(\U) = \{\rho| \rho: V \to [0, \infty), L_{\rho}(A) \ge 1\,\, \text{for any $A \in \U$}\}.
\]
Moreover for $\rho: V \to [0, \infty)$, define
\[
M_p(\rho) = \sum_{x \in V} \rho(x)^p
\]
and
\[
Mod_p(\U) = \inf_{\rho \in \A(\U)} M_p(\rho).
\]
Note that if $\U = \emptyset$, then $\A(\U) =[0, \infty)^{V}$ and $Mod_p(\U) = 0$.

\lemma\label{UFD.lemma10}
Assume that $\U$ consists of finite sets. Then there exists $\rho_*  \in \A(\U)$ such that
\[
Mod_p(\U)  = M_p(\rho_*).
\]
\endlemma

\demo
Choose $\{\rho_i \}_{i \ge 1} \subseteq \A(\U)$ such that $M_p(\rho_i) \to Mod_p(\U)$ as $i \to \infty$. Since $V$ is countable, there exists a subsequence $\{\rho_{n_j}\}_{j \ge 1}$ such that, for any $v \in V$, $\rho_{n_j}(v)$ is convergent as $j \to \infty$. Set $\rho_*(p) = \lim_{j \to \infty} \rho_{n_j}(p)$.  For any $A \in \U$, since $A$ is a finite set, it follows that $L_{\rho_*}(A) \ge 1$. Hence $\rho_* \in \A(\U)$. For any $\e > 0$, there exists a finite set $X_{\e}$ such that $\sum_{v \in X_{\e}} \rho_*(v)^p \ge  M_p(\rho_*) - \e$. As
\[
Mod_p(\U) = \lim_{j \to \infty} M_p(\rho_{n_j})  \ge \lim_{j \to \infty} \sum_{v \in X_{\e}} \rho_{n_j}(v)^p,
\]
we obtain $Mod_p(\U) \ge M_p(\rho_*) - \e$ for any $\e > 0$. Hence $Mod_p(\U) \ge M_p(\rho_*)$. On the other hand, since $\rho_* \in \A(\U)$, we see $M_p(\rho_*) \ge Mod_p(\U)$. Therefore $M_p(\rho_*) = Mod_p(\U)$.
\enddemo

\lemma\label{UFD.lemma20}
Assume that $\U$ consists of finite sets. For $v \in V$, define $\U_v = \{A| A \in \U, v \in A\}$. Then
\[
\rho_*(v)^p \le Mod_p(\U_v)
\]
for any $\rho_* \in \A(\U)$ with $M_p(\rho_*) = Mod_p(\U)$. In particular, if $\U_v = \emptyset$, then $\rho_*(v) = 0$.
\endlemma

\demo
Suppose that $\rho_* \in \A(\U)$ and $M_p(\rho_*) = Mod_p(\U)$. Assume that $\U_v = \emptyset$ and $\rho_*(v) > 0$. Define $\rho_*'$ by
\[
\rho_*'(u) = \begin{cases} 
\rho_*(u)&\quad\text{if $u \neq v$,}\\
0&\quad\text{if $u = v$.}
\end{cases}\]
Then $\rho_*' \in \A(\U)$ and $M_p(\rho_*') < M_p(\rho_*)$. This contradicts the fact that $M_p(\rho_*) = Mod_p(\U)$. Thus if $\U_v = \emptyset$, then $\rho_*(v) = 0$. Next assume that $\U_v \neq \emptyset$. Let $\rho_v \in \A(\U_v)$ with $M_p(\rho_v) = Mod_p(\U_v)$. Note that such a $\rho_v$ does exist by Lemma~\ref{UFD.lemma10}. Define
\[
\tilde{\rho}(u) = \begin{cases}
\max\{\rho_*(u), \rho_v(u)\}\quad&\text{if $u \neq v$,}\\
\rho_v(v) \quad&\text{if $u = v$.}
\end{cases}
\]
Let $A \in \U$. If $v \notin A$, then $\tilde{\rho} \ge \rho_*$ on $A$, so that $\tilde{\rho} \in \A(A)$. If $v \in A$, then $\tilde{\rho} \ge \rho_v$ on $A$ and hence $\tilde{\rho} \in \A(A)$. Thus we see that $\tilde{\rho} \in \A(\U)$. Therefore,
\[
Mod_p(\U) \le M_p(\tilde{\rho}) \le \sum_{u \neq v} \rho_*(u)^p + \sum_{u \in V} \rho_v(u)^p = Mod_p(\U) - \rho_*(v)^p + Mod_p(\U_v).
\]
\enddemo

Define $\ell_+(V) = \{f| f: V \to [0, \infty)\}$.

\lemma\label{UFD.lemma40}
Let $V_1$ and $V_2$ be finite sets. Let $\U_i \subseteq \P(V_i)$ for $i = 1, 2$. If there exist maps $\xi: \U_2 \to \U_1$, $F:\ell_+(V_1) \to \ell_+(V_2)$ and constants $C_1, C_2 > 0$ such that
\[
C_1L_{F(\rho)}(\c) \ge L_{\rho}(\xi(\c)) \quad\text{and}\quad M_p(F(\rho)) \le C_2M_p(\rho)
\]
for any $\rho \in \ell_+(V_1)$ and $\c \in \U_2$, then
\[
Mod_p(\U_2) \le (C_1)^pC_2Mod_p(\U_1)
\]
for any $p > 0$.
\endlemma

\demo
Note that $C_1F(\rho) \in \A(\U_2)$ for any $\rho \in \A(\U_1)$. Hence if $F'(\rho) = C_1F(\rho)$, then
\begin{multline*}
Mod_p(\U_2) = \min_{\rho \in \A(\U_2)} M_p(\rho) \le \min_{\rho \in \A(\U_1)} M_p(F'(\rho))\\
 \le (C_1)^PC_2\min_{\rho \in \A(\U_1)}M_p(\rho) = (C_1)^PC_2Mod_p(\U_1).
\end{multline*}
\enddemo

\lemma\label{UFD.lemma50}
Let $V_1$ and $V_2$ be countable sets and let $\U_i \subseteq \P(V_i)$ for $i = 1, 2$. Assume that $H_v \subseteq V_1$ and $\#(H_v) < \infty$ for any $v \in V_2$.  Furthermore, assume that, for any $B \in \U_2$, there exists $A \in \U_1$ such that $A \subseteq \cup_{v \in B} H_v$. Then
\[
Mod_p(\U_2) \le \sup_{v \in V_2} \#(H_v)^p \sup_{u \in V_1} \#(\{v| v \in V_2, u \in H_v\})Mod_p(\U_1)
\]
for any $p > 0$.
\endlemma

\demo
For $\rho: V_1 \to \BbR$, define
\[
F(\rho)(v) = \max_{u \in H_v} \rho(u)
\]
for any $v \in V_2$. Then $F: \ell_+(V_1) \to \ell_+(V_2)$ and
\begin{multline*}
M_p(F(\rho)) = \sum_{v \in V_2} \max_{u \in H_v}\rho(u)^p \le \sum_{v \in V_2} \sum_{u \in H_v} \rho(u)^p \\
\le  \sup_{u \in V_1} \#(\{v| v \in V_2, u \in H_v\})M_p(\rho).
\end{multline*}
On the other hand, for $B \in \U_2$, choose $\xi(B) \in \U_1$ such that $\xi(B) \subseteq \cup_{v \in B} H_v$. Then for any $\rho \in \ell_+(V_1)$ and $B \in \U_2$,

\begin{multline*}
\sup_{u \in V_2}\#(H_u)L_{B}(F(\rho)) \ge \sum_{u \in B} \#(H_u)F(\rho)(u)  \ge \sum_{u \in B}\sum_{v \in H_u} \rho(v)\\
 = \sum_{v \in \cup_{u \in B}H_u} \#(\{u| v \in H_u\})\rho(v) \ge \sum_{v \in \xi(B)} \rho(v) = L_{\xi(B)}(\rho).
\end{multline*}
Hence by Lemma~\ref{UFD.lemma40}, we have the desired conclusion.
\enddemo

\section{An Arzel{\`a}-Ascoli theorem for discontinuous functions}\label{APP1}

The following lemma is a version of Arzel{\`a}-Ascoli theorem showing the existence of a uniformly convergent subsequence of a sequence of functions. The difference between the original version and the current one is that it can handle a sequence of discontinuous functions.

\lemma[Extension of Arzel{\`a}-Ascoli]\label{GHT.lemma10}
Let $(X, d_X)$ be a totally bounded metric space and let $(Y, d_Y)$ be a metric space. Let $u_i: X \to Y$ for any $i \ge 1$. Assume that there exist a monotonically increasing function $\eta: [0, \infty) \to [0, \infty)$ and  a sequence $\{\delta_i\}_{i \ge 1} \in [0, \infty)$ such that $\eta(t) \to 0$ as $t \downarrow 0$, $\delta_i \to 0$ as $i \to \infty$ and
\begin{equation}\label{GHT.eq100}
d_Y(u_i(x_1), u_i(x_2)) \le \eta(d_X(x_1, x_2)) + \delta_i
\end{equation}
for any $i \ge 1$ and $x_1, x_2 \in X$. If $\overline{\cup_{i \ge 1} u_i(X)}$ is compact, then there exists a subsequence $\{u_{n_j}\}_{j \ge 1}$ such that $\{u_{n_j}\}_{j \ge 1}$ converges uniformly to a continuous function $u: X \to Y$ as $j \to \infty$ satisfying $d_Y(u(x_1),u(x_2)) \le \eta(d_X(x_1, x_2))$ for any $x_1, x_2 \in X$.
\endlemma

\demo
Since $X$ is totally bounded, there exists a countable subset $A \subseteq X$ which is dense in $X$ and contains a finite $\tau$-net $A_{\tau}$ of $X$ for any $\tau > 0$. Let $K = \overline{\cup_{i \ge 1} u_i(X)}$. Since $K$ is compact and $\{u_i(x)\}_{i \ge 1} \subseteq K$ is bounded for any $x \in A$, there exists a subsequence $\{u_{m_k}(x)\}_{k \ge 1}$ converging as $k \to \infty$. By the standard diagonal argument, we may find a subsequence $\{u_{n_j}\}_{j \ge 1}$ such that $\{u_{n_j}(x)\}_{j \ge 1}$ converges as $j \to \infty$ for any $x \in A$. Set $v_j = u_{n_j}$ and $\a_j = \delta_{n_j}$. Define $v(x) = \lim_{j \to \infty} v(x)$ for any $x \in A$.  By \eqref{GHT.eq100}, 
\[
d_Y(v_j(x_1), v_j(x_2)) \le \eta(d_X(x_1, x_2)) + \a_j
\]
for any $x_1, x_2 \in A$. Letting $j \to \infty$, we see that
\begin{equation}\label{GHT.eq110}
d_Y(v(x_1), v(x_2)) \le \eta(d_X(x_1, x_2))
\end{equation}
for any $x_1, x_2 \in A$. Since $A$ is dense in $X$, $v$ is extended to a continuous function on $X$ satisfying \eqref{GHT.eq110} for any $x_1, x_2 \in X$. Fix $\e > 0$. Choose $\tau > 0$ such that $\eta(\tau) < \e/3$. Since the $\tau$-net $A_{\tau}$ is a finite set, there exists $k_0$ such that if $k \ge k_0$, then $\a_k < \e/3$ and $d_Y(v(z), v_k(z)) < \e$ for any $z \in A_{\tau}$. Let $x \in X$ and choose $z \in A_{\tau}$ such that $d_X(x, z) < \tau$.  If $k \ge k_0$, then
\begin{multline*}
d_Y(v_k(x), v(x)) \le d_Y(v_k(x), v_k(z)) + d_Y(v_k(z), v(z)) + d_Y(v(z), v(x))
 \\ \le 2\eta(d_X(x, z)) + \a_k + d_Y(v_k(z), v(z)) < 2\e.
\end{multline*}
Thus $\{v_j\}_{j \ge 1}$ converges uniformly to $v$ as $j \to \infty$.
\enddemo

\setcounter{equation}{0}
\section{Geometric properties of strongly symmetric self-similar sets}\label{SNS}
\def\ws{\widetilde{\s}}
In this appendix, we will give proofs  of claims on topological and geometric properties of self-similar sets treated in Section~\ref{OEX}. Namely we will give proofs of Propositions~\ref{OEX.prop10} and \ref{OEX.prop20}.   First we recall the setting of Section~\ref{OEX}.  Let $S$ be a finite subset of $\BbR^L$ and let $\rho \in (0, 1)$.  Let $U_q \in O(L)$ for any $q \in S$. Define $f_q: \BbR^L \to \BbR^L$ by
\[
f_q(x) = \rho{U_q}(x - q) + q
\]
for $x \in \BbR^L$.  Let $K$ be the self-similar set with respect to $\{f_q\}_{q \in S}$, i.e. $K$ is the unique non-empty compact set $K$ satisfying
\[
K = \bigcup_{q \in S} f_q(K).
\]
The triple $(K, S, \{f_q\}_{q \in S})$ is know to be a self-similar structure defined in Definition~\ref{SSF.def00} and the map $\chi: S^{\BbN} \to K$ is given by 
\[
\{\chi(q_1q_2\ldots)\} = \bigcap_{m \ge 0} f_{\word qm}(K)
\]
as we have seen in Section~\ref{SSF}.
\definition\label{SNS.def10}
(1)\,\,Define $\ws: S^{\BbN} \to S^{\BbN}$ by $\ws(q_1q_2\ldots) = q_2q_3\ldots$ for $q_1q_2\ldots \in S^{\BbN}$.\\
(2)\,\,Define 
\[
C_K = \bigcup_{i \neq j \in S} K_i \cap K_j,\quad
\C = \chi^{-1}(C_K),\quad
\P = \bigcup_{k \ge 1} \ws^k(\C),
\]
and $V_0 = \chi(\P)$. $\C$ and $\P$ are called the critical set and the post critical set of $(K, S, \{f_q\}_{q \in S})$ respectively. A self-similar structure $(K, S, \{f_q\}_{q \in S})$ is said to be post critically finite (p.c.f. for short) if $\P$ is a finite set.
\enddefinition

By \cite[Theorem~1.2.3]{AOF}, we have the following proposition.

\prop\label{SNS.prop10}
The map $\chi$ is continuous and surjective. Moreover,
\begin{equation}\label{SNS.eq10}
\chi(q_1q_2\ldots) = f_{q_1}(\chi(\ws(q_1q_2\ldots)))
\end{equation}
for any $q_1q_2\ldots \in S^{\BbN}$.
\endprop

In this appendix, we suppose that Assumption~\ref{OEX.ass10} holds. \par
The next lemma gives a proof of Proposition~\ref{OEX.prop10}.

\lemma\label{SNS.lemma10}
Under Assumption~\ref{OEX.ass10}, we have\\
{\rm (1)}\,\,
For any $i = 1, \ldots, N$, $\chi^{-1}(q) = \overline{q}$,where $\overline{q} = qqq\ldots \in S^{\BbN}$.\\
{\rm (2)}\,\,
$\P = \{\overline{q}| q \in U \}$. In particular, a self-similar structure $(K, S, \{f_q\}_{q \in S})$ is post critically finite and $V_0 = U$.
\endlemma
\demo
(1)\,\,
Suppose $\chi(\tau_1\tau_2\ldots) = q$. Then by \eqref{SNS.eq10},
\[
q = \chi(\tau_1\tau_2\ldots) = f_{\tau_1}(\chi(\tau_2\tau_3\ldots)) \in K_{\tau_1}.
\]
By Assumption~\ref{OEX.ass10}-(1), it follows that $\tau_1 = q$. Since $f_q$ is invertible, we see that $\chi(\tau_2\tau_3\ldots) = q$. Using the same argument as above, we see that $\tau_2 = q$ as well. Thus we deduce that $\tau_k = q$ for any $k \in \BbN$ inductively.\\
(2)\,\,
Suppose that $\chi(\tau_1\tau_2\ldots) \in f_{\tau_1}(K)  \cap f_q(K)$ for some $q \neq \tau_1$. By \eqref{SNS.eq10}, it follows that $\chi(\tau_1\tau_2\ldots) = f_{\tau_1}(\chi(\tau_2\tau_3\ldots))$. Hence by Assumption~\ref{OEX.ass10}-(2),
\[
\chi(\tau_2\tau_3\ldots) \in (f_{\tau_1})^{-1}(f_{\tau_1}(K) \cap f_q(K)) \subseteq U.
\]
Thus $\tau_2\tau_3\ldots = \overline{q'}$ for some $q' \in U$.  Therefore, $\P \subseteq U$.\\
Conversely, again by Assumption~\ref{OEX.ass10}-(2), for any $q \in U$, there exist $p_1, p_2 \in S$ with $p_1 \neq p_2$ such that $\chi(p_1\overline{q}) \in f_{p_1}(K) \cap f_{p_2}(K)$. This shows that $p_1\overline{q} \in \C$ and hence $\overline{q} \in \P$.
\enddemo

In the next two lemmas, we are going to show a sufficient condition for Assumption~\ref{OEX.ass30}.

\lemma\label{SNS.lemma20}
Suppose that  Assumption~\ref{OEX.ass10} holds and that $U_q$ is the identity map for any $q \in V_0$. Let $q = f_{p_1}(q_1) = f_{p_2}(q_2)$ for some $p_1, p_2 \in S$  with $p_1 \neq p_2$ and $q_1, q_2 \in V_0$. Then there exists $\c = \c(p_1, p_2, q_1, q_2) > 0$ such that
\[
d(\overline{\sd{K_{p_1}}{K_{p_1(q_1)^{m - 1}}}}, K_{p_2}) \ge \c\rho^m
\]
for any $m \ge 1$, where $d(A, B) = \inf_{x \in A, y \in B} |x - y|$ and $(q)^k = \underset{\text{$k$-times}}{q\ldots{q}} \in T_k$.
\endlemma

In the following proof, we assume that 
\[
\#(f_{p_1}(K) \cap f_{p_2}(K)) \le 1.
\]
to avoid a non-essential complication of arguments. Without this assumption, the lemma is still true with a technical modification of the proof.
\demo
Set $c_m = \inf\{d(K_w, K_v)| w, v \in T_m, K_w \cap K_v = \emptyset\}$. Define
\[
X_m = \overline{\sd{K_{p_1}}{K_{p_1(q_1)^{m - 1}}}}\,\,\text{and}\,\,Y_m = \overline{\sd{K_{p_1q_1}}{K_{p_1(q_1)^{m - 1}}}}
\]
for $m \ge 1$. Then $X_m = Y_m \cup (\cup_{q \neq q_1} K_{p_1q})$ and $K_{p_2} = K_{p_2q_2} \cup (\cup_{q \neq q_2} K_{p_2q})$. This implies that
\[
d(X_m, K_{p_2}) \ge \min\{d(Y_m, K_{p_2q_2}), c_2\}.
\]
On the other hand, letting $f(x) = \rho(x - q) + q$, we see that $Y_m \cup K_{p_2q_2} = f(X_{m - 1} \cup K_{p_2})$. This yields $d(Y_m, K_{p_2q_2}) = {\rho}d(X_{m - 1}, K_{p_2})$. Consequently, we have
\[
d(X_m, K_{p_2}) \ge \min\{{\rho}d(X_{m - 1}, K_{p_2}), c_2\}.
\]
Now inductive argument suffices.
\enddemo

\lemma\label{SNS.lemma30}
Suppose that  Assumption~\ref{OEX.ass10} holds and that $U_q$ is the identity map for any $q \in V_0$. Then Assumption~\ref{OEX.ass30} holds.\endlemma

\remark
According to the notation in the proof of Lemma~\ref{SNS.lemma20}, this lemma claims $c_m \ge c\rho^m$ for any $m \ge 1$.
\endremark

\demo
Suppose that $w, v \in T_m$ and $K_w \cap K_v = \emptyset$. Let $w = \word wm$ and let $v = \word vm $. In case $w_1 = w_2$, then
\[
d(K_w, K_v) = {\rho}d(K_{w_2\ldots{w_m}}, K_{v_2\ldots{v_m}}) \ge c_{m - 1}\rho^m.
\]
Otherwise, assume that $w_1 \neq v_1$. If $K_{w_1} \cap K_{v_1} = \emptyset$, then $d(K_w, K_v) \ge c_1$. So, the remaining possibility is that $w_1 = v_1$ and $K_{w_1} \cap K_{v_1} \neq \emptyset$. In this case let $q = K_{w_1} \cap K_{v_1}$. Then $q = f_{w_1}(p_{j_1}) = f_{w_2}(p_{j_2})$ for some $j_1, j_2 \in \{1, \ldots, L\}$. By Lemma~\ref{SNS.lemma20}, it follows that
\[
d(K_w, K_v) \ge \overline{\c}\rho^m,
\]
where $\overline{\c} = \min\{\c(p_1, p_2, q_1, q_2)| p_1, p_2 \in S, q_1, q_2 \in V_0, f_{p_2}(q_1) = f_{p_1}(q_2)\}$. Combining all the cases, we see that
\[
c_m \ge \min\{c_{m - 1}, \c, c_1\rho^{-m}\} \ge \min\{c_1, \c\}
\]
for any $m \ge 1$.
\enddemo

Now we start showing Proposition~\ref{OEX.prop20}, that is,  Assumption~\ref{ALFR} hold under Assumptions~\ref{OEX.ass10} and \ref{OEX.ass30}.

\lemma\label{SNS.lemma40}
Under Assumptions~\ref{OEX.ass10} and \ref{OEX.ass30},  Assumption~\ref{ALFR}-(2) holds with $r = \rho$, $M_* = 1$, and $d = d_*$, where $d_*$ is the restriction of the Euclidean metric.
\endlemma

\demo
The condition (2A) is obvious. Set
\[
\GG_{1, n}(x) = \bigcup_{\substack{w \in T_n\\x \in K_w}} \GG_1(w).
\]
for $x \in K$ and $n \ge 1$.
Then for any $v \in \sd{T_n}{\GG_{1, n}(x)}$, there exists $w \in T_n$ such that $x \in K_w$ and $K_w \cap K_v = \emptyset$. By Lemma~\ref{SNS.lemma30}, we see that $d(K_w, x) \ge c\rho^n$ and hence $B_{d_*}(x, cr^n) \cap K_v = \emptyset$. Thus we have 
\begin{equation}\label{SNS.eq20}
B_{d_*}(x, c\rho^n) \subseteq U_1(x: n).
\end{equation}
On the other hand, by (2A), there exists $c' > 0$ such that $\diam{K_w, d_*} \le c'\rho^{|w|}$ for any $w \in T$. This implies
\begin{equation}\label{SNS.eq30}
U_1(x: n) \subseteq B_{d_*}(x, 3c'\rho^n).
\end{equation}
So we have (2B). Choose $x_0 \in \sd{K}{V_0}$ and choose $m_0 \in \BbN$ such that $2\rho^{m_0} < d(x_0, V_0)$. Let $w \in T_n$ and let $u \in  \GG_{1, m_0 + n}(f_w(x_0))$. Suppose that $u \in T(v)$ for some $v \in T_n$ with $v \neq w$. Since $u \in \GG_{1, m_0 + n}(f_w(x_0))$, there exists $u_0 \in T_{n + m_0}$ such that $f_w(x_0) \in K_{u_0}$ and $K_{u_0} \cap K_u \neq \emptyset$. Let $y \in K_u$. Since $K$ is connected (and hence arcwise connected by \cite[Theorem 1.6.2]{AOF}), there exists a continuous curve $\zeta: [0, 1] \to K_{u_0} \cup K_u$ such that $\zeta(0) = x$ and $\zeta(1) = y$. Note that $x \in K_w$ and $y \in K_v$. By \eqref{OEX.eq100}, the curve $\zeta$ intersects with $f_w(V_0)$. Therefore, $(K_u \cup K_{u_0}) \cap f_w(V_0) \neq \emptyset$.  However, since $\diam{K_u, d_*} = \diam{K_{u_0}, d_*} = \rho^{m_0 + n}$, it follows
\[
d(f_w(x_0), K_u \cup K_{u_0}) \le 2\rho^{m_0 + n} < d(f_w(x_0), f_w(V_0)),
\]
so that $(K_{u_0} \cup K_u) \cap f_w(V_0) = \emptyset$. This contradiction shows that $u \in T(w)$ and hence $U_1(f_w(x_0): m_0 + n) \subseteq K_w$. By \eqref{SNS.eq20}, we see that
\[
B_{d_*}(f_w(x_0), c\rho^{m_0 + n}) \subseteq U_1(f_w(x_0): m_0 + n) \subseteq K_w.
\]
This shows (2C).
\enddemo

Next set $\a_H = -\log N/\log{\rho}$. Note that $\rho^{\a_H} = N^{-1}$. Let $\mu$ be the self-similar measure on $K$ with weight $(\rho^{\a_H}, \ldots, \rho^{\a_H})$. By \cite[Theorem~1.2.7]{Ki13}, we see that $\mu(K_w) = \rho^{|w|}$ for any $w \in T$ and consequently $\mu(\{x\}) = 0$ for any $x \in K_w$. These facts show that $\mu$ satisfies Assumption~\ref{ASS.30}. Moreover, we have the following proposition.

\prop\label{SNS.prop20}
Under Assumptions~\ref{OEX.ass10} and \ref{OEX.ass30}, there exist $c_1, c_2 > 0$ such that
\begin{equation}\label{SNS.eq40}
c_1s^{\a_H} \le \mu(B_{d_*}(x, s)) \le c_1s^{\a_H}.
\end{equation}
for any $s \in [0, 1]$. In particular, $\mu$ is $\a_H$-Ahlfors regular with respect to $d_*$ and the Hausdorff dimension of $(K, d_*)$ equals $\a_H$.
\endprop

\demo
By \eqref{SNS.eq30}, for any $x \in K$ and $n \ge 1$, if $w \in \GG_{1, n}(x)$, then
\begin{equation}\label{SNS.eq50}
(\rho^{n})^{\a_H} = \mu(K_w) \le \mu(B_{d_*}(x, 3c'\rho^n).
\end{equation}
On the other hand, by \cite[Proposition~1.6.11]{Ki13}, there exists $J_* \in \BbN$ such that
\begin{equation}\label{SNS.eq55}
\#(\GG_{1, n}(x)) \le J_*
\end{equation}
for any $x \in T$ and $n \ge 0$. (Note that $\LL^1_{\rho^n, x}$ defined in \cite[Definition~1.3.3]{Ki13} equals $\GG_{1, n}(x)$.) Therefore by \eqref{SNS.eq20}, 
\begin{equation}\label{SNS.eq60}
\mu(B_{d_*}(x, c\rho^n)) \le \sum_{v \in \GG_{1, n}(x)} \mu(K_v) \le J_*(\rho^n)^{\a_H}.
\end{equation}
Combining \eqref{SNS.eq50} and \eqref{SNS.eq60}, we obtain \eqref{SNS.eq40}.
\enddemo

The following proposition is immediately deduced from the previous propositions and lemmas. Note that $\GG_1(w) \subseteq \GG_{1, n}(x)$ for any $w \in T$ and $x \in K_w$. Hence by \eqref{SNS.eq55}, we see that the partition  $\{K_w\}_{w \in T}$ is uniformly finite.

\prop\label{SNS.prop30}{\rm [Proposition~\ref{OEX.prop20}]}
Under Assumptions~\ref{OEX.ass10} and \ref{OEX.ass30}, Assumption~\ref{ALFR} holds with $r = \rho$, $d = d_*$ and $M_* = M_0 = 1$.
\endprop

The fact that $M_0 = 1$ is due to the second remark after Assumption~\ref{ASS.10}.

\newpage
\section{List of Definitions and notations}

{\bf Definitions}\par\vspace{3pt}\noindent
adjacency matrix -- Definition~\ref{TWR.def10}\\
Ahlfors regular -- \eqref{ASS.eq50}\\
Ahlfors regular conformal dimension -- \eqref{INT.eq00}\\
Arzel{\`a}-Ascoli -- Appendix~\ref{APP1}\\
chipped Sierpinski carpet -- Example~\ref{CYP.ex10}\\
conductance constant -- Definition~\ref{MAC.def10}\\
conductively homogeneous(conductive homogeneity) -- Definition~\ref{CPE.def05}\\
critical set -- Definition~\ref{SNS.def10}\\
exponential -- Lemma~\ref{ASS.lemma20}\\
folding map --  Definition~\ref{CYP.def10}-(2)\\
graph -- Definition~\ref{TWR.def10}\\
hyperoctahedral group -- Definition~\ref{SHD.def00}\\
geodesic -- Definition~\ref{TWR.def10}-(3)\\
locally finite -- Definition~\ref{TWR.def10}-(1)\\
locally symmetric --  Definition~\ref{CYP.def10}-(4)\\
Markov property -- Theorem~\ref{CPE.thm10}-(c)\\
minimal -- Definition~\ref{ASS.def10}-(1)\\
Moulin -- Example~\ref{CYP.ex30}\\
neighbor disparity constant -- Definition~\ref{REL.def20}\\
nested fractal -- Definition~\ref{OEX.ex10}\\
non-degenerate -- Definition~\ref{CYP.def10}-(1)\\
$m$-walk -- Definition~\ref{OEX.ass20}\\
partition --Definition~\ref{PAS.def20}\\
path -- Definition~\ref{TWR.def10}-(2)\\
$p$-energy -- Theorem~\ref{CPE.thm10}\\
pentakun -- Example~\ref{OEX.ex10}\\
pinwheel -- Example~\ref{CYP.ex30}\\
Poincar{\'e} constant -- Definition~\ref{MAC.def120}\\
post critical set -- Definition~\ref{SNS.def10}\\
post critically finite -- Definition~\ref{SNS.def10}\\
p.c.f -- Definition~\ref{SNS.def10}\\
Sierpinski cross -- Section~\ref{SCR}\\
ray -- Definition~\ref{TWR.def20}\\
reference point-- Definition~\ref{TWR.def20}\\
root -- Definition~\ref{TWR.def20}\\
self-similar set -- \eqref{SSF.eq210}\\
self-similar structure -- Definition~\ref{SSF.def00}\\
simple -- Definition~\ref{TWR.def10}-(2)\\
snowflake -- Example~\ref{OEX.ex20}\\
strict $0$-walk -- Definition~\ref{OEX.ass20}\\
strongly connected --  Definition~\ref{CYP.def10}-(3)\\
strongly symmetric -- Definition~\ref{OEX.ass20}\\
sub-multiplicative inequality(conductance) -- Corollary~\ref{CMS.cor10}\\
sub-multiplicative inequality(Modulus) -- Theorem~\ref{CMS.thm10}\\
symmetry -- Definition~\ref{CHS.def10}\\
sub-multiplicative inequality(Neighbor disparity) --  Lemma~\ref{REL.lemma40}\\
subsystem of cubic tiling -- Definition~\ref{CYP.def10}\\
super-exponential -- Assumption~\ref{ASS.30}\\
tree -- Definition~\ref{TWR.def10}-(3)\\
uniformly finite -- Definition~\ref{ASS.def10}-(3)\\

\par

\vspace{5pt}\noindent
{\bf Notations}\par\vspace{3pt}\noindent
$\A_m^{(M)}(A_1, A_2, A)$ -- Definition~\ref{CMS.def10}-(2)\\
$\A_{N, m}^{(M)}(w)$ -- Definition~\ref{CMS.def10}-(3)\\
$A_s$ -- Definition~\ref{CYP.def10}\\
$B_d(x, r)$ -- Assumption~\ref{ALFR}\\
$B_{j, i}$ -- Definition~\ref{SHD.def00}\\
$\BB_L$ -- Definition~\ref{SHD.def00}\\
$B_{M, k}(w)$ -- Definition~\ref{MAC.def20}\\
$B_w$ -- Definition~\ref{ASS.def10}\\
$c_s^{L, N}$ -- Definition~\ref{SHD.def00}\\
$\ulc_{\E}(L, N, p), \olc_{\E}(L, N, p)$ -- Definition~\ref{FIN.def20}\\
$\ulc_{\lambda}(p, L, N), \olc_{\lambda}(p, L, N)$ -- Definition~\ref{FIN.def40}\\
$\ulc_{\s}(L, N, \kappa), \olc_{\s}(L, N, \kappa)$ -- Definition~\ref{FIN.def60}\\
$C_*^L$ -- Definition~\ref{SHD.def00}\\
$C_s^{L, N}$ -- Definition~\ref{SHD.def00}\\
$\C_m^{(M)}(A_1, A_2, A)$ -- Definition~\ref{CMS.def10}-(2)\\
$\C_{N, m}^{(M)}(w)$ -- Definition~\ref{CMS.def10}-(3)\\
$\diam{K, d}$ -- Assumption~\ref{ALFR}\\
$\dim_{AR}(K, d)$ -- \eqref{INT.eq00}\\
$\overline{D}_k$ -- Lemma~\ref{REX.lemma10}\\
$E_n^*$ -- Proposition~\ref{FRW.prop10}\\
$E_{M, n}^*$ -- Definition~\ref{CMS.def10}\\
$E_n^{\ell}$ --  Definition~\ref{CYP.def10}-(3)\\
$\E_{p, A}^n(\cdot)$ -- Definition~\ref{MAC.def10}\\
$\wE_p^m(\cdot)$ -- \eqref{CPE.eq200}, \eqref{SSF.eq00}\\
$\hE_p(\cdot)$ -- Theorem~\ref{CPE.thm10}\\
$\E_{p, m}(A_1, A_2, A)$ -- Definition~\ref{MAC.def10}\\
$\E_{M, p, m, n}$ -- Definition~\ref{CPE.def00}\\
$\E_{M, p, m}(w, A)$ -- Definition~\ref{MAC.def10}\\
$\overline{f}$ -- Definition~\ref{CPE.def100}\\
$g(w)$ -- \eqref{SSF.eq100}\\
$\G(L, N) $ -- Definition~\ref{FIN.def30}\\
$\G_{\E}(L, N)$ -- Definition~\ref{FIN.def10}\\
$\G_{\s}(L, N)$ -- Definition~\ref{FIN.def50}\\
$\G_{(K, T)}$ -- Definition~\ref{CHS.def10}\\
$h_{M, w, m}^*$ -- Definition~\ref{MAC.def30}\\
$h_{M_*, w}^*$ -- Lemma~\ref{CPE.lemma75}\\
$\H_{j_1, j_2}^i$ -- Definition~\ref{SHD.def10}\\
$I_{A, k, m}$ -- Lemma~\ref{MAC.lemma50}\\
$\hat{I}_{A, m}$ -- Lemma~\ref{MAC.lemma20}\\
$\wI_{A, k}$ -- Lemma~\ref{MAC.lemma40}\\
$\IT(K, T)$ -- Definition~\ref{CHS.def10}\\
$j(w)$ -- \eqref{SSF.eq100}\\
$J_n$ -- \eqref{CPE.eq500}\\
$K(\cdot)$ -- \eqref{SCR.eq10}\\
$K_T, K_B, K_R, K_L$ -- \eqref{SCR.eq20}\\
$\ell(\cdot)$ -- \eqref{CON.eq10}\\
$\ell_{w, v}$ -- \eqref{CYP.eq25}\\
$\ell_T, \ell_B, \ell_R, \ell_L$ -- Definition~\ref{SCR.def00}\\
$L_*$ -- \eqref{ASS.eq00}\\
$M_0$ -- Assumption~\ref{ASS.10}-(3), Assumption~\ref{ALFR}-(4)\\
$M_*$ -- Assumption~\ref{ASS.10}-(2), Assumption~\ref{ALFR}-(2)\\
$\M_{p, m}(A_1, A_2, A)$ -- Definition~\ref{CMS.def10}-(2)\\
$\M_{N, p, m}^{(M)}(w)$ -- Definition~\ref{CMS.def10}-(3)\\
$n_L(\cdot, \cdot)$ -- Definition~\ref{CPE.def10}\\
$\N_p(\cdot)$ -- Lemma~\ref{CPE.lemma200}\\
$N_*$ -- \eqref{ASS.eq05}\\
$O_w$ -- Definition~\ref{ASS.def10}\\
$P_n$ -- Definition~\ref{CPE.def20}\\
$P_{n, m}$ -- Lemma~\ref{REL.lemma20}\\
$\P(V, E)$ -- Definition~\ref{FIN.def30}\\
$\P(V, E, V_1, V_2, \kappa)$ -- Definition~\ref{FIN.def50}\\
$Q_n$ -- \eqref{CPE.eq600}\\
$R_j, R_{j_1, j_2}^i$ -- Definition~\ref{SHD.def10}\\
$R_{i, jk}, R^*_{i, jk}$ -- Definition~\ref{SCR.def10}\\
$S(w), S^m(w)$ -- Definition~\ref{TWR.def20}-(1)\\
$T_m$ -- Definition~\ref{TWR.def20}-(2)\\
$T_n^n, T_n^{n + 1}$ -- Lemma~\ref{SCR.lemma05}\\
$T(w)$ -- Definition~\ref{TWR.def20}-(3)\\
$U_M(w)$ -- Lemma~\ref{CPE.lemma75}\\
$U_M(x: n)$ -- Assumption~\ref{ALFR}\\
$|w|$ -- Definition~\ref{TWR.def20}-(2)\\
$\overline{wv}$ -- Definition~\ref{TWR.def10}-(3)\\
$\W^p$ -- Lemma~\ref{CPE.lemma200}\\
$X(e)$-- Definition~\ref{CHS.def10}\\
$\b_*$ -- Theorem~\ref{SUC.thm30}\\
$\c$ -- Assumption~\ref{ASS.30}\\
$\GG^A_M(w), \GG_M(w)$ -- Definition~\ref{ASS.def10}\\
$\delta_L(\cdot, \cdot)$ -- Definition~\ref{CPE.def10}\\
$\partial{S^m(w)}$ -- Definition~\ref{FRW.def10}\\
$\kappa$ -- Assumption~\ref{ASS.30}\\
$\lambda_{p, m}(A), \tla_{p, m}(A)$ -- Definition~\ref{MAC.def120}\\
$\ola_{p, m}$ -- Definition~\ref{REX.def10}\\
$\LL_{r^n}^g$ -- \eqref{SSF.eq110}\\
$\theta_m(\cdot, \cdot)$ -- Definition~\ref{CMS.def10}\\
$\Theta_{\frac {\pi}2}$ --  Theorem~\ref{CYP.thm20}\\
$\xi_n$ -- Lemma~\ref{REX.lemma20}\\
$\xi_n(w)$ -- Definition~\ref{POC.def10}\\
$\pi$ -- Definition~\ref{TWR.def20}\\
$\s$ -- Theorem~\ref{SUC.thm10}\\
$\s_{p, m}(w, v), \s_{p, m, n}, \s_{p, m}$ -- Definition~\ref{REL.def20}\\
$\tau$ -- Lemma~\ref{CPE.lemma20}\\
$\tau_p$ -- Lemma~\ref{SUC.lemma20}\\
$\tau_*$ -- Theorem~\ref{SUC.thm30}\\
$\Phi_s$ -- Definition~\ref{CYP.def10}\\
$\vp_e$ -- Definition~\ref{CHS.def10}\\
$\vp_{M, w, m}^*$ -- Definition~\ref{MAC.def30}\\
$\vp_{M_*, w}^*$ -- Lemma~\ref{CPE.lemma75}\\
$\psi_n$ -- Definition~\ref{CHS.def10}\\
$\psi^*_{n, m}$ -- Definition~\ref{SCR.def20}-(1)\\
$\Sigma$ -- Definition~\ref{TWR.def20}-(4)\\
$\#(\cdot)$ -- Definition~\ref{ASS.def10}\\
$\norm{\cdot}_{p, \mu}$ -- Lemma~\ref{CPE.lemma200}, Definition~\ref{MAC.def120}\\

\bibliography{papers,aof,book}

\providecommand{\bysame}{\leavevmode\hbox to3em{\hrulefill}\thinspace}
\providecommand{\MR}{\relax\ifhmode\unskip\space\fi MR }
\providecommand{\MRhref}[2]{%
  \href{http://www.ams.org/mathscinet-getitem?mr=#1}{#2}
}
\providecommand{\href}[2]{#2}
\begin{thebibliography}{10}

\bibitem{BB1}
M.~T. Barlow and R.~F. Bass, \emph{The construction of {Brownian} motion on the
  {Sierpinski} carpet}, Ann. Inst. Henri Poincar{\'e} \textbf{25} (1989),
  225--257.

\bibitem{BB2}
\bysame, \emph{Local time for {Brownian} motion on the {Sierpinski} carpet},
  Probab. Theory Related Fields \textbf{85} (1990), 91--104.

\bibitem{BB3}
\bysame, \emph{On the resistance of the {Sierpinski} carpet}, Proc. R. Soc.
  London A \textbf{431} (1990), 354--360.

\bibitem{BB4}
\bysame, \emph{Transition densities for {Brownian} motion on the {Sierpinski}
  carpet}, Probab. Theory Related Fields \textbf{91} (1992), 307--330.

\bibitem{BB5}
\bysame, \emph{Coupling and {Harnack} inequalities for {Sierpinski} carpets},
  Bull. Amer. Math. Soc. (N. S.) \textbf{29} (1993), 208--212.

\bibitem{BB6}
\bysame, \emph{{Brownian} motion and harmonic analysis on {Sierpinski}
  carpets}, Canad. J. Math. \textbf{51} (1999), 673--744.

\bibitem{BBKT}
M.~T. Barlow, R.~F. Bass, T.~Kumagai, and A.~Teplyaev, \emph{Uniqueness of
  {Brownian} motion on {Sierpinski} carpets}, J. Eur. Math. Soc. \textbf{12}
  (2010), 665--701.

\bibitem{BP}
M.~T. Barlow and E.~A. Perkins, \emph{{Brownian} motion on the {Sierpinski}
  gasket}, Probab. Theory Related Fields \textbf{79} (1988), 542--624.

\bibitem{Bjorn2}
A.~Bj{\"o}rn and J.~Bj{\"o}rn, \emph{Nonlinear potential theory on metric
  spaces}, Tracts in Math., vol.~17, European Math. Soc., 2011.

\bibitem{BonkSaks}
M.~Bonk and E.~Saksman, \emph{Sobolev spaces and hyperbolic fillings}, J. Reine
  Angew. Math. \textbf{737} (2018), 161--187.

\bibitem{BouKleiner}
M.~Bourdon and B.~Kleiner, \emph{Combinatorial modulus, the combinatorial
  {Loewner} property, and {Coxeter} groups}, Groups, Geometry, and Dynamics
  \textbf{7} (2013), 39--107.

\bibitem{Braides1}
A.~Braides, \emph{A handbook of {$\Gamma$}-convergence}, Handbook of
  Differential Equations: Stationary Partial Differential Equations, vol.~3,
  Elsevier/North-Holland, 2006, pp.~101--213.

\bibitem{CaoQui}
S.~Cao and H.~Qui, \emph{Dirichlet forms on unconstrained {Sierpinski}
  carpets}, preprint.

\bibitem{CarPiag}
M.~{Carrasco Piaggio}, \emph{On the conformal gauge of a compact metric space},
  Ann. Sci. Ecole. Norm. Sup. \textbf{46} (2013), 495--548.

\bibitem{Cheeger}
J.~Cheeger, \emph{Differentiability of {Lipschitz} functions on metric measure
  spaces}, Geom. Funct. Anal. \textbf{3} (1999), 428--517.

\bibitem{ChenFuku}
Z.-Q. Chen and M.~Fukushima, \emph{{Symmetric Markov processes, Time Change,
  and Boundary Theory}}, London Math. Soc. Monographs, vol. 35, Princeton Univ.
  Press, 2012.

\bibitem{Dav2}
E.~B. Davies, \emph{{Spectral Theory and Differential Operators}}, Cambridge
  Studies in Advanced Math. vol. 42, Cambridge University Press, 1995.

\bibitem{DS}
P.~G. Doyle and J.~L. Snell, \emph{{Random Walks and Electrical Networks}},
  Math. Assoc. Amer., Washington, 1984.

\bibitem{FOT}
M.~Fukushima, Y.~Oshima, and M.~Takeda, \emph{{Dirichlet Forms and Symmetric
  Markov Processes}}, de Gruyter Studies in Math. vol. 19, de Gruyter, Berlin,
  1994.

\bibitem{Gri1}
A.~Grigor'yan, \emph{The heat equation on noncompact {Riemannian} manifolds.
  (in {Russian})}, Mat. Sb. \textbf{182} (1991), 55--87, English translation in
  Math. USSR-Sb. 72(1992), 47--77.

\bibitem{Griyang}
A.~Grigor'yan and M.~Yang, \emph{Local and non-local {Dirichlet} forms on the
  {Sierpinski} carpet}, Trans. A.M.S. \textbf{372} (2019), 3985--4030.

\bibitem{Haj1}
P.~Haj{\l}asz, \emph{Sobolev spaces on an arbitrary metric spaces}, Pontential
  Analysis \textbf{5} (1996), 403--415.

\bibitem{HeiKoShTy}
J.~Heinonen, P.~Koskela, N.~Shanmugalingam, and J.~Tyson, \emph{{Sobolev Spaces
  on Metric Measure Spaces - An Approach Based on Upper Gradients}}, New
  Mathematical Monographs, Cambridge University Press, 2015.

\bibitem{HerPeiStr}
P.~E. Herman, R.~Peirone, and R.~S. Strichartz, \emph{$p$-energy and
  $p$-harmonic fufnction on {Sierpinski} gasket type fractals}, Potential Anal.
  \textbf{20} (2004), 125--148.

\bibitem{Hino0}
M.~Hino, \emph{On short time asymptotic behavior of some symmetric diffusions
  on general state spaces}, Potential Anal. \textbf{16} (2002), 249--264.

\bibitem{KajMur3}
N.~Kajino and M.~Murugan, \emph{On the conformal walk dimension {II}:
  Non-attainment for some {Sierpi\'{n}ski} carpets}, in preparation.

\bibitem{KajMur2}
\bysame, \emph{On the conformal walk dimension: Quasisymmetric uniformization
  for symmetric diffusions}, preprint.

\bibitem{Ki10}
J.~Kigami, \emph{Markov property of {Kusuoka-Zhou's Dirichlet} forms on
  self-similar sets}, J. Math. Sci. Univ. Tokyo \textbf{7} (2000), 27--33.

\bibitem{AOF}
\bysame, \emph{{Analysis on Fractals}}, Cambridge Tracts in Math. vol. 143,
  Cambridge University Press, 2001.

\bibitem{Ki15}
\bysame, \emph{{Measurable Riemannian geometry on the Sierpinski gasket: the
  Kusuoka measure and the Gaussian heat kernel estimate}}, Math. Ann.
  \textbf{340} (2008), 781--804.

\bibitem{Ki13}
\bysame, \emph{Volume doubling measures and heat kernel estimates on
  self-similar sets}, Memoirs of the American Mathematical Society \textbf{199}
  (2009), no.~932.

\bibitem{Ki16}
\bysame, \emph{Resistance forms, quasisymmetric maps and heat kernel
  estimates}, Memoirs of the American Mathematical Society \textbf{216} (2012),
  no.~1015.

\bibitem{Ki21}
\bysame, \emph{Time changes of the {Brownian} motion: Poincar{\'e} inequality,
  heat kernel estimate and protodistance}, Memoirs of the American Mathematical
  Society \textbf{259} (2019), no.~1250.

\bibitem{GAMS}
\bysame, \emph{{Geometry and Analysis of Metric Spaces via Weighted
  Partitions}}, Lecture Notes in Math. vol. 2265, Springer, 2020.

\bibitem{Kum3}
T.~Kumagai, \emph{Estimates of the transition densities for {Brownian} motion
  on nested fractals}, Probab. Theory Related Fields \textbf{96} (1993),
  205--224.

\bibitem{KusZ1}
S.~Kusuoka and X.~Y. Zhou, \emph{{Dirichlet} forms on fractals: {Poincar{\'e}}
  constant and resistance}, Probab. Theory Related Fields \textbf{93} (1992),
  169--196.

\bibitem{Li1}
T.~Lindstr{\o}m, \emph{{Brownian} motion on nested fractals}, Mem. Amer. Math.
  Soc. \textbf{420} (1990).

\bibitem{PPStos}
K.~Pietruska-Paluba and A.~St{\'o}s, \emph{Poincar{\'e} inequality and
  {Haj{\l}asz-Sobolev} spaces on nested fractals}, Studia Math. \textbf{218}
  (2013), 1--26.

\bibitem{Saloff1}
L.~Saloff-Coste, \emph{A note on {Poincar{\'e}, Sobolev, and Harnack}
  inequalities}, Internat. Math. Res. Notices (1992), 27--38.

\bibitem{Shanm1}
N.~Shanmugalingam, \emph{Newtonian spaces: an extension of {Sobolev} spaces to
  metric measure spaces}, Rev. Mat. Iberoamer. \textbf{16} (2000), 243--279.

\bibitem{Shimizu2}
R.~Shimizu, \emph{Construction of $p$-energy and associated energy measure on
  the {Sierpinski} carpet}, preprint.

\bibitem{Sturm1}
K.T. Sturm, \emph{Analysis on local {Dirichlet} spaces {III}. the parabolic
  {Harnack} inequality}, J. Math. Pures Appl. \textbf{75} (1996), 275--297.

\bibitem{TysonWu1}
J.~T. Tyson and J.-M. Wu, \emph{On quasiconformal dimensions of self-similar
  fractals}, Rev. Mat. Iberoamericana \textbf{22} (2006), 205--258.

\bibitem{Woess3}
W.~Woess, \emph{{Denumerable Markov Chains}}, European Math. Soc., 2009.

\bibitem{Yosida}
K.~Yosida, \emph{{Functional Analysis, sixth ed.}}, Classics in Math.,
  Springer, 1995, originally published in 1980 as Grundlehren der
  mathematischen Wissenschaften band 123.

\end{thebibliography}
\bibliographystyle{amsplain}

\end{document}